\documentclass[11pt]{amsart}

\usepackage{amsfonts,amssymb,mathtools,mathrsfs} \usepackage{microtype} \usepackage{enumitem} \usepackage{booktabs,array,longtable} \usepackage{graphicx,float} \usepackage[hmargin=1.15in,vmargin=1.10in]{geometry} \IfFileExists{xurl.sty}{\usepackage{xurl}}{} \usepackage{hyperref}

\allowdisplaybreaks[3]  

\newtheorem{theorem}{Theorem}[section] \newtheorem{lemma}[theorem]{Lemma} \newtheorem{proposition}[theorem]{Proposition} \newtheorem{corollary}[theorem]{Corollary} \newtheorem{appendixBlemma}{Lemma}  \theoremstyle{definition} 

\newcommand{\proofpart}{front}          

\newcommand{\TheTitle}{A Complete Resolution of Forsythe's Conjecture\\for Restarted Conjugate Gradients} \newcommand{\ThePDFTitle}{A Complete Resolution of Forsythe's Conjecture for Restarted Conjugate Gradients} \newcommand{\TheAuthors}{Matthew J. Colbrook, George Stepaniants, and Alex Townsend} \title[Complete resolution of Forsythe's conjecture]{\TheTitle} \author[M. J. Colbrook]{Matthew J. Colbrook} \address{Department of Applied Mathematics and Theoretical Physics\\ University of Cambridge\\Cambridge CB3 0WA, United Kingdom} \email{m.colbrook@damtp.cam.ac.uk} \author[G. Stepaniants]{George Stepaniants} \address{Department of Computing and Mathematical Sciences\\ California Institute of Technology\\Pasadena, California 91125, USA} \email{gstepan@caltech.edu} \author[A. Townsend]{Alex Townsend} \address{Department of Mathematics\\Cornell University\\ Ithaca, New York 14853, USA} \email{townsend@cornell.edu}

\date{5 September 2026}


\hypersetup{hidelinks,pdftitle={\ThePDFTitle},pdfauthor={\TheAuthors}} 

\begin{document}

\begin{abstract} Forsythe's conjecture, published in 1968, asserts that for each restart length $s$, every exact-arithmetic restarted conjugate-gradient iteration on a real symmetric positive definite problem either terminates or has normalised residuals that converge separately along the even and odd restart subsequences. Apart from the classical steepest-descent case, this asymptotic question remained unresolved in full generality for nearly six decades. We give a complete classification by restart length in the original finite-dimensional setting and identify a sharp threshold. For $s=2$ and $s=3$, every problem either terminates or has convergent even and odd residual directions. For every $s\ge4$, there is a diagonal positive definite counterexample of dimension $s+4$ which never terminates and whose even residual directions do not converge. Together with Akaike's theorem for $s=1$, this shows that the conjectured universal conclusion is true precisely for $s\in\{1,2,3\}$ and false for every $s\ge4$. The positive results follow from a degree-independent double-orthogonality identity and an analysis of the low-degree fixed-point sets. At restart length four, rational interval arithmetic and Sturm sequences certify a transverse Hopf point of the leading vector field of the rescaled squared-weight map. Analytic periodic-orbit and shadowing arguments yield the counterexample at restart length four, and degree elevation extends the construction to every larger restart length. The classification for all $s\ge2$ is also formally verified in Lean. \end{abstract}

\maketitle

\setcounter{tocdepth}{2} \tableofcontents

\section{Introduction}\label{sec:introduction}

Forsythe's conjecture, first formulated in a 1967 report and published in 1968 \cite{Forsythe1967,Forsythe1968}, asks whether every nonterminating sequence of normalised residual directions at the restart points approaches a two-cycle---that is, whether its even and odd subsequences converge separately. The classification by restart length remained unresolved for nearly six decades. We complete it and identify a sharp threshold at restart length four. Together with Akaike's theorem for $s=1$, our results show that the termination-or-two-cycle alternative holds universally for $s\in\{1,2,3\}$, whereas for every $s\ge4$ it fails already for a diagonal positive definite matrix of dimension $s+4$.

The conjugate gradient method, introduced by Hestenes and Stiefel in 1952, is a standard Krylov method for real symmetric positive definite (SPD) linear systems \cite{HestenesStiefel1952,LiesenStrakos2013}. Classical convergence estimates quantify the decay of the error and residual norm, but do not in general determine the asymptotic direction of the residual. Under fixed-length restarting, the normalised residual evolves by a nonlinear map of the unit sphere. The conjecture therefore concerns directional dynamics not captured by the classical norm bounds.

Let $A\in\mathbb R^{n\times n}$ be symmetric positive definite, let $x_*=A^{-1}b$, and define $r_k=b-Ax_k$. Restarted CG with restart length $s$ is the exact iteration $$ x_{k+1}\in x_k+\mathcal K_s(A,r_k),\qquad \|x_*-x_{k+1}\|_A =\min_{z\in x_k+\mathcal K_s(A,r_k)}\|x_*-z\|_A, $$ where $\mathcal K_s(A,r)=\operatorname{span}\{r,Ar,\ldots,A^{s-1}r\}$. Here $k$ counts restart blocks and $\|v\|_A^2=v^\top Av$. This is Forsythe's $s$-dimensional optimum gradient method.

The precise statement for $s\ge2$ is as follows.

\begin{theorem}[Sharp classification]\label{thm:sharp} Fix an integer $s\ge2$ and carry out restarted CG in exact arithmetic. For every instance, let $r_k=b-Ax_k$ and, whenever $r_k\ne0$, let $y_k=r_k/\|r_k\|$. \begin{enumerate}[label=\textup{(\roman*)}] \item If $s\in\{2,3\}$, then for every $n\ge1$, every real symmetric positive definite matrix $A\in\mathbb R^{n\times n}$, and all $b,x_0\in\mathbb R^n$, the iteration either reaches $r_K=0$ for some finite $K$, or else there are unit vectors $y_{\rm e},y_{\rm o}$ such that $$ y_{2k}\longrightarrow y_{\rm e},\qquad y_{2k+1}\longrightarrow y_{\rm o}. $$ \item If $s\ge4$, there are $b,x_0\in\mathbb R^{s+4}$ and a diagonal positive definite matrix $A\in\mathbb R^{(s+4)\times(s+4)}$ for which the iteration never terminates and $(y_{2k})$ does not converge. \end{enumerate} \end{theorem}

Theorem~\ref{thm:sharp} has also been formalised in Lean, including both the positive cases and the diagonal counterexamples in exactly dimension $s+4$. All numerical claims needed for the formal proof are established in Lean. Comparator, a separate verification tool, checked that the results proved in Lean match the specified formal statements, using Lean's core proof checker to verify the proofs and all the results they rely on. \hyperref[app:lean]{Appendix D} lists the relevant statements and source files, explains how to review their mathematical meaning, and gives instructions for repeating the verification.

When $s=1$, the iteration is the method of steepest descent. Forsythe and Motzkin formulated its two-direction asymptotics in 1951, and Akaike proved them in 1959 \cite{ForsytheMotzkin1951,Akaike1959}. Motivated by this theorem and by computations for $s=2$, Forsythe conjectured the same behaviour for every fixed finite $s$. In the formulation considered here, the assertion is universal: for each $s$, it ranges over every dimension, every SPD matrix, and every initial residual.

A persistent difficulty in the earlier literature is that convergence of iteration parameters need not imply convergence of residual directions. For $s=2$, Zhuk and Bondarenko established convergence of the associated polynomial coefficients and claimed convergence of the residual directions \cite{ZhukBondarenko1984}. Faber, Liesen, and Tich\'y observed that the result cited from Zabolotskaya \cite[p.~238]{Zabolotskaya1979} does not supply the final implication and consequently left the first nontrivial restarted case open \cite[Section~4.1]{FaberLiesenTichy2023}. Two contemporaneous manuscripts treating the case $s=2$ appeared in 2026 \cite{Liesen2026,ColbrookStepaniantsTownsend2026}; Part A gives a revised, self-contained treatment of the argument in the latter, adapted to the unified framework of this paper.

For general restart length, Zhuk proved convergence of the iteration parameters in a related Hilbert-space setting \cite{Zhuk1995}; this again does not determine a limiting residual direction. He also constructed direction-invariant states and a sufficient local attraction condition \cite{Zhuk1982}. The probability-measure dynamics, including invariant measures with small support, were studied by Pronzato, Wynn, and Zhigljavsky \cite{PronzatoWynnZhigljavsky2009}. For a survey of delayed gradient methods for symmetric positive definite linear systems, including their spectral interpretation, see \cite{ZouMagoules2022}. A recent workshop survey describes the general problem and discusses the contemporaneous work on the case $s=2$ \cite[Section~2.7]{AmselEtAl2026}. These results describe several aspects of the dynamics but do not decide the universal two-cycle question for arbitrary restart length.

The threshold is dynamical. A useful way to view the proof is that each restart reweights the spectral mass at every eigenvalue by the square of an orthogonal polynomial. Over two restarts, a monotone energy controls the change in these polynomials, drawing every nonterminating orbit toward a family of exact two-cycles. For $s=2,3$, the low-degree algebra is rigid: sign-definite drift and logarithmic balance laws rule out sustained motion along the family, so each residual parity selects a single limit. At $s=4$, two neutral external spectral modes make persistent oscillatory motion possible in the effective slow dynamics. A transverse Hopf point produces a periodic profile in the leading rescaled dynamics, and an exact restarted-CG orbit shadows it rather than settling. Adding distant eigenvalues preserves this mechanism for every $s>4$.

The proof separates a degree-independent spectral core from the mechanisms on the two sides of the threshold. In an eigenbasis of $A$, the squared coordinates of a normalised residual form a probability vector $w$, and one restart is governed by the monic degree-$s$ polynomial orthogonal to lower-degree polynomials for the discrete measure $\sum_iw_i\delta_{\lambda_i}$. A monotone two-block energy and a double-orthogonality identity imply convergence of the even and odd orthogonal factors. They also show that every limit of a nonterminating orbit is supported on between $s+1$ and $2s$ eigenvalues.

The remaining geometry depends on $s$. For $s=2$, the only possibilities are three- and four-node limits. The four-node two-cycles form a compact one-parameter family; a logarithmic cocycle and a retraction that remains uniform at its endpoints give finite total variation of the parameter. For $s=3$, the limiting support has four, five, or six nodes. The product of the two limiting orthogonal cubics determines a degree-six completion polynomial. The proof controls motion transverse to the associated family of two-cycles before treating stable, unstable, and neutral external spectral modes. The neutral case is decided by logarithmic cocycles and an analysis of the strata on which a weight vanishes.

\begin{figure}[t]
\centering
\begin{minipage}[t]{0.49\textwidth}
  \centering
  \includegraphics[width=\linewidth]{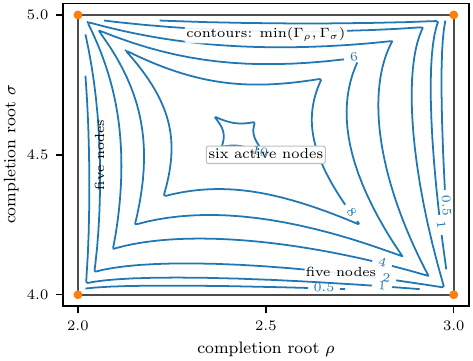}
\end{minipage}\hfill
\begin{minipage}[t]{0.49\textwidth}
  \centering
  \includegraphics[width=\linewidth]{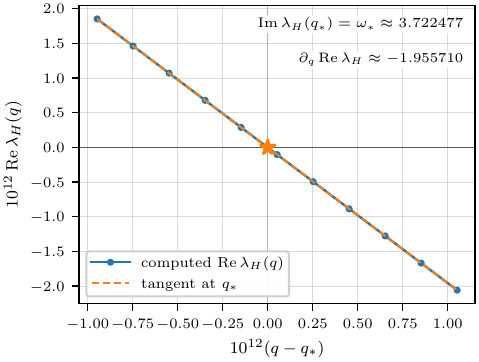}
\end{minipage}
\caption{The mechanism at the threshold. Left: at restart length three, the limiting two-cycles form a parameter rectangle whose interior, edges, and corners correspond to six-, five-, and four-node states; the contours show the positive interior control used in Part B. This is fixed-set geometry, not an orbit. Right: a numerical reconstruction of the certified transverse Hopf crossing at restart length four. The solid curve is the real part of the critical eigenvalue branch and the dashed curve is its tangent; the zero crossing with nonzero slope shows transversality.}
\label{fig:threshold-mechanism}
\end{figure}

At $s=4$, a transverse Hopf point in the leading vector field of the rescaled squared-weight map is used to construct a nonconstant periodic slow profile. An exact restarted-CG orbit shadows this profile on harmonic time scales, and a degree-elevation argument carries the construction to every larger restart length. The finite numerical input is confined to restart length four; the periodic-orbit, shadowing, and degree-elevation arguments are analytic. Figure~\ref{fig:threshold-mechanism} illustrates the two mechanisms at the threshold.

The proof begins with the degree-independent argument. Parts A and B treat restart lengths two and three, respectively; Part C constructs the counterexamples and carries out degree elevation. Appendix B contains the local calculations at the six-node boundary in Part B, Appendix C describes the certificate for restart length four, and Appendix D documents the Lean formalisation. Equation and theorem numbers carry the letter of the relevant part.

\section*{AI-assistance}
\label{sec:ai-process}

This manuscript can be read, in part, as an experiment by three numerical analysts, asking how far a frontier language model could be pushed when its suggestions were assessed by experts. During this project, all three authors used ChatGPT with OpenAI's GPT-5.6 model, sometimes independently and sometimes while discussing the responses together. We repeatedly worked through promising ideas on our own and checked arguments. GPT-5.6 and GPT-6 were used to formalise all proofs in this paper in Lean. The process lasted several weeks.

We began with $s=2$. The proof was developed by alternating conversations with the model and discussions among the authors: once a possible route appeared, we asked for simplifications and explanations, and then thoroughly checked the resulting mathematics ourselves. This led to the proof first circulated in our separate $s=2$ preprint \cite{ColbrookStepaniantsTownsend2026}. Part A grew out of that proof, but has since been revised and reorganised to fit the unified framework of the present paper. By the time we were ready to circulate that preprint, we learned that two independent $s=2$ manuscripts were already public, by Jarek Liesen and Richard Peng; a subsequent status update records all three manuscripts \cite[Section~2.7]{AmselEtAl2026}.

After that we continued to $s=3$. This was qualitatively harder. Several days of exchanges produced many incomplete arguments before the interaction eventually uncovered the structural route that became Part B. Once $s=3$ had been settled, we asked GPT-5.6 a different question. Rather than merely asking it to push the same proof one degree further, we asked it to explain why the arguments for $s=2$ and $s=3$ worked and what qualitatively new behaviour should become possible in higher degree. In that discussion, persistent oscillatory motion in an effective slow dynamics emerged as a possible obstruction to convergence. This suggested a much more concrete next question: could the slow dynamics near a family of two-cycles possess a nonconstant periodic orbit for larger $s$?

Over roughly two further weeks, candidate mechanisms were proposed, computed, rejected, and repaired. The idea was eventually realised at restart length four: a transverse Hopf point in the leading vector field of the rescaled squared-weight map is used to construct a nonconstant periodic orbit of that field, and a shadowing argument produces an exact restarted-CG orbit that follows this periodic profile on harmonic time scales. The degree-elevation construction then carries the mechanism to every $s\ge4$. The counterexample became visible only after the positive cases had been understood and the model had been repeatedly pressed to explain the structure behind them.

Finally, GPT-5.6 and GPT-6 enabled us to formalise and verify all proofs in this paper in Lean, a process which took two additional weeks to complete.

Forsythe's conjecture had remained open for decades, but resolving the remaining cases required a large amount of technically specialised work relative to the likely mathematical reward. Without AI assistance, we doubt that we would have invested the time needed to pursue the classification through $s=3$, discover the mechanism at $s=4$, and turn it into a proof for every $s\ge4$. The cost-reward ratio would simply have been unfavourable, and the conjecture might plausibly have remained unresolved for considerably longer. GPT-5.6 changed that ratio for us on this conjecture.

\section*{The degree-independent spectral core}

Throughout, $\mathcal P_j$ denotes the real polynomials of degree at most $j$. For a polynomial $R$ and any displayed indeterminate $z$, the symbol $[z^m]R$ denotes the coefficient of $z^m$, and $\|R\|_{\rm coeff}$ denotes the Euclidean norm of its coefficient vector in the monomial basis.

This preliminary argument confines every nonterminating orbit to a finite resonant support. The geometry of the resulting fixed set depends on the restart length and is treated in Parts A--C.

Fix $s\ge1$. If $r_0=0$, the iteration has already terminated; hence suppose $r_0\ne0$. We orthogonally diagonalise $A$ and retain only those distinct eigenvalues for which the spectral projection of $r_0$ is nonzero. Index these eigenvalues so that $$ 0<\lambda_1<\cdots<\lambda_N, $$ and define \begin{equation*} v_i=\frac{\operatorname{proj}_{\ker(A-\lambda_iI)}r_0} {\|\operatorname{proj}_{\ker(A-\lambda_iI)}r_0\|}. \end{equation*} The $v_i$ are orthonormal. Polynomial functional calculus shows that every nonzero residual direction has the unique form \begin{equation*} y=\sum_{i=1}^N\eta_i v_i,\qquad \sum_i\eta_i^2=1, \end{equation*} and no other direction in a repeated eigenspace can occur. Define \begin{equation*} w_i=\eta_i^2,\qquad \langle f,g\rangle_w=\sum_iw_if(\lambda_i)g(\lambda_i). \end{equation*} We identify $y$ with its coefficient vector when no confusion can arise, denote its $i$th coordinate by $y_i=\eta_i$, and set $w(y)=(y_i^2)_i$; along an orbit, $w_k=w(y_k)$. If a node $\xi$ equals $\lambda_i$, a subscript $\xi$ denotes the same coordinate as the subscript $i$, so $w_\xi=w_i$ and $w_{k,\xi}=w_{k,i}$. The vector $w$ is a probability vector on the fixed set of active spectral nodes. An initially zero coordinate remains zero under every later block. Ratios of signed coordinates occur only on punctured charts where the denominator is nonzero. If a coordinate vanishes, the dynamics continues on the corresponding invariant face of lower support.

If $w$ has more than $s$ active nodes, its moment Gram matrix on $\mathcal P_{s-1}$ is positive definite, so there is a unique monic degree-$s$ polynomial $P_w$ orthogonal to $\mathcal P_{s-1}$. Define \begin{equation*} H(w)=\langle P_w,P_w\rangle_w. \end{equation*} Then $H(w)>0$. Moreover, $P_w$ has $s$ real, simple roots in the open convex hull of the active positive nodes. To see this, suppose that $P_w$ has fewer than $s$ sign changes between the ordered active nodes. The product of the corresponding sign-change factors, with a suitable overall sign, is a polynomial $R\in\mathcal P_{s-1}$ for which $P_wR$ is nonnegative on the support and does not vanish there identically. This contradicts $\langle P_w,R\rangle_w=0$. All $s$ roots must therefore be real and simple, and must lie strictly between the extreme active nodes. In particular, $\operatorname{sign}P_w(0)=(-1)^s$. The Galerkin equations for one exact restart block consequently give \begin{equation*} (Tw)_i=\frac{w_iP_w(\lambda_i)^2}{H(w)},\qquad F(y)_i=(-1)^s\frac{P_w(\lambda_i)\eta_i}{\sqrt{H(w)}}. \tag{G.5}\label{eq:common-block-map} \end{equation*} The unnormalised residual polynomial is $P_w/P_w(0)$; the sign in \eqref{eq:common-block-map} is the sign of $P_w(0)$. If the support has at most $s$ nodes, a polynomial of degree at most $s$ vanishing on that support annihilates the residual, so the block terminates. Conversely, termination implies that the grade, hence the number of active distinct eigenvalues, is at most $s$. This includes termination at the initial state.

Assume for the remainder of this section that the orbit never terminates. At every block boundary, define \begin{equation*} P_k=P_{w_k},\qquad H_k=H(w_k),\qquad q_k=P_kP_{k+1},\qquad \Delta_k=q_k-H_{k+1}. \end{equation*} The signed two-block recurrence is \begin{equation*} y_{k+2,i}=\frac{q_k(\lambda_i)}{\sqrt{H_kH_{k+1}}}\,y_{k,i}. \tag{G.7}\label{eq:common-two-block} \end{equation*}

\begin{lemma}[Energy and same-parity chords]\label{lem:common-energy} For every $k$, \begin{equation*} H_{k+1}-H_k=\frac1{H_{k+1}} \sum_iw_{k,i}\Delta_k(\lambda_i)^2, \tag{G.8}\label{eq:common-variance} \end{equation*} and \begin{equation*} \sqrt{H_{k+1}}-\sqrt{H_k} =\frac{\sqrt{H_{k+1}}}{2}\|y_{k+2}-y_k\|^2. \tag{G.9}\label{eq:common-chord} \end{equation*} Then $H_k\uparrow h$ for some $h>0$ and \begin{equation*} \sum_{k=0}^{\infty}(H_{k+1}-H_k)<\infty, \qquad \sum_{k=0}^{\infty}\|y_{k+2}-y_k\|^2<\infty. \tag{G.10}\label{eq:common-sums} \end{equation*} \end{lemma}

\begin{proof} Because $P_{k+1}-P_k\in\mathcal P_{s-1}$, orthogonality of $P_k$ gives \begin{equation*} \sum_iw_{k,i}q_k(\lambda_i) =\langle P_k,P_{k+1}\rangle_{w_k}=H_k. \tag{G.11}\label{eq:common-first-moment} \end{equation*} Using the weight update once gives \begin{equation*} \sum_iw_{k,i}q_k(\lambda_i)^2 =H_k\sum_iw_{k+1,i}P_{k+1}(\lambda_i)^2 =H_kH_{k+1}. \tag{G.12}\label{eq:common-second-moment} \end{equation*} Expanding $\sum_iw_{k,i}(q_k(\lambda_i)-H_{k+1})^2$ and substituting \eqref{eq:common-first-moment} and \eqref{eq:common-second-moment} gives $H_{k+1}(H_{k+1}-H_k)$, which is \eqref{eq:common-variance}. Equation \eqref{eq:common-two-block} then gives $$ \langle y_{k+2},y_k\rangle =\frac{\sum_iw_{k,i}q_k(\lambda_i)} {\sqrt{H_kH_{k+1}}} =\sqrt{H_k/H_{k+1}}. $$ Both vectors are unit vectors, so $\|y_{k+2}-y_k\|^2=2-2\sqrt{H_k/H_{k+1}}$, which is \eqref{eq:common-chord}. The fixed monic competitor $t^s$ gives $$ H_k\le\sum_iw_{k,i}\lambda_i^{2s}\le\lambda_N^{2s}. $$ Equation \eqref{eq:common-variance} shows that $H_k$ increases to a finite limit $h\ge H_0>0$. The first sum in \eqref{eq:common-sums} telescopes. Since $\sqrt{H_{k+1}}\ge\sqrt{H_0}$, summing \eqref{eq:common-chord} proves the second. \end{proof}

The following identity replaces the separate coefficient-drift arguments for restart lengths two and three.

\begin{lemma}[Degree-independent double orthogonality] \label{lem:common-double-orthogonality} For every $\varphi\in\mathcal P_{s-1}$, \begin{equation*} \left\langle P_{k+2}-P_k,\varphi\right\rangle_{w_{k+2}} =-\frac1{H_kH_{k+1}} \left\langle P_k\varphi,\Delta_k^2\right\rangle_{w_k}. \tag{G.13}\label{eq:common-double-orthogonality} \end{equation*} \end{lemma}

\begin{proof} Two applications of \eqref{eq:common-block-map} give \begin{equation*} w_{k+2,i}=\frac{w_{k,i}q_k(\lambda_i)^2}{H_kH_{k+1}}. \end{equation*} Since $P_{k+2}\perp\mathcal P_{s-1}$ under $w_{k+2}$, \begin{equation*} \left\langle P_{k+2}-P_k,\varphi\right\rangle_{w_{k+2}} =-\frac1{H_kH_{k+1}} \langle P_k\varphi,q_k^2\rangle_{w_k}. \tag{G.15}\label{eq:common-double-expand} \end{equation*} Since $q_k=H_{k+1}+\Delta_k$, the constant term in the square on the right of \eqref{eq:common-double-expand} vanishes because $P_k\perp\mathcal P_{s-1}$. The linear term vanishes as well, because $$ \begin{aligned} \langle P_k\varphi,\Delta_k\rangle_{w_k} &=\sum_iw_{k,i}P_k(\lambda_i)^2 P_{k+1}(\lambda_i)\varphi(\lambda_i)-H_{k+1}\langle P_k,\varphi\rangle_{w_k}=H_k\langle P_{k+1},\varphi\rangle_{w_{k+1}}=0. \end{aligned} $$ Only the $\Delta_k^2$ term remains, and \eqref{eq:common-double-orthogonality} follows. \end{proof}

The moment identity yields coefficient convergence once the moment Gram matrices are uniformly invertible. For a subset $J$ of the original active nodes with $|J|\le s$, define \begin{equation*} R_J(t)=t^{s-|J|}\prod_{j\in J}(t-\lambda_j),\qquad M_J=\max_{i\notin J}R_J(\lambda_i)^2. \end{equation*} The original orbit has more than $s$ active nodes, so $M_J>0$. Since $R_J$ is an admissible monic competitor, \begin{equation*} H_0\le H_k\le\sum_iw_{k,i}R_J(\lambda_i)^2 \le M_J\sum_{i\notin J}w_{k,i}. \tag{G.17}\label{eq:common-face-separation} \end{equation*} Equation \eqref{eq:common-face-separation} shows that no $\omega$-limit weight is supported on at most $s$ nodes. At every $\omega$-limit point, the moment Gram matrix of $1,t,\ldots,t^{s-1}$ is therefore positive definite. Compactness of the $\omega$-limit set gives a common positive lower bound for its least eigenvalue, and all sufficiently late moment Gram matrices have uniformly bounded inverses. Enlarging the constants accounts for the finitely many earlier indices.

\begin{proposition}[Factor convergence and spectral localisation] \label{prop:common-factor-convergence} There are monic degree-$s$ polynomials $P_{\rm e},P_{\rm o}$ and a constant $C<\infty$ such that \begin{equation*} \|P_{k+2}-P_k\|_{\rm coeff} \le C(H_{k+1}-H_k), \tag{G.18}\label{eq:common-factor-step} \end{equation*} \begin{equation*} P_{2n}\longrightarrow P_{\rm e},\qquad P_{2n+1}\longrightarrow P_{\rm o},\qquad q_k\longrightarrow q_\infty:=P_{\rm e}P_{\rm o}, \tag{G.19}\label{eq:common-factor-limits} \end{equation*} and \begin{equation*} \|q_k-q_\infty\|_{\rm coeff}\le C(h-H_k). \tag{G.20}\label{eq:common-product-tail} \end{equation*} Every $\omega$-limit state is fixed by the signed two-block map. Its support has between $s+1$ and $2s$ nodes and is contained in \begin{equation*} E=\{\lambda_i:q_\infty(\lambda_i)=h\}. \tag{G.21}\label{eq:common-resonant-set} \end{equation*} Moreover, every nonresonant coordinate has summable squared mass: \begin{equation*} \sum_{k=0}^{\infty}w_{k,i}<\infty \qquad(\lambda_i\notin E). \tag{G.22}\label{eq:common-external-summable} \end{equation*} With $\mathscr D_k=q_k-\sqrt{H_kH_{k+1}}$ and $\mathscr D_*=q_\infty-h$, \begin{equation*} \|\mathscr D_k-\mathscr D_*\|_{\rm coeff} \le C(h-H_k). \tag{G.23}\label{eq:common-completion-tail} \end{equation*} \end{proposition}

\begin{proof} Since the nodes lie in a fixed finite interval, the right-hand sides of the moment equations are uniformly bounded. The uniform bounds for the inverse Gram matrices therefore give a uniform coefficient bound for $P_k$. Define $D_k=P_{k+2}-P_k\in\mathcal P_{s-1}$. For each monomial $\varphi=t^j$, $0\le j<s$, Lemma~\ref{lem:common-double-orthogonality} and \eqref{eq:common-variance} give \begin{equation*} \begin{aligned} |\langle D_k,\varphi\rangle_{w_{k+2}}| &\le\frac{\max_i|P_k(\lambda_i)\varphi(\lambda_i)|} {H_kH_{k+1}}\sum_iw_{k,i}\Delta_k(\lambda_i)^2\le C(H_{k+1}-H_k). \end{aligned} \tag{G.24}\label{eq:common-moment-bound} \end{equation*} Uniform inversion of the $w_{k+2}$ moment Gram matrix converts the $s$ bounds in \eqref{eq:common-moment-bound} into \eqref{eq:common-factor-step}. Since \begin{equation*} q_{k+1}-q_k=P_{k+1}(P_{k+2}-P_k), \tag{G.25}\label{eq:common-product-increment} \end{equation*} the same coefficient bound gives $\|q_{k+1}-q_k\|_{\rm coeff}\le C(H_{k+1}-H_k)$. The summability in \eqref{eq:common-sums} makes both parity factor sequences and the full product sequence Cauchy, giving \eqref{eq:common-factor-limits}. Summing \eqref{eq:common-product-increment} from $k$ to infinity gives \eqref{eq:common-product-tail}.

It remains to identify the accumulation points and to control the nonresonant coordinates. Let $y_{k_j}\to y_*$ along one parity. Inequality \eqref{eq:common-face-separation} keeps the limit in the regular part of the simplex, where the solution of the moment equations and the block map are continuous, so $F$ is defined and continuous at $y_*$. Therefore $y_{k_j+1}=F(y_{k_j})\to F(y_*)$. The latter is an opposite-parity $\omega$-limit state and hence, again by \eqref{eq:common-face-separation}, is regular with more than $s$ active nodes. Thus $F$ is defined and continuous at $F(y_*)$, so $F^2(y_*)$ is defined and $y_{k_j+2}\to F^2(y_*)$; by \eqref{eq:common-sums} the same sequence tends to $y_*$. Thus $F^2(y_*)=y_*$. Cancelling every nonzero coordinate in \eqref{eq:common-two-block} and using \eqref{eq:common-factor-limits} gives $$ q_\infty(\lambda_i)=h $$ on the support of $y_*$. The lower support bound is \eqref{eq:common-face-separation}. The upper bound is $2s$, because the support is contained in the zero set of the nonzero monic degree-$2s$ polynomial $q_\infty-h$, which gives \eqref{eq:common-resonant-set}.

For the nonresonant coordinates, \eqref{eq:common-two-block} gives the exact identity \begin{equation*} \|y_{k+2}-y_k\|^2 =\sum_iw_{k,i}\left( \frac{q_k(\lambda_i)}{\sqrt{H_kH_{k+1}}}-1\right)^2. \tag{G.26}\label{eq:common-coordinate-chord} \end{equation*} If $\lambda_i\notin E$, the scalar in parentheses converges to the nonzero number $q_\infty(\lambda_i)/h-1$. Its absolute value is therefore bounded below for all sufficiently large $k$. Summing \eqref{eq:common-coordinate-chord} and applying \eqref{eq:common-sums} proves \eqref{eq:common-external-summable}. Equation \eqref{eq:common-completion-tail} follows from \eqref{eq:common-product-tail} and $0\le h-\sqrt{H_kH_{k+1}}\le h-H_k$. \end{proof}

Thus the support of every accumulation point is contained in the resonant set $E$. Convergence is reduced to motion within the corresponding fixed set of the two-block map.

\section*{Part A. Convergence at restart length two} \gdef\proofpart{A} \renewcommand{\thesection}{A.\arabic{section}} \setcounter{section}{0} \setcounter{theorem}{0}

Part A specialises the preceding argument to restart length two, retaining the quadratic notation needed for the endpoint analysis; see also \cite{ColbrookStepaniantsTownsend2026}. The support bound leaves only three- and four-node accumulation points. The former are rigid, whereas the latter lie on a compact one-parameter family of exact two-cycles, so convergence reduces to controlling motion along that family.

\section{The theorem for restart length two and the spectral map}\label{s2:sec:map}

\begin{theorem}[Forsythe's conjecture for restart length two]\label{s2:thm:main} Let $n\ge1$, let $A\in\mathbb R^{n\times n}$ be symmetric positive definite, and let $b,x_0\in\mathbb R^n$. Define $x_*=A^{-1}b$ and $r_k=b-Ax_k$. Whenever $r_k\ne0$, let $x_{k+1}$ be the unique minimiser of $$ \|x_*-z\|_A, \qquad z\in x_k+\mathcal K_2(A,r_k), \qquad \mathcal K_2(A,r)=\operatorname{span}\{r,Ar\}. $$ Then either $r_K=0$ for some finite $K$, or there are unit vectors $y_{\rm e},y_{\rm o}\in\mathbb R^n$ such that, with $y_k=r_k/\|r_k\|$, $$ y_{2k}\longrightarrow y_{\rm e}, \qquad y_{2k+1}\longrightarrow y_{\rm o} $$ in Euclidean norm. \end{theorem}

\begin{lemma}[Spectral reduction and one block at restart length two]\label{s2:lem:spectral} Suppose $r_0\ne0$. Let $0<\lambda_1<\cdots<\lambda_N$ be the distinct eigenvalues for which the spectral projection of $r_0$ is nonzero, and define $$ v_i=\frac{\operatorname{proj}_{\ker(A-\lambda_iI)}r_0} {\|\operatorname{proj}_{\ker(A-\lambda_iI)}r_0\|}. $$ Then the $v_i$ are orthonormal and every normalised nonzero residual has a unique representation $$ y=\sum_{i=1}^N\eta_i v_i, \qquad \sum_i\eta_i^2=1. $$ No other direction inside a repeated eigenspace can appear. The associated weights and inner product are $$ w_i=\eta_i^2, \qquad \langle f,g\rangle_w=\sum_iw_if(\lambda_i)g(\lambda_i). $$ If the support of $w$ has at least two nodes, there is a unique monic quadratic $$ P_w(t)=t^2+a(w)t+b(w) $$ orthogonal to $\mathcal P_1$. Define $H(w)=\|P_w\|_w^2$. Whenever $H(w)>0$, define \begin{equation*} (Tw)_i=\frac{w_iP_w(\lambda_i)^2}{H(w)}, \qquad F(y)_i=\frac{P_w(\lambda_i)\eta_i}{\sqrt{H(w)}}. \tag{A.1}\label{s2:eq:blockmap} \end{equation*} For a state on at least three nodes, $H(w)>0$, $P_w(0)>0$, and one block with restart length two satisfies \begin{equation*} y^+=F(y),\qquad w^+=Tw. \tag{A.2}\label{s2:eq:exactblock} \end{equation*} If the support has at most two nodes, the next block terminates. Conversely, a block with restart length two can terminate only on a support of at most two nodes. A zero coordinate remains zero forever. \end{lemma}

\begin{proof} We orthogonally diagonalise $A$ and translate by $x_*$. This reduces the iteration to the active spectral lines without changing the residual directions. Polynomial functional calculus keeps the component in each repeated eigenspace on the single line generated by its initial projection. The Vandermonde argument identifies the grade of a residual with the number of its active nodes.

For $P_w=t^2+at+b$, orthogonality to $1,t$ is the system \begin{equation*} \Gamma(w)\binom b a =-\binom{\sum_iw_i\lambda_i^2}{\sum_iw_i\lambda_i^3}, \qquad \Gamma(w)= \begin{pmatrix} 1&\sum_iw_i\lambda_i\\ \sum_iw_i\lambda_i&\sum_iw_i\lambda_i^2 \end{pmatrix}. \tag{A.3}\label{s2:eq:gram} \end{equation*} If at least two distinct nodes are active, then $$ (u,v)\Gamma(w)(u,v)^T=\sum_iw_i(u+v\lambda_i)^2>0 $$ for $(u,v)\ne(0,0)$, so the solution is unique. Orthogonality gives the least-norm identity \begin{equation*} \|P_w+\ell\|_w^2=H(w)+\|\ell\|_w^2 \qquad(\ell\in\mathcal P_1). \tag{A.4}\label{s2:eq:least} \end{equation*} This identity shows that $H(w)=0$ precisely when the monic quadratic vanishes on the entire support. This happens on two nodes and is impossible on three or more.

For at least three active nodes, the ordered values $P_w(\lambda_i)$ have at least two sign changes. Otherwise a linear polynomial could be chosen with the same sign as all nonzero numbers $w_iP_w(\lambda_i)$, contradicting their annihilation of $\mathcal P_1$. The two roots of $P_w$ are distinct and lie between positive spectral nodes. In particular $P_w(0)>0$.

The residual after one block is $\pi(A)r$, where $\pi\in\mathcal P_2$, $\pi(0)=1$, and the Galerkin equations say $\pi\perp\mathcal P_1$ for the current spectral measure. A polynomial of degree at most one satisfying those equations would be zero by positive definiteness of $\Gamma(w)$, contradicting $\pi(0)=1$. We therefore have $\pi=P_w/P_w(0)$. Since $P_w(0)>0$, normalisation of the residual gives \eqref{s2:eq:blockmap}--\eqref{s2:eq:exactblock}. The grade statement gives the termination assertions, while \eqref{s2:eq:blockmap} gives face invariance. \qedhere \end{proof}

For a nonterminating orbit, the relevant evolution is the two-block map, since the theorem concerns the two parity subsequences. If $r_0=0$, termination occurs at index zero. Henceforth suppose that the orbit does not terminate; every state then has at least three active nodes. Let $$ P_k=P_{w_k},\qquad H_k=H(w_k),\qquad q_k=P_kP_{k+1}. $$ The signed two-block recurrence is \begin{equation*} y_{k+2,i}=\frac{q_k(\lambda_i)}{\sqrt{H_kH_{k+1}}}\,y_{k,i}. \tag{A.5}\label{s2:eq:twoblock} \end{equation*}

\section{Energy, spectral localisation, and polynomial convergence}\label{s2:sec:global}

For $s=2$, Lemma~\ref{lem:common-energy} and Proposition~\ref{prop:common-factor-convergence} take the following quadratic form, involving the quartic and endpoint faces used throughout Part A. The block energy controls the motion on each parity, and the summability of its increments gives convergence of both orthogonal factors.

\begin{lemma}[Energy and same-parity displacement]\label{s2:lem:energy} For every $k$, \begin{equation*} H_{k+1}-H_k =\frac1{H_{k+1}}\sum_iw_{k,i} \bigl(q_k(\lambda_i)-H_{k+1}\bigr)^2, \tag{A.6}\label{s2:eq:variance} \end{equation*} and \begin{equation*} \sqrt{H_{k+1}}-\sqrt{H_k} =\frac{\sqrt{H_{k+1}}}{2}\|y_{k+2}-y_k\|^2. \tag{A.7}\label{s2:eq:chord} \end{equation*} Then $H_k\uparrow h$ for some $h>0$, \begin{equation*} \sum_{k=0}^\infty(H_{k+1}-H_k)<\infty, \qquad \sum_{k=0}^\infty\|y_{k+2}-y_k\|^2<\infty. \tag{A.8}\label{s2:eq:energytail} \end{equation*} \end{lemma}

\begin{proof} Since $P_{k+1}-P_k\in\mathcal P_1$, orthogonality gives $\langle P_k,P_{k+1}\rangle_{w_k}=H_k$, so $$ \sum_iw_{k,i}q_k(\lambda_i)=H_k. $$ Using the weight update once gives $$ \sum_iw_{k,i}q_k(\lambda_i)^2 =H_k\sum_iw_{k+1,i}P_{k+1}(\lambda_i)^2 =H_kH_{k+1}. $$ Expansion of the square proves \eqref{s2:eq:variance}. Moreover, \eqref{s2:eq:twoblock} gives $$ \langle y_{k+2},y_k\rangle=\sqrt{H_k/H_{k+1}}, $$ and the unit-vector chord identity proves \eqref{s2:eq:chord}. The monic competitor $t^2$ bounds $H_k\le\lambda_N^4$, while \eqref{s2:eq:variance} makes $H_k$ nondecreasing. It therefore converges to $h\ge H_0>0$, and the first sum in \eqref{s2:eq:energytail} telescopes. For the second, \eqref{s2:eq:chord} gives, for every $m$, $$ \sum_{k=0}^{m}\|y_{k+2}-y_k\|^2 \le \frac{2}{\sqrt{H_0}} \bigl(\sqrt{H_{m+1}}-\sqrt{H_0}\bigr). $$ Letting $m\to\infty$ gives the assertion. \qedhere \end{proof}

\begin{lemma}[Regular $\omega$-limits]\label{s2:lem:omega} The even and odd $\omega$-limit sets of the signed states, and of their squared weights, are nonempty, compact, and connected. Every $\omega$-limit state $y^*$ satisfies \begin{equation*} F^2(y^*)=y^*. \tag{A.9}\label{s2:eq:fixedtwo} \end{equation*} Let $w^*=((y_i^*)^2)_i$. Its support $S$ has three or four nodes. If $P=P_{w^*}$ and $Q=P_{Tw^*}$, then \begin{equation*} H(w^*)=H(Tw^*)=h, \qquad P(\lambda)Q(\lambda)=h\quad(\lambda\in S). \tag{A.10}\label{s2:eq:cycleproduct} \end{equation*} \end{lemma}

\begin{proof} A sequence in a compact metric space whose successive distances tend to zero has a connected $\omega$-limit set. Its distance from the $\omega$-limit set tends to zero, since otherwise a subsequence remaining a fixed positive distance away would have a cluster point outside that set. If the $\omega$-limit set were the disjoint union of two nonempty compact sets, choose disjoint neighbourhoods whose closures are still disjoint. Every sufficiently late term lies in their union, but sufficiently small successive steps prevent passage from one neighbourhood to the other. The sequence would eventually be confined to one neighbourhood, contrary to the presence of the other compact set in the $\omega$-limit set. By \eqref{s2:eq:energytail}, this observation applies to each parity subsequence. The $\omega$-limit set of the squared weights is precisely the coordinatewise-square image of the signed $\omega$-limit set. One inclusion follows from continuity; for the other, take a further convergent signed subsequence of any weight-convergent subsequence. The weight $\omega$-limit sets are therefore connected as well.

A limit weight cannot have support of size at most two. If it did, a fixed monic quadratic vanishing on that support would, by \eqref{s2:eq:least}, force a convergent subsequence of $H_k$ to tend to zero, contrary to $H_k\ge H_0$.

Let $y_{k_j}\to y^*$ along one parity. Moment inversion is continuous near a state with at least three positive coordinates, so $F$ and $H$ are continuous there. We have $y_{k_j+1}\to F(y^*)$, and the intermediate limit also has at least three active nodes. Applying continuity once more and $\|y_{k_j+2}-y_{k_j}\|\to0$ gives \eqref{s2:eq:fixedtwo}. Passing to the limit in the energies gives their common value $h$. Applying the block map twice and cancelling each active coordinate yields \eqref{s2:eq:cycleproduct}. The monic quartic $PQ-h$ vanishes on the support, so that support has at most four nodes. \qedhere \end{proof}

Define $$ \Delta_k=q_k-H_{k+1}, \qquad d_k=H_{k+1}-H_k. $$

\begin{lemma}[Double orthogonality and coefficient convergence]\label{s2:lem:factors} For every $\varphi\in\mathcal P_1$, \begin{equation*} \langle P_{k+2}-P_k,\varphi\rangle_{w_{k+2}} =-\frac{1}{H_kH_{k+1}} \langle P_k\varphi,\Delta_k^2\rangle_{w_k}. \tag{A.11}\label{s2:eq:doubleorth} \end{equation*} There are monic quadratics $P_{\rm e},P_{\rm o}$, a monic quartic $q_\infty=P_{\rm e}P_{\rm o}$, and a constant $C$ such that, for all large $k$, \begin{equation*} \|P_{k+2}-P_k\|_{\rm coeff}\le Cd_k, \qquad \|q_{k+1}-q_k\|_{\rm coeff}\le Cd_k, \tag{A.12}\label{s2:eq:coefstep} \end{equation*} \begin{equation*} P_{2n}\to P_{\rm e}, \qquad P_{2n+1}\to P_{\rm o}, \qquad q_k\to q_\infty, \tag{A.13}\label{s2:eq:factorlimits} \end{equation*} and \begin{equation*} \|q_k-q_\infty\|_{\rm coeff}\le C(h-H_k). \tag{A.14}\label{s2:eq:quartictail} \end{equation*} \end{lemma}

\begin{proof} Orthogonality of $P_{k+2}$ with respect to $w_{k+2}$, followed by the two-block weight recurrence, gives $$ \langle P_{k+2}-P_k,\varphi\rangle_{w_{k+2}} =-\frac1{H_kH_{k+1}} \langle P_k\varphi,(H_{k+1}+\Delta_k)^2\rangle_{w_k}. $$ The constant term vanishes because $P_k\perp\mathcal P_1$. The linear term vanishes as well, since $$ \langle P_k\varphi,\Delta_k\rangle_{w_k} =\sum_iw_{k,i}P_k(\lambda_i)^2P_{k+1}(\lambda_i)\varphi(\lambda_i) =H_k\langle P_{k+1},\varphi\rangle_{w_{k+1}}=0. $$ This gives \eqref{s2:eq:doubleorth}.

On the compact union of the two weight $\omega$-limit sets, every moment Gram matrix $\Gamma(w)$ is positive definite by Lemma~\ref{s2:lem:omega}. Its least eigenvalue therefore has a positive minimum. All sufficiently late $\Gamma(w_k)$ have uniformly bounded inverses, and \eqref{s2:eq:gram} gives a uniform coefficient bound for $P_k$. Let $D_k=P_{k+2}-P_k\in\mathcal P_1$. For $\varphi\in\{1,t\}$, equations \eqref{s2:eq:doubleorth} and \eqref{s2:eq:variance} give $$ |\langle D_k,\varphi\rangle_{w_{k+2}}| \le C d_k. $$ Uniform inversion of $\Gamma(w_{k+2})$ converts these two moment bounds into the first inequality in \eqref{s2:eq:coefstep}. The factorisation $$ q_{k+1}-q_k=P_{k+1}(P_{k+2}-P_k) $$ and the coefficient bound for $P_{k+1}$ give the second. Since $\sum_kd_k<\infty$, the even and odd factors and the full quartic are Cauchy. Summing the second inequality from $k$ to infinity gives \eqref{s2:eq:quartictail}. \qedhere \end{proof}

Define the resonant spectral set \begin{equation*} E=\{\lambda_i:q_\infty(\lambda_i)=h\}. \end{equation*} Every $\omega$-limit support is contained in $E$. Along a convergent subsequence, \eqref{s2:eq:cycleproduct} and \eqref{s2:eq:factorlimits} give $q_\infty(\lambda)=h$ at each active node. Since $q_\infty-h$ is a monic quartic, $|E|\le4$.

The support bound leaves two alternatives. A three-node accumulation point is isolated by the moment equations; four-node accumulation points instead fill a one-parameter family of exact two-cycles. Both conclusions follow from the interpolation identity below.

\section{Three- and four-node limiting geometry}\label{s2:sec:geometry}

For a finite set $S$ of distinct nodes, define $$ R_S(t)=\prod_{\xi\in S}(t-\xi). $$ Lagrange interpolation gives \begin{equation*} \sum_{\xi\in S}\frac{f(\xi)}{R_S'(\xi)}=[t^{|S|-1}]f \qquad(\deg f\le |S|-1). \tag{A.16}\label{s2:eq:lagrange} \end{equation*}

\begin{lemma}[The three-node case]\label{s2:lem:three} Let $w$ be a three-node $\omega$-limit weight, let $S=\operatorname{supp}w$, and let $P=P_w$. Then \begin{equation*} w_\xi P(\xi)=\frac{h}{R_S'(\xi)} \qquad(\xi\in S). \tag{A.17}\label{s2:eq:threeweights} \end{equation*} The pair $S,P$ determines $w$ uniquely. If every $\omega$-limit weight has three active nodes, both parity weight sequences converge. \end{lemma}

\begin{proof} The vector $(w_\xi P(\xi))_{\xi\in S}$ annihilates $1,t$. By \eqref{s2:eq:lagrange}, its one-dimensional nullspace consists of $\chi/R_S'(\xi)$. Multiplication by $P(\xi)$ and summation gives $$ h=H(w)=\chi\sum_{\xi\in S}\frac{P(\xi)}{R_S'(\xi)}=\chi, $$ because $P$ is monic quadratic. This is \eqref{s2:eq:threeweights}; its denominators and $P(\xi)$ are nonzero. There are only finitely many possible supports, and each support determines at most one point in a given parity $\omega$-limit set because the factor limit is fixed. A finite connected compact set is a singleton. \qedhere \end{proof}

The four-node alternative is not rigid: after normalisation, the annihilator of $\mathcal P_1$ retains one scalar degree of freedom. The following construction identifies that parameter and its admissible interval. Suppose that a four-node $\omega$-limit point exists. Then $E=\{e_1<e_2<e_3<e_4\}$, and \begin{equation*} R(t)=\prod_{i=1}^4(t-e_i), \qquad q_\infty=R+h, \qquad P_{\rm e}P_{\rm o}=R+h. \end{equation*} Let $\mathcal S_E$ be the probability simplex supported on $E$. For $w\in\mathcal S_E$ with at least three positive coordinates, let $P=P_w$ and $H=H(w)$. The annihilator of $1,t$ on four nodes has the form \begin{equation*} w_{e_i}P(e_i)=H\frac{e_i-\alpha(w)}{R'(e_i)}. \tag{A.19}\label{s2:eq:alphaweights} \end{equation*} The general annihilator is $(\zeta_1e_i+\zeta_0)/R'(e_i)$, and multiplication by $P(e_i)$ and \eqref{s2:eq:lagrange} give $\zeta_1=H$. If one weight vanishes, its node is $\alpha(w)$. The formula \begin{equation*} \alpha(w)= \sum_{i=1}^4e_i- \frac{\sum_{i=1}^4e_i^3w_{e_i}P_w(e_i)}{H(w)} \tag{A.20}\label{s2:eq:alphaextension} \end{equation*} agrees with the preceding parameter in the interior and extends it smoothly through the three-node endpoints, because $\sum e_i^3/R'(e_i)=1$ and $\sum e_i^4/R'(e_i)=\sum e_i$.

For every admissible value of $\alpha$, the endpoint sign relations are \begin{equation*} \begin{array}{c|c} \text{quantity}&\text{sign or endpoint behaviour}\\ \hline h&>0\\ P_{\rm e}(e_i)P_{\rm o}(e_i)&=h>0\\ (e_i-\alpha)/(R'(e_i)P_{\rm e}(e_i))&\ge0\\ (e_i-\alpha)/(R'(e_i)P_{\rm o}(e_i))&\ge0. \end{array} \end{equation*} At either endpoint exactly the row corresponding to $e_i=\alpha$ has zero numerator; every factor $P_{\rm e}(e_i),P_{\rm o}(e_i)$ remains nonzero. The signed amplitudes $\eta_i$ are the analytic face coordinates, while the moment polynomials, energies, and weight map are analytic even functions of them through $w_i=\eta_i^2$. Away from the zero face, signed-amplitude ratios are ordinary quotients; at the face, the corresponding analytic multiplier extensions are used.

For $\alpha\in\mathbb R$, define \begin{equation*} w^{\rm e}_{e_i}(\alpha)= \frac{h(e_i-\alpha)}{R'(e_i)P_{\rm e}(e_i)}, \qquad w^{\rm o}_{e_i}(\alpha)= \frac{h(e_i-\alpha)}{R'(e_i)P_{\rm o}(e_i)}. \tag{A.21}\label{s2:eq:familyweights} \end{equation*} All other coordinates are zero.

\begin{lemma}[The closed four-node two-cycle family]\label{s2:lem:family} The vectors in \eqref{s2:eq:familyweights} have total mass one for every $\alpha$. The values for which they are nonnegative form one nondegenerate compact interval $I=[e_j,e_{j+1}]$ for some $j\in\{1,2,3\}$. In its interior all four weights are positive, while exactly one weight vanishes at either endpoint. For every $\alpha\in I$, \begin{equation*} P_{w^{\rm e}(\alpha)}=P_{\rm e}, \quad H(w^{\rm e}(\alpha))=h, \quad Tw^{\rm e}(\alpha)=w^{\rm o}(\alpha), \quad Tw^{\rm o}(\alpha)=w^{\rm e}(\alpha). \tag{A.22}\label{s2:eq:familycycle} \end{equation*} Moreover, \begin{equation*} \ell_0:=\min_{\alpha\in I} \min\{|P_{\rm e}(\alpha)|,|P_{\rm o}(\alpha)|\}>0. \tag{A.23}\label{s2:eq:factoraway} \end{equation*} Every even $\omega$-limit point lies on the even segment and every odd $\omega$-limit point on the odd segment. \end{lemma}

\begin{proof} The partial-fraction expansion of $P_{\rm o}/R$, followed by comparison at infinity, gives $$ \sum_i\frac{P_{\rm o}(e_i)}{R'(e_i)}=0, \qquad \sum_i\frac{e_iP_{\rm o}(e_i)}{R'(e_i)}=1. $$ Since $h/P_{\rm e}(e_i)=P_{\rm o}(e_i)$, the even weights sum to one; the odd case is symmetric. Moreover, $$ w^{\rm o}_{e_i}(\alpha) =\frac{P_{\rm e}(e_i)^2}{h}w^{\rm e}_{e_i}(\alpha), $$ so both families have the same signs and zeros.

Each coordinate is a nonconstant affine function of $\alpha$. The four-node $\omega$-limit state supplies an interior point of the intersection of the four nonnegativity half-lines. Slopes of both signs occur because their sum is zero, so the intersection is compact and nondegenerate. A coordinate can vanish only at $\alpha=e_i$, so each endpoint is one node and the endpoints are consecutive.

For $\alpha$ in the interior, positivity gives $$ \operatorname{sign}P_{\rm e}(e_i) =\operatorname{sign}P_{\rm o}(e_i) =\operatorname{sign}\frac{e_i-\alpha}{R'(e_i)}. $$ The two sign changes occur in the gaps disjoint from $I$, so each limiting quadratic has one root in each of those gaps and no root in $I$. The node values are also nonzero because their products equal $h$, proving \eqref{s2:eq:factoraway}.

Equation \eqref{s2:eq:alphaweights} shows directly that $P_{\rm e}$ is orthogonal for $w^{\rm e}(\alpha)$. Multiplication by $P_{\rm e}(e_i)$ and \eqref{s2:eq:lagrange} give energy $h$. The weight update gives the first map identity in \eqref{s2:eq:familycycle}; the second is symmetric. Factor and energy convergence, together with \eqref{s2:eq:alphaweights}, place all $\omega$-limit points on the stated segments, including three-node endpoints. \qedhere \end{proof}

Define $$ \mathcal W_{\rm e}=\{w^{\rm e}(\alpha):\alpha\in I\}, \quad \mathcal W_{\rm o}=\{w^{\rm o}(\alpha):\alpha\in I\}, \quad \mathcal W_{\rm cyc}=\mathcal W_{\rm e}\cup\mathcal W_{\rm o}. $$ Then the parity sequences approach their respective segments.

The three-node alternative is rigid. In the four-node case, convergence means that the orbit approaches a single point of the limiting segment. The defect coordinate below measures transverse motion and the resulting tangential drift.

\section{Dynamics on the four-node family}\label{s2:sec:scalar}

Consider an exact orbit in $\mathcal S_E$, with at least three positive weights at every step. Let $$ \alpha_k=\alpha(w_k), \qquad K_\alpha(t)=\frac{R(t)-R(\alpha)}{t-\alpha}, \qquad c(w)=[t^3](P_wP_{Tw}-R),\qquad c_k=c(w_k). $$

\begin{lemma}[Defect and drift identities]\label{s2:lem:defect} Along such an orbit, \begin{equation*} q_k=R+H_{k+1}+c_kK_{\alpha_k}, \qquad \alpha_{k+1}-\alpha_k= \frac{c_kR(\alpha_k)}{H_{k+1}}, \tag{A.24}\label{s2:eq:defectdrift} \end{equation*} and \begin{equation*} d_k=\frac{c_k^2}{H_{k+1}} \sum_{i=1}^4w_{k,e_i}K_{\alpha_k}(e_i)^2. \tag{A.25}\label{s2:eq:defectenergy} \end{equation*} \end{lemma}

\begin{proof} Applying \eqref{s2:eq:alphaweights} at times $k$ and $k+1$, together with the weight update, gives $$ (e_i-\alpha_k)q_k(e_i) =H_{k+1}(e_i-\alpha_{k+1}) $$ at every active node. If a node of $E$ is absent, face invariance keeps it absent and both parameters equal that node, so the same identity holds there as well. The polynomial $$ (t-\alpha_k)q_k-H_{k+1}(t-\alpha_{k+1}) $$ vanishes at all four roots of $R$. Its quotient by $R$ is the monic linear polynomial $t-\alpha_k+c_k$. Hence \begin{equation*} (t-\alpha_k)(q_k-R) =c_kR+H_{k+1}(t-\alpha_{k+1}). \tag{A.26}\label{s2:eq:keypoly} \end{equation*} Evaluation at $t=\alpha_k$ yields the drift formula. Substituting $R(t)-R(\alpha_k)=(t-\alpha_k)K_{\alpha_k}(t)$ into \eqref{s2:eq:keypoly} and cancelling $t-\alpha_k$ gives the first formula in \eqref{s2:eq:defectdrift}. At the nodes of $E$, $q_k-H_{k+1}=c_kK_{\alpha_k}$; insertion into \eqref{s2:eq:variance} gives \eqref{s2:eq:defectenergy}. \qedhere \end{proof}

The drift identity alone does not control the sign of the defect $c_k$. A divided difference of $K_\alpha$ supplies the multiplicative transport law needed for this purpose. For a monic quadratic $P(t)=t^2+p_1t+p_0$, define \begin{equation*} L(\alpha,P)= p_1^2-p_0-(R_3+\alpha)p_1+R_2+R_3\alpha+\alpha^2, \tag{A.27}\label{s2:eq:L} \end{equation*} where $R=t^4+R_3t^3+R_2t^2+R_1t+R_0$. If $z_1,z_2$ are distinct roots of $P$, this equals $$ \frac{K_\alpha(z_1)-K_\alpha(z_2)}{z_1-z_2}. $$ The coefficient formula defines it without choosing roots and remains valid for a double root.

\begin{lemma}[Sign transport]\label{s2:lem:transport} For every exact orbit supported on $E$, \begin{equation*} c_{k+1}L(\alpha_{k+1},P_{k+1}) =c_kL(\alpha_k,P_{k+1}). \tag{A.28}\label{s2:eq:transport} \end{equation*} On the limiting family, \begin{equation*} L(\alpha,P_{\rm o})=P_{\rm e}(\alpha), \qquad L(\alpha,P_{\rm e})=P_{\rm o}(\alpha). \tag{A.29}\label{s2:eq:Lfamily} \end{equation*} These identities give a relative neighbourhood $\mathcal U_E$ of $\mathcal W_{\rm cyc}$ in $\mathcal S_E$ on which \begin{equation*} c(Tu)=m(u)c(u), \qquad m(u)=\frac{L(\alpha(u),P_{Tu})} {L(\alpha(Tu),P_{Tu})}>0, \end{equation*} and $m$, $m^{-1}$, and their first derivatives are uniformly bounded. \end{lemma}

\begin{proof} Subtract the first identity in \eqref{s2:eq:defectdrift} at two consecutive times: $$ c_{k+1}K_{\alpha_{k+1}}-c_kK_{\alpha_k} +(H_{k+2}-H_{k+1}) =P_{k+1}(P_{k+2}-P_k). $$ Evaluate at the two roots of $P_{k+1}$, subtract, and divide by the root difference. This gives \eqref{s2:eq:transport}; the coefficient identity extends it through a double root.

Let $P_{\rm e}=t^2+a_{\rm e}t+b_{\rm e}$, with analogous notation for the odd factor. Comparing the cubic and quadratic coefficients in $P_{\rm e}P_{\rm o}=R+h$ and substituting in \eqref{s2:eq:L} gives \eqref{s2:eq:Lfamily}. The lower bound \eqref{s2:eq:factoraway}, continuity, and compactness now provide the stated neighbourhood and positivity. \qedhere \end{proof}

These identities control the exact map on the resonant face. Small external coordinates perturb the dynamics on this face. The estimates below control their influence uniformly at the three-node endpoints while preserving the drift orientation.

\section{External modes and convergence of the weights}\label{s2:sec:convergence}

For an external node $\lambda\notin E$, define $$ \rho_\lambda=\frac{q_\infty(\lambda)}h. $$ Every external weight tends to zero, since a positive subsequential lower bound would produce an $\omega$-limit support not contained in $E$. If the coordinate is never annihilated, then \begin{equation*} \frac{w_{k+2,\lambda}}{w_{k,\lambda}} =\frac{q_k(\lambda)^2}{H_kH_{k+1}} \longrightarrow \rho_\lambda^2. \tag{A.31}\label{s2:eq:externalratio} \end{equation*} The limit excludes $|\rho_\lambda|>1$, and $\rho_\lambda=1$ would force $\lambda\in E$. A surviving external mode is therefore contractive, unless $\rho_\lambda=-1$.

\begin{lemma}[Exclusion of the neutral external multiplier]\label{s2:lem:neutral} Assume a four-node $\omega$-limit point exists. Every external coordinate with $\rho_\lambda=-1$ is annihilated after finitely many blocks. \end{lemma}

\begin{proof} Here $R(\lambda)=-2h$. Define $$ \Lambda_i(\lambda)= \frac{R(\lambda)}{(\lambda-e_i)R'(e_i)}. $$ These are the degree-three Lagrange cardinal values at $\lambda$: $\sum_i\Lambda_i=1$ and $f(\lambda)=\sum_i\Lambda_i f(e_i)$ for $f\in\mathcal P_3$.

After interchanging parity if necessary, choose an interior four-node even $\omega$-limit point. All four weights on $E$ have then been positive from the start, because zero weights are absorbing. Whenever the external coordinate is also positive, define \begin{equation*} J_{\lambda,n}= \log w_{2n,\lambda}+ \sum_{i=1}^4\Lambda_i(\lambda)\log w_{2n,e_i}. \tag{A.32}\label{s2:eq:J} \end{equation*} Let $$ H_n^{\rm e}=H_{2n},\quad H_n^{\rm o}=H_{2n+1},\quad \varepsilon_n=q_{2n}-(R+H_n^{\rm o}). $$ Then $\varepsilon_n\in\mathcal P_3$, and the exact two-block weight recurrence, followed by Lagrange interpolation, gives \begin{equation*} J_{\lambda,n+1}-J_{\lambda,n} =2\log\frac{2h-H_n^{\rm o}}{H_n^{\rm e}}+2\mathfrak r_n, \tag{A.33}\label{s2:eq:Jincrement} \end{equation*} where $$ \mathfrak r_n= \log\left(1- \frac{\varepsilon_n(\lambda)}{2h-H_n^{\rm o}}\right) +\sum_i\Lambda_i(\lambda) \log\left(1+ \frac{\varepsilon_n(e_i)}{H_n^{\rm o}}\right). $$ All denominators are uniformly separated from zero for all sufficiently large $n$. The linear part of $\mathfrak r_n$ is $$ \frac{2(h-H_n^{\rm o})} {H_n^{\rm o}(2h-H_n^{\rm o})}\,\varepsilon_n(\lambda), $$ because interpolation cancels the other first-order terms. Taylor's theorem and \eqref{s2:eq:quartictail} give $$ |\mathfrak r_n|\le C(h-H_n^{\rm e})^2. $$ On the other hand, $$ 2\log\frac{2h-H_n^{\rm o}}{H_n^{\rm e}} \ge 2\log\frac h{H_n^{\rm e}} \ge \frac{2(h-H_n^{\rm e})}{h}. $$ Let $\delta_n=h-H_n^{\rm e}$. The increment in \eqref{s2:eq:Jincrement} is at least $2\delta_n/h-2C\delta_n^2$. Since $\delta_n\to0$, it is nonnegative whenever $n$ is sufficiently large that $C\delta_n\le1/h$.

If the external coordinate were never annihilated, it would tend to zero by \eqref{s2:eq:externalratio} and the preceding support argument. Along a subsequence converging to the interior four-node $\omega$-limit point, all four weights on $E$ stay bounded below, so \eqref{s2:eq:J} tends to $-\infty$. This contradicts eventual monotonicity of $J_{\lambda,n}$. \qedhere \end{proof}

The total external mass $$ U(w)=\sum_{\lambda_i\notin E}w_i $$ therefore decays geometrically: there are $C<\infty$ and $0<\theta<1$ such that \begin{equation*} U(w_k)\le C\theta^k. \tag{A.34}\label{s2:eq:externaldecay} \end{equation*} Every external coordinate that survives forever has $|\rho_\lambda|<1$, and \eqref{s2:eq:externalratio} gives a uniform geometric estimate after taking the maximum over the finitely many nodes. Coordinates deleted at a finite time remain zero and already satisfy this estimate.

The parity orbit is an exponentially small perturbation of its retraction to the resonant face. The orientation lemma below shows that a transverse defect transported without reversal causes only finite motion along the compact cycle family.

\begin{lemma}[Orientation under summable perturbations]\label{s2:lem:abstract} Let $\alpha_n$ lie in a compact interval and suppose \begin{equation*} c_{n+1}=M_nc_n+\xi_n, \qquad \alpha_{n+1}-\alpha_n=\mathcal A_nc_n+\eta_n, \tag{A.35}\label{s2:eq:abstractrec} \end{equation*} where, for some $C<\infty$ and $0<\theta<1$, $$ M_n>0, \quad M_n\to1, \quad |\xi_n|+|\eta_n|\le C\theta^n, \quad |\mathcal A_n|\le C, $$ and all nonzero $\mathcal A_n$ have one common sign. Then $\sum_n|\alpha_{n+1}-\alpha_n|<\infty$. \end{lemma}

\begin{proof} If $c_n$ eventually has one sign, then $\mathcal A_nc_n$ also has one sign. Summing the second relation in \eqref{s2:eq:abstractrec} shows that the monotone partial sums of $\mathcal A_nc_n$ are bounded. They are therefore absolutely summable, as are the $\alpha$-increments.

Otherwise there are arbitrarily late crossings $m$ with $c_mc_{m+1}\le0$. Choose $m_0\in(\theta,1)$ so that $M_n\ge m_0$ late in the sequence. At a crossing, $$ |c_m|+|c_{m+1}|\le C\theta^m. $$ Solving the first recurrence backward from such a crossing gives, for $n<m$, $$ |c_n| \le Cm_0^{-(m-n)}\theta^m +C\sum_{j=n}^{m-1}m_0^{-(j-n+1)}\theta^j \le C'\theta^n. $$ Letting $m\to\infty$ through crossings proves $\sum_n|c_n|<\infty$, and the second recurrence again gives finite total variation. \qedhere \end{proof}

To apply Lemma~\ref{s2:lem:abstract} to the original orbit, compare each two-block step with its retraction onto the resonant face. Let $\mathcal G=T^2$.

\begin{lemma}[Endpoint-uniform regularity and retraction]\label{s2:lem:retraction} There are a relative neighbourhood $\mathcal N$ of $\mathcal W_{\rm cyc}$ in the probability simplex and $\nu>0$ such that, for every $w\in\mathcal N$, $$ \det\Gamma(w),\ H(w),\ \det\Gamma(Tw),\ H(Tw)\ge\nu. $$ The maps $w\mapsto P_w,H,T,\mathcal G$ and the coordinates $\alpha$ from \eqref{s2:eq:alphaextension} and $c$ defined above extend as $C^1$ maps with uniformly bounded first derivatives to an open neighbourhood of $\overline{\mathcal N}$ in the affine hyperplane $\sum_iw_i=1$. For $U(w)<1$, define $$ (\operatorname{ret}_Ew)_i= \begin{cases} w_i/(1-U(w)),&\lambda_i\in E,\\ 0,&\lambda_i\notin E. \end{cases} $$ After shrinking $\mathcal N$, \begin{equation*} \|w-\operatorname{ret}_Ew\|_1=2U(w), \qquad \|\operatorname{ret}_E\mathcal G(w) -\mathcal G(\operatorname{ret}_Ew)\|_1 \le CU(w). \tag{A.36}\label{s2:eq:retraction} \end{equation*} Moreover, every $u\in\mathcal N\cap\mathcal S_E$ satisfies \begin{equation*} \alpha(u)\in I, \qquad \alpha(Tu)\in I. \tag{A.37}\label{s2:eq:alphagap} \end{equation*} \end{lemma}

\begin{proof} Every point of $\mathcal W_{\rm cyc}$ has at least three positive coordinates on $E$. Its moment Gram determinant and energy are positive; the same is true after applying $T$, because the cycle family is $T$-invariant. Continuity and compactness give a relative neighbourhood $\mathcal V$ of the family and a common positive lower bound for the four displayed quantities. The formulas obtained by solving \eqref{s2:eq:gram}, followed by \eqref{s2:eq:blockmap}, are rational functions whose only denominators are those determinants and energies. It follows that $P,H,T$, and $\mathcal G$ extend as $C^1$ maps on $\mathcal V$. Formula \eqref{s2:eq:alphaextension} supplies the endpoint-regular extension of $\alpha$, while the coefficient definition of $c$ supplies its extension. Choose $\mathcal N$ with compact closure in $\mathcal V$, with $\mathcal N\cap\mathcal S_E\subset\mathcal U_E$, and small enough that $T(\overline{\mathcal N})\cup\mathcal G(\overline{\mathcal N})\subset \mathcal V$. Compactness gives the uniform derivative bounds.

Since every cycle state is supported on $E$, the total mass on $E$ of both $w$ and $\mathcal G(w)$ is bounded below on a smaller neighbourhood. The first equality in \eqref{s2:eq:retraction} follows directly: the total change on $E$ equals $U(w)$, and the total change off $E$ also equals $U(w)$. Let $v=\operatorname{ret}_Ew$. The segment from $w$ to $v$, and its image under $\mathcal G$, remain in a fixed compact regular neighbourhood after one further shrink. The map $\Psi=\operatorname{ret}_E\circ\mathcal G$ has bounded derivative there, and $\mathcal G(v)\in\mathcal S_E$. We therefore have $$ \|\operatorname{ret}_E\mathcal G(w)-\mathcal G(v)\|_1 =\|\Psi(w)-\Psi(v)\|_1 \le C\|w-v\|_1=2CU(w), $$ which proves the second estimate.

The node values of $P_{\rm e}$ and $P_{\rm o}$ are nonzero and have the common sign pattern that defines $I$. By continuity, the values of $P_u$ and $P_{Tu}$ retain this pattern on a smaller neighbourhood. For $u\in\mathcal S_E$, the nonnegativity of the weights in \eqref{s2:eq:alphaweights} then forces both $\alpha(u)$ and $\alpha(Tu)$ into the same closed interval $I$. \qedhere \end{proof}

\begin{lemma}[Convergence in the four-node case]\label{s2:lem:fourconv} If a four-node $\omega$-limit point exists, both parity weight sequences converge. \end{lemma}

\begin{proof} For the even subsequence, let $$ \widehat w_n=\operatorname{ret}_Ew_{2n}, \qquad \alpha_n=\alpha(\widehat w_n), \qquad c_n=c(\widehat w_n). $$ The approach to $\mathcal W_{\rm e}$, \eqref{s2:eq:externaldecay}, and \eqref{s2:eq:retraction} give \begin{equation*} \|\widehat w_{n+1}-\mathcal G(\widehat w_n)\|_1 \le C\theta^n. \tag{A.38}\label{s2:eq:pseudoorbit} \end{equation*} For the exact map restricted to $E$, Lemma~\ref{s2:lem:transport} gives $$ c(\mathcal G(u))=M(u)c(u), \qquad M(u)=m(Tu)m(u)>0. $$ Applying the drift formula in Lemma~\ref{s2:lem:defect} twice gives $$ \alpha(\mathcal G(u))-\alpha(u)=\mathcal A(u)c(u), $$ where $$ \mathcal A(u)= \frac{R(\alpha(u))}{H(Tu)} +m(u)\frac{R(\alpha(Tu))}{H(\mathcal G(u))}. $$ By \eqref{s2:eq:alphagap}, all nonzero values of $\mathcal A$ have the fixed sign of $R$ on the interior of $I$, and $\mathcal A$ is bounded.

The map $w\mapsto P_wP_{Tw}$ is $C^1$. Since $\widehat w_n-w_{2n}\to0$, we have $P_{\widehat w_n}P_{T\widehat w_n}\to q_\infty$, and hence $c_n\to0$. On the exact family $m=1$. The drift identity and the uniform bounds on $L$ imply $$ |m(u)-1|\le C|c(u)|, \qquad |m(Tu)-1|\le C|c(u)|, $$ and therefore $M(\widehat w_n)\to1$. For each $n$, let $$ M_n=M(\widehat w_n),\qquad \mathcal A_n=\mathcal A(\widehat w_n), $$ and $$ \xi_n=c(\widehat w_{n+1})-c(\mathcal G(\widehat w_n)),\qquad \eta_n=\alpha(\widehat w_{n+1})-\alpha(\mathcal G(\widehat w_n)). $$ The uniform $C^1$ bounds in Lemma~\ref{s2:lem:retraction} and \eqref{s2:eq:pseudoorbit} give $|\xi_n|+|\eta_n|\le C\theta^n$. For all sufficiently large $n$, $M_n>0$, $M_n\to1$, $\mathcal A_n$ is bounded, every nonzero $\mathcal A_n$ has the fixed sign of $R$ on $I$, and $\alpha_n\in I$. These quantities satisfy \eqref{s2:eq:abstractrec}, so Lemma~\ref{s2:lem:abstract} gives $\alpha_n\to\alpha_\infty$.

Equation \eqref{s2:eq:alphaweights}, together with factor and energy convergence, shows that every weight of $\widehat w_n$ on $E$ converges. The denominators $P_{\rm e}(e_i)$ are nonzero because $P_{\rm e}(e_i)P_{\rm o}(e_i)=h$. Since the external mass tends to zero, $w_{2n}$ converges. Its limit has at least three active nodes, so continuity of $T$ gives convergence of $w_{2n+1}=Tw_{2n}$. \qedhere \end{proof}

\section{Signed recovery and proof of the theorem for restart length two}\label{s2:sec:signed}

The preceding sections give convergence of the squared spectral coordinates. It remains to recover their signs and return to the original eigenbasis.

\begin{proof}[Proof of Theorem~\ref{s2:thm:main}] If no four-node $\omega$-limit point exists, Lemma~\ref{s2:lem:three} gives convergence of both parity weight sequences. Otherwise Lemma~\ref{s2:lem:fourconv} does. Fix one parity and let $k$ range over it. If the limiting weight at node $\lambda_i$ is positive, then $\lambda_i\in E$ and $q_\infty(\lambda_i)=h$. The signed recurrence \eqref{s2:eq:twoblock} gives $$ \frac{y_{k+2,i}}{y_{k,i}} =\frac{q_k(\lambda_i)}{\sqrt{H_kH_{k+1}}} \longrightarrow1. $$ The sign of that coordinate is eventually constant, while its magnitude converges to the square root of the limiting weight. Coordinates whose limiting weights vanish tend to zero regardless of sign. Since there are finitely many coordinates, each parity of the signed spectral vector converges in Euclidean norm.

The spectral vectors $v_i$ form a fixed orthonormal family, so convergence of their coefficient vectors is exactly convergence of the original normalised residual directions. Lemma~\ref{s2:lem:spectral} covers repeated eigenvalues and initially absent spectral components; later support loss is permanent and hence does not affect the argument. \qedhere \end{proof}

\section*{Part B. Convergence at restart length three} \gdef\proofpart{B} \renewcommand{\thesection}{B.\arabic{section}} \setcounter{section}{0} \setcounter{theorem}{0}

The common argument shows that the two parity cubic sequences converge and that every limiting support contains four, five, or six spectral nodes. For restart length three this does not by itself imply convergence of the weights. The limiting cubics determine a monic sextic. Every limiting spectral node is one of its roots; any further roots are completion roots rather than spectral nodes. The external coordinates are classified by their limiting linear multipliers. Stable coordinates are handled by a shadowing comparison while remaining in the exact state, unstable coordinates disappear after finitely many returns, and neutral coordinates are governed by logarithmic cocycles on the five- and six-node boundary strata. The collision coordinates, support-loss strata, and endpoint normal forms needed for this analysis are developed below and in Appendix B.

The notation used after the limiting completion has been selected is summarised here. A tilde marks data of the exact full state, while untilded data belong to the normalised limiting core. A plus on an untilded coordinate denotes the normalised core induced by the exact return; setting the external masses to zero gives the restricted-core return.
\begin{center}
\begingroup\small
\begin{tabular}{@{}>{\raggedright\arraybackslash}p{0.19\textwidth}>{\raggedright\arraybackslash}p{0.38\textwidth}>{\raggedright\arraybackslash}p{0.35\textwidth}@{}}
 & exact full state & normalised core\\ \hline
squared weights & $Y_{n,i}^2$; retained $W_{n,i}$; external $E_{h,n}$ & $w_{n,i}=W_{n,i}/\sum_{j\in S_*}W_{n,j}$\\
monic factors & $\widetilde P_n,\widetilde Q_n$ & $P_n,Q_n$\\
one-block norms & $\widetilde a_n,\widetilde b_n$ & $a_n,b_n$\\
two-block scale & $\widetilde C_n=\widetilde a_n\widetilde b_n$ & $C_n=a_nb_n$\\
completion and normal & $\widetilde{\mathscr D}_n,\widetilde U_n$ & $\mathscr D_n,U_n$
\end{tabular}
\endgroup
\end{center}
Starred symbols are fixed limiting data rather than orbit quantities. In particular, $a_n,b_n$ are one-block norms, whereas $C_n=a_nb_n$ is their two-block product. The global proof dependency is
\[
 \begin{gathered}
 \text{\ref{sec:2}--\ref{sec:3}}\longrightarrow\text{\ref{sec:4}}\longrightarrow\text{\ref{sec:5}--\ref{sec:6}}
 \longrightarrow\text{\ref{sec:7}/\ref{sec:8}/\ref{sec:9}}\longrightarrow\text{\ref{sec:10}},\\
 \text{Appendix B}\longrightarrow\text{local calculations for Sections~\ref{sec:8}--\ref{sec:9}}.
 \end{gathered}
\]
These stages give, respectively, the spectral dynamics, fixed-set geometry, limiting completion and external-mode reduction, the support cases, and signed recovery. Within Sections~\ref{sec:8}--\ref{sec:9} the local order is the exact full return, the induced normalised-core return, the analytic extension of a specified log-return, and only then its stable-weight comparison. No estimate for a logarithm is invoked before the logarithm of the corresponding positive return ratio has been shown to extend real-analytically to the relevant signed boundary charts.

\section{The theorem for restart length three}\label{sec:1}

\begin{theorem}[Forsythe's conjecture for restart length three]\label{thm:main}
Let $n\ge 1$, let $A\in\mathbb R^{n\times n}$ be symmetric positive definite, and let $b,x_0\in\mathbb R^n$. Define $x_*=A^{-1}b$ and $r_k=b-Ax_k$. At every restart step for which $r_k\ne0$, let $x_{k+1}$ be the unique minimiser of
$$ \|x_*-z\|_A,\qquad z\in x_k+\mathcal K_3(A,r_k), \qquad \mathcal K_3(A,r)=\operatorname{span}\{r,Ar,A^2r\}. $$
Then either $r_K=0$ for some finite $K$, or there are unit vectors $y_{\rm e},y_{\rm o}\in\mathbb R^n$ such that, with $y_k=r_k/\|r_k\|$,
$$ y_{2k}\longrightarrow y_{\rm e}, \qquad y_{2k+1}\longrightarrow y_{\rm o} $$
in Euclidean norm.
\end{theorem}

\begin{figure}[H]
\centering
\includegraphics[width=0.72\textwidth]{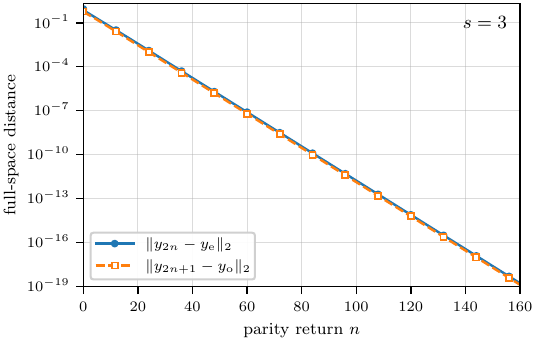}
\caption{Direct restarted-CG computation at restart length three. We take $A=\operatorname{diag}(1,\ldots,10)$ and $y_0=10^{-1/2}(1,\ldots,1)^T$. At each block, the degree-three residual polynomial is obtained by solving the Galerkin equations \eqref{eq:B.2.3} in 120-digit arithmetic. The curves are the Euclidean distances of the even and odd residual directions from separate reference parity limits, computed by continuing the same iteration to block 1600. Both distances tend to zero. The plot is not a phase portrait: the even and odd sequences converge to two fixed points of $F^2$, equivalently the one-block orbit approaches a two-cycle of $F$.}
\label{fig:s3-convergence}
\end{figure}

Figure~\ref{fig:s3-convergence} shows a representative direct computation of this two-cycle convergence.

\section{Spectral reduction and the signed residual map}\label{sec:2}

A restart-three block has the following spectral representation. The squared coordinates determine the monic orthogonal cubic, while the coordinate signs recover the residual direction.

For a polynomial $f$, we abbreviate $f_i=f(\lambda_i)$ whenever the spectral nodes are understood.

\begin{lemma}[Spectral reduction, exact block polynomial, and grade]\label{lem:2.1}
Suppose that $r_0\ne0$. Let $\lambda_1<\cdots<\lambda_N$ be the distinct eigenvalues for which the spectral projection of $r_0$ is nonzero, and let
$$ v_i=\frac{\operatorname{proj}_{\ker(A-\lambda_i I)}r_0} {\|\operatorname{proj}_{\ker(A-\lambda_i I)}r_0\|}. $$
Then the $v_i$ are orthonormal, every later residual is uniquely of the form
\begin{equation*} r_k=\rho_k\sum_{i=1}^N\eta_{k,i}v_i,\qquad \rho_k=\|r_k\|,\qquad \sum_i\eta_{k,i}^2=1, \end{equation*}
and no direction orthogonal to $v_i$ inside a repeated eigenspace can ever appear. In particular, a repeated eigenvalue is one spectral node and an eigenvalue with zero initial projection is absent forever.

For a current nonzero residual $r=\rho y$, define
\begin{equation*} y=\sum_i\eta_i v_i,\qquad w_i=\eta_i^2,\qquad S=\{i:w_i>0\}. \end{equation*}
If $|S|>3$, there is a unique $\pi_w\in\mathcal P_3$ satisfying
\begin{equation*} \pi_w(0)=1,\qquad \sum_iw_i\pi_w(\lambda_i)q(\lambda_i)=0 \quad(q\in\mathcal P_2), \tag{B.2.3}\label{eq:B.2.3} \end{equation*}
and the resulting normalised state after one exact restarted block is
\begin{equation*} F(y)= \frac{\sum_i\eta_i\pi_w(\lambda_i)v_i} {\left(\sum_iw_i\pi_w(\lambda_i)^2\right)^{1/2}}. \tag{B.2.4}\label{eq:B.2.4} \end{equation*}
The squared coordinates obey
\begin{equation*} (Tw)_i= \frac{w_i\pi_w(\lambda_i)^2} {\sum_jw_j\pi_w(\lambda_j)^2}. \tag{B.2.5}\label{eq:B.2.5} \end{equation*}
If instead $|S|\le3$, the next block terminates. Conversely, a block can terminate only when $|S|\le3$. A coordinate which becomes zero remains zero, so later support loss is permanent and there are only finitely many exact support losses.
\end{lemma}

\begin{proof}
Translation by $x_*$ changes the problem to $b=0$ without changing any residual. Orthogonal diagonalisation of $A$ preserves both Euclidean and $A$-norms. In a $\lambda_i$-eigenspace every polynomial in $A$ acts as a scalar; hence all Krylov corrections and residuals stay on the one oriented line spanned by the initial projection $v_i$. The vectors $r,Ar,\ldots,A^{m-1}r$, where $m=|S|$, are independent because their coordinate matrix is a nonzero diagonal matrix times the Vandermonde matrix $(\lambda_i^j)$. The grade of $r$ is exactly $|S|$.

Every block candidate is $x+q(A)r$, $q\in\mathcal P_2$, and its residual is
$$ r^+=(1-Aq(A))r=\pi(A)r,\qquad \pi(t)=1-tq(t). $$
The Galerkin condition $x_*-x^+\perp_A\mathcal K_3(A,r)$ is equivalent to
$$ (r^+)^TA^jr=0,\qquad j=0,1,2, $$
which is exactly \eqref{eq:B.2.3}. If $|S|>3$, the bilinear form
$$ B_w(p,q)=\sum_iw_i\lambda_i p(\lambda_i)q(\lambda_i), \qquad p,q\in\mathcal P_2, $$
is positive definite: a nonzero quadratic cannot vanish at all of more than three distinct active nodes. Solving $B_w(q,p)=\sum_iw_ip(\lambda_i)$ gives a unique $q$, hence a unique $\pi_w=1-tq$. Its residual norm is nonzero, since a nonzero cubic normalised at zero cannot vanish at more than three active nodes. Spectral functional calculus now gives \eqref{eq:B.2.4} and \eqref{eq:B.2.5}.

If $|S|\le3$, the polynomial
$$ \prod_{i\in S}\left(1-\frac{t}{\lambda_i}\right) $$
has degree at most three and annihilates $r$; the feasible block contains $x_*$, which is the unique minimiser. Conversely, termination would require a nonzero polynomial of degree at most three, equal to one at zero, to vanish on every active node; this is impossible when $|S|>3$. Finally \eqref{eq:B.2.5} proves coordinate-face invariance and permanent deletion.
\end{proof}

If $r_0=0$, the iteration terminates at restart index zero; henceforth assume $r_0\ne0$.

\begin{lemma}[Monic form and signs]\label{lem:2.2}
For a probability weight $w$ on at least four distinct positive nodes, let $P_w$ be the monic cubic orthogonal to $\mathcal P_2$, and define
\begin{equation*} \sigma(w)=\left(\sum_iw_iP_w(\lambda_i)^2\right)^{1/2}. \end{equation*}
For a unit spectral state $y=\sum_i\eta_iv_i$, let $w(y)=(\eta_i^2)_i$ and abbreviate $$ P_y=P_{w(y)},\qquad \sigma(y)=\sigma(w(y)). $$
Then $P_w$ has three simple positive roots, its values at the ordered active nodes have at least three sign changes, and $P_w(0)<0$. Defining
\begin{equation*} L(y)=\frac{P_w(A)y}{\sigma(w)} \end{equation*}
gives
\begin{equation*} F(y)=-L(y),\qquad F^2(y)=L^2(y). \tag{B.2.8}\label{eq:B.2.8} \end{equation*}
\end{lemma}

\begin{proof}
The coefficients $w_iP_w(\lambda_i)$ annihilate every quadratic. If their ordered signs had at most two changes, a polynomial of degree at most two with the same sign at every active node would make their scalar product strictly positive, a contradiction. Hence $P_w(\lambda_i)$ changes sign at least three times. The intermediate value theorem supplies three distinct roots between positive nodes; a cubic has no others. In particular, $P_w(0)<0$. Both $P_w/P_w(0)$ and $\pi_w$ satisfy \eqref{eq:B.2.3}, so they are equal. Normalising their residuals proves \eqref{eq:B.2.8}.
\end{proof}

Henceforth, $w_k=w(y_k)$ denotes the squared-coordinate weight at the $k$th block boundary. The parity return of the signed residual map is $L^2$. We therefore analyse the squared spectral weights first and recover the signs after proving convergence.

\section{Parity energy and \texorpdfstring{$\omega$}{omega}-limit sets}\label{sec:3}

We group two consecutive restart blocks into one parity return. This section restates the $s=3$ energy and localisation conclusions of the common core in the two-phase notation used throughout Part B. The phase-rigidity lemma at the end supplies the additional cubic-specific ingredient.

Define
\begin{equation*} Y_n=y_{2n},\qquad Z_n=L(Y_n)=-y_{2n+1},\qquad Y_{n+1}=L(Z_n). \tag{B.3.1}\label{eq:B.3.1} \end{equation*} We use $Y_n^2:=w(Y_n)=(Y_{n,i}^2)_i$ for the coordinatewise square.

Let
\begin{equation*} P_n=P_{Y_n},\quad Q_n=P_{Z_n},\quad a_n=\sigma(Y_n),\quad b_n=\sigma(Z_n). \end{equation*}

\begin{lemma}[Cross projection, interlacing, and finite quadratic energy]\label{lem:3.1}
For every nonterminating orbit,
\begin{equation*} 0<a_0\le a_n\le b_n\le a_{n+1}\le\|A^3\|, \tag{B.3.3}\label{eq:B.3.3} \end{equation*}
and $a_n,b_n$ converge to one number $\tau>0$. Moreover
\begin{equation*} \frac12\|Y_{n+1}-Y_n\|^2=1-\frac{a_n}{b_n}, \qquad \frac12\|Z_{n+1}-Z_n\|^2=1-\frac{b_n}{a_{n+1}}, \tag{B.3.4}\label{eq:B.3.4} \end{equation*}
and
\begin{equation*} \sum_{n=0}^{\infty} \bigl(\|Y_{n+1}-Y_n\|^2+\|Z_{n+1}-Z_n\|^2\bigr)<\infty. \tag{B.3.5}\label{eq:B.3.5} \end{equation*}
\end{lemma}

\begin{proof}
Let $v$ be a unit state, $z=L(v)$, $v^+=L(z)$, and set $z^+=L(v^+)$, $P=P_v$, $Q=P_z$, $P^+=P_{v^+}$, $a=\sigma(v)$, $b=\sigma(z)$, and $a^+=\sigma(v^+)$. Since two monic cubics differ by a quadratic and $P(A)v\perp\mathcal K_3(A,v)$,
$$ a^2=\langle P(A)v,Q(A)v\rangle. $$
Symmetry of $A$ and $P(A)v=az$, $Q(A)z=bv^+$ give
$$ a^2=ab\langle v,v^+\rangle. $$
We have $a=b\langle v,v^+\rangle\le b$, and the unit-vector chord identity gives the first equation in \eqref{eq:B.3.4}. At the next phase, $$ b^2=\langle Q(A)z,P^+(A)z\rangle =ba^+\langle z,z^+\rangle, $$ so $b=a^+\langle z,z^+\rangle\le a^+$, and the chord identity gives the second equation in \eqref{eq:B.3.4}. This gives the full interlacing \eqref{eq:B.3.3}. The monic competitor $t^3$ bounds every norm by $\|A^3\|$, so the two interlaced sequences have a common limit $\tau\ge a_0>0$.

Since every denominator is at least $a_0$,
$$ \left(1-\frac{a_n}{b_n}\right) +\left(1-\frac{b_n}{a_{n+1}}\right) \le\frac{a_{n+1}-a_n}{a_0}. $$
Summation telescopes and proves \eqref{eq:B.3.5}.
\end{proof}

\begin{lemma}[Uniform separation from supports of size at most three]\label{lem:3.2}
For every subset $J$ of the original active nodes with $|J|\le3$, define
\begin{equation*} R_J(t)=t^{3-|J|}\prod_{j\in J}(t-\lambda_j),\qquad M_J=\max_{i\notin J}R_J(\lambda_i)^2. \end{equation*}
Then $0<M_J<\infty$, and every block-boundary weight satisfies
\begin{equation*} \sum_{i\notin J}w_{k,i}\ge\frac{\sigma(w_0)^2}{M_J}>0. \tag{B.3.7}\label{eq:B.3.7} \end{equation*}
Every $\omega$-limit state has at least four positive coordinates, and the degree-three moment Gram matrix is uniformly nonsingular on the orbit closure.
\end{lemma}

\begin{proof}
The monic cubic $R_J$ is an admissible competitor for the least-norm monic cubic, so
$$ \sigma(w)^2\le\sum_iw_iR_J(\lambda_i)^2 \le M_J\sum_{i\notin J}w_i. $$
Lemma~\ref{lem:3.1}, applied at consecutive single blocks, makes the monic norm nondecreasing and at least $\sigma(w_0)>0$. This is \eqref{eq:B.3.7}. A limit on a face with at most three nodes would contradict it. Finally, failure of uniform Gram nonsingularity would produce at a limit a nonzero quadratic vanishing on at least four distinct positive-weight nodes, which is impossible.
\end{proof}

\begin{lemma}[Localisation to a two-cycle on four to six nodes]\label{lem:3.3}
Every $\omega$-limit state of either parity has between four and six positive coordinates. At every paired accumulation point $(Y,Z)$,
\begin{equation*} Z=L(Y),\qquad Y=L(Z), \tag{B.3.8}\label{eq:B.3.8} \end{equation*}
so the two phases have the same support and form an exact signed two-cycle. The $\omega$-limit set of either parity, and the paired $\omega$-limit set, is compact and connected.
\end{lemma}

\begin{proof}
The solution of the coefficient equations defining $P_w$, the norm, and $L$ is continuous at every state having at least four positive coordinates, by Lemma~\ref{lem:3.2}. If $Y_{n_j}\to Y$, then \eqref{eq:B.3.5} gives $Y_{n_j+1}\to Y$; continuity and the exact iteration give $L^2(Y)=Y$. Likewise, along paired convergence $(Y_{n_j},Z_{n_j})\to(Y,Z)$, continuity gives \eqref{eq:B.3.8}, and polynomial multiplication in both directions forces equality of supports.

Let $P,Q$ be the two monic cubics and let $a,b>0$ be their norms at such a point. Equality in the cross-projection chain gives $a=b$, and on every positive coordinate
\begin{equation*} P(\lambda_i)Q(\lambda_i)=ab. \end{equation*}
The nonzero monic sextic $P Q-ab$ has every support node as a root. It has degree six, so the support has at most six nodes; Lemma~\ref{lem:3.2} gives the lower bound four.

Finally, in a compact metric space the $\omega$-limit set of a sequence whose successive distances tend to zero is connected. Suppose it split into two nonempty compact pieces. Choose disjoint open neighbourhoods of those pieces whose closures are still disjoint. Every sufficiently late term lies in their union, since otherwise a subsequence outside the union would have an $\omega$-limit point there. Once successive distances are smaller than the positive distance between the two closures, the sequence cannot pass from one neighbourhood to the other. It is eventually confined to one, contradicting the presence of the other piece in the $\omega$-limit set. Apply this to \eqref{eq:B.3.5}, also in the product space.
\end{proof}

\begin{lemma}[Phase rigidity]\label{lem:3.4}
Let $w,w'$ be positive weights on the same finite support, with $w'=Tw$ and $w=Tw'$. If $p,q$ are their normalised-at-zero block polynomials and $D,D'>0$ are the corresponding squared residual normalisations, then
\begin{equation*} p(\lambda_i)q(\lambda_i)=\sqrt{DD'}>0 \tag{B.3.10}\label{eq:B.3.10} \end{equation*}
at every support node. Hence every signed lift of a regular weight two-cycle is a signed two-cycle, not merely a projective one.
\end{lemma}

\begin{proof}
The weight equations first give signs $\varepsilon_i\in\{\pm1\}$ with
$$ p(\lambda_i)q(\lambda_i)=\varepsilon_i\sqrt{DD'}. $$
Let $S_-=\{i:\varepsilon_i=-1\}$. The moment equations of the first and second block, after substitution of $w'_i=w_ip(\lambda_i)^2/D$, imply
\begin{equation*} \sum_{i\in S_-}w_ip(\lambda_i)s(\lambda_i)=0 \quad(s\in\mathcal P_2). \tag{B.3.11}\label{eq:B.3.11} \end{equation*}
If nonzero coefficients $z_i$ on ordered nodes annihilate $\mathcal P_2$, their signs have at least three changes: otherwise a quadratic with one root between each sign change can be chosen with the same sign as all $z_i$, contradicting annihilation.

Each cubic $p,q$ has exactly three simple positive roots. Their product is positive at zero and at positive infinity. Count the boundary incidences of the components of $\{t>0:p(t)q(t)<0\}$ with the roots of $pq$, using multiplicity. Every component has two boundary incidences. A common root of $p$ and $q$, although omitted from the strict set and possibly separating two negative components, is a double root of $pq$ and supplies the corresponding two incidences. With common roots counted twice, the total incidence count is at most six, so the strict negative set has at most three components. On each such component $p$ has constant sign. The ordered signs of $w_ip(\lambda_i)$, restricted to $S_-$, have at most three constant blocks and hence at most two changes, contradicting \eqref{eq:B.3.11}. Hence $S_-=\varnothing$, proving \eqref{eq:B.3.10}.
\end{proof}

Thus every accumulation pair belongs to the signed fixed set of $G=L^2$. To pass from this qualitative localisation to convergence, we need coordinates on that fixed set which remain regular at common-root collisions and support loss.

\section{Algebra of fixed two-cycles and collision coordinates}\label{sec:4}

Throughout Sections~\ref{sec:4}--\ref{sec:9}, a superscript $+$ denotes one two-block parity return $G=L^2$; intermediate one-block quantities carry a tilde or a phase subscript.

The fixed-point set of the parity return has the following description. A fixed pair of monic cubics $P,Q$ satisfies a factorisation in which the spectral nodes are completed, with multiplicity, to a monic sextic. The resulting coordinates remain regular when a node is lost or when $P$ and $Q$ acquire a common root, and will be used to convert decay of the parity energy into convergence of the completion polynomial. Throughout Part B, a \emph{fixed component} means a local analytic component of this fixed-point set. A \emph{simultaneous collision} means that two or more common-root collisions occur at the same critical value.

The four lemmas of this section have distinct roles. Lemma~\ref{lem:4.1} describes the fixed components and their support faces. Lemma~\ref{lem:4.2} identifies an affine normal and its first nonzero drift on the coprime six-node stratum, while Lemma~\ref{lem:4.3} extends this structure across common-root collisions and support losses. Lemma~\ref{lem:4.4} then compares completion drift with the squared parity chord.

For a finite set $J$ of distinct nodes define
\begin{equation*} D_J(t)=\prod_{\lambda_i\in J}(t-\lambda_i), \qquad c_i^J=\frac1{D_J'(\lambda_i)}. \end{equation*}
The Lagrange leading-coefficient identities used repeatedly below are
\begin{equation*} \begin{aligned} \sum_{i\in J}c_i^JR(\lambda_i)&=0 &&(\deg R\le |J|-2),\\ \sum_{i\in J}c_i^JR(\lambda_i)&=[t^{|J|-1}]R &&(\deg R\le |J|-1). \end{aligned} \tag{B.4.2}\label{eq:B.4.2} \end{equation*}
They follow by comparing the highest coefficient in the Lagrange interpolant of $R$.

\begin{lemma}[Completion polynomial and the six-node parameter rectangle]\label{lem:4.1}
Let $Y$ be a regular $L^2$-fixed state with support $J$, $4\le |J|\le6$. Let $P,Q$ be its consecutive monic cubics and let $C>0$ be their common squared monic norm. There are unique monic polynomials $K_J,\Lambda_J$, of degrees
\begin{equation*} \deg K_J=|J|-4,\qquad \deg\Lambda_J=6-|J|, \end{equation*}
such that, for the two phase weights,
\begin{equation*} w_i^P=c_i^JK_J(\lambda_i)Q(\lambda_i),\qquad w_i^Q=c_i^JK_J(\lambda_i)P(\lambda_i), \tag{B.4.4}\label{eq:B.4.4} \end{equation*}
and
\begin{equation*} P Q-C=D_J\Lambda_J. \tag{B.4.5}\label{eq:B.4.5} \end{equation*}
The polynomial $\Omega_J=K_J\Lambda_J$ is always monic quadratic. If a node $\alpha$ is lost in passing from $J\cup\{\alpha\}$ to $J$, then
\begin{equation*} K_{J\cup\{\alpha\}}=(t-\alpha)K_J,\qquad \Lambda_J=(t-\alpha)\Lambda_{J\cup\{\alpha\}}, \tag{B.4.6}\label{eq:B.4.6} \end{equation*}
and $\Omega_J$ is unchanged.

For six fixed ordered nodes $S=\{\lambda_1<\cdots<\lambda_6\}$, let $D=D_S$ and $c_i=1/D'(\lambda_i)$. If $P Q=D+C$, every fixed state in the component incident to a positive six-node state, including its five-node edges and four-node corners, is
\begin{equation*} w_i^P=c_i\Omega(\lambda_i)Q(\lambda_i),\qquad w_i^Q=c_i\Omega(\lambda_i)P(\lambda_i), \tag{B.4.7}\label{eq:B.4.7} \end{equation*}
where
\begin{equation*} \Omega=(t-\rho)(t-\sigma),\qquad (\rho,\sigma)\in [\lambda_2,\lambda_3]\times[\lambda_4,\lambda_5]. \tag{B.4.8}\label{eq:B.4.8} \end{equation*}
The rectangle interior has six positive weights, an open edge has five, and a corner has four.
\end{lemma}

\begin{proof}
Signed two-block fixedness and Lemma~\ref{lem:3.4} give
\begin{equation*} P(\lambda_i)Q(\lambda_i)=C,\qquad i\in J. \tag{B.4.9}\label{eq:B.4.9} \end{equation*}
Orthogonality says that $(w_i^PP(\lambda_i))_{i\in J}$ annihilates $\mathcal P_2$. The nullspace of the three evaluation rows is, by \eqref{eq:B.4.2},
$$ \{(c_i^JK(\lambda_i))_{i\in J}:\deg K\le |J|-4\}. $$
Its highest coefficient fixes the scalar and gives the first equation in \eqref{eq:B.4.4}, with $K_J$ monic; \eqref{eq:B.4.9} gives the second. The monic sextic $PQ-C$ vanishes on $J$, so division by $D_J$ gives \eqref{eq:B.4.5}. If the $\alpha$-weight vanishes, \eqref{eq:B.4.4} and $C>0$ imply $K_{J\cup\{\alpha\}}(\alpha)=0$, proving the first equation in \eqref{eq:B.4.6}; comparison in \eqref{eq:B.4.5} proves the second.

For six nodes, \eqref{eq:B.4.2} shows directly that the weights in \eqref{eq:B.4.7} sum to one, that $P,Q$ satisfy the three moment equations, and that both squared monic norms equal $C$. Conversely the preceding argument recovers every incident fixed state.

Both $P$ and $Q$ are monic orthogonal cubics for positive measures on the six distinct nodes, so each has three simple real roots in the open convex hull of the support. Since $PQ=D+C$ and $C>0$, these six roots can occur only in the three node gaps on which $D<0$. Rolle's theorem gives exactly one zero of $D'$ in each of the five node gaps, and hence $D+C$ has at most two roots in any one negative gap. It has the six real roots of $P Q$, so it has exactly two in each negative gap. Orthogonality of either monic cubic for a positive six-node measure forces three sign changes on the ordered nodes. Hence $P$ and $Q$ each take one root from each of the three gaps $(\lambda_1,\lambda_2)$, $(\lambda_3,\lambda_4)$, and $(\lambda_5,\lambda_6)$. Hence
$$ \operatorname{sign}P(\lambda_i) =\operatorname{sign}Q(\lambda_i)=(-,+,+,-,-,+), $$
whereas $\operatorname{sign}c_i=(-,+,-,+,-,+)$. Positivity in \eqref{eq:B.4.7} is equivalent to
$$ \Omega(\lambda_1),\Omega(\lambda_2)\ge0,\quad \Omega(\lambda_3),\Omega(\lambda_4)\le0,\quad \Omega(\lambda_5),\Omega(\lambda_6)\ge0. $$
A monic quadratic has this sign pattern exactly when its roots lie in the two intervals in \eqref{eq:B.4.8}.
\end{proof}

\begin{lemma}[Affine normal and positive quadratic root drift]\label{lem:4.2}
On a regular neighbourhood of a positive six-node fixed state on which every $P(\lambda_i)$ and $Q(\lambda_i)$ is separated from zero, let $H_0$ be the current squared monic norm, $P$ its monic cubic, $Q$ the next monic cubic, and $H_1$ the next squared monic norm. There are unique monic quadratics $K,K_1$ satisfying the general recovery identities \begin{equation*} w_i=\frac{H_0c_iK(\lambda_i)}{P(\lambda_i)},\qquad \widetilde w_i=\frac{H_1c_iK_1(\lambda_i)}{Q(\lambda_i)}. \tag{B.4.9a}\label{eq:B.4.9a} \end{equation*} There is a unique monic quadratic $\Sigma$ satisfying
\begin{equation*} K P Q=D\Sigma+H_1K_1. \tag{B.4.10}\label{eq:B.4.10} \end{equation*}
Define the affine normal polynomial
\begin{equation*} U=\Sigma-K=U_1t+U_0. \end{equation*}
Then $U=0$ is exactly the fixed-cycle equation. Near a coprime fixed component,
\begin{equation*} P^+-P=O(\|U\|^2),\qquad U^+-U=B_{P,K}(U,U)+O(\|U\|^3), \tag{B.4.12}\label{eq:B.4.12} \end{equation*}
\begin{equation*} K^+-K=-\frac2{H_0}\operatorname{rem}_K(DU)+O(\|U\|^2). \tag{B.4.13}\label{eq:B.4.13} \end{equation*}
If $K=(t-\rho)(t-\sigma)$, the two quadratic root-normal coefficients at a fixed state are
\begin{equation*} \Gamma_\rho= -\frac{D(\sigma)}{H_0} \left(\frac{P(\sigma)}{P(\rho)} +\frac{Q(\sigma)}{Q(\rho)}\right)>0, \tag{B.4.14}\label{eq:B.4.14} \end{equation*}
\begin{equation*} \Gamma_\sigma= -\frac{D(\rho)}{H_0} \left(\frac{P(\rho)}{P(\sigma)} +\frac{Q(\rho)}{Q(\sigma)}\right)>0. \tag{B.4.15}\label{eq:B.4.15} \end{equation*}
They are respectively $B_{P,K}(t-\rho,t-\rho)(\rho)$ and $B_{P,K}(t-\sigma,t-\sigma)(\sigma)$. Their displayed rational expressions have finite positive limits as any off-node common-root collision of $P,Q$ is approached.
\end{lemma}

\begin{proof}
Orthogonality says that $(w_iP(\lambda_i))_{i=1}^6$ annihilates $\mathcal P_2$. By \eqref{eq:B.4.2}, it has the form $$ w_iP(\lambda_i)=c_iR(\lambda_i),\qquad \deg R\le2. $$ Moreover $$ H_0=\sum_iw_iP(\lambda_i)^2 =\sum_i c_iR(\lambda_i)P(\lambda_i) =[t^5]\{R(t)P(t)\}. $$ The leading coefficient of $R$ is $H_0$, and $R=H_0K$ for a unique monic quadratic $K$. This gives the first identity in \eqref{eq:B.4.9a}. Applying the result at the next phase gives the second identity. The exact weight update is \begin{equation*} \widetilde w_i=\frac{w_iP(\lambda_i)^2}{H_0} =c_iK(\lambda_i)P(\lambda_i). \tag{B.4.16a}\label{eq:B.4.16a} \end{equation*} Hence $KPQ-H_1K_1$ vanishes at all six nodes. It is monic of degree eight, so division by the monic sextic $D$ gives the unique monic quadratic $\Sigma$ in \eqref{eq:B.4.10}.

If $U=0$, then \begin{equation*} K(PQ-D)=H_1K_1. \end{equation*} Reduction modulo $K$, together with $H_1>0$, shows that $K\mid K_1$. Both polynomials are monic quadratics, hence $K_1=K$, and $PQ=D+H_1$. At a node $H_1/P(\lambda_i)=Q(\lambda_i)$. The mass normalisation and \eqref{eq:B.4.2} give \begin{equation*} 1=\sum_i\frac{H_0c_iK(\lambda_i)}{P(\lambda_i)} =\frac{H_0}{H_1}\sum_i c_iK(\lambda_i)Q(\lambda_i) =\frac{H_0}{H_1}, \end{equation*} because $KQ$ is monic of degree five. Thus $H_0=H_1$, and \eqref{eq:B.4.16a} at the two phases gives exact two-block return.

Conversely, a signed two-block return has squared monic norms $H_0\le H_1\le H_2=H_0$. Equality in the cross-projection chord identity \eqref{eq:B.3.4} gives $$ P(\lambda_i)Q(\lambda_i)=H_0 $$ at every positive-weight node. The two monic sextics $PQ-H_0$ and $D$ agree at six distinct nodes and have the same leading coefficient, so $PQ=D+H_0$. Substitution in \eqref{eq:B.4.10} gives $U=0$. This establishes the exact fixed equivalence.

Near a coprime fixed point, set $C=H_0=H_1$ at the fixed state. The following local analytic coordinates are available. No $P(\lambda_i)$ vanishes, and mass normalisation recovers \begin{equation*} H_0^{-1}=\sum_i\frac{c_iK(\lambda_i)}{P(\lambda_i)}. \tag{B.4.16d}\label{eq:B.4.16d} \end{equation*} The next cubic is obtained from an invertible $3\times3$ moment Gram matrix, and division by the monic $D$ is polynomial in coefficients. Hence \begin{equation*} \Phi(P,K)=(H_0,K,U) \tag{B.4.16e}\label{eq:B.4.16e} \end{equation*} is analytic. Its derivative is injective at a coprime fixed point. To see this, suppose $dH_0=0,dK=0,dU=0$, and denote the variations of $P,Q,K_1,H_1$ by $p,q,k_1,\eta$. Differentiating \eqref{eq:B.4.10} gives \begin{equation*} K(Qp+Pq)=\eta K+Ck_1. \end{equation*} Here $k_1$ has degree at most one because $K_1$ remains monic. Reduction modulo $K$ forces $k_1=0$. The function $H_1-H_0$ is nonnegative near the fixed state and vanishes there, so its differential is zero; hence $\eta=dH_1=dH_0=0$. Thus $Qp+Pq=0$. Coprimality and $\deg p,\deg q\le2$ imply $p=q=0$. Both sides of \eqref{eq:B.4.16e} have dimension five, so the analytic inverse-function theorem makes $(H_0,K,U)$ a local coordinate system.

For the full derivative of the return in these coordinates, hold $H_0=C$ and $K$ fixed, let $u=dU$ be affine, and divide \begin{equation*} D u=KR+r,\qquad \deg r\le1. \tag{B.4.16g}\label{eq:B.4.16g} \end{equation*} There is a unique pair $p,q\in\mathcal P_2$ satisfying \begin{equation*} Qp+Pq=R: \tag{B.4.16h}\label{eq:B.4.16h} \end{equation*} a kernel pair would satisfy $Qp=-Pq$, and coprimality would force it to vanish; dimensions then give existence.

The perturbation $p$ preserves the current norm coordinate. To verify this directly, the derivative of \eqref{eq:B.4.16d} is proportional to $\sum_i c_iK_i p_i/P_i^2$. Since $P_iQ_i=C$ on the fixed component, $$ C^2\sum_i\frac{c_iK_ip_i}{P_i^2} =[t^5]\operatorname{rem}_D(KpQ^2). $$ Using \eqref{eq:B.4.16g}--\eqref{eq:B.4.16h}, $$ KpQ^2=KQ(R-Pq) =Q(Du-r)-K(D+C)q \equiv-Qr-CKq\pmod D. $$ The final polynomial has degree at most four, so its coefficient of $t^5$ is zero.

Let $$ (A_0,B_0)=(P,Q),\quad(A_1,B_1)=(Q,P),\quad (A_2,B_2)=(P,Q), $$ $$ p_0=p,\quad p_1=q,\quad p_2=p,\quad p_3=q,\qquad \kappa_j=-\frac{j}{C}r. $$ For $j=0,1,2$, direct substitution gives the exact coefficient identity \begin{equation*} \begin{aligned} C\kappa_j+K(p_jB_j+A_jp_{j+1}) &=-jr+KR\\ &=Du-(j+1)r =Du+C\kappa_{j+1}. \end{aligned} \end{equation*} This is precisely the coefficientwise differential of the phase factor equation $$ K_jA_jB_j=D(K_j+U_j)+H_{j+1}K_{j+1}, $$ where $K_j,U_j,H_j$ denote the recovery quadratic, affine normal, and squared monic norm at phase $j$. The decomposition into affine normal, scalar norm, affine monic-factor variation, and degree-at-most-two cubic variation is unique. After two blocks, $$ dP^+=dP,\qquad dU^+=dU,\qquad dH_0^+=dH_0, $$ \begin{equation*} dK^+-dK=-\frac2C\operatorname{rem}_K(D\,dU). \end{equation*}

Tangential variations $U=0$ are returned exactly. If $f$ is any coefficient of $P^+-P$, or of $K^+-K+(2/H_0)\operatorname{rem}_K(DU)$, then $f(z,0)=D_Uf(z,0)=0$, where $z=(H_0,K)$. The exact Hadamard formula \begin{equation*} f(z,U)=\int_0^1(1-s)D_U^2f(z,sU)[U,U]\,ds \end{equation*} gives the two $O(\|U\|^2)$ statements uniformly on compact coprime charts. Applying the third-order integral formula after subtracting one half of the $U$-Hessian at zero defines a symmetric bilinear tensor $B_{P,K}$ and proves $$ U^+-U=B_{P,K}(U,U)+O(\|U\|^3). $$

For these two root evaluations, affinely normalise $\rho=0,\sigma=1$, so $K=t(t-1)$, and take the critical normal $u=t$. Define $$ D_1(t)=\frac{D(t)-D(1)}{t-1},\qquad \delta_D=-\frac{D(1)}C, $$ and take the unique $p,q\in\mathcal P_2$ with \begin{equation*} Qp+Pq=D_1. \end{equation*} Because $Dt=KD_1-C\delta_Du$, the first-order affine $K$-displacement at phase $j$ is $j\delta_Du$. Expanding the factor equation through second order gives the exact forcing \begin{equation*} W_j=Kpq+j\delta_DuD_1. \tag{B.4.16m}\label{eq:B.4.16m} \end{equation*} The static term $Kpq$ has zero two-block normal return. Let its unique phase-zero decomposition be \begin{equation*} KPr+Kpq=Da+\eta K+Ck. \tag{B.4.16n}\label{eq:B.4.16n} \end{equation*} At the next phase the incoming cubic and $K$-corrections are $r,k$. Choosing outgoing cubic correction $0$ and outgoing $K$-correction $2k$ gives $$ Ck+KrP+Kpq=Da+\eta K+2Ck $$ by \eqref{eq:B.4.16n}. The affine normal output is again $a$; uniqueness of the decomposition proves that its net contribution is zero.

Factor $\delta_D$ from the remaining forcing. Its phase-one response is the unique decomposition \begin{equation*} KQr+uD_1=Da+\eta K+Ck. \tag{B.4.16o}\label{eq:B.4.16o} \end{equation*} Let $J_Q=(Q-Q(1))/(t-1)$. Because $$ (t-1)D_1=D-D(1),\qquad (t-1)J_Q=Q-Q(1), $$ and $t-1$ is invertible modulo $Q$, \begin{equation*} \operatorname{rem}_Q(uD_1)=P(1)\{J_Q-Q(1)\}. \tag{B.4.16p}\label{eq:B.4.16p} \end{equation*} Its value at zero is $-P(1)Q(0)$. Evaluation of \eqref{eq:B.4.16o} at zero and reduction modulo $Q$, followed by evaluation at zero, give $$ D(0)a(0)+Ck(0)=0, $$ $$ -P(1)Q(0)=-Ca(0)+Ck(0). $$ Using $P(0)Q(0)=D(0)+C$ and solving gives \begin{equation*} a(0)=\frac{P(1)}{P(0)},\qquad k(0)=-\frac{D(0)P(1)}{CP(0)}. \tag{B.4.16q}\label{eq:B.4.16q} \end{equation*}

For the next phase, let $s_{\rm cub}\in\mathcal P_2$ denote the outgoing cubic correction. The corresponding identity is \begin{equation*} Ck+KrQ+KP s_{\rm cub}+2uD_1=Db+\zeta K+C\ell. \tag{B.4.16r}\label{eq:B.4.16r} \end{equation*} Substitute \eqref{eq:B.4.16o} into \eqref{eq:B.4.16r} and reduce modulo $P$. Since the remaining representatives have degree at most two, the congruence is an identity: \begin{equation*} -Ca+\eta K+2Ck+\operatorname{rem}_P(uD_1) =-Cb+\zeta K+C\ell. \end{equation*} The $P$-version of \eqref{eq:B.4.16p} has value $$ \operatorname{rem}_P(uD_1)(0)=-Q(1)P(0). $$ Evaluation of \eqref{eq:B.4.16r} at zero also gives $$ Ck(0)=D(0)b(0)+C\ell(0). $$ Eliminating $\ell(0)$ from these two equations and using \eqref{eq:B.4.16q} produces $$ P(0)Q(0)b(0)=P(1)Q(0)+Q(1)P(0), $$ so \begin{equation*} b(0)=\frac{P(1)}{P(0)}+\frac{Q(1)}{Q(0)}. \end{equation*} In the second-order expansion, $a$ and $b$ are respectively the phase-one and phase-two affine normal coefficients caused by the unit forcing $uD_1$. The input coordinate $U=\varepsilon u$ has no quadratic term, the $Kpq$ contribution was proved above to have zero two-block normal return, and \eqref{eq:B.4.16m} multiplies the remaining response by $\delta_D$. Hence the returned quadratic normal coefficient is $\delta_Db$. Multiplication by $\delta_D=-D(1)/C$ proves \begin{equation*} B_{P,K}(t,t)(0)= -\frac{D(1)}C \left\{\frac{P(1)}{P(0)}+\frac{Q(1)}{Q(0)}\right\}. \tag{B.4.16u}\label{eq:B.4.16u} \end{equation*}

For $g=\sigma-\rho$, let $x=(t-\rho)/g$ and $$ \bar P(x)=g^{-3}P(\rho+gx),\quad \bar Q(x)=g^{-3}Q(\rho+gx),\quad \bar D(x)=g^{-6}D(\rho+gx), $$ $$ \bar C=g^{-6}C,\qquad \bar U(x)=g^{-2}U(\rho+gx),\qquad \bar K(x)=g^{-2}K(\rho+gx)=x(x-1). $$ Substitution in the phase identity gives \begin{equation*} B_{\bar P,\bar K}(\bar U,\bar U)(x) =g^{-2}B_{P,K}(U,U)(\rho+gx). \tag{B.4.16v}\label{eq:B.4.16v} \end{equation*} Taking $\bar U=x$ corresponds to $U=g(t-\rho)$; bilinearity cancels the factor $g^2$ in \eqref{eq:B.4.16v}. Formula \eqref{eq:B.4.16u} is exactly \eqref{eq:B.4.14}. Exchanging $\rho,\sigma$ gives \eqref{eq:B.4.15}.

Each cubic has one root in every negative spectral gap. Hence $$ P(\rho),Q(\rho)>0,\qquad P(\sigma),Q(\sigma)<0, $$ whereas $D(\rho),D(\sigma)>0$. The brackets in \eqref{eq:B.4.14}--\eqref{eq:B.4.15} are strictly negative and their prefactors are negative, so both $\Gamma$'s are positive. At a common-factor collision the common root remains in a negative gap while $\rho,\sigma$ remain in the two complementary positive gaps. None of the four factor evaluations in the denominators vanishes. The formulas have finite strictly positive limits.
\end{proof}

\begin{lemma}[Simple and simultaneous collision ideals]\label{lem:4.3}
Every common-root collision of $P,Q$ on a positive six-node fixed set is an off-node collision at a nondegenerate critical point of $D$. In the regular local analytic state ring the fixed ideal
\begin{equation*} \mathcal I=(U_0,U_1) \end{equation*}
is a radical complete intersection and
\begin{equation*} \mathcal I^{(m)}=\mathcal I^m\qquad(m\ge1). \end{equation*}
The formulae \eqref{eq:B.4.12}--\eqref{eq:B.4.15} and the positive quadratic norm gain extend through every simple or simultaneous collision.

At one omitted $+1$-multiplier node $\alpha$, let $x$ be its signed amplitude and define
\begin{equation*} \bar K_\alpha(t)=\frac{K(t)-K(\alpha)}{t-\alpha},\qquad A=[U,\bar K_\alpha],\qquad \eta=\frac{U(\alpha)}{\bar K_\alpha(\alpha)}, \tag{B.4.19}\label{eq:B.4.19} \end{equation*}
where $[f,g]=f_1g_0-f_0g_1$ for affine polynomials. Then $\mathcal I=(A,\eta)$, and the union of the adjacent five- and six-node fixed strata has radical complete-intersection ideal
\begin{equation*} \mathcal I_1=(A,x\eta). \end{equation*}

At two such nodes $\alpha,\beta$, with signed amplitudes $x,y$, define $\bar K_\gamma(t)=\{K(t)-K(\gamma)\}/(t-\gamma)$ for $\gamma\in\{\alpha,\beta\}$, and
\begin{equation*} \eta_\alpha=[U,\bar K_\beta],\qquad \eta_\beta=[U,\bar K_\alpha]. \tag{B.4.21}\label{eq:B.4.21} \end{equation*}
Then $\mathcal I=(\eta_\alpha,\eta_\beta)$, and the union of the four-, five-, and six-node fixed strata has radical complete-intersection ideal
\begin{equation*} \mathcal I_2=(x\eta_\alpha,y\eta_\beta). \end{equation*}
Both union ideals also have symbolic powers equal to ordinary powers.
\end{lemma}

\begin{proof}
A common root $h$ of $P,Q$ satisfies $$ D(h)+C=P(h)Q(h)=0 $$ and, after differentiating $PQ=D+C$, $$ D'(h)=P'(h)Q(h)+P(h)Q'(h)=0. $$ It is not a spectral node because $C>0$. Rolle's theorem gives at least one zero of $D'$ in each of the five open gaps between the six nodes. Since $\deg D'=5$, these are all its zeros and each is simple. Hence $D''(h)\ne0$; in a negative component $h$ is the unique strict minimum, so $D''(h)>0$.

Complexify the real-analytic state chart and work in the regular local ring \begin{equation*} R=\mathbb C\{p_0,p_1,p_2,k_0,k_1\}, \end{equation*} which has dimension five. Lemma~\ref{lem:4.2} identifies $V(\mathcal I)$, $\mathcal I=(U_0,U_1)$, with the factorisations $PQ=D+C$, with the two coefficients of $K$ free.

We analyse the ramification at a simultaneous critical level via a finite normalisation. Let $\delta=C-C_*$. At the nondegenerate critical point in the $j$th participating component of $\{D<0\}$, analytic Morse preparation gives analytic coordinates and a unit $u_j$ such that \begin{equation*} D(t)+C=u_j(t,\delta)\{s_j(t,\delta)^2+\delta\}. \tag{B.4.17b}\label{eq:B.4.17b} \end{equation*} Pass now to the finite normalisation $\delta=-q^2$. On this cover the two local roots are the analytic branches $s_j=\pm q$, after absorbing an analytic square root of the unit. Every remaining root is analytic in $\delta$ by the implicit-function theorem. Assigning the two roots in each participating component, and the remaining roots, between the monic factors $P$ and $Q$ gives finitely many smooth reduced $q$-branches. Changing $q$ to $-q$ reverses every participating allocation without changing the image branch, so one, two, or three simultaneous collisions give respectively $1,2$, or $4$ factor-allocation branches. Conversely, Weierstrass division shows that every nearby monic factorisation occurs on one of these branches. The cover is finite and surjective, so dimension is preserved on descent. Each factor branch has one parameter $q$ and the free $K$-plane has dimension two; every irreducible component of $V(\mathcal I)$ in the original ring has dimension three and height exactly two. Since $R$ is regular and $\mathcal I$ has two generators, $U_0,U_1$ form a regular sequence.

Here is the local algebra in the maximally ramified case. If the three negative components collide at $h_1,h_2,h_3$, set $S(t)=\prod_{j=1}^3(t-h_j)$. The monic sextics $D+C_*$ and $S^2$ have the same three double roots, and hence $D+C_*=S^2$. The three roots of $P$ near the simple roots of $S$ are analytic functions of its coefficients. Applying the parameterised Morse coordinates in \eqref{eq:B.4.17b} at those roots gives analytic factor coordinates $z_1,z_2,z_3$, with the two coefficients of $K$ as spectator variables, such that the selected root in the $j$th pair satisfies $z_j^2+\delta=0$. Eliminating $\delta$ gives the triple-collision germ
\begin{equation*} \mathcal I_{\rm tri}=(z_2^2-z_1^2,z_3^2-z_1^2) =\bigcap_{\epsilon_2,\epsilon_3\in\{\pm1\}} (z_2-\epsilon_2z_1,z_3-\epsilon_3z_1). \tag{B.4.17b1}\label{eq:B.4.17b1} \end{equation*}
The four displayed primes are exactly the four factor-allocation branches, and each has height two. The two quadrics therefore have height two in the regular factor-coefficient ring and form a regular sequence. Each branch is smooth away from the central point, and at its points with $z_1\ne0$ the two gradients have rank two, so the complete intersection is generically reduced on every minimal component. The same $S_1$ argument used below then shows directly that \eqref{eq:B.4.17b1} is radical. For comparison, after analytic changes of coordinates the one- and two-collision germs are respectively
\begin{equation*} (z_2,z_3),\qquad (z_3,z_2^2-z_1^2)=(z_3,z_2-z_1)\cap(z_3,z_2+z_1). \tag{B.4.17b2}\label{eq:B.4.17b2} \end{equation*}
At the central triple point the ordinary differentials of the two quadratic generators vanish. No independence of those differentials is asserted or needed: height two gives the regular sequence, equivalently the two generator classes form a free basis of $\mathcal I_{\rm tri}/\mathcal I_{\rm tri}^2$.

Every normalised branch contains points with $q\ne0$. There the allocated roots are distinct, $P$ and $Q$ are coprime, and the coprime coordinate calculation of Lemma~\ref{lem:4.2} gives $\operatorname{rank}dU=2$. The localisation of $R/\mathcal I$ at each minimal prime is generically a field and is reduced. A complete-intersection quotient is Cohen--Macaulay, hence satisfies $S_1$ and has no embedded associated prime. A nonzero nilradical would have an associated prime that was either embedded or a minimal prime with nonreduced localisation; both alternatives have just been excluded. Hence $R/\mathcal I$ is reduced and $\mathcal I$ is radical. This reducedness holds on every branch of the finite normalisation.

The symbolic-power assertion follows from the linear-type property of a regular sequence: \begin{equation*} \operatorname{gr}_{\mathcal I}R =\bigoplus_{m\ge0}\mathcal I^m/\mathcal I^{m+1} \cong(R/\mathcal I)[T_0,T_1]. \tag{B.4.17c}\label{eq:B.4.17c} \end{equation*} Each graded piece is a free $R/\mathcal I$-module. Since the latter is reduced and Cohen--Macaulay, its associated primes are exactly $\operatorname{Min}(\mathcal I)$. The exact sequences \begin{equation*} 0\longrightarrow\mathcal I^{m-1}/\mathcal I^m \longrightarrow R/\mathcal I^m \longrightarrow R/\mathcal I^{m-1}\longrightarrow0 \end{equation*} show inductively that $\operatorname{Ass}(R/\mathcal I^m)\subset \operatorname{Min}(\mathcal I)$; the reverse inclusion holds because these are the minimal primes of $\mathcal I^m$. Hence $\mathcal I^m$ has no embedded primary component, and \begin{equation*} \mathcal I^m =\bigcap_{\mathfrak p\in\operatorname{Min}\mathcal I} (\mathcal I^mR_{\mathfrak p}\cap R) =\mathcal I^{(m)}. \tag{B.4.17e}\label{eq:B.4.17e} \end{equation*}

For completeness, let $R_{\mathbb R}$ be the real convergent-power-series state ring and $R_{\mathbb C}=R_{\mathbb R}\otimes_{\mathbb R}\mathbb C=R_{\mathbb R}[i]$ its complexification. This extension is finite free and faithfully flat, so $\mathcal I_{\mathbb C}\cap R_{\mathbb R}=\mathcal I_{\mathbb R}$; in particular, the radicality of $\mathcal I_{\mathbb C}$ proved above first gives the radicality of $\mathcal I_{\mathbb R}$. Every minimal prime $\mathfrak q$ of $\mathcal I_{\mathbb C}$ lies over a minimal prime $\mathfrak p$ of $\mathcal I_{\mathbb R}$, and for each $\mathfrak p$ the semilocal map
\begin{equation*} (R_{\mathbb R})_{\mathfrak p}\longrightarrow \prod_{\mathfrak q\mid\mathfrak p}(R_{\mathbb C})_{\mathfrak q} \tag{B.4.17e1}\label{eq:B.4.17e1} \end{equation*}
is faithfully flat. Consequently, for $f\in R_{\mathbb R}$,
\begin{equation*} f\in\mathcal I_{\mathbb R}^{(m)} \Longleftrightarrow f\otimes1\in\mathcal I_{\mathbb C}^{(m)} \Longleftrightarrow f\otimes1\in\mathcal I_{\mathbb C}^m \Longleftrightarrow f\in\mathcal I_{\mathbb R}^m. \tag{B.4.17e2}\label{eq:B.4.17e2} \end{equation*}
The first equivalence is the preceding localisation test over every prime above $\mathfrak p$, and the last is contraction of an extended ideal under a faithfully flat map. Regularity of the two-generator sequence descends because multiplication by either successive generator is injective after the faithfully flat base change. This proves the real radical, complete-intersection, and symbolic-power assertions without assuming that a real minimal prime remains prime after complexification.

At the generic point $\mathfrak p$ of every factor branch, Lemma~\ref{lem:4.2} gives \begin{equation*} P^+-P,\quad U^+-U,\quad K^+-K+\frac2{H_0}\operatorname{rem}_K(DU) \ \in(\mathcal IR_{\mathfrak p})^2. \end{equation*} Their coefficients belong to the symbolic square, and \eqref{eq:B.4.17e} promotes them to $\mathcal I^2$.

We compute the norm-gain quadratic explicitly. Define $$ \Delta K=-\frac2C\operatorname{rem}_K(DU). $$ At the two fixed phases, $$ w_i^P=c_iK_iQ_i,\qquad w_i^Q=c_iK_iP_i. $$ Their returned first weight variations are respectively $\delta w_i^P=c_iQ_i\Delta K_i$ and $\delta w_i^Q=c_iP_i\Delta K_i$. Passing to signed coordinates gives $$ \|\delta Y^P\|^2 =\frac14\sum_i\frac{c_iQ_i}{K_i}\Delta K_i^2,\qquad \|\delta Y^Q\|^2 =\frac14\sum_i\frac{c_iP_i}{K_i}\Delta K_i^2. $$ Let $H_0,H_1,H_2$ be the three consecutive squared monic norms. The exact chord identities \eqref{eq:B.3.4} read $$ \sqrt{H_0/H_1}=1-\frac12\|Y^{P,+}-Y^P\|^2,\qquad \sqrt{H_1/H_2}=1-\frac12\|Y^{Q,+}-Y^Q\|^2. $$ Both squared chords are $O(\|U\|^2)$, while $H_j=C+O(\|U\|)$. Squaring the two identities and adding yields $$ H_2-H_0=\mathcal Q(U)+O(\|U\|^3), $$ where \begin{equation*} \mathcal Q(U)=\frac C4\sum_i \frac{c_i\{P(\lambda_i)+Q(\lambda_i)\}}{K(\lambda_i)} \left\{-\frac2C \operatorname{rem}_K(DU)(\lambda_i)\right\}^2. \end{equation*} At a generic branch point the displayed Taylor remainder is in $(\mathcal IR_{\mathfrak p})^3$ on every branch and hence, by \eqref{eq:B.4.17e}, in $\mathcal I^3$. Moreover \begin{equation*} \frac{c_iQ(\lambda_i)}{K(\lambda_i)} =\frac{w_i^P}{K(\lambda_i)^2}>0,\qquad \frac{c_iP(\lambda_i)}{K(\lambda_i)} =\frac{w_i^Q}{K(\lambda_i)^2}>0. \tag{B.4.17h}\label{eq:B.4.17h} \end{equation*} Neither root of $K$ is a spectral node, so $D$ is a unit in $R[t]/(K)$; multiplication by $D$, followed by remainder modulo $K$, is invertible on affine polynomials. The form $\mathcal Q$ is positive definite. Compactness of the collision set gives a uniform lower eigenvalue, while the $\mathcal I^3$ remainder is bounded by $M\|U\|^3$. After shrinking the chart, both the norm gain and the same-parity chord are coercive in $\|U\|^2$. Coefficient functions that agree on the dense coprime part agree in the reduced ring $R/\mathcal I$; hence the quadratic-drift formulas \eqref{eq:B.4.14}--\eqref{eq:B.4.15} have the stated collision values.

At the signed support boundaries, first suppose that one $+1$-multiplier node $\alpha$ is omitted. Neither $P(\alpha)$ nor $Q(\alpha)$ vanishes because their product is $C$. On the fixed quotient, \begin{equation*} x^2=w_\alpha=c_\alpha K(\alpha)Q(\alpha), \tag{B.4.19a}\label{eq:B.4.19a} \end{equation*} and $c_\alpha Q(\alpha)$ is a unit. Evaluation $K\mapsto K(\alpha)$, together with one complementary coefficient of $K$, is a coordinate on the free $K$-plane. Substitution from \eqref{eq:B.4.19a} identifies the signed fixed quotient with $B\{x,\kappa\}$, where $B$ is the reduced factor-curve ring. It is reduced. With two omitted nodes $\alpha,\beta$, the evaluation map $K\mapsto(K(\alpha),K(\beta))$ is invertible because $\alpha\ne\beta$, and the signed quotient is $B\{x,y\}$, again reduced.

After one support loss, define $\bar K_\alpha=(K-K(\alpha))/(t-\alpha)$. On $x=0$, $$ K=(t-\alpha)(t-\theta). $$ A five-node fixed state has one remaining completion root $\xi$ and $$ PQ=D_J(t)(t-\xi)+C. $$ Multiplication by $K$, with $D=D_J(t)(t-\alpha)$, gives $$ KPQ=D(t)(t-\theta)(t-\xi)+CK. $$ Comparison with $KPQ=D(K+U)+CK$ yields \begin{equation*} U=(\alpha-\xi)(t-\theta) =(\alpha-\xi)\bar K_\alpha. \tag{B.4.19b}\label{eq:B.4.19b} \end{equation*} Hence $A=[U,\bar K_\alpha]=0$ is its fixed equation and $\eta=U(\alpha)/\bar K_\alpha(\alpha)=\alpha-\xi$ is its free shear. If $U=u_1t+u_0$ and $\bar K_\alpha=t-\theta$, then \begin{equation*} A=-\theta u_1-u_0,\qquad \eta=\frac{\alpha u_1+u_0}{\alpha-\theta}. \tag{B.4.19c}\label{eq:B.4.19c} \end{equation*} The determinant of this change from $(u_0,u_1)$ is $-1$, so $(U_0,U_1)=(A,\eta)$. The lower and higher fixed components have ideals $(x,A)$ and $(A,\eta)$; their union is \begin{equation*} V(A,x\eta). \tag{B.4.19d}\label{eq:B.4.19d} \end{equation*}

After two support losses, define $\eta_\alpha,\eta_\beta$ as in \eqref{eq:B.4.21}. Since $\bar K_\gamma=t+\gamma+k_1$, their constant coefficients differ by $\alpha-\beta$, so they are an invertible linear change of the two coefficients of $U$. If $\beta$ is absent, then $K=(t-\beta)(t-\rho)$, and \eqref{eq:B.4.19b}, with $\beta$ in place of $\alpha$, gives \begin{equation*} U=(\beta-\xi)(t-\rho),\qquad \eta_\alpha=0,\qquad \eta_\beta=(\alpha-\beta)(\beta-\xi). \end{equation*} Exchanging the two nodes gives the other face. The four strata have ideals \begin{equation*} (x,y),\quad(y,\eta_\alpha),\quad (x,\eta_\beta),\quad(\eta_\alpha,\eta_\beta), \tag{B.4.21b}\label{eq:B.4.21b} \end{equation*} and their union is \begin{equation*} V(x\eta_\alpha,y\eta_\beta). \end{equation*}

We calculate the height in the actual signed state ring, without treating the coefficients of $U$ as independent state coordinates at a collision. This ring has dimension five: near a two-loss point one may use three independent retained weights and the two signed omitted amplitudes. The retained weights determine the remaining weight by normalisation; the moment Gram matrices remain invertible. Near a one-loss point one uses four independent retained weights and one signed amplitude.

A common parameter count describes all the support components. Let $J$ be an active support of size $m\in\{4,5,6\}$, and write $D_J(t)=\prod_{i\in J}(t-\lambda_i)$ and $c_i^J=1/D_J'(\lambda_i)$. Every nearby fixed state on this support has
\begin{equation*}
PQ=D_J R+C,\qquad
w_i=c_i^J L(\lambda_i)Q(\lambda_i),\qquad i\in J,
\tag{B.4.21c}\label{eq:B.4.21c}
\end{equation*}
where $R$ is monic of degree $6-m$ and $L$ is monic of degree $m-4$. To verify the converse, the degree-$(m-1)$ Lagrange coefficient identity gives $\sum_iw_i=1$. At the nodes, $P Q=C$, so
\[
\sum_i w_iP(\lambda_i)\lambda_i^j
 =C\sum_i c_i^J L(\lambda_i)\lambda_i^j=0
 \quad(0\le j\le2),\qquad
\sum_iw_iP(\lambda_i)^2=C.
\]
The returned weights are $c_i^J L(\lambda_i)P(\lambda_i)$; the same calculation identifies their monic cubic as $Q$, their norm as $C$, and their return as $w$. Conversely, the fixed-cycle identity gives $PQ-C=D_JR$, and recovery gives $L$ and the displayed weights. All these identifications use only nonzero node evaluations and invertible moment systems, so they are analytic on the chosen charts.

For fixed $(C,R)$ there are finitely many monic cubic factor allocations of $D_JR+C$. Adjoin the roots successively using the monic root equation and the successive monic quotients. This gives a finite free ordered-root cover of the $(C,R)$ parameter space; assigning three roots to each factor gives all the local factor allocations. No minimal component of this cover is supported only over the discriminant: a nonzero parameter function acts injectively on the finite free algebra, whereas an element of a minimal prime is a zero divisor. Consequently each factor-allocation branch has $1+(6-m)$ parameters before the free $L$ coefficients are included. Thus every support component has dimension
\[
1+(6-m)+(m-4)=3.
\]
This count also holds at the boundary of a component. At an additional omitted weight, the recovery coefficient multiplying $L(\lambda_i)$ is a unit. Replacing that coefficient of $L$ by the signed amplitude via $w_i=x_i^2$ preserves the dimension. With two omitted weights the two evaluation functionals on the quadratic recovery polynomial are independent, as in \eqref{eq:B.4.19a}. At factor collisions the finite root covers specify the different factor-allocation branches. They do not make the coefficients of $U$ ordinary independent coordinates.

The exact divisions \eqref{eq:B.4.19b}--\eqref{eq:B.4.21b} identify the zero sets of $\mathcal I_1=(A,x\eta)$ and $\mathcal I_2=(x\eta_\alpha,y\eta_\beta)$ with these unions. Every irreducible component therefore has height two in the regular signed state ring. Both two-generated ideals are complete intersections.

We verify generic reducedness, including the face equations. On a five-node face with omitted node $\gamma$, write $K=(t-\gamma)(t-\rho)$ and $K_1=(t-\gamma)(t-\rho_1)$. Divide $U$ as $U=U_1(t-\rho)+U(\rho)$ and put $\xi=\gamma-U_1$. Cancelling $t-\gamma$ in \eqref{eq:B.4.10} and then evaluating at $\rho$ gives the exact identity
\begin{equation*}
PQ-D_J(t)(t-\xi)-H_1
 =U(\rho)\frac{D_J(t)-D_J(\rho)}{t-\rho}.
\tag{B.4.21d}\label{eq:B.4.21d}
\end{equation*}
The recovery formula makes the three coefficients of $P$ and the pivot $\rho$ ordinary coordinates on this face. At a coprime fixed point, the differential $(\dot P,\dot Q)\mapsto Q\dot P+P\dot Q$ is invertible on the six monic-factor coefficients. Hence the factor family $PQ=D_J(t)(t-\xi)+C$ has two parameters $(\xi,C)$ and projects immersively into the three $P$ coefficients. Indeed, if $\dot P=0$, then $P\dot Q=\dot C-D_J\dot\xi$. If $\dot\xi\ne0$, reduction modulo $P$ would make $D_J$ constant modulo $P$; the fixed identity would then make $d(t-\xi)+C$ divisible by the cubic $P$ for a constant $d$. Degree comparison forces $d=C=0$, contrary to $C>0$. Thus $\dot\xi=\dot C=\dot Q=0$.

If $d(U(\rho))=0$ on a variation of the five-node state, differentiation of \eqref{eq:B.4.21d} puts its $\dot P$ in this two-dimensional factor tangent space. The free pivot supplies just one further direction. The kernel of $d(U(\rho))$ therefore has dimension at most three in the four-dimensional face, so this differential is nonzero. This proves that the absent-amplitude equation and its single fixed normal define a reduced smooth component at the generic coprime point.

On a six-node component, choose a coprime factor allocation and nonzero omitted amplitudes; the normal calculation of Lemma~\ref{lem:4.2} gives two independent defining equations. On the four-node component, choose $R$ in \eqref{eq:B.4.21c} nonzero at the omitted nodes. The bracket normals are then units, and the union ideal is simply $(x,y)$. Such points are dense on every relevant support component. The finite free root cover above also shows that the coprime points are dense in every factor-allocation branch. Thus $\mathcal I_1$ and $\mathcal I_2$ are generically reduced at every minimal prime.

A complete-intersection quotient has no embedded associated prime. Componentwise generic reducedness therefore proves that $\mathcal I_1,\mathcal I_2$ are radical. The associated-graded and localisation arguments \eqref{eq:B.4.17c}--\eqref{eq:B.4.17e}, applied to these regular sequences, give $\mathcal I_j^{(m)}=\mathcal I_j^m$ for every $m$. The faithfully flat argument \eqref{eq:B.4.17e1}--\eqref{eq:B.4.17e2} gives the same conclusions over the real analytic ring. Throughout, minimal primes mean the individual irreducible factor-allocation branches of the support components. In particular, the last ideal in \eqref{eq:B.4.21b} has two minimal primes at a proper double collision and four at a triple collision. No independence of the ordinary differentials of the collision generators is used.
\end{proof}

\begin{lemma}[Energy control of the completion polynomial]\label{lem:4.4}
For a regular state $Y$, let $Q_Y=P_{L(Y)}$ and define the monic completion sextic
\begin{equation*} \mathscr D_Y=P_YQ_Y-\sigma(Y)\sigma(L(Y)). \end{equation*}
On every compact regular neighbourhood of a four-, five-, or six-node fixed point, including every simple or simultaneous off-node factor collision and every one- or two-weight $+1$ incidence, there is $C_{\mathscr D}<\infty$ such that
\begin{equation*} \|\mathscr D_{L^2(Y)}-\mathscr D_Y\|_{\rm coeff} \le C_{\mathscr D}\|L^2(Y)-Y\|^2. \tag{B.4.24}\label{eq:B.4.24} \end{equation*}
\end{lemma}

\begin{proof}
Let $G=L^2$ and define $$ g(Y)=\mathscr D_{GY}-\mathscr D_Y. $$ All coefficient functions below are real analytic on the regular, pivot-separated charts because the moment Gram matrices are uniformly invertible. All norms on the finite-dimensional space of sextic coefficients are equivalent.

First consider six positive nodes. On the fixed set, $\mathscr D=PQ-C=D$, independently of the fixed parameter $K$. At a generic coprime fixed point, Lemma~\ref{lem:4.2} shows that the derivative of $G-I$ is the $K$-shear \begin{equation*} \delta K=-\frac2C\operatorname{rem}_K(DU), \tag{B.4.24a}\label{eq:B.4.24a} \end{equation*} which is tangent to that fixed set. Every coefficient of $g$ has zero value and zero first normal derivative. It belongs to the symbolic square of $(U_0,U_1)$, hence to its ordinary square by Lemma~\ref{lem:4.3}: \begin{equation*} \|g(Y)\|_{\rm coeff}\le C_1\|U\|^2. \tag{B.4.24b}\label{eq:B.4.24b} \end{equation*} At a fixed state, the first returned same-phase weight variation is $$ \delta w_i=c_iQ(\lambda_i)\delta K(\lambda_i). $$ For signed coordinates, \begin{equation*} \|\delta Y\|^2 =\frac14\sum_i \frac{c_iQ(\lambda_i)}{K(\lambda_i)} \delta K(\lambda_i)^2. \tag{B.4.24c}\label{eq:B.4.24c} \end{equation*} The coefficients are positive by \eqref{eq:B.4.17h}, and multiplication by $D$ modulo $K$ is invertible. Hence \eqref{eq:B.4.24a}--\eqref{eq:B.4.24c} define a positive definite quadratic form in $U$. After subtracting this linear chord, the remainder of $GY-Y$ lies in the ordinary square of the fixed ideal. Shrinking a compact chart gives \begin{equation*} \|GY-Y\|^2\ge c_1\|U\|^2. \tag{B.4.24d}\label{eq:B.4.24d} \end{equation*} The equality of symbolic and ordinary powers in Lemma~\ref{lem:4.3} makes \eqref{eq:B.4.24b}--\eqref{eq:B.4.24d} valid at simultaneous collisions as well.

Next consider a positive five-node chart. Let $J$ be its node set, let $\theta$ be the pivot, and define \begin{equation*} x_i=\lambda_i-\theta,\quad \mathfrak A=\sum_i\frac{w_i}{x_i},\quad B=\sum_i\frac{w_i}{x_i^2},\quad \delta=\frac{\mathfrak A}{B},\quad \mathcal N=1-\frac{\mathfrak A^2}{B}. \end{equation*} We now derive the exact signed return. With $c_i^J=1/D_J'(\lambda_i)$, the three moment equations and \eqref{eq:B.4.2} give $$ w_iP(\lambda_i)=Hc_i^J(\lambda_i-\theta). $$ The coefficient $H$ is the squared monic norm by the degree-four leading-coefficient identity. Moreover $\sum_i c_i^JP(\lambda_i)=0$, so substitution of this formula gives $$ \theta= \frac{\sum_i(c_i^J)^2\lambda_i/w_i} {\sum_i(c_i^J)^2/w_i}. $$ After one block, $$ \widetilde w_i =\frac{w_iP(\lambda_i)^2}{H} =\frac{H(c_i^J)^2x_i^2}{w_i}. $$ The preceding dual formula gives the partner pivot $$ \widetilde\theta =\theta+ \frac{\sum_iw_i/x_i}{\sum_iw_i/x_i^2} =\theta+\delta. $$ Applying the recovery identity to the partner and then the second block gives $$ w_i^+=\frac{\widetilde H}{H}\, w_i(1-\delta/x_i)^2. $$ Summing this identity and using $$ \sum_iw_i(1-\delta/x_i)^2 =1-2\delta\mathfrak A+\delta^2B =1-\frac{\mathfrak A^2}{B}=\mathcal N $$ shows $\widetilde H/H=1/\mathcal N$. Taking the signed square root near the fixed component, where every factor is positive, proves the exact formula \begin{equation*} (GY)_i=Y_i\frac{1-\delta/x_i}{\sqrt{\mathcal N}}. \tag{B.4.24f}\label{eq:B.4.24f} \end{equation*} Its inner product with $Y$ is $$ \langle GY,Y\rangle =\frac{1-\delta\mathfrak A}{\sqrt{\mathcal N}} =\sqrt{\mathcal N}. $$ It follows that \begin{equation*} \|GY-Y\|^2 =2(1-\sqrt{\mathcal N}) =\frac{2\mathfrak A^2}{B(1+\sqrt{\mathcal N})} \ge\frac{\mathfrak A^2}{B}. \tag{B.4.24g}\label{eq:B.4.24g} \end{equation*} On a compact pivot-separated chart $B$ is bounded above, so the chord controls $\mathfrak A^2$ throughout the simultaneous factor-collision strata.

The fixed set is $V(\mathfrak A)$. It is a reduced principal complete intersection: each irreducible component contains dense coprime points at which \eqref{eq:B.4.24f} displays a reduced fixed hypersurface, and a principal Cohen--Macaulay quotient has no embedded primes. On that component, \begin{equation*} PQ-C=D_J(t)(t-\xi),\qquad w_i=c_i^J(\lambda_i-\theta)Q(\lambda_i), \tag{B.4.24h}\label{eq:B.4.24h} \end{equation*} where the remaining completion root $\xi$ is fixed. Holding $C,\xi$, and the discrete factor allocation fixed while differentiating $\theta$ gives $$ \partial_\theta w_i=-c_i^JQ(\lambda_i). $$ Differentiation of \eqref{eq:B.4.24f} in the normal $\mathfrak A$-direction gives the weight shear \begin{equation*} V_i^{\rm wt}=-\frac{2w_i}{B(\lambda_i-\theta)} =-\frac2B c_i^JQ(\lambda_i) =\frac2B\,\partial_\theta w_i. \tag{B.4.24i}\label{eq:B.4.24i} \end{equation*} The first return derivative changes only $\theta$ and leaves the completion \eqref{eq:B.4.24h} fixed. Every coefficient of $g$ has zero first normal derivative, so \begin{equation*} \|g(Y)\|_{\rm coeff}\le C_2\mathfrak A^2. \tag{B.4.24j}\label{eq:B.4.24j} \end{equation*} On four nodes $G$ is exactly the identity: the one-dimensional moment annihilator is multiplied by its reciprocal and returned by the second block. Hence $g=0$ there.

The strata fit together as follows. At one $+1$ loss use $\mathcal I_1=(A,x\eta)$ from Lemma~\ref{lem:4.3}. At a generic point of the five-node component, evenness in the absent signed amplitude $x$ forces its derivative to vanish, and the $A$-normal returns as the $\theta$-shear \eqref{eq:B.4.24i}, which fixes the completion. At a generic point of the six-node component, \eqref{eq:B.4.24a} is tangent and the completion is $D$. Hence each coefficient of $g$ belongs locally to the square at every minimal prime. Lemma~\ref{lem:4.3} gives $$ g\in\mathcal I_1^{(2)}=\mathcal I_1^2 =(A^2,Ax\eta,x^2\eta^2). $$ Since $2|Ax\eta|\le A^2+x^2\eta^2$, \begin{equation*} \|g(Y)\|_{\rm coeff}\le C_3(A^2+x^2\eta^2). \tag{B.4.24k}\label{eq:B.4.24k} \end{equation*}

Every component of $\Xi=GY-Y$ vanishes on $V(\mathcal I_1)$. Here the ambient signed chart is regular, $\mathcal I_1=(A,x\eta)$ is the radical complete intersection proved in Lemma~\ref{lem:4.3}, and the vector-valued analytic Hadamard lemma is applied componentwise to its two generators. It gives \begin{equation*} \Xi=AV+x\eta W \end{equation*} with analytic vector coefficients. The retained signed-coordinate coefficient from \eqref{eq:B.4.24f} is \begin{equation*} V_i=-\frac{Y_i}{B(\lambda_i-\theta)}, \end{equation*} so $V$ has a nonzero retained part.

The unit needed to normalise the omitted multiplier is explicit. On the five-node face, $$ P(t)Q(t)=D_J(t)(t-\xi)+C,\qquad \eta=\alpha-\xi. $$ The analytic extension $M_\alpha=Y_\alpha^+/Y_\alpha$ satisfies \begin{equation*} M_\alpha-1 =\frac{P(\alpha)Q(\alpha)-C}{C} =\frac{D_J(\alpha)}{C}(\alpha-\xi) =\frac{D_J(\alpha)}{C}\eta. \tag{B.4.24m1}\label{eq:B.4.24m1} \end{equation*} Because $\alpha\notin J$ and $C>0$, $D_J(\alpha)/C$ is a nonzero analytic unit. Multiplying $\eta$ by this unit, which changes neither $\mathcal I_1$ nor any two-sided estimate, therefore makes $M_\alpha-1=\eta$ and $W_\alpha=1$ at the incidence. The vectors $V,W$ are linearly independent there. Their Gram matrix is continuous and positive definite, so on a smaller compact neighbourhood \begin{equation*} \|\Xi\|^2\ge c_3(A^2+x^2\eta^2). \tag{B.4.24n}\label{eq:B.4.24n} \end{equation*}

At two $+1$ losses define $$ g_1=x\eta_\alpha,\qquad g_2=y\eta_\beta. $$ The equality of symbolic and ordinary squares gives \begin{equation*} \|g(Y)\|_{\rm coeff}\le C_4(g_1^2+g_2^2). \tag{B.4.24o}\label{eq:B.4.24o} \end{equation*} To prove chord coercivity uniformly when $x/y$ is arbitrarily small or large, let $$ M_\alpha=\frac{Y_\alpha^+}{Y_\alpha},\qquad M_\beta=\frac{Y_\beta^+}{Y_\beta}, $$ using their analytic extensions to $x=0$ or $y=0$. The map is even in $x,y$. Since $M_\alpha-1$ vanishes when $\eta_\alpha=\eta_\beta=0$, Hadamard division gives $$ M_\alpha-1=a_{11}\eta_\alpha+b_{12}\eta_\beta. $$ On $y=0$, $\eta_\alpha=0$ is the exact five-node fixed component for every $\eta_\beta$. Hence $b_{12}$ vanishes at $y=0$, and evenness makes it divisible by $y^2$: $b_{12}=y^2a_{12}$. Exchanging the nodes gives $$ M_\beta-1=x^2a_{21}\eta_\alpha+a_{22}\eta_\beta. $$ The one-boundary computation \eqref{eq:B.4.24m1}, applied on the two adjacent faces, shows that $a_{11}$ and $a_{22}$ are nonzero analytic units at the four-node incidence. We may rescale the normal coordinates by these units. Multiplying by the signed amplitudes gives the exact matrix identity \begin{equation*} \binom{\Xi_\alpha}{\Xi_\beta} =\mathbf M\binom{g_1}{g_2},\qquad \mathbf M=\begin{pmatrix} a_{11}&xy\,a_{12}\\ xy\,a_{21}&a_{22} \end{pmatrix}. \tag{B.4.24p}\label{eq:B.4.24p} \end{equation*} The two leading signed-coordinate vectors are the columns $$ v_1=(a_{11},xy a_{21})^T,\qquad v_2=(xy a_{12},a_{22})^T, $$ whose Gram matrix is \begin{equation*} \begin{pmatrix} a_{11}^2+x^2y^2a_{21}^2 & xy(a_{11}a_{12}+a_{21}a_{22})\\ xy(a_{11}a_{12}+a_{21}a_{22}) & a_{22}^2+x^2y^2a_{12}^2 \end{pmatrix}. \tag{B.4.24p1}\label{eq:B.4.24p1} \end{equation*} At $x=y=0$ it is the positive diagonal matrix $\operatorname{diag}(a_{11}^2,a_{22}^2)$. Shrink the neighbourhood so that both diagonal units have modulus at least $d>0$ and the off-diagonal part of $\mathbf M$ has operator norm at most $d/2$. Then the least singular value of $\mathbf M$ is at least $d/2$, equivalently the least eigenvalue of \eqref{eq:B.4.24p1} is at least $d^2/4$. Hence \begin{equation*} \|GY-Y\|^2\ge \|\bigl(\Xi_\alpha,\Xi_\beta\bigr)\|^2 \ge\frac{d^2}{4} (x^2\eta_\alpha^2+y^2\eta_\beta^2). \tag{B.4.24q}\label{eq:B.4.24q} \end{equation*} Only the product $|xy|$ enters the perturbation, so this bound is uniform as $x^2/y^2\to0$ or $+\infty$.

The symbolic-square assertion in \eqref{eq:B.4.24o} holds on all four components. At the four-node prime $(x,y)$, completion coefficients are even in $x,y$, so their first derivatives vanish. At $(y,\eta_\alpha)$ and $(x,\eta_\beta)$, the absent-amplitude derivative vanishes and the other normal returns as the support-five $\theta$-shear, fixing the remaining completion root. At $(\eta_\alpha,\eta_\beta)$, the six-node return is tangent and its completion is $D$. Hence the coefficients lie in the localised square at every minimal prime, hence in $\mathcal I_2^{(2)}=\mathcal I_2^2$. A common off-node factor collision changes the bounded coefficients in \eqref{eq:B.4.24p}, but not its invertible diagonal limit; Lemma~\ref{lem:4.3} and compactness keep the constants uniform.

For reference, the finite signed-chart ledger used in this proof is as follows. The displayed variables are the normal variables together with the factor variables needed at collisions; all unlisted analytic variables are tangential, and $\Xi=GY-Y$.
\begingroup
\small
\renewcommand{\arraystretch}{1.18}
\begin{center}
\begin{tabular}{@{}>{\raggedright\arraybackslash}p{0.18\textwidth}>{\raggedright\arraybackslash}p{0.22\textwidth}>{\raggedright\arraybackslash}p{0.26\textwidth}>{\raggedright\arraybackslash}p{0.24\textwidth}@{}}
\toprule
Chart and ambient variables & Fixed ideal and minimal components & Units and first-order return & Drift and chord control \\
\midrule
Six nodes: $(p_0,p_1,p_2,k_0,k_1)$; use $(z_1,z_2,z_3,k_0,k_1)$ at a collision
& $\mathcal I=(U_0,U_1)$; the factor-allocation primes are those in \eqref{eq:B.4.17b1}--\eqref{eq:B.4.17b2}
& $D$ is a unit modulo $K$, the coefficients in \eqref{eq:B.4.17h} are positive, and $dK^+-dK=-2\operatorname{rem}_K(D\,dU)/C$
& $g\in\mathcal I^2$ and $\|\Xi\|^2\ge c\|U\|^2$ by \eqref{eq:B.4.24b}--\eqref{eq:B.4.24d} \\
Five nodes: $(\mathfrak A,\theta,\xi,\vartheta)$
& $(\mathfrak A)$ on each factor-allocation component
& $B$ and every retained $\lambda_i-\theta$ are units; $\delta=\mathfrak A/B$, and the first return is the tangential $\theta$-shear \eqref{eq:B.4.24i}
& $g=O(\mathfrak A^2)$ and $\|\Xi\|^2\ge c\mathfrak A^2$ by \eqref{eq:B.4.24g}--\eqref{eq:B.4.24j} \\
One $+1$ loss: $(x,A,\eta,\vartheta)$
& $\mathcal I_1=(A,x\eta)$, with primes $(A,x)$ and $(A,\eta)$
& $D_J(\alpha)/C$ is a unit; after its absorption, $\Xi=AV+x\eta W$ with $V,W$ independent
& $g\in\mathcal I_1^2$ and $\|\Xi\|^2\ge c(A^2+x^2\eta^2)$ by \eqref{eq:B.4.24k}--\eqref{eq:B.4.24n} \\
Two $+1$ losses: $(x,y,\eta_\alpha,\eta_\beta,\vartheta)$
& $\mathcal I_2=(x\eta_\alpha,y\eta_\beta)$, with the four primes in \eqref{eq:B.4.21b}
& $a_{11},a_{22}$ are units and the matrix $\mathbf M$ in \eqref{eq:B.4.24p} has uniformly invertible diagonal limit
& $g\in\mathcal I_2^2$ and $\|\Xi\|^2\ge c(x^2\eta_\alpha^2+y^2\eta_\beta^2)$ by \eqref{eq:B.4.24o}--\eqref{eq:B.4.24q} \\
\bottomrule
\end{tabular}
\end{center}
\endgroup

On every overlap, the displayed generator changes are analytic with unit determinant. After the charts are shrunk, those changes and their inverses have uniformly bounded derivatives, so the drift and chord estimates in the last column transfer in both directions with uniform constants, including on the simultaneous-collision branches.

If an omitted signed multiplier $m_0\ne+1$, its chord is $$ (m_0-1)x+O(x\,\operatorname{dist}(Y,\operatorname{Fix}G)), $$ and hence has modulus at least $c|x|$ in a small neighbourhood. Every retained polynomial and completion coefficient is even in $x$, so its change is $O(x^2)$, controlled by the squared chord. If a positive active coordinate is deleted by a block polynomial, the point is not fixed; continuity gives $\|GY-Y\|\ge\delta>0$ on a small neighbourhood while the completion drift is bounded.

All constants above are continuous local bounds. The relevant closed union of the collision and boundary strata is compact and has finitely many node patterns. A finite subcover, the largest upper constant, and the smallest positive lower constant prove \eqref{eq:B.4.24}.
\end{proof}

\section{Convergence of the completion polynomial and external modes}\label{sec:5}

For the cubic orbit, the degree-independent factor estimate produces a single limiting completion polynomial and hence at most six limiting spectral nodes. The estimates of Section~\ref{sec:4} control motion near the two-dimensional $(\rho,\sigma)$ fixed rectangles at fixed limiting factor data and node allocation, uniformly across simultaneous factor collisions and support-loss strata. The limiting multiplier of the parity return then classifies every remaining external coordinate as stable, unstable, or neutral.

From this point onward, a tilde distinguishes an object formed from the exact full spectral state. Thus set
\begin{equation*} \widetilde P_n=P_{Y_n},\qquad \widetilde Q_n=P_{Z_n},\qquad
 \widetilde a_n=\sigma(Y_n),\qquad \widetilde b_n=\sigma(Z_n), \end{equation*}
and define \begin{equation*} \widetilde{\mathscr D}_n=\widetilde P_n\widetilde Q_n-\widetilde a_n\widetilde b_n. \end{equation*}
Untilded orbit quantities introduced after \eqref{eq:B.5.19} are formed from the normalised $S_*$-core and its restricted one-block partner. Starred quantities always denote the common limiting fixed data.

\begin{lemma}[Convergence of the completion polynomial]\label{lem:5.1}
There are a monic sextic $\mathscr D_*$, monic cubics $P_*,Q_*$, and $C_*=\tau^2>0$ such that
\begin{equation*} \widetilde{\mathscr D}_n\longrightarrow\mathscr D_*, \qquad \widetilde P_n\longrightarrow P_*, \qquad \widetilde Q_n\longrightarrow Q_*, \qquad P_*Q_*=\mathscr D_*+C_*. \tag{B.5.2}\label{eq:B.5.2} \end{equation*}
Moreover,
\begin{equation*} \|\widetilde{\mathscr D}_n-\mathscr D_*\|_{\rm coeff} +|\widetilde a_n\widetilde b_n-C_*| +\|\widetilde P_n-P_*\|_{\rm coeff}+\|\widetilde Q_n-Q_*\|_{\rm coeff} \le C(\tau-\widetilde a_n). \tag{B.5.3}\label{eq:B.5.3} \end{equation*}
If
\begin{equation*} S_*=\{\lambda_i:\mathscr D_*(\lambda_i)=0\}, \end{equation*}
then
\begin{equation*} 4\le |S_*|\le6, \tag{B.5.5}\label{eq:B.5.5} \end{equation*}
every $\omega$-limit support is contained in $S_*$, and
\begin{equation*} \sum_{n=0}^{\infty}\sum_{\lambda_j\notin S_*}Y_{n,j}^{\,2}<\infty. \tag{B.5.6}\label{eq:B.5.6} \end{equation*}
\end{lemma}

\begin{proof}
In the notation of the common core, \begin{equation*} \widetilde P_n=P_{2n},\qquad \widetilde Q_n=P_{2n+1},\qquad H_{2n}=\widetilde a_n^2,\qquad H_{2n+1}=\widetilde b_n^2, \qquad H_\infty:=\tau^2=C_*. \end{equation*} Proposition~\ref{prop:common-factor-convergence} gives monic cubics $P_*,Q_*$ with $\widetilde P_n\to P_*$ and $\widetilde Q_n\to Q_*$. More precisely, summing \eqref{eq:common-factor-step} separately on the two parities gives \begin{equation*} \begin{aligned} \|\widetilde P_n-P_*\|_{\rm coeff} &\le C\sum_{j=n}^{\infty}(H_{2j+1}-H_{2j}) \le C(H_\infty-H_{2n}),\\ \|\widetilde Q_n-Q_*\|_{\rm coeff} &\le C\sum_{j=n}^{\infty}(H_{2j+2}-H_{2j+1}) \le C(H_\infty-H_{2n+1}). \end{aligned} \tag{B.5.7a}\label{eq:B.5.7a} \end{equation*} Since $\widetilde a_n\le \widetilde b_n\le\tau$, \begin{equation*} H_\infty-H_{2n}=\tau^2-\widetilde a_n^2\le2\tau(\tau-\widetilde a_n), \qquad H_\infty-H_{2n+1}\le H_\infty-H_{2n}, \end{equation*} so both factor tails in \eqref{eq:B.5.3} have the required bound. The same interlacing gives \begin{equation*} 0\le C_*-\widetilde a_n\widetilde b_n =\tau(\tau-\widetilde a_n)+\widetilde a_n(\tau-\widetilde b_n) \le2\tau(\tau-\widetilde a_n). \tag{B.5.7c}\label{eq:B.5.7c} \end{equation*} Define \begin{equation*} \mathscr D_*=P_*Q_*-C_*. \end{equation*} The factor coefficients are uniformly bounded, so the product difference identity $$ \widetilde P_n\widetilde Q_n-P_*Q_*=(\widetilde P_n-P_*)\widetilde Q_n+P_*(\widetilde Q_n-Q_*) $$ together with \eqref{eq:B.5.7a}--\eqref{eq:B.5.7c} proves all the convergences in \eqref{eq:B.5.2} and the quantitative completion tail \begin{equation*} \|\widetilde{\mathscr D}_n-\mathscr D_*\|_{\rm coeff} \le C(\tau-\widetilde a_n). \tag{B.5.8}\label{eq:B.5.8} \end{equation*}

Let $Y$ be an $\omega$-limit point with support $J$. The fixed-cycle identity of Lemma~\ref{lem:4.1}, or equivalently \eqref{eq:common-resonant-set}, gives \begin{equation*} \mathscr D_*=D_J\Lambda_J. \end{equation*} Every active node of an $\omega$-limit state is a root of $\mathscr D_*$. Lemma~\ref{lem:3.2} supplies at least four such nodes, while the nonzero monic sextic $\mathscr D_*$ has at most six distinct roots. This establishes \eqref{eq:B.5.5} and the containment.

For an index $j$ with $\lambda_j\notin S_*$, the exact signed two-block coordinate multiplier is
\begin{equation*} \widetilde M_{j,n} =\frac{\widetilde P_n(\lambda_j)\widetilde Q_n(\lambda_j)}{\widetilde a_n\widetilde b_n} =1+\frac{\widetilde{\mathscr D}_n(\lambda_j)}{\widetilde a_n\widetilde b_n}. \tag{B.5.10}\label{eq:B.5.10} \end{equation*}
Its distance from $1$ is eventually bounded below. Since
\begin{equation*} Y_{n+1,j}-Y_{n,j}=(\widetilde M_{j,n}-1)Y_{n,j}. \tag{B.5.11}\label{eq:B.5.11} \end{equation*}
the coordinate form \eqref{eq:common-coordinate-chord} of the summable same-parity chord, or directly \eqref{eq:B.3.5}, implies \eqref{eq:B.5.6}.
\end{proof}

Once $\mathscr D_*$ is fixed, each off-completion node is governed to first order by a scalar multiplier. The resulting trichotomy determines which external coordinates can influence the subsequent dynamics.

\begin{lemma}[Stable, unstable, and persistent neutral external modes]\label{lem:5.2}
For an index $j$ with $\lambda_j\notin S_*$, define
\begin{equation*} m_j=1+\frac{\mathscr D_*(\lambda_j)}{C_*}\ne1, \qquad E_{j,n}=Y_{n,j}^{\,2}, \qquad X_n=\tau-\widetilde a_n. \end{equation*} When $h=\lambda_j$ is used as a node label, abbreviate $m_h:=m_j$, $E_{h,n}:=E_{j,n}$, and $\widetilde M_{h,n}:=\widetilde M_{j,n}$.

Exactly one of the following holds.
1. If $|m_j|<1$, then $E_{j,n}$ decays exponentially. 2. If $|m_j|>1$, the coordinate becomes exactly zero after finitely many two-block returns. 3. If $|m_j|=1$, then $m_j=-1$. If the coordinate is never deleted, then \begin{equation*} E_{j,n}=O(X_n^2). \tag{B.5.13}\label{eq:B.5.13} \end{equation*}

If at least one neutral coordinate persists, the total stable weight satisfies
\begin{equation*} E_{{\rm st},n}=\sum_{\substack{\lambda_j\notin S_*\\ |m_j|<1}}E_{j,n},\qquad E_{{\rm neu},n}=\sum_{\substack{\lambda_j\notin S_*\\ |m_j|=1}}E_{j,n},\qquad E_{{\rm st},n}=o(E_{{\rm neu},n}). \tag{B.5.14}\label{eq:B.5.14} \end{equation*}
\end{lemma}

\begin{proof}
Equation \eqref{eq:B.5.10} gives $\widetilde M_{j,n}\to m_j$ and
\begin{equation*} E_{j,n+1}=\widetilde M_{j,n}^{\,2}E_{j,n}. \tag{B.5.15}\label{eq:B.5.15} \end{equation*}
This gives exponential decay in the stable case. By \eqref{eq:B.5.6}, $E_{j,n}\to0$. If $|m_j|>1$ and the coordinate were never zero, then \eqref{eq:B.5.15} would eventually make it grow geometrically, a contradiction. Hence some block polynomial vanishes at $\lambda_j$, after which the coordinate is permanently zero by Lemma~\ref{lem:2.1}. Since $m_j\ne1$, neutrality $|m_j|=1$ means $m_j=-1$.

For a neutral index $j$, define
$$ \delta_{j,n}=\widetilde{\mathscr D}_n(\lambda_j)+2\widetilde a_n\widetilde b_n. $$
Equations \eqref{eq:B.5.8} and $\widetilde a_n,\widetilde b_n\to\tau$, together with $\widetilde a_n\le \widetilde b_n\le\tau$, give
\begin{equation*} |\delta_{j,n}|\le C X_n,\qquad \widetilde M_{j,n}=-1+\frac{\delta_{j,n}}{\widetilde a_n\widetilde b_n}. \end{equation*}
Hence
\begin{equation*} \left|\log\frac{E_{j,n+1}}{E_{j,n}}\right|\le C X_n \tag{B.5.17}\label{eq:B.5.17} \end{equation*}
whenever the coordinate is nonzero and $n$ is large. On the other hand, the $j$th coordinate of \eqref{eq:B.5.11}, with $\widetilde M_{j,n}\to-1$, and \eqref{eq:B.3.4} imply
\begin{equation*} X_n-X_{n+1}=\widetilde a_{n+1}-\widetilde a_n \ge c\sum_{|m_g|=1}E_{g,n}. \tag{B.5.18}\label{eq:B.5.18} \end{equation*}

Fix a persistent $j$, and let $m>n$ be the first index for which $E_{j,m}\le e^{-1}E_{j,n}$; it exists by \eqref{eq:B.5.6}. Since $X_k\le X_n$, \eqref{eq:B.5.17} gives $1\le C(m-n)X_n$. Until that index, \eqref{eq:B.5.18} gives
$$ X_n\ge c\sum_{k=n}^{m-1}E_{j,k} \ge ce^{-1}(m-n)E_{j,n}. $$
Eliminating $m-n$ proves \eqref{eq:B.5.13}. There are finitely many external nodes, so the same estimate holds for their sum.

Finally, $E_{{\rm st},n}\le Cq^n$ for some $q<1$. A persistent neutral weight has successive ratio tending to one, so for every $q<r<1$ it is eventually bounded below by a positive multiple of $r^n$. This is \eqref{eq:B.5.14}.
\end{proof}

The neutral regime is studied on the normalised $S_*$-core. The next lemma compares its radial and intrinsic normal coordinates with the parity-energy tail $X_n$.

\begin{lemma}[Normalisation of the limiting core]\label{lem:5.3}
Suppose that at least one neutral external mode persists, and pass beyond the last deletion of an unstable mode. Let $\Pi_{S_*}$ denote coordinate projection onto the $S_*$-coordinates, and normalise the projection of $Y_n$ onto $S_*$:
\begin{equation*} \widehat Y_n=\frac{\Pi_{S_*}Y_n}{\|\Pi_{S_*}Y_n\|}, \tag{B.5.19}\label{eq:B.5.19} \end{equation*}
put $\widehat Z_n=L(\widehat Y_n)$, and let
$$ P_n=P_{\widehat Y_n},\qquad Q_n=P_{\widehat Z_n},\qquad
 a_n=\sigma(\widehat Y_n),\qquad b_n=\sigma(\widehat Z_n). $$
Thus $Q_n$ and $b_n$ are obtained from the restricted core block, not by silently replacing the exact full partner by its projection. Let
\begin{equation*} C_n=a_nb_n,\qquad \mathscr D_n=P_nQ_n-C_n,\qquad x_n=C_n-C_*. \end{equation*}
and let $U_n$ be any of the affine normal coordinates for the completion polynomial used in Sections~\ref{sec:4},~\ref{sec:6},~\ref{sec:8}, and~\ref{sec:9}. On a nonstationary tail on which a neutral mode persists,
\begin{equation*} x_n<0,\qquad cX_n\le -x_n\le CX_n,\qquad |x_n|+\|U_n\|\le CX_n. \end{equation*}

More generally, write $\mathcal E=\sum_hE_h$, reconstruct a nearby full state from a core signed-chart coordinate $z$ by
\begin{equation*} \widetilde\mu(z,E)=(1-\mathcal E)\mu(z)+\sum_hE_h\delta_h, \tag{B.5.24}\label{eq:B.5.24} \end{equation*}
and form $\widetilde P,\widetilde Q,\widetilde C,\widetilde{\mathscr D},\widetilde U$ from its two exact blocks. On every compact regular signed chart used in Sections~\ref{sec:8}--\ref{sec:9}, including its support incidences and common-factor collisions,
\begin{equation*} \begin{aligned}
 &\max_{|\gamma|\le2}\bigl\|\partial_z^\gamma(\widetilde P-P)\bigr\|_{\rm coeff}
 +\max_{|\gamma|\le2}\bigl\|\partial_z^\gamma(\widetilde Q-Q)\bigr\|_{\rm coeff}
 +\max_{|\gamma|\le2}\left|\partial_z^\gamma(\widetilde C-C)\right|\\
 &\qquad +\max_{|\gamma|\le2}\bigl\|\partial_z^\gamma(\widetilde{\mathscr D}-\mathscr D)\bigr\|_{\rm coeff}
 +\max_{|\gamma|\le2}\bigl\|\partial_z^\gamma(\widetilde U-U)\bigr\|
 \le C_{\rm fc}\mathcal E.
\end{aligned} \tag{B.5.25}\label{eq:B.5.25} \end{equation*}
The same estimate holds for the normalised retained core after either one block or one parity return, compared with the corresponding restricted-core map. Derivatives involving an $E_h$ are uniformly bounded, but are not asserted to be $O(\mathcal E)$; their generally nonzero values at $E=0$ are the Uvarov derivatives used below.

If $W_i$ are the full retained weights, $S=\sum_{i\in S_*}W_i=1-\mathcal E$, and $w_i=W_i/S$, the exact full signed multiplier and the induced normalised-core recurrence are
\begin{equation*} \widetilde M(t)=\frac{\widetilde P(t)\widetilde Q(t)}{\widetilde C}
 =1+\frac{\widetilde{\mathscr D}(t)}{\widetilde C},\qquad
 E_h^+=E_h\widetilde M(h)^2,\qquad
 \frac{w_i^+}{w_i}=\widetilde M(\lambda_i)^2\frac S{S^+}. \tag{B.5.26}\label{eq:B.5.26} \end{equation*}
In particular, along the orbit,
\begin{equation*} \widetilde C_n-C_n=O(\mathcal E_n),\qquad
 \|\widetilde U_n-U_n\|=O(\mathcal E_n),\qquad
 \mathcal E_n=O(X_n^2). \tag{B.5.27}\label{eq:B.5.27} \end{equation*}
\end{lemma}

\begin{proof}
For a probability measure $\nu$, write $m_j(\nu)=\int t^j\,d\nu$. If
$$ P_\nu(t)=t^3+p_2t^2+p_1t+p_0, $$
then its lower coefficient vector is the unique solution of
\begin{equation*}
 \bigl(m_{i+j}(\nu)\bigr)_{0\le i,j\le2}(p_0,p_1,p_2)^T
 =-\bigl(m_3(\nu),m_4(\nu),m_5(\nu)\bigr)^T. \tag{B.5.28}\label{eq:B.5.28}
\end{equation*}
Every core on the compact cover has at least four active distinct nodes, so this moment Gram matrix is positive definite. Its least eigenvalue has a positive lower bound on a sufficiently small finite cover. The difference between the measures in \eqref{eq:B.5.24} has total variation at most $2\mathcal E$; hence all moments through degree six differ by $O(\mathcal E)$. Equation \eqref{eq:B.5.28}, the inverse-matrix identity, and compactness prove the estimate for $\widetilde P-P$.

The squared monic norm is a polynomial in the moments and the coefficients of $P$. It is bounded away from zero on the same cover. Consequently the one-block weight map, its normalisation on $S_*$, the partner polynomial $Q$, and the second norm are analytic functions of $(z,E)$. Repeating the preceding argument after the first block proves the assertions for $\widetilde Q,\widetilde C,\widetilde{\mathscr D}$ and for the returned normalised core after one and two blocks. The coefficient divisions and signed-chart changes defining $U$ are analytic on the regular cover. They remain so at a factor collision because the moment systems and the polynomial coefficient maps do not divide by a factor separation; at a support incidence the remaining core Gram matrix is still uniformly positive definite.

For any one of these analytic functions $F$, the multivariable Hadamard formula gives
$$ F(z,E)-F(z,0)=\sum_hE_h\int_0^1
 \partial_{E_h}F(z,E_1,\ldots,E_{h-1},sE_h,0,\ldots,0)\,ds. $$
Differentiating this identity at most twice in $z$ and taking the finite-cover maximum proves \eqref{eq:B.5.25}. It also explains why external derivatives are bounded rather than $O(\mathcal E)$. Formula \eqref{eq:B.5.26} follows directly from the exact full two-block coordinate recurrence and the identity $w_i=W_i/S$.

Every persistent neutral weight is $O(X_n^2)$ by Lemma~\ref{lem:5.2}, unstable weights have been deleted, and the total stable weight is smaller by \eqref{eq:B.5.14}. Thus $\mathcal E_n=O(X_n^2)$, and \eqref{eq:B.5.25} gives
\begin{equation*} |a_n-\widetilde a_n|+|b_n-\widetilde b_n|=O(X_n^2). \tag{B.5.22}\label{eq:B.5.22} \end{equation*}
The completion coefficients determine the normal uniformly, including at collisions. The factor identity $$ KPQ=D(K+U)+H_1K_1 $$ and the division $DU=Kq+r$ imply that the quadratic $r+H_1K_1$ is divisible by the monic quadratic $K$, hence equals $H_1K$. Thus \begin{equation*} PQ=D+\operatorname{quo}_K(DU)+H_1. \tag{B.5.22a}\label{eq:B.5.22a} \end{equation*} If $U=u_0+u_1t$, the coefficients of degrees five and four of $\operatorname{quo}_K(DU)$ are $u_1$ and $u_0$ plus a bounded multiple of $u_1$. This triangular recovery is uniform on the closed parameter rectangle, including along the collision locus. On a five-node face, the same polynomial division is its exact quotient at the support boundary and the bracket coordinates \eqref{eq:B.4.19} give the identical uniform conclusion. In particular,
\begin{equation*} a_n-\tau=-X_n+O(X_n^2). \tag{B.5.23}\label{eq:B.5.23} \end{equation*}
Since $\widetilde b_n\le\tau$, \eqref{eq:B.5.22}--\eqref{eq:B.5.23} imply
$$ C_n-C_* =-\tau\{(\tau-a_n)+(\tau-b_n)\} +O(X_n^2), $$
which is negative and has magnitude between positive multiples of $X_n$ on a sufficiently late nonstationary tail. The completion tail bound \eqref{eq:B.5.8}, transferred to the normalised core using \eqref{eq:B.5.25}, gives $\|U_n\|=O(X_n)$. If $X_n=0$ at a finite index, \eqref{eq:B.3.3} forces equality through the cross-projection chain, so the parity state is already fixed and the conclusion follows directly.
\end{proof}

\section{Intrinsic convergence and removal of stable modes}\label{sec:6}

Two reductions separate the intrinsic dynamics from the external forcing. First, an orbit confined to the six-node limiting face converges to one fixed point, uniformly through collisions and asymptotic support loss. Second, summable stable coordinates can be removed by shadowing. The same estimates show that stable forcing is negligible whenever neutral modes persist.

\begin{lemma}[Convergence near the limiting six-node fixed set]\label{lem:6.1}
Fix six nodes $S=\{\lambda_1<\cdots<\lambda_6\}$. Suppose that on a tail all six coordinates are positive, every $\omega$-limit point lies on the closed fixed set
\begin{equation*} P_*Q_*=D_S+C_*, \end{equation*}
and no external forcing is present. Then the squared parity state converges to one point. The conclusion remains valid when the limiting cubics have any simple or simultaneous off-node common roots, and when one or two core weights tend to zero without becoming exactly zero.
\end{lemma}

\begin{proof}
Let $G=L^2$, abbreviate $D=D_S$, and let $C=H_0$ denote the current squared monic norm; on the limiting fixed set, $C=C_*$. After the last exact deletion, let $$ K(t)=(t-\rho)(t-\sigma),\qquad (\rho,\sigma)\in [\lambda_2,\lambda_3]\times[\lambda_4,\lambda_5], $$ and define $$ z_\rho=U(\rho),\qquad z_\sigma=U(\sigma), $$ $$ d_\rho=(\rho-\lambda_2)(\lambda_3-\rho),\qquad d_\sigma=(\sigma-\lambda_4)(\lambda_5-\sigma). $$ The root gap is uniformly separated: $$ \sigma-\rho\ge\lambda_4-\lambda_3>0. $$ If $U=U_0+U_1t$, then $$ \binom{U_0}{U_1} = \frac1{\sigma-\rho} \begin{pmatrix}\sigma&-\rho\\-1&1\end{pmatrix} \binom{z_\rho}{z_\sigma}. $$ Hence $(z_\rho,z_\sigma)$ is a uniformly invertible set of generators for the intrinsic fixed ideal \begin{equation*} \mathfrak n=(U_0,U_1). \end{equation*} This remains true at a simultaneous factor collision because it uses only the separated roots of $K$, not a local factor-splitting branch of $P,Q$.

Let $$ F_\rho=U^+(\rho^+),\qquad F_\sigma=U^+(\sigma^+). $$ Lemmas~\ref{lem:4.2}--\ref{lem:4.3} give the exact ideal-valued jets \begin{equation*} U^+-U\in\mathfrak n^2,\qquad K^+-K+\frac2C\operatorname{rem}_K(DU)\in\mathfrak n^2. \tag{B.6.17}\label{eq:B.6.17} \end{equation*} The roots of $K^+$ are analytic because they remain separated, and the second relation makes their increments elements of $\mathfrak n$. Since $U_1^+\in\mathfrak n$, $$ U^+(\rho^+)-U^+(\rho) =U_1^+(\rho^+-\rho)\in\mathfrak n^2, $$ with the corresponding identity at $\sigma$. Hence \begin{equation*} F_\rho-z_\rho\in\mathfrak n^2,\qquad F_\sigma-z_\sigma\in\mathfrak n^2. \tag{B.6.18}\label{eq:B.6.18} \end{equation*} Because $z_\rho,z_\sigma$ form a regular sequence, there are analytic coefficients for which \begin{equation*} F_\rho=z_\rho\widetilde R_\rho+z_\sigma^2c_\rho,\qquad \widetilde R_\rho|_{\mathfrak n=0}=1. \tag{B.6.19}\label{eq:B.6.19} \end{equation*} An element of $\mathfrak n^2$ is $z_\rho^2a+z_\rho z_\sigma b+z_\sigma^2c_\rho$; the first two terms are absorbed into $\widetilde R_\rho=1+z_\rho a+z_\sigma b$.

We divide the last coefficient uniformly at the two horizontal faces. Let $\beta\in\{\lambda_4,\lambda_5\}$, and let $y$ be the signed amplitude of the omitted coordinate. The exact recovery identity gives \begin{equation*} y^2=u_\beta d_\sigma \tag{B.6.20}\label{eq:B.6.20} \end{equation*} with $u_\beta>0$ an analytic unit. On $y=0$, $$ \bar K_\beta(t)=\frac{K(t)-K(\beta)}{t-\beta}=t-\rho,\qquad [U,\bar K_\beta]=-U(\rho)=-z_\rho. $$ The lower-support fixed component has ideal $(y,z_\rho)$, while $z_\sigma$ is its free shear coordinate. Its quotient is reduced, and $z_\sigma$ belongs to no minimal prime. Hence it is a non-zero-divisor and \begin{equation*} ((y,z_\rho):z_\sigma^2)=(y,z_\rho). \tag{B.6.21}\label{eq:B.6.21} \end{equation*} Setting $y=z_\rho=0$ in \eqref{eq:B.6.19}, exact fixedness gives $z_\sigma^2c_\rho=0$. Equation \eqref{eq:B.6.21} yields $c_\rho\in(y,z_\rho)$. Every retained quantity is even in $y$; expanding an even representative gives $$ c_\rho=z_\rho c_{\rho,1}+y^2c_{\rho,2} =z_\rho c_{\rho,1}+d_\sigma q_\rho. $$ Absorb $z_\rho z_\sigma^2c_{\rho,1}$ into $\widetilde R_\rho$. Hence \begin{equation*} F_\rho=z_\rho R_\rho+d_\sigma z_\sigma^2q_\rho. \tag{B.6.22}\label{eq:B.6.22} \end{equation*} Exchanging $$ (\rho,\lambda_2,\lambda_3,z_\rho) \longleftrightarrow (\sigma,\lambda_4,\lambda_5,z_\sigma) $$ gives \begin{equation*} F_\sigma=z_\sigma R_\sigma+d_\rho z_\rho^2q_\sigma. \tag{B.6.23}\label{eq:B.6.23} \end{equation*}

We calculate all four signs. On $\mathfrak n=0$, $R_\rho=R_\sigma=1$. On the additional face $z_\rho=0$, $$ U(t)=\frac{z_\sigma}{\sigma-\rho}(t-\rho). $$ The quadratic tensor in Lemma~\ref{lem:4.2} gives \begin{equation*} q_\rho\big|_{\mathfrak n=0} =\frac{\Gamma_\rho}{(\sigma-\rho)^2d_\sigma}, \tag{B.6.24}\label{eq:B.6.24} \end{equation*} where \begin{equation*} \Gamma_\rho =-\frac{D(\sigma)}C \left\{\frac{P(\sigma)}{P(\rho)} +\frac{Q(\sigma)}{Q(\rho)}\right\}. \tag{B.6.25}\label{eq:B.6.25} \end{equation*} Interlacing gives $$ P(\rho),Q(\rho)>0,\qquad P(\sigma),Q(\sigma)<0, $$ and $D(\sigma)>0$; hence $\Gamma_\rho>0$. At either endpoint $\sigma=\lambda_4,\lambda_5$, $$ \frac{D(\sigma)}{d_\sigma} =-\!\!\prod_{j\notin\{4,5\}}(\sigma-\lambda_j) $$ has a finite strictly positive limit, while the bracket in \eqref{eq:B.6.25} has a finite strictly negative limit. Hence \eqref{eq:B.6.24} extends to a finite positive endpoint value. At a common-factor collision, the common root of $P,Q$ lies in a negative gap, whereas $\rho,\sigma$ lie in the complementary positive gaps; none of the four displayed evaluations vanishes, so the strict sign persists. Compactness gives, after shrinking the tube, \begin{equation*} R_\rho,R_\sigma\ge R_{\min}^{\rm core}>0,\qquad q_\rho,q_\sigma\ge q_{\min}>0. \tag{B.6.26}\label{eq:B.6.26} \end{equation*} It follows from \eqref{eq:B.6.22}--\eqref{eq:B.6.26} that \begin{equation*} z_\rho\ge0\Longrightarrow z_\rho^+\ge0,\qquad z_\sigma\ge0\Longrightarrow z_\sigma^+\ge0. \tag{B.6.27}\label{eq:B.6.27} \end{equation*} Each normal either stays negative or crosses once into its nonnegative half-line.

Evaluating the second identity in \eqref{eq:B.6.17} at $\rho$ gives $$ K^+(\rho)=-\frac{2D(\rho)}Cz_\rho+\mathfrak n^2. $$ The integral identity $$ 0=K^+(\rho^+) =K^+(\rho)+(\rho^+-\rho) \int_0^1(K^+)'(\rho+s(\rho^+-\rho))\,ds $$ has an integral equal to $\rho-\sigma\ne0$ on the fixed set. Division by this analytic unit gives \begin{equation*} \rho^+-\rho =\frac{2D(\rho)}{C(\rho-\sigma)}z_\rho+\mathfrak n^2. \tag{B.6.28}\label{eq:B.6.28} \end{equation*} Evaluation at $\sigma$ gives \begin{equation*} \sigma^+-\sigma =\frac{2D(\sigma)}{C(\sigma-\rho)}z_\sigma+\mathfrak n^2. \end{equation*}

Let $\gamma\in\{\lambda_2,\lambda_3\}$ be an endpoint of the $\rho$-interval and $x$ its signed amplitude. Recovery gives \begin{equation*} x^2=\frac{H_0c_\gamma K(\gamma)}{P(\gamma)} =u_\gamma d_\rho \end{equation*} with $u_\gamma>0$ analytic. On $x=0$, coordinate-factor invariance keeps that weight zero and the nearby returned root is exactly $\gamma$. The increment is even in $x$ and vanishes at $x=0$, so it is divisible by $x^2$, equivalently by $d_\rho$.

This division preserves the $\mathfrak n^2$ remainder. Lemma~\ref{lem:4.3} makes $\mathfrak n^2$ unmixed with the same minimal primes as $\mathfrak n$, and $x$ belongs to none of them. Hence \begin{equation*} (\mathfrak n^2:x^2)=\mathfrak n^2. \end{equation*} After dividing \eqref{eq:B.6.28} by $d_\rho$, denote its quadratic remainder by $z_\rho a_1+z_\sigma^2b_1$. At a horizontal endpoint, imposing $z_\rho=0$ gives the exact lower-support fixed component, so the increment vanishes. The colon identity \eqref{eq:B.6.21}, evenness in the horizontal amplitude, and \eqref{eq:B.6.20} give $$ b_1=z_\rho b_2+d_\sigma B_\rho. $$ Absorbing the first term into the linear coefficient yields \begin{equation*} \rho^+-\rho =d_\rho\left(z_\rho A_\rho +d_\sigma z_\sigma^2B_\rho\right). \tag{B.6.32}\label{eq:B.6.32} \end{equation*} The four endpoint divisions for $\sigma$ give \begin{equation*} \sigma^+-\sigma =d_\sigma\left(z_\sigma A_\sigma +d_\rho z_\rho^2B_\sigma\right). \tag{B.6.33}\label{eq:B.6.33} \end{equation*} On the fixed set, \begin{equation*} A_\rho=\frac{2D(\rho)} {C(\rho-\sigma)d_\rho}<0,\qquad A_\sigma=\frac{2D(\sigma)} {C(\sigma-\rho)d_\sigma}>0. \end{equation*} The quotients $D(\rho)/d_\rho$ and $D(\sigma)/d_\sigma$ extend as strictly positive functions to their closed intervals. Compactness of the explicit analytic coefficients gives \begin{equation*} A_\rho\le-c_\rho<0,\qquad A_\sigma\ge c_\sigma>0,\qquad |B_\rho|+|B_\sigma|\le B_0. \tag{B.6.35}\label{eq:B.6.35} \end{equation*}

We transfer the chord normal form of Lemma~\ref{lem:4.4} to these root coordinates. Near a corner $(\rho,\sigma)=(\gamma,\beta)$, let $x,y$ be the two signed omitted amplitudes. Then $$ x^2=u_\gamma d_\rho,\qquad y^2=u_\beta d_\sigma, $$ and the bracket generators are exactly \begin{equation*} \eta_\gamma =-z_\rho-\frac{\sigma-\beta}{\sigma-\rho} (z_\sigma-z_\rho), \quad \eta_\beta =-z_\sigma-\frac{\rho-\gamma}{\sigma-\rho} (z_\sigma-z_\rho). \tag{B.6.36}\label{eq:B.6.36} \end{equation*} The first off-diagonal coefficient is $O(d_\sigma)$, and the second is $O(d_\rho)$. Conjugating this two-by-two transformation by $\operatorname{diag}(\sqrt{d_\rho},\sqrt{d_\sigma})$ gives diagonal entries $-1+O(d_\rho+d_\sigma)$ and off-diagonal entries $O(\sqrt{d_\rho d_\sigma})$. It and its inverse are uniformly bounded near the corner. We have \begin{equation*} x^2\eta_\gamma^2+y^2\eta_\beta^2 \asymp d_\rho z_\rho^2+d_\sigma z_\sigma^2, \tag{B.6.37}\label{eq:B.6.37} \end{equation*} with constants independent of $d_\rho/d_\sigma$. On an open edge the corresponding one-boundary restriction of \eqref{eq:B.6.36} gives the equivalence, and away from all faces both $d$'s are bounded below. The chord-coercivity estimates \eqref{eq:B.4.24d}, \eqref{eq:B.4.24g}, \eqref{eq:B.4.24n}, and \eqref{eq:B.4.24q}, proved in Lemma~\ref{lem:4.4}, together with a finite cover now give \begin{equation*} \|Y^+-Y\|^2 \ge c\{d_\rho z_\rho^2+d_\sigma z_\sigma^2\} \tag{B.6.38}\label{eq:B.6.38} \end{equation*} through simple collisions, simultaneous collisions, edges, and corners. Lemma~\ref{lem:3.1} implies \begin{equation*} \sum_n\{d_{\rho,n}z_{\rho,n}^2+ d_{\sigma,n}z_{\sigma,n}^2\}<\infty. \tag{B.6.39}\label{eq:B.6.39} \end{equation*}

After the possible crossing in \eqref{eq:B.6.27}, suppose first that $z_\rho\ge0$. Equations \eqref{eq:B.6.32}, \eqref{eq:B.6.35} give $$ [\rho^+-\rho]_+ \le B_0d_\rho d_\sigma z_\sigma^2 \le C d_\sigma z_\sigma^2. $$ If $z_\rho<0$ forever, they give the identical bound for $[\rho-\rho^+]_+$. One directional variation is summable. Since $\rho_n$ is bounded, the other is summable as well: for $m>N$, $$ \sum_{n=N}^{m-1}[\rho_n-\rho_{n+1}]_+ = \sum_{n=N}^{m-1}[\rho_{n+1}-\rho_n]_+ -(\rho_m-\rho_N). $$ Hence \begin{equation*} \sum_n|\rho_{n+1}-\rho_n|<\infty. \tag{B.6.40}\label{eq:B.6.40} \end{equation*} For $\sigma$, \eqref{eq:B.6.33} makes the negative variation summable when $z_\sigma\ge0$, and the positive variation summable when $z_\sigma<0$. Hence \begin{equation*} \sum_n|\sigma_{n+1}-\sigma_n|<\infty. \end{equation*}

It follows that $K_n\to K_*$. Lemma~\ref{lem:5.1} gives $a_n\to\tau>0$ and $P_n\to P_*$. The recovery identity \begin{equation*} w_{n,i} =\frac{a_n^2c_iK_n(\lambda_i)}{P_n(\lambda_i)} \tag{B.6.42}\label{eq:B.6.42} \end{equation*} gives coordinatewise convergence; every denominator has nonzero limit because $P_*(\lambda_i)Q_*(\lambda_i)=C_*>0$. Equations \eqref{eq:B.6.37}--\eqref{eq:B.6.42} are uniform when one or two weights tend to zero at arbitrary relative rates. If a weight becomes exactly zero, the remaining tail lies on an invariant five- or four-node face; those lower-support cases are proved independently in Section~\ref{sec:7}.
\end{proof}

The intrinsic dynamics on the limiting face are therefore settled. The stable external coordinates are next removed by a shadowing construction, which also quantifies their effect when neutral modes remain.

\begin{lemma}[Shadowing after removal of stable modes]\label{lem:6.2}
Let a tail approach a compact fixed set determined by the completion polynomial. If every surviving external mode is stable, the full tail is exponentially close to an exact orbit of positive weights on the invariant coordinate face determined by the limiting nodes. If a neutral mode persists, stable forcing is $o(E_{{\rm neu},n})$. Deleting the stable weights changes an arbitrary quantity analytic on the signed charts by $O(\sum_g|\xi_g|)=O(E_{{\rm st},n}^{1/2})$, where the sum is over the stable signed amplitudes. For every retained quantity even in each deleted signed amplitude, equivalently analytic in the deleted squared masses, the sharper bound is $O(E_{{\rm st},n})$; this applies in particular to moments, inverse Gram matrices, orthogonal-polynomial coefficients, monic norms, completion coefficients, and returned core coordinates. In the persistent-neutral case both comparison errors are $o(E_{{\rm neu},n})$.
\end{lemma}

\begin{proof}
Work after the last exact deletion. First suppose that every surviving external mode is stable. This is the only case in which the shadowing orbit is constructed. Let $z$ be the vector of stable squared weights, and let $\chi$ contain the retained normalised-core variables and intrinsic completion coordinates. Thus, on $z=0$, the map $f$ has no retained neutral external coordinate. Coordinate-factor invariance and regularity of the degree-three moment Gram matrix give analytic equations \begin{equation*} \begin{aligned} \chi^+&=f(\chi)+a(\chi)z+R(\chi,z),\\ z^+&=M(\chi)z+S(\chi,z), \end{aligned} \tag{B.6.43}\label{eq:B.6.43} \end{equation*} where $$ R(\chi,z)=O(\|z\|^2),\qquad S(\chi,z)=O(\|z\|^2), $$ and $z=0$ is invariant. The finitely many stable squared multipliers have limiting modulus below one. Shrink the tube so that \begin{equation*} \|M(\chi)\|\le q_{\rm st}<1. \end{equation*}

The shadowing argument rests on the following weak-centre derivative bound. On a six-node fixed component, Lemma~\ref{lem:4.2} gives \begin{equation*} Df=I+N,\qquad N(\delta P,\delta K)\in\{0\}\oplus\mathcal P_1,\qquad N(0,\delta K)=0. \end{equation*} The first assertion is the coefficient identity $$ dP^+=dP,\qquad dK^+-dK=-\frac2{H_0}\operatorname{rem}_K(D\,dU), $$ and the second holds because a pure $K$-variation is tangent to the fixed rectangle and has $dU=0$. Hence $N^2=0$.

On a five-node interior component with support $J$, define $$ c_i=\frac1{D_J'(\lambda_i)},\qquad x_i=\lambda_i-\theta,\qquad \mathfrak A=\sum_i\frac{w_i}{x_i},\qquad B=\sum_i\frac{w_i}{x_i^2}. $$ The exact return is \begin{equation*} G(w)_i=w_i \frac{(1-\mathfrak A/(Bx_i))^2}{1-\mathfrak A^2/B}. \tag{B.6.46}\label{eq:B.6.46} \end{equation*} Differentiating at $\mathfrak A=0$ gives \begin{equation*} DG_w=I+V\otimes d\mathfrak A_w,\qquad V_i=-\frac{2w_i}{Bx_i}. \end{equation*} To verify nilpotence, define $$ \gamma_i=c_i^2,\qquad b_i=\frac{\gamma_i}{w_i},\qquad B_{\rm dual}=\sum_i b_i. $$ For a simplex tangent vector $h$, differentiation of the dual mean $$ \theta=\frac{\sum_i\gamma_i\lambda_i/w_i} {\sum_i\gamma_i/w_i} $$ gives \begin{equation*} d\theta[h]=-\frac1{B_{\rm dual}} \sum_i\frac{\gamma_ix_i}{w_i^2}h_i. \end{equation*} Substitution of $h=V$ yields $d\theta[V]=2/B$. Also $$ d\mathfrak A[h]=\sum_i\frac{h_i}{x_i}+B\,d\theta[h], $$ so \begin{equation*} d\mathfrak A[V] =-\frac2B\sum_i\frac{w_i}{x_i^2}+2=0. \end{equation*} Hence $(V\otimes d\mathfrak A)^2=0$: the nilpotent part sends the normal to the fixed tangent and kills that tangent.

On a four-node component the exact squared return is the identity. At an endpoint of a five-node interval, let $\omega$ be the fifth weight. We verify the endpoint derivative here. If the endpoint node is $\alpha=\lambda_h$, use the parametrisation $$ w_h=\omega,\qquad w_j=(1-\omega)v_j,\qquad \gamma_i=c_i^2,\qquad D_0=\sum_{j\ne h}\frac{\gamma_j(\lambda_j-\alpha)}{v_j}. $$ The dual-pivot formula derived before \eqref{eq:B.4.24f} gives exactly $$ \theta-\alpha =\frac{\omega D_0} {\gamma_h(1-\omega)+\omega\sum_{j\ne h}\gamma_j/v_j}. $$ At the fixed endpoint $D_0\ne0$, because the fixed recovery weight is a nonzero multiple of $\theta-\alpha$. With $x_i=\lambda_i-\theta$, $$ B=\frac{\gamma_h^2}{D_0^2\omega}+O(1),\qquad \delta=\frac{\mathfrak A}{B}=O(\omega\mathfrak A), $$ $$ \delta/x_h=O(\mathfrak A),\qquad \delta/x_j=O(\omega\mathfrak A)\quad(j\ne h). $$ Substitution in the exact dual-pivot quotient $$ \theta^+-\theta= \frac{\sum_i(\gamma_i/w_i)x_i^3/(x_i-\delta)^2} {\sum_i(\gamma_i/w_i)x_i^2/(x_i-\delta)^2} $$ and the identity $$ (1-r)^{-2}=1+2r+ r^2\frac{3-2r}{(1-r)^2} $$ give $\theta^+-\theta=2\delta\{1+O(|\mathfrak A|)\}$. Finally the elementary identity $$ \frac{(x-\delta)^2}{x^2(x-t)} =\frac1x+\frac{t-2\delta}{x^2} +\frac{(t-\delta)^2}{x^2(x-t)} $$ summed against $w_i$, with $t=\theta^+-\theta$, gives $$ \left(1-\frac{\mathfrak A^2}{B}\right)\mathfrak A^+ =\mathfrak A+(t-2\delta)B +(t-\delta)^2\sum_i\frac{w_i}{x_i^2(x_i-t)}. $$ The displayed endpoint bounds make the last two terms $O(\mathfrak A^2)$. Formula \eqref{eq:B.6.46} gives the corresponding fifth weight estimate. Hence \begin{equation*} \omega^+-\omega=O(\omega\mathfrak A),\qquad \mathfrak A^+-\mathfrak A=O(\mathfrak A^2). \end{equation*} The derivative is again the identity plus a map from the $\mathfrak A$-normal into the fixed tangent, and its tangent coefficient has a finite endpoint limit.

The derivative bound remains available without choosing subbundles across the support-loss strata. Work in the ambient smooth core charts on which the restricted return is analytic, and let $\mathscr F$ be the compact fixed set under consideration. For $p\in\mathscr F$, define $$ N_p=Df(p)-I:T_p\mathscr M\longrightarrow T_p\mathscr M. $$ The calculations above show $N_p^2=0$ on each regular six-node, five-node, and endpoint stratum; on a four-node component $N_p=0$. The same identity holds at support incidences and factor collisions. Indeed, in every signed core chart the entries of $N_p^2$ are analytic, and the regular coprime points at which the calculation applies are dense in each fixed component. Continuity therefore gives \begin{equation*} N_p^2=0\quad\hbox{and hence}\quad N_p^3=0 \qquad(p\in\mathscr F). \tag{B.6.51}\label{eq:B.6.51} \end{equation*} Thus the third level in the metric below is redundant but harmless.

Choose a smooth background inner product $\langle\cdot,\cdot\rangle_0$ on a neighbourhood of $\mathscr F$. For $0<\varepsilon<1$ and $p\in\mathscr F$ define \begin{equation*} \langle v,w\rangle_{\varepsilon,p} =\varepsilon^4\langle v,w\rangle_0 +\varepsilon^2\langle N_pv,N_pw\rangle_0 +\langle N_p^2v,N_p^2w\rangle_0. \tag{B.6.52}\label{eq:B.6.52} \end{equation*} This form is positive definite. With $$ J_{\varepsilon,p}v=(\varepsilon^2v,\varepsilon N_pv,N_p^2v) $$ and $S(u_0,u_1,u_2)=(u_1,u_2,0)$, equation \eqref{eq:B.6.51} gives $$ J_{\varepsilon,p}Df(p)=(I+\varepsilon S)J_{\varepsilon,p}, \quad J_{\varepsilon,p}Df(p)^{-1} =(I-\varepsilon S+\varepsilon^2S^2)J_{\varepsilon,p}. $$ The two operator norms on $\mathscr F$ are at most $1+\varepsilon$ and $1+\varepsilon+\varepsilon^2$, respectively.

The metric \eqref{eq:B.6.52} extends to a neighbourhood of the compact closed set $\mathscr F$: extend its matrix entries in local trivialisations, shrink the neighbourhood so that they remain positive definite, and patch them with a partition of unity. Choose $\kappa>1$ with $q_{\rm st}\kappa<1$, and then choose $\varepsilon$ so that $1+\varepsilon+\varepsilon^2<\kappa$. Continuity, compactness, and the inverse-function theorem give a smaller neighbourhood on which \begin{equation*} \|Df(\chi)\|_{\varepsilon,\chi\to f(\chi)},\quad \|Df(\chi)^{-1}\|_{\varepsilon,f(\chi)\to\chi} \le\kappa. \tag{B.6.53}\label{eq:B.6.53} \end{equation*} Fix $q_{\rm st}<q_{\rm sh}<1/\kappa$. After one further shrink, the second equation of \eqref{eq:B.6.43}, with $S(\chi,0)=D_zS(\chi,0)=0$, gives \begin{equation*} \|z_n\|\le q_{\rm sh}^{\,n-N}\|z_N\|. \end{equation*}

Using the exponential map of the extended metric as a local addition, seek the orbit on the invariant coordinate face in the form $$ \bar\chi_n=\exp_{\chi_n}(d_n),\qquad \bar\chi_{n+1}=f(\bar\chi_n). $$ For a sufficiently late tail define $$ \Psi_n(d)=\exp_{\chi_{n+1}}^{-1} \bigl(f(\exp_{\chi_n}d)\bigr),\qquad A_n=D\Psi_n(0). $$ Equations \eqref{eq:B.6.43} and Taylor's formula give \begin{equation*} d_{n+1}=A_nd_n-a(\chi_n)z_n+\mathcal R_n(d_n,z_n), \tag{B.6.55}\label{eq:B.6.55} \end{equation*} where the $O(\|z_n\|^2)$ difference between $\Psi_n(0)$ and $-a(\chi_n)z_n$ is included in $\mathcal R_n$. Since $\chi_{n+1}-f(\chi_n)=O(\|z_n\|)$, continuity and \eqref{eq:B.6.53} allow the tail to be chosen so that $$ \|A_n\|,\ \|A_n^{-1}\|\le\kappa. $$ Uniform equivalence on the compact neighbourhood gives one constant $L$ such that \begin{equation*} \|\mathcal R_n(d,z)\|_{\varepsilon,\chi_{n+1}} \le L\{\|d\|_{\varepsilon,\chi_n}^2 +\|d\|_{\varepsilon,\chi_n}\|z\|+\|z\|^2\}. \tag{B.6.56}\label{eq:B.6.56} \end{equation*} Let $$ \Phi(n,k+1)=A_n^{-1}\cdots A_k^{-1}. $$ These maps act between the indicated tangent spaces, and \begin{equation*} \|\Phi(n,k+1)\|_{\varepsilon,\chi_{k+1}\to\chi_n} \le\kappa^{\,k+1-n}. \tag{B.6.57}\label{eq:B.6.57} \end{equation*} A solution with $d_n\to0$ is equivalent, by backward iteration, to \begin{equation*} d_n=\sum_{k=n}^\infty\Phi(n,k+1) \{a(\chi_k)z_k-\mathcal R_k(d_k,z_k)\}. \tag{B.6.58}\label{eq:B.6.58} \end{equation*}

Use the Banach norm $$ \|d\|_*=\sup_{n\ge N}q_{\rm sh}^{-(n-N)} \|d_n\|_{\varepsilon,\chi_n}. $$ The linear forcing is bounded because \begin{equation*} \sum_{k=n}^\infty \kappa^{k+1-n}q_{\rm sh}^{k-N} =\frac{\kappa}{1-\kappa q_{\rm sh}}q_{\rm sh}^{n-N}. \end{equation*} On the ball $\|d\|_*\le R_0\|z_N\|$, \eqref{eq:B.6.56} bounds the nonlinear summand by a constant times $$ R_0^2\|z_N\|^2q_{\rm sh}^{2(k-N)}. $$ Combining this with \eqref{eq:B.6.57} gives a multiple of $$ \frac{\kappa}{1-\kappa q_{\rm sh}^2} R_0^2\|z_N\|^2q_{\rm sh}^{2(n-N)} \le C R_0^2\|z_N\|^2q_{\rm sh}^{n-N}. $$ The difference estimate has contraction constant $O(R_0\|z_N\|)$. Taking $N$ late enough makes it less than $1/2$; \eqref{eq:B.6.58} has a unique fixed point and \begin{equation*} \|d_n\|_{\varepsilon,\chi_n}+\|z_n\| \le Cq_{\rm sh}^{\,n-N}\|z_N\|. \tag{B.6.60}\label{eq:B.6.60} \end{equation*} Substitution back into \eqref{eq:B.6.55}, equivalently into the defining equation for $\Psi_n$, proves $\bar\chi_{n+1}=f(\bar\chi_n)$ exactly.

We verify that the shadowing orbit has positive weights. Let $u_{n,i}>0$ be a retained normalised squared coordinate after the last exact deletion. Every retained normalised-core squared weight has limiting two-block ratio one. Hence $$ \frac{u_{n+1,i}}{u_{n,i}}\longrightarrow1. $$ Choose $r$ with $q_{\rm sh}<r<1$. Increasing $N$ gives \begin{equation*} u_{n,i}\ge r^{\,n-N}u_{N,i}. \end{equation*} Work on the finite signed cover. On each member of this cover, every retained positive weight is an analytic function of the local coefficient variables, with uniformly bounded derivative on the compact chart. If a coordinate tends to a face, its signed amplitude is the chart variable and its weight is its square; the same derivative bound follows by differentiation and remains uniform at zero. The mean-value theorem in that signed chart, followed by the finite-cover maximum, gives $$ |\bar u_{n,i}-u_{n,i}| \le Cq_{\rm sh}^{\,n-N}\|z_N\|. $$ The stable coordinates decay at a fixed rate smaller than $r$, whereas each retained coordinate has successive ratio eventually at least $r$. It follows that $\|z_N\|/u_{N,i}\to0$ as $N\to\infty$. Choose $N$ so that $$ C\frac{\|z_N\|}{u_{N,i}}<\frac12 \quad\hbox{for every retained }i. $$ Then \begin{equation*} \frac{|\bar u_{n,i}-u_{n,i}|}{u_{n,i}} \le C\left(\frac{q_{\rm sh}}{r}\right)^{n-N} \frac{\|z_N\|}{u_{N,i}}<\frac12. \tag{B.6.62}\label{eq:B.6.62} \end{equation*} Every $\bar u_{n,i}$ is positive. The coordinate chart imposes total mass one, so $\bar u_n$ is a point of the probability simplex. Coordinate faces are invariant under the restricted return map, and \eqref{eq:B.6.62} keeps the shadowing orbit away from supports on which the iteration terminates. This yields a positive orbit on the invariant coordinate face.

The case in which a neutral external weight persists does not use the preceding shadowing construction. Suppose now that at least one neutral external weight persists. For every stable node $g$, $$ E_{g,n+1}\le q_sE_{g,n} $$ eventually, for some $q_s<1$. For a persistent neutral node $h$, $E_{h,n+1}/E_{h,n}\to1$. Choose $r\in(\sqrt{q_s},1)$. Then $$ E_{h,n}\ge c_hr^{\,n-N},\qquad E_{{\rm st},n}\le Cq_s^{\,n-N},\qquad \sum_{\text{stable }g}|\xi_{g,n}|\le Cq_s^{(n-N)/2}, $$ and hence \begin{equation*} E_{{\rm st},n}=o(E_{h,n}),\qquad \sum_{\text{stable }g}|\xi_{g,n}|=o(E_{h,n})=o\!\left(\sum_{\text{neutral }j}E_{j,n}\right). \tag{B.6.63}\label{eq:B.6.63} \end{equation*} Uniform invertibility of the moment matrices and the mean-value theorem in, respectively, the squared-mass and signed charts give the two asserted comparison bounds. Equation \eqref{eq:B.6.63} makes both little-oh of total neutral forcing. No assertion is made here about a logarithm of a vanishing coordinate. The only logarithmic quantities used below are specified log-increments, and their stable-weight comparisons are made after their analytic boundary extensions have been proved.
\end{proof}

Lemmas~\ref{lem:5.2},~\ref{lem:6.1}, and~\ref{lem:6.2} reduce the problem to the normalised core on the invariant $S_*$-face with persistent neutral modes retained explicitly. Unstable modes have disappeared; stable modes remain in the exact recurrence but contribute $o(E_{{\rm neu},n})$ to every analytic returned coordinate. For coordinates even in the stable signed amplitudes, the sharper comparison is $O(E_{{\rm st},n})$. The remaining cases are organised by support size and by whether a neutral mode persists.

\section{Four-node limits and five-node limits without persistent neutral modes}\label{sec:7}

The four-node case requires no further division according to external modes: the parity return on the core face is the identity, and all outside squared mass is summable. For five nodes without a persistent neutral mode, the fixed states form an interval parametrised by a pivot $\theta$, subsequently identified with a root of the common second-kind polynomial. An orientation argument yields convergence along this interval.

\begin{lemma}[Four-node completions]\label{lem:7.1}
If $|S_*|=4$, then $Y_n^2$ converges coordinatewise, even when coordinates outside $S_*$ remain positive at every finite index.
\end{lemma}

\begin{proof}
First let $J$ be any four-node support and $D_J(t)=\prod_{i\in J}(t-\lambda_i)$, $c_i=1/D_J'(\lambda_i)$. For a positive weight $w$ on $J$, orthogonality of its monic cubic $P$ says
\begin{equation*} w_iP(\lambda_i)=\kappa c_i. \end{equation*}
The leading-coefficient identity \eqref{eq:B.4.2} gives
\begin{equation*} \kappa=\sum_iw_iP(\lambda_i)^2=:a^2. \end{equation*}
The next weight is
\begin{equation*} \widetilde w_i=\frac{w_iP(\lambda_i)^2}{a^2} =\frac{a^2c_i^2}{w_i}. \tag{B.7.3}\label{eq:B.7.3} \end{equation*}
If $Q$ and $b^2$ are the next monic cubic and norm, its orthogonality gives $$ \widetilde w_iQ(\lambda_i)=\widetilde\kappa c_i,\qquad \widetilde\kappa =\sum_i\widetilde w_iQ(\lambda_i)^2=b^2, $$ where the second identity is again the leading-coefficient identity \eqref{eq:B.4.2}. Substitution of \eqref{eq:B.7.3} gives
\begin{equation*} Q(\lambda_i)=\frac{b^2c_i}{\widetilde w_i} =\frac{b^2w_i}{a^2c_i}. \end{equation*}
The degree-three leading-coefficient identity gives $$ 1=[t^3]Q=\sum_i c_iQ(\lambda_i) =\frac{b^2}{a^2}\sum_iw_i=\frac{b^2}{a^2}. $$ Hence $b=a$ and $P(\lambda_i)Q(\lambda_i)/(ab)=1$ for every $i$. The exact signed, and therefore squared, two-block return on four nodes is the identity.

Now take $J=S_*$, let $\Pi_J$ be coordinate projection, and define
\begin{equation*} E_n=\|(\operatorname{Id}-\Pi_J)Y_n\|^2,\qquad \widehat Y_n=\frac{\Pi_JY_n}{\sqrt{1-E_n}}. \end{equation*}
By \eqref{eq:B.5.6}, $\sum_nE_n<\infty$. Uniform moment-Gram regularity and analytic dependence on squared coordinates give
\begin{equation*} \widehat Y_{n+1}=G_J(\widehat Y_n)+\delta_n^{\rm sh},\qquad \|\delta_n^{\rm sh}\|\le C_{\rm sh}(E_n+E_{n+1}), \end{equation*}
where $G_J=L^2|_J$. Since $G_J$ is the identity, $\sum_n\|\widehat Y_{n+1}-\widehat Y_n\|<\infty$. Hence $\widehat Y_n$ and its squared coordinates converge; the outside coordinates tend to zero.
\end{proof}

The four-node case is now complete. On five nodes the fixed set has one degree of freedom, described by the pivot below.

\begin{lemma}[Five-node recovery and the fixed pivot interval]\label{lem:7.2}
Let $J=\{\lambda_1<\cdots<\lambda_5\}$, let $w_i>0$, and let $P$ be its monic orthogonal cubic. Define
\begin{equation*} c_i=\frac{1}{D_J'(\lambda_i)},\qquad H=\sum_iw_iP(\lambda_i)^2,\qquad \theta= \frac{\sum_i c_i^2\lambda_i/w_i}{\sum_i c_i^2/w_i}. \tag{B.7.7}\label{eq:B.7.7} \end{equation*}
Then
\begin{equation*} \boxed{\quad w_iP(\lambda_i)=Hc_i(\lambda_i-\theta). \quad} \tag{B.7.8}\label{eq:B.7.8} \end{equation*}
In particular, when $P(\lambda_i)\ne0$, this identity may be divided by $P(\lambda_i)$ to recover $w_i$.

At a fixed two-block pair $P,Q$, where $P(\lambda_i)Q(\lambda_i)=C>0$,
\begin{equation*} w_i^P=c_i(\lambda_i-\theta)Q(\lambda_i),\qquad w_i^Q=c_i(\lambda_i-\theta)P(\lambda_i). \tag{B.7.9}\label{eq:B.7.9} \end{equation*}
The simultaneous strict positivity conditions in \eqref{eq:B.7.9} define one bounded open interval $I=(\ell,r)$. At either endpoint exactly one weight vanishes and the other four stay positive.
\end{lemma}

\begin{proof}
The vector $(w_iP(\lambda_i))$ annihilates all quadratics. By \eqref{eq:B.4.2}, the nullspace of the $3\times5$ evaluation matrix is
$$ \{(c_i(a\lambda_i+b))_i:a,b\in\mathbb R\}. $$
Hence
\begin{equation*} w_iP(\lambda_i)=\kappa c_i(\lambda_i-\vartheta). \tag{B.7.10}\label{eq:B.7.10} \end{equation*}
Substitution in the reciprocal-weight expression in \eqref{eq:B.7.7}, together with $\sum_ic_iP(\lambda_i)=0$, gives $\vartheta=\theta$. Also
$$ \sum_ic_i(\lambda_i-\theta)P(\lambda_i)=1 $$
by the degree-four leading-coefficient identity. Multiplying \eqref{eq:B.7.10} by $P(\lambda_i)$ and summing gives $\kappa=H$, proving \eqref{eq:B.7.8}. At a fixed pair $H=C$ and $P(\lambda_i)Q(\lambda_i)=C>0$, so division is legitimate and gives \eqref{eq:B.7.9}.

Each inequality in \eqref{eq:B.7.9} is an affine half-line in $\theta$, and the two phase coefficients have the same signs. Their intersection is an interval. It is bounded because $\sum_ic_iQ(\lambda_i)=0$, so the five coefficients cannot all have one sign. An endpoint can equal only one of the distinct nodes, and its one linear factor vanishes simply; the remaining strict inequalities persist.
\end{proof}

\begin{lemma}[Convergence of an exact five-node tail]\label{lem:7.3}
Suppose that, after a finite index, a nonterminating parity orbit is supported on at most five fixed nodes and has no external forcing. Then its squared parity state converges.
\end{lemma}

\begin{proof}
The four-node case is Lemma~\ref{lem:7.1}. If a five-node block deletes a coordinate, deletion is permanent and the resulting tail is again covered by Lemma~\ref{lem:7.1}. Otherwise pass beyond the last such block. The one-step update is $\widetilde w_i=w_iP(\lambda_i)^2/H$; hence $P(\lambda_i)\ne0$ at every retained coordinate on this tail. Formula \eqref{eq:B.7.8} now gives $\theta\ne\lambda_i$ for every retained $i$. Define \begin{equation*} x_i=\lambda_i-\theta,\qquad \mathfrak A=\sum_i\frac{w_i}{x_i},\qquad B=\sum_i\frac{w_i}{x_i^2},\qquad \delta=\frac{\mathfrak A}{B},\qquad \mathcal N=1-\frac{\mathfrak A^2}{B}. \tag{B.7.11}\label{eq:B.7.11} \end{equation*} On the fixed interval $\mathfrak A=0$; thus $\theta$ is the tangential coordinate and $\mathfrak A$ measures transverse displacement. The strict Cauchy--Schwarz inequality gives $\mathfrak A^2\le B$. Equality would force $1/x_i$ to be constant on the support, which is impossible because the retained nodes are distinct. Hence $0<\mathcal N\le1$. Applying \eqref{eq:B.7.8} at the two consecutive phases gives the exact return \begin{equation*} w_i^+=w_i\frac{(1-\delta/x_i)^2}{\mathcal N}. \tag{B.7.12}\label{eq:B.7.12} \end{equation*} To verify this, the partner pivot computed from its reciprocal weights is $$ \frac{\sum_i c_i^2\lambda_i/\widetilde w_i} {\sum_i c_i^2/\widetilde w_i} =\theta+\frac{\sum_iw_i/x_i}{\sum_iw_i/x_i^2} =\theta+\delta, $$ and substitution in the second application of \eqref{eq:B.7.8} yields \eqref{eq:B.7.12}; summing its numerator gives $\mathcal N=1-\mathfrak A^2/B$.

Let $$ \gamma_i=c_i^2,\qquad b_i=\frac{\gamma_i}{w_i},\qquad B_{\rm dual}=\sum_i b_i. $$ The dual-mean definition of $\theta$ gives $$ \sum_i b_ix_i=0. $$ If $t=\theta^+-\theta$, substitution of \eqref{eq:B.7.12} into that definition gives the exact formula \begin{equation*} t= \frac{\displaystyle\sum_i b_ix_i^3/(x_i-\delta)^2} {\displaystyle\sum_i b_ix_i^2/(x_i-\delta)^2}. \tag{B.7.13}\label{eq:B.7.13} \end{equation*} Here $$ \frac{\gamma_i}{w_i^+} =b_i\mathcal N\frac{x_i^2}{(x_i-\delta)^2}, $$ and the common factor $\mathcal N$ cancels from the dual mean.

The following estimates are uniform at either endpoint. Let $\alpha=\lambda_h$ be the vanishing node and use the parametrisation $$ w_h=\varepsilon,\qquad w_j=(1-\varepsilon)v_j\quad(j\ne h),\qquad \sum_{j\ne h}v_j=1. $$ For this parametrisation, define $$ \Delta_j=\lambda_j-\alpha,\qquad D_0(v)=\sum_{j\ne h}\frac{\gamma_j\Delta_j}{v_j},\qquad G_0(v)=\sum_{j\ne h}\frac{\gamma_j}{v_j}. $$ Subtracting $\alpha$ from the exact dual mean gives \begin{equation*} q:=\theta-\alpha =\frac{\varepsilon D_0(v)} {\gamma_h(1-\varepsilon)+\varepsilon G_0(v)}. \tag{B.7.14}\label{eq:B.7.14} \end{equation*} At the selected fixed endpoint, \eqref{eq:B.7.8} reads $$ w_h=-c_hQ_*(\alpha)(\theta-\alpha). $$ Since $c_hQ_*(\alpha)\ne0$, comparison with \eqref{eq:B.7.14} shows $D_0(v)\ne0$ there. Shrink the endpoint chart so that it is bounded away from zero.

Define $$ H_\alpha=\gamma_h(1-\varepsilon)+\varepsilon G_0(v). $$ Then $x_h=-q$, and direct substitution in \eqref{eq:B.7.11} gives the analytic expression \begin{equation*} \mathfrak A =-\frac{H_\alpha}{D_0(v)} +(1-\varepsilon)\sum_{j\ne h}\frac{v_j}{\Delta_j-q}. \tag{B.7.15}\label{eq:B.7.15} \end{equation*} Also \begin{equation*} \begin{aligned} B&=\frac{\varepsilon}{q^2} +(1-\varepsilon)\sum_{j\ne h}\frac{v_j}{(\Delta_j-q)^2}\\ &=\frac{\gamma_h^2}{D_0(v)^2\varepsilon}+O(1). \end{aligned} \tag{B.7.16}\label{eq:B.7.16} \end{equation*} Uniformly for arbitrary relative rates of $\varepsilon\to0$ and $\mathfrak A\to0$, \begin{equation*} \delta= \frac{D_0(v)^2}{\gamma_h^2}\varepsilon\mathfrak A\{1+O(\varepsilon)\}, \end{equation*} \begin{equation*} \left|\frac{\delta}{x_h}\right|\le K|\mathfrak A|, \qquad \left|\frac{\delta}{x_j}\right| \le K\varepsilon|\mathfrak A|\quad(j\ne h). \tag{B.7.18}\label{eq:B.7.18} \end{equation*} Choose the chart so that every ratio in \eqref{eq:B.7.18} has modulus at most $1/2$. The exact scalar identity \begin{equation*} (1-r)^{-2} =1+2r+r^2\frac{3-2r}{(1-r)^2} \tag{B.7.19}\label{eq:B.7.19} \end{equation*} applied termwise in \eqref{eq:B.7.13}, with $\sum_i b_ix_i=0$, gives \begin{equation*} \begin{aligned} \operatorname{num}\eqref{eq:B.7.13} &=2\delta B_{\rm dual}+\delta^2 \sum_i\frac{b_i}{x_i} \frac{3-2\delta/x_i}{(1-\delta/x_i)^2},\\ \operatorname{den}\eqref{eq:B.7.13} &=B_{\rm dual}+2\delta\sum_i\frac{b_i}{x_i} +\delta^2\sum_i\frac{b_i}{x_i^2} \frac{3-2\delta/x_i}{(1-\delta/x_i)^2}. \end{aligned} \end{equation*} At the endpoint, $$ B_{\rm dual}\asymp\varepsilon^{-1},\qquad \frac{\sum_i b_i/|x_i|}{B_{\rm dual}}\le K\varepsilon^{-1},\qquad \frac{\sum_i b_i/|x_i|^2}{B_{\rm dual}}\le K\varepsilon^{-2}. $$ Together with $|\delta|\le K\varepsilon|\mathfrak A|$, these estimates give $$ \operatorname{num}\eqref{eq:B.7.13} =2\delta B_{\rm dual}\{1+O(|\mathfrak A|)\},\qquad \operatorname{den}\eqref{eq:B.7.13} =B_{\rm dual}\{1+O(|\mathfrak A|)\}. $$ Hence \begin{equation*} \boxed{\theta^+-\theta =2\delta\{1+O(|\mathfrak A|)\}.} \tag{B.7.21}\label{eq:B.7.21} \end{equation*}

We derive the normal return from an exact rational identity. For $x\ne0,t$, \begin{equation*} \frac{(x-\delta)^2}{x^2(x-t)} =\frac1x+\frac{t-2\delta}{x^2} +\frac{(t-\delta)^2}{x^2(x-t)}. \tag{B.7.22}\label{eq:B.7.22} \end{equation*} Multiply by $w_i$, sum, and use \eqref{eq:B.7.12}. This gives \begin{equation*} \mathcal N\mathfrak A^+ =\mathfrak A+(t-2\delta)B +(t-\delta)^2\sum_i\frac{w_i}{x_i^2(x_i-t)}. \tag{B.7.23}\label{eq:B.7.23} \end{equation*} Equation \eqref{eq:B.7.21} and $\delta B=\mathfrak A$ imply $$ (t-2\delta)B=O(\mathfrak A^2). $$ Equations \eqref{eq:B.7.18}, \eqref{eq:B.7.21} also give $|x_i-t|\ge|x_i|/2$. Hence $$ \left|\sum_i\frac{w_i}{x_i^2(x_i-t)}\right| \le K\varepsilon^{-2}, $$ while $\delta^2\le K\varepsilon^2\mathfrak A^2$. The last term of \eqref{eq:B.7.23} is $O(\mathfrak A^2)$. Finally, $$ 1-\mathcal N=\frac{\mathfrak A^2}{B}=O(\mathfrak A^2), $$ because $B$ has a uniform positive lower bound. Hence \begin{equation*} \boxed{\mathfrak A^+ =\mathfrak A+O(\mathfrak A^2) =\mathfrak A\{1+O(|\mathfrak A|)\}.} \tag{B.7.24}\label{eq:B.7.24} \end{equation*} For the vanishing endpoint weight, \eqref{eq:B.7.12}, \eqref{eq:B.7.18}, and $\mathcal N=1+O(\mathfrak A^2)$ also give \begin{equation*} \varepsilon^+-\varepsilon=O(\varepsilon|\mathfrak A|). \tag{B.7.25}\label{eq:B.7.25} \end{equation*}

All constants in \eqref{eq:B.7.21}, \eqref{eq:B.7.24}, and \eqref{eq:B.7.25} are uniform up to the endpoint and impose no relation between the rates at which $\varepsilon$ and $\mathfrak A$ vanish. Away from the two endpoints every $x_i$ is uniformly separated from zero, and the estimates follow directly from \eqref{eq:B.7.13}, \eqref{eq:B.7.19}, and \eqref{eq:B.7.22}. A finite cover of the closed fixed interval makes them uniform everywhere, including off-node factor collisions because no factor resultant occurs in these weight identities.

Every $\omega$-limit point is fixed by Lemma~\ref{lem:3.3}, so $\mathfrak A_n\to0$. Equation \eqref{eq:B.7.24} implies that $\mathfrak A_n$ is either zero from some index onward or has one eventual sign. Since $B_n>0$, \eqref{eq:B.7.21} gives the same eventual sign to $\theta_{n+1}-\theta_n$. The pivot is a positive weighted average of the five nodes, hence bounded; it converges.

Finally, Lemma~\ref{lem:5.1} gives convergence of the norm and cubic factors, and \eqref{eq:B.7.8} gives convergence of each squared weight. If $\theta=\lambda_i$ at a finite block, the exact block factor at that node is zero; the coordinate is permanently deleted and Lemma~\ref{lem:7.1} applies. Endpoint convergence without exact deletion gives one zero limiting weight and the same four-node limit.
\end{proof}

\begin{corollary}[Absence of persistent neutral modes]\label{cor:7.4}
If $|S_*|\le5$ and no persistent neutral external mode remains, then $Y_n^2$ converges.
\end{corollary}

\begin{proof}
Unstable modes have been deleted by Lemma~\ref{lem:5.2}. Lemma~\ref{lem:6.2} supplies an exact orbit of positive weights on the $S_*$-face exponentially close to the original tail. Apply Lemma~\ref{lem:7.1} or Lemma~\ref{lem:7.3} to that restricted orbit and transfer its limit through \eqref{eq:B.6.60}.
\end{proof}

\section{Neutral external modes at a five-node limit}\label{sec:8}

The remaining five-node case has at least one persistent neutral external mode. The missing sixth root of the completion polynomial supplies a second-kind root outside the pivot interval. Lemmas~\ref{lem:8.1}--\ref{lem:8.4} establish the geometry of the two roots and settle all configurations with one-sided forcing by a single-crossing argument. Mixed configurations require logarithmic cocycles: Lemma~\ref{lem:8.5} treats a distinct completion root and Lemma~\ref{lem:8.6} its confluent counterpart.

Let $S_*=J=\{\lambda_1<\cdots<\lambda_5\}$, let $D_J(t)=\prod_{i=1}^5(t-\lambda_i)$ and $c_i=1/D_J'(\lambda_i)$, and factor the limiting completion, with multiplicity, as \begin{equation*} \mathscr D_*(t)=D_J(t)(t-\alpha),\qquad P_*Q_*=\mathscr D_*+C_*. \tag{B.8.1}\label{eq:B.8.1} \end{equation*} The sixth root $\alpha$ need not be spectral and may coincide with one retained non-pivot node. Work with the normalised core of Lemma~\ref{lem:5.3}. From this point through Section~\ref{sec:9}, $P_n,Q_n,a_n,b_n,C_n,\mathscr D_n,U_n$ and $w_{n,i}$ are untilded core objects. A superscript $+$ on a core coordinate denotes the normalised retained core induced by the exact full parity return; setting all external masses to zero gives the restricted-core return. The corresponding exact full polynomials, norms, completion, and multiplier carry tildes as in \eqref{eq:B.5.24}--\eqref{eq:B.5.26}. The core fixed component is parametrised by the pivot $\theta\in I=(\ell,r)$ from Lemma~\ref{lem:7.2}. In the analytic completion chart \begin{equation*} K=(t-\theta)(t-\alpha) \end{equation*} and departure from the fixed component is an affine polynomial $U$. Define \begin{equation*} g=\alpha-\theta,\quad s=\operatorname{sign}g,\quad a=-U(\theta),\quad e=-\frac{U(\alpha)}{\alpha-\theta},\quad q=se, \end{equation*} and \begin{equation*} d_\theta=(\theta-\ell)(r-\theta). \end{equation*}

\begin{lemma}[Geometry of the common second-kind polynomial]\label{lem:8.1}
For every positive five-node fixed pair \eqref{eq:B.8.1}, the common monic second-kind quadratic of its two phase measures is \begin{equation*} (t-\theta)(t-\alpha). \end{equation*} Exactly one of the following orders occurs: \begin{equation*} \begin{array}{c|c|c} I&\theta&\alpha\\ \hline (\lambda_2,\lambda_3)&\text{lower second-kind root} &\lambda_3<\alpha<\lambda_5,\\ (\lambda_3,\lambda_4)&\text{upper second-kind root} &\lambda_1<\alpha<\lambda_3. \end{array} \tag{B.8.6}\label{eq:B.8.6} \end{equation*} Hence \begin{equation*} \alpha\in(\lambda_1,\lambda_5)\setminus[\ell,r], \qquad |\alpha-\theta|\ge g_0>0, \tag{B.8.7}\label{eq:B.8.7} \end{equation*} and no neutral external node lies in $I$.
\end{lemma}

\begin{proof}
Using \eqref{eq:B.7.9} and matching residues and behaviour at infinity gives \begin{equation*} G_P(z)=\sum_i\frac{w_i^P}{z-\lambda_i} =\frac{(z-\theta)Q(z)}{D_J(z)} \end{equation*} and \begin{equation*} \sum_i\frac{w_i^PP(\lambda_i)}{z-\lambda_i} =C_*\frac{z-\theta}{D_J(z)}. \end{equation*} The monic degree-two second-kind polynomial is \begin{equation*} P(z)G_P(z)- \sum_i\frac{w_i^PP(\lambda_i)}{z-\lambda_i} =(z-\theta)(z-\alpha). \end{equation*} The same calculation at the partner phase gives the identical quadratic. For a positive measure this polynomial is the characteristic polynomial of the irreducible trailing $2\times2$ principal minor of the $3\times3$ Jacobi matrix. Its roots are real and simple and strictly interlace the three roots of the orthogonal cubic \cite[sec.~3.3]{Szego1975}; see also \cite[chaps.~2--3]{GolubMeurant2010}.

The signs in \eqref{eq:B.7.9} force both cubics to have one root in each node gap except the pivot gap, and no root in that gap. Strict interlacing leaves exactly \eqref{eq:B.8.6}, which proves \eqref{eq:B.8.7}. On $I$, neither cubic vanishes and $P(h)Q(h)>0$. A neutral external node satisfies $P(h)Q(h)=-C_*<0$, so it cannot belong to $I$.
\end{proof}

\begin{lemma}[Transverse curvature at collisions]\label{lem:8.2}
Define $$ g=\alpha-\theta,\qquad s=\operatorname{sign}g,\qquad a=-U(\theta),\qquad e=-\frac{U(\alpha)}g,\qquad q=se. $$ The evaluation map from affine polynomials to $(a,e)$ is invertible because \begin{equation*} U(t)=\frac{a(t-\alpha)}g-e(t-\theta). \tag{B.8.11a}\label{eq:B.8.11a} \end{equation*} On the unforced five-node return, \begin{equation*} q^+-q=-c_\theta d_\theta a^2+O(d_\theta|a|^3), \tag{B.8.11b}\label{eq:B.8.11b} \end{equation*} where $c_\theta$ extends over the closed pivot interval as a positive analytic unit. Formula \eqref{eq:B.8.11b} remains valid when $\alpha\in J$ and through every simple or simultaneous off-node common-root collision of $P,Q$.
\end{lemma}

\begin{proof}
First suppose that $\alpha\notin J$ and $P,Q$ are coprime. Adjoin $\alpha$ as an additional spectral node with zero weight. The added coordinate is invariantly zero, so the residual polynomials and the five weights on $J$ do not change. In this six-node face, $$ \mathscr D_*(t)=D_J(t)(t-\alpha),\qquad K(t)=(t-\theta)(t-\alpha). $$ Moreover $P(\alpha)Q(\alpha)=C_*$, so the virtual signed multiplier is $+1$. The one-loss chart and the symbolic-to-ordinary power passage in Lemma~\ref{lem:4.3} extend the quadratic tensor of Lemma~\ref{lem:4.2} to this invariant face.

We compute the quadratic coefficient at $\alpha$. Affinely send $\alpha$ to $0$ and $\theta$ to $1$, including the associated monic rescaling. Then $$ K=t(t-1),\qquad u=t. $$ Define \begin{equation*} D_1(t)=\frac{\mathscr D_*(t)-\mathscr D_*(1)}{t-1}, \qquad \delta_{\mathscr D}=-\frac{\mathscr D_*(1)}{C_*}, \end{equation*} and let $p_{\rm aux},q_{\rm aux}\in\mathcal P_2$ be the unique solution of \begin{equation*} Qp_{\rm aux}+Pq_{\rm aux}=D_1. \end{equation*} The Sylvester map $$ \mathcal P_2\oplus\mathcal P_2\longrightarrow\mathcal P_5, \qquad(\phi,\psi)\longmapsto Q\phi+P\psi $$ has zero kernel: $Q\phi=-P\psi$, coprimality, and $\deg\phi,\deg\psi\le2<3$ force $\phi=\psi=0$. Its domain and codomain both have dimension six, proving existence and uniqueness.

The first-order lift of $u$ through phase $j$ has affine factor displacement $j\delta_{\mathscr D}u$. Expanding the exact factor equation to second order shows that its phase-$j$ forcing is \begin{equation*} Kp_{\rm aux}q_{\rm aux}+j\delta_{\mathscr D}uD_1. \tag{B.8.12c}\label{eq:B.8.12c} \end{equation*} The first summand cancels as follows. If its phase-zero decomposition is \begin{equation*} KPr+Kp_{\rm aux}q_{\rm aux}=\mathscr D_*v+\eta K+C_*k, \tag{B.8.12d}\label{eq:B.8.12d} \end{equation*} with $r\in\mathcal P_2$, $k\in\mathcal P_1$, $v\in\mathcal P_1$, and $\eta$ scalar, then at the next phase the incoming cubic and $K$-corrections are $r,k$. Choosing outgoing cubic correction $0$ and outgoing $K$-correction $2k$ gives the exact identity $$ C_*k+KrP+Kp_{\rm aux}q_{\rm aux}=\mathscr D_*v+\eta K+2C_*k $$ by \eqref{eq:B.8.12d}. The decomposition into a multiple of $\mathscr D_*$, a multiple of $K$, an affine monic-factor correction, and a degree-at-most-two cubic correction is unique. The phase-two and phase-zero affine normal outputs are both $v$, so $Kp_{\rm aux}q_{\rm aux}$ has zero two-block normal return.

We compute the response to $uD_1$. Its phase-one decomposition is \begin{equation*} KQr+uD_1=\mathscr D_*v_1+\eta K+C_*k. \tag{B.8.12e}\label{eq:B.8.12e} \end{equation*} With $$ J_Q(t)=\frac{Q(t)-Q(1)}{t-1}, $$ the identities $$ (t-1)D_1=\mathscr D_*-\mathscr D_*(1),\qquad (t-1)J_Q=Q-Q(1) $$ give, after reduction modulo $Q$, \begin{equation*} \operatorname{rem}_Q(uD_1) =P(1)\{J_Q-Q(1)\}. \end{equation*} Its value at zero is $-P(1)Q(0)$. Evaluating \eqref{eq:B.8.12e} at zero and also reducing it modulo $Q$, then evaluating, gives \begin{equation*} \mathscr D_*(0)v_1(0)+C_*k(0)=0, \tag{B.8.12g}\label{eq:B.8.12g} \end{equation*} \begin{equation*} -P(1)Q(0)=-C_*v_1(0)+C_*k(0). \end{equation*} Because $\mathscr D_*(0)+C_*=P(0)Q(0)$, solving this two-by-two system yields \begin{equation*} v_1(0)=\frac{P(1)}{P(0)},\qquad k(0)=-\frac{\mathscr D_*(0)P(1)}{C_*P(0)}. \tag{B.8.12i}\label{eq:B.8.12i} \end{equation*}

Let $s_{\rm cub}\in\mathcal P_2$ denote the outgoing cubic correction. The phase-two affine output $B$ satisfies the exact factor equation \begin{equation*} C_*k+KrQ+KP s_{\rm cub}+2uD_1 =\mathscr D_*B+\zeta K+C_*\ell. \tag{B.8.12i1}\label{eq:B.8.12i1} \end{equation*} Substitution of \eqref{eq:B.8.12e}, followed by reduction modulo $P$, gives \begin{equation*} -C_*v_1+\eta K+2C_*k+\operatorname{rem}_P(uD_1) =-C_*B+\zeta K+C_*\ell. \tag{B.8.12i2}\label{eq:B.8.12i2} \end{equation*} The remainder identity is $$ \operatorname{rem}_P(uD_1) =Q(1)\left\{\frac{P-P(1)}{t-1}-P(1)\right\}, $$ whose value at zero is $-Q(1)P(0)$. Evaluation of \eqref{eq:B.8.12i1} at zero also gives $$ C_*k(0)=\mathscr D_*(0)B(0)+C_*\ell(0). $$ Evaluate \eqref{eq:B.8.12i2} at zero, eliminate $\ell(0)$ with this identity, and use \eqref{eq:B.8.12g}--\eqref{eq:B.8.12i}. The result is $$ P(0)Q(0)B(0)=P(1)Q(0)+Q(1)P(0), $$ and hence \begin{equation*} B(0)=\frac{P(1)}{P(0)}+\frac{Q(1)}{Q(0)}. \end{equation*} The input $U=\varepsilon u$ has no quadratic coefficient, the $Kp_{\rm aux}q_{\rm aux}$ contribution has zero return, and the remaining forcing in \eqref{eq:B.8.12c} is $\delta_{\mathscr D}$ times the response whose phase-two affine output is $B$. Hence the returned quadratic coefficient at the zero root is $\delta_{\mathscr D}B(0)$. Restoring $\delta_{\mathscr D}=-\mathscr D_*(1)/C_*$ and undoing the affine normalisation proves \begin{equation*} \Gamma_\alpha =-\frac{\mathscr D_*(\theta)}{C_*} \left\{\frac{P(\theta)}{P(\alpha)} +\frac{Q(\theta)}{Q(\alpha)}\right\}. \tag{B.8.12k}\label{eq:B.8.12k} \end{equation*} The denominators in \eqref{eq:B.8.12k} are $C_*,P(\alpha),Q(\alpha)$, all nonzero.

In either order in \eqref{eq:B.8.6}, $$ \mathscr D_*(\theta)=D_J(\theta)(\theta-\alpha)>0. $$ Strict second-kind interlacing places exactly one root of each of $P,Q$ between $\theta$ and $\alpha$. Hence $$ \frac{P(\theta)}{P(\alpha)}<0,\qquad \frac{Q(\theta)}{Q(\alpha)}<0, $$ and \eqref{eq:B.8.12k} gives $\Gamma_\alpha>0$. At a pivot endpoint $j\in\{\ell,r\}$, $\mathscr D_*(\theta)$ has the same simple zero as $d_\theta$. The four factor evaluations in \eqref{eq:B.8.12k} stay nonzero because $$ P(j)Q(j)=P(\alpha)Q(\alpha)=C_*>0. $$ Hence $\Gamma_\alpha/d_\theta$ extends as a positive analytic unit.

For a pure canonical normal, $e=0$ and $$ U(t)=\frac{a(t-\alpha)}g. $$ Bilinearity of the quadratic return tensor and \eqref{eq:B.8.12k} give $$ U^+(\alpha)-U(\alpha) =\frac{a^2}{g^2}\Gamma_\alpha+O(d_\theta|a|^3). $$ Since $e=-U(\alpha)/g$, $$ e^+-e=-\frac{\Gamma_\alpha}{g^3}a^2 +O(d_\theta|a|^3). $$ Multiplication by $s$, together with $sg^{-3}=|g|^{-3}$, proves \eqref{eq:B.8.11b} with $$ c_\theta=\frac{\Gamma_\alpha}{d_\theta|g|^3}>0. $$

Suppose next that $\alpha=c\in J$ and that $P,Q$ are coprime. Work in coefficient space and replace the zero-weight virtual node $c$ by $\alpha'$. We solve \begin{equation*} P_{\alpha'}Q_{\alpha'}=D_J(t)(t-\alpha')+C_* \end{equation*} for monic cubics $P_{\alpha'},Q_{\alpha'}$. The derivative with respect to their lower coefficients is $$ (p,q)\longmapsto Qp+Pq: \mathcal P_2\oplus\mathcal P_2\longrightarrow\mathcal P_5. $$ It is injective: if $Qp=-Pq$, coprimality implies $P\mid p$, and the degree bound forces $p=q=0$. The spaces have the same dimension, so the derivative is an isomorphism. The analytic implicit-function theorem applied directly at $\alpha'=c$ now gives analytic coefficient branches.

The recovered five-node weights extend analytically to $\alpha'=c$. If $S_{\alpha'}=J\cup\{\alpha'\}$ and $K_{\alpha'}=(t-\theta)(t-\alpha')$, then for $\lambda_i\in J$ \begin{equation*} \frac{K_{\alpha'}(\lambda_i)Q_{\alpha'}(\lambda_i)} {D_{S_{\alpha'}}'(\lambda_i)} =\frac{(\lambda_i-\theta)Q_{\alpha'}(\lambda_i)} {D_J'(\lambda_i)}. \end{equation*} The factor $\lambda_i-\alpha'$ has cancelled. The right-hand side is analytic at $\alpha'=c$ and is precisely the original recovery formula; the same argument applies in the other phase with $P$ and $Q$ interchanged. The coefficients, weights, and their first two derivatives converge. Formula \eqref{eq:B.8.12k} has the same limit because its only displayed denominators are $C_*,P(\alpha'),Q(\alpha')$, and $P(c)Q(c)=C_*>0$.

At a common-factor collision, Lemma~\ref{lem:4.3} proves that the fixed ideal is a radical complete intersection and that its symbolic and ordinary powers agree. The coprime-locus quadratic coefficient extends in the ordinary square on every component. The corresponding endpoint factorisation is as follows. In a regular signed pivot chart let $y$ be the omitted amplitude. Recovery gives $$ y^2=u\,d_\theta $$ with $u>0$ an analytic unit. At a fixed five-node state, differentiation of the dual-pivot formula gives \begin{equation*} d\mathfrak A=\frac{B D_J(\theta)}{C_*}\,dU(\theta) =-\frac{B D_J(\theta)}{C_*}\,da. \end{equation*} The product $B D_J(\theta)$ has a finite nonzero limit at either pivot endpoint, as follows from \eqref{eq:B.7.14}--\eqref{eq:B.7.16}. The coordinates $a$ and $\mathfrak A$ are equivalent local conormals there. The four-node component has ideal $(y)$ and the five-node fixed component has ideal $(a)$; their union therefore has the radical principal ideal \begin{equation*} \mathcal I_{\rm ep}=(ya). \end{equation*} On the component $(y)$ the return is exactly the identity. On the component $(a)$ its value and first normal derivative vanish by the coprime calculation above. The localised remainder lies in $(\mathcal I_{\rm ep}R_{\mathfrak p})^2$ at both minimal primes $\mathfrak p=(y),(a)$. Since a radical principal complete intersection has $\mathcal I_{\rm ep}^{(2)}=\mathcal I_{\rm ep}^2$, the remainder is $y^2a^2$ times an analytic coefficient. Substituting $y^2=u d_\theta$ and its nonzero coprime-locus coefficient factors it as $$ d_\theta a^2\{c_\theta+O(a)\}. $$ The resulting coefficient belongs to the ordinary square and is uniform at the endpoints and collisions.
\end{proof}

Lemma~\ref{lem:8.2} gives the unforced quadratic drift. The next calculation identifies the term contributed by each external neutral mass.

\begin{lemma}[External kernel derivative and endpoint signs]\label{lem:8.3}
Let $h\notin J$ be neutral: $$ P(h)Q(h)=-C_*. $$ For either fixed phase let $k_P(\,\cdot\,,h)$ and $k_Q(\,\cdot\,,h)$ denote the degree-at-most-two reproducing kernels. The external affine completion response is \begin{equation*} A_h(t)=4P(h)\operatorname{quo}_{D_J} \left((t-\theta) \{P(t)k_Q(t,h)-Q(t)k_P(t,h)\}\right). \tag{B.8.14}\label{eq:B.8.14} \end{equation*} If $\rho<\sigma$ are the two roots $\theta,\alpha$, then \begin{equation*} \operatorname{sign}A_h(\rho)=\operatorname{sign}(h-\rho),\qquad \operatorname{sign}A_h(\sigma)=\operatorname{sign}(\sigma-h). \tag{B.8.15}\label{eq:B.8.15} \end{equation*} Hence \begin{equation*} \operatorname{sign} \left.\partial_{E_h}(e^+-e)\right|_{\rm fixed} =\operatorname{sign}(h-\alpha). \tag{B.8.16}\label{eq:B.8.16} \end{equation*} Define the distinct-root moment normal \begin{equation*} \mathfrak A=\sum_{i\in J}\frac{w_i}{\lambda_i-\theta}, \qquad B=\sum_{i\in J}\frac{w_i}{(\lambda_i-\theta)^2}. \tag{B.8.17a}\label{eq:B.8.17a} \end{equation*} Then \begin{equation*} \left.\partial_{E_h}\mathfrak A^+\right|_{\rm fixed} =\frac{B\,D_J(\theta)}{C_*}A_h(\theta), \tag{B.8.17b}\label{eq:B.8.17b} \end{equation*} and this coefficient is a finite nonzero signed unit at both pivot endpoints.
\end{lemma}

\begin{proof}
Here $P,Q$ are the untilded polynomials of the fixed core. Insert squared mass $E$ at $h$: $$ \mu_E=(1-E)\mu_P+E\delta_h. $$ The factor $1-E$ does not affect monic orthogonality. If the full polynomial is $\widetilde P_E=P+Ep_h+O(E^2)$, differentiation of $$ \int \widetilde P_Ef\,d\mu_E=0\qquad(f\in\mathcal P_2) $$ gives $$ \int p_hf\,d\mu_P=-P(h)f(h), $$ hence \begin{equation*} p_h=-P(h)k_P(\,\cdot\,,h). \tag{B.8.18a}\label{eq:B.8.18a} \end{equation*} Let $\widetilde Q_E=Q+Eq_h+O(E^2)$ be the next full monic cubic. Its orthogonality is $$ \int \widetilde P_E^2\widetilde Q_Ef\,d\mu_E=0. $$ At a core node $PQ=C_*$, whereas $P(h)Q(h)=-C_*$. Dividing the differentiated equation by $C_*$ gives $$ \begin{aligned} \langle q_h,f\rangle_{\mu_Q} &=-\frac2{C_*}\int Pp_hQf\,d\mu_P -\frac{P(h)^2Q(h)}{C_*}f(h)\\ &=2P(h)f(h)+P(h)f(h). \end{aligned} $$ Hence \begin{equation*} q_h=3P(h)k_Q(\,\cdot\,,h). \tag{B.8.18b}\label{eq:B.8.18b} \end{equation*}

After two blocks, normalise the returned core to total mass one. A common constant in the logarithmic core-weight derivative disappears, and at the retained nodes \begin{equation*} \delta\log w_i^+ =\frac{2p_h(\lambda_i)}{P(\lambda_i)} +\frac{2q_h(\lambda_i)}{Q(\lambda_i)} +\text{a constant independent of }i. \tag{B.8.18c}\label{eq:B.8.18c} \end{equation*} If $p_{\rm out}$ is the returned monic-cubic derivative, then $$ \begin{aligned} \langle p_{\rm out},f\rangle_{\mu_P} &=-\sum_iw_iP_i f_i \left(\frac{2(p_h)_i}{P_i}+\frac{2(q_h)_i}{Q_i}\right)\\ &=-2\langle p_h,f\rangle_{\mu_P} -2\langle q_h,f\rangle_{\mu_Q}\\ &=-4P(h)f(h), \end{aligned} $$ so \begin{equation*} p_{\rm out}=-4P(h)k_P. \tag{B.8.18d}\label{eq:B.8.18d} \end{equation*} Differentiate the returned partner orthogonality $$ \int P_{\rm out}^2Q_{\rm out}f\,d\mu_P^{\rm out}=0. $$ Using \eqref{eq:B.8.18c}, its derivative is $$ \begin{aligned} 0={}&C_*\langle q_{\rm out},f\rangle_{\mu_Q} +2C_*\langle p_{\rm out},f\rangle_{\mu_P}\\ &+2\int P Q p_h f\,d\mu_P +2\int P^2q_hf\,d\mu_P. \end{aligned} $$ Since $PQ=C_*$ on the core and $d\mu_Q=P^2d\mu_P/C_*$, equations \eqref{eq:B.8.18a}--\eqref{eq:B.8.18b}, \eqref{eq:B.8.18d} reduce this to $$ 0=C_*\langle q_{\rm out},f\rangle_{\mu_Q} -8C_*P(h)f(h)-2C_*P(h)f(h)+6C_*P(h)f(h). $$ Hence \begin{equation*} q_{\rm out}=4P(h)k_Q. \tag{B.8.18e}\label{eq:B.8.18e} \end{equation*} More generally, if $k$ is any external node at a nearby fixed core and $m_k=P(k)Q(k)/C$, where $C$ is that core's product, the same differentiation without the neutral substitution gives $$ p^{[k]}=-P(k)k_P,\qquad q^{[k]}=P(k)(2-m_k)k_Q, $$ $$ p_{\rm out}^{[k]}=-2(1-m_k)P(k)k_P,\qquad q_{\rm out}^{[k]}=2(1-m_k)P(k)k_Q. $$ The displayed neutral formulas are the special case $C=C_*$ and $m_k=-1$; the general coefficients are analytic and bounded, but no endpoint sign is asserted for them. Differentiate the virtual six-node factor identity at the fixed state: $$ \delta K\,PQ+K(p_{\rm out}Q+Pq_{\rm out}) =\mathscr D_*(\delta K+\delta U)+\delta H\,K+C_*\delta K_1. $$ Because $PQ=\mathscr D_*+C_*$, cancellation of $\mathscr D_*\delta K$ leaves $$ K(p_{\rm out}Q+Pq_{\rm out}) =\mathscr D_*\delta U+ \{\delta H\,K+C_*(\delta K_1-\delta K)\}. $$ The braces have degree at most two. Polynomial division by $\mathscr D_*$ makes $\delta U$ the quotient of the left-hand side. Cancelling the common virtual factor $t-\alpha$ shows equivalently that the five-node affine response is the quotient by $D_J$ of $$ (t-\theta)\{Pq_{\rm out}+Qp_{\rm out}\}, $$ which proves \eqref{eq:B.8.14}.

For the endpoint signs, affinely normalise $\rho=0,\sigma=1$. The common second-kind polynomial is $K=t(t-1)$. The last two monic Jacobi recurrences for the $P$-phase measure have the form \begin{equation*} P=(t-a_0)K-b(t-c),\qquad b>0,\quad0<c<1, \tag{B.8.19a}\label{eq:B.8.19a} \end{equation*} and the preceding monic quadratic and its squared norm are \begin{equation*} R_P=(t-1+c)(t-a_0)-b,\qquad H_P=bc(1-c). \end{equation*} If the last two Jacobi diagonal entries are $a_1,a_2$ and the last off-diagonal square is $\beta_2>0$, the second-kind identity $$ (t-a_2)(t-a_1)-\beta_2=t(t-1) $$ gives $a_1+a_2=1$ and $\beta_2=a_2(1-a_2)$, with $c=a_2$. For the $Q$-phase, \begin{equation*} Q=(t-A_0)K-B_0(t-f),\qquad B_0>0,\quad0<f<1, \end{equation*} \begin{equation*} R_Q=(t-1+f)(t-A_0)-B_0,\qquad H_Q=B_0f(1-f). \tag{B.8.19d}\label{eq:B.8.19d} \end{equation*} Christoffel--Darboux gives \begin{equation*} k_P(t,h)= \frac{P(t)R_P(h)-R_P(t)P(h)}{H_P(t-h)} \tag{B.8.19e}\label{eq:B.8.19e} \end{equation*} and the displayed formula with $P,R_P,H_P$ replaced by $Q,R_Q,H_Q$.

Define $$ S(t)=P(t)k_Q(t,h)-Q(t)k_P(t,h),\qquad s_5=[t^5]S,\quad s_4=[t^4]S, $$ where $p_2=[t^2]P$ and $q_2=[t^2]Q$. Taking the two leading coefficients in \eqref{eq:B.8.19e} gives $$ s_5=\frac{R_Q(h)}{H_Q}-\frac{R_P(h)}{H_P}, $$ $$ s_4=(h+p_2+q_2)s_5+ \frac{P(h)}{H_P}-\frac{Q(h)}{H_Q}. $$ Let $\Theta=\operatorname{quo}_{\mathscr D_*}(KS)$. Since $$ \mathscr D_*=t^6+(p_2+q_2)t^5+\cdots,\qquad K=t^2-t, $$ comparison of the coefficients of $t^7,t^6$ in $KS=\mathscr D_*\Theta+\operatorname{rem}_{\mathscr D_*}(KS)$ gives $$ \Theta(0)=(h-1)s_5+\frac{P(h)}{H_P}-\frac{Q(h)}{H_Q}, $$ $$ \Theta(1)=hs_5+\frac{P(h)}{H_P}-\frac{Q(h)}{H_Q}. $$ Substitution of \eqref{eq:B.8.19a}--\eqref{eq:B.8.19d}, followed by the two elementary identities $$ (h-a_0)(h-1)-b=\frac{P(h)-bc}{h}, $$ $$ h(h-a_0)-b=\frac{P(h)+b(1-c)}{h-1}, $$ and their $Q$-phase versions, gives \begin{equation*} \frac{A_h(0)}{4P(h)} =\frac1h \left(\frac{P(h)}{bc}-\frac{Q(h)}{B_0f}\right), \tag{B.8.19f}\label{eq:B.8.19f} \end{equation*} \begin{equation*} \frac{A_h(1)}{4P(h)} =\frac1{h-1} \left(-\frac{P(h)}{b(1-c)} +\frac{Q(h)}{B_0(1-f)}\right). \tag{B.8.19g}\label{eq:B.8.19g} \end{equation*} Under these substitutions, the constant terms from the two phases cancel. Using $P(h)Q(h)=-C_*$ in \eqref{eq:B.8.19f}--\eqref{eq:B.8.19g} yields \begin{equation*} A_h(0)=\frac4h \left(\frac{P(h)^2}{bc}+\frac{C_*}{B_0f}\right), \end{equation*} \begin{equation*} A_h(1)=-\frac4{h-1} \left(\frac{P(h)^2}{b(1-c)} +\frac{C_*}{B_0(1-f)}\right). \end{equation*} The parenthesised quantities are strictly positive, proving \eqref{eq:B.8.15}.

At a fixed state, $U=0$, so \begin{equation*} \left.\partial_{E_h}(e^+-e)\right|_{\rm fixed} =-\frac{A_h(\alpha)}g. \end{equation*} If $\theta<\alpha$, then $\alpha=\sigma$; if $\alpha<\theta$, then $\alpha=\rho$. Formula \eqref{eq:B.8.15} in the appropriate case proves \eqref{eq:B.8.16}.

To obtain \eqref{eq:B.8.17b}, consider a five-node fixed state, where \begin{equation*} w_i=C_*c_i\frac{\lambda_i-\theta}{P_i}. \end{equation*} For the external insertion direction, \eqref{eq:B.8.18c} gives \begin{equation*} \delta\mathfrak A^+ =\sum_i\frac{w_i\,\delta\log w_i^+}{\lambda_i-\theta} +B\,\delta\theta. \tag{B.8.21b}\label{eq:B.8.21b} \end{equation*} Differentiate the returned fixed recovery identity $$ w_i^+=C_*c_i\frac{\lambda_i-\theta^+}{P_{{\rm out},i}}. $$ In the derivative of $\mathfrak A^+$, the two $\delta\theta$-terms cancel, while the common norm term is multiplied by $\mathfrak A=0$. Hence $$ \delta\mathfrak A^+ =-C_*\sum_i\frac{c_ip_{{\rm out},i}}{P_i^2}. $$ Since $p_{\rm out}=4p_h$, the $p$-part of the first sum in \eqref{eq:B.8.21b} equals $$ 2C_*\sum_i\frac{c_i(p_h)_i}{P_i^2} =-\frac12\,\delta\mathfrak A^+. $$ The $q$-part is $$ 2\sum_i c_i(q_h)_i=0 $$ by \eqref{eq:B.4.2}, since $\deg q_h\le2$. Hence \begin{equation*} \delta\mathfrak A^+=\frac23B\,\delta\theta. \tag{B.8.21c}\label{eq:B.8.21c} \end{equation*} For the general source above the corresponding relation is $$ B\,\delta\theta=\frac{2-m_k}{1-m_k}\,\delta\mathfrak A^+, $$ whenever $m_k\ne1$. Thus \eqref{eq:B.8.21c}, and the cancellation for which it is used below, are neutral-source identities.

For an arbitrary intrinsic perturbation of the five-node core, the exact return formula of Lemma~\ref{lem:7.3} gives \begin{equation*} \delta\theta=\frac{2}{B}\,d\mathfrak A. \tag{B.8.21d}\label{eq:B.8.21d} \end{equation*} On the other hand, evaluation of $$ dK^+-dK=-\frac2{C_*}\operatorname{rem}_K(\mathscr D_*\,dU) $$ at $\theta$, followed by differentiation of $K^+(\theta^+)=0$, gives \begin{equation*} \delta\theta=\frac{2D_J(\theta)}{C_*}\,dU(\theta). \tag{B.8.21e}\label{eq:B.8.21e} \end{equation*} Both $d\mathfrak A$ and $dU(\theta)$ vanish on the tangent space of the fixed component. Equating \eqref{eq:B.8.21d} and \eqref{eq:B.8.21e} gives the conormal identity \begin{equation*} d\mathfrak A =\frac{B D_J(\theta)}{C_*}\,dU(\theta). \tag{B.8.21f}\label{eq:B.8.21f} \end{equation*} For external insertion the returned completion response is $dU^+(\theta)=A_h(\theta)dE_h$; applying \eqref{eq:B.8.21f} at the returned fixed state proves \eqref{eq:B.8.17b}.

Near a pivot endpoint $j$, the vanishing fixed weight is a positive unit times $|\theta-j|$. Hence $$ B\asymp|\theta-j|^{-1},\qquad D_J(\theta)\asymp\theta-j, $$ so $BD_J(\theta)$ has a finite nonzero limit. Formula \eqref{eq:B.8.15} keeps $A_h(\theta)$ nonzero there.
\end{proof}

\begin{lemma}[Single-crossing estimate for the transverse recurrence]\label{lem:8.4}
Let $\mathcal N$ be the persistent neutral external nodes, let $E_{\rm neu}=\sum_{h\in\mathcal N}E_h$, and put $\mathcal E=E_{\rm neu}+E_{\rm st}=\sum_hE_h$. Let $\zeta$ denote all remaining tangential deviations from the limiting fixed component. Uniformly on the closed neighbourhood of the five-node face, \begin{equation*} \begin{aligned} q^+-q={}&-c_\theta d_\theta a^2 +\sum_{h\in\mathcal N}\varepsilon_h\rho_h^{\rm ext}E_h+\mathcal R_{\rm st}\\ &+O\!\left(d_\theta|a|^3+|a|\mathcal E+\mathcal E^2 +\|\zeta\|\mathcal E\right), \end{aligned} \tag{B.8.22}\label{eq:B.8.22} \end{equation*} where the $\rho_h^{\rm ext}$ are positive analytic units, \begin{equation*} \varepsilon_h=s\,\operatorname{sign}(h-\alpha), \end{equation*} and $\mathcal R_{\rm st}$ contains the exact stable first jets and satisfies $|\mathcal R_{\rm st}|\le C E_{\rm st}$. The pivot side of $\alpha$ means the component of $\mathbb R\setminus\{\alpha\}$ containing the pivot interval $I$. If every persistent neutral mode lies on the pivot side of $\alpha$, then eventually \begin{equation*} q^+-q\le-c(d_\theta a^2+E_{\rm neu}),\qquad q\ge0. \tag{B.8.24}\label{eq:B.8.24} \end{equation*} If all persistent neutral modes lie on one side of the pivot interval, then $\mathfrak A$ crosses zero at most once and $\theta_n$ converges. Stable modes remain in the exact recurrence in both conclusions.
\end{lemma}

\begin{proof}
Throughout this proof, $\rho_h^{\rm ext}$ in \eqref{eq:B.8.22}, for $h\in\mathcal N$, is evaluated at the current tangential base point on the fixed component. Its dependence on $C$, the free completion shear, the pivot, and the remaining tangential chart variables is suppressed. The same convention is used for $R$ and $\kappa_h$ in \eqref{eq:B.8.25}. Thus these coefficients vary with the current tangential base point.

Lemma~\ref{lem:8.2} supplies the first term of \eqref{eq:B.8.22}, and Lemma~\ref{lem:8.3} supplies the value and sign of every neutral coefficient at the limiting fixed component. The general calculation following \eqref{eq:B.8.18e} gives the actual coefficient at the current tangential base point; continuity preserves its sign, while the difference from its limiting value is included in $O(\|\zeta\|\mathcal E)$. The remaining external first jets are analytic and uniformly bounded; their exact contribution is $\mathcal R_{\rm st}=O(E_{\rm st})$. By \eqref{eq:B.5.25}, the map being expanded is the normalised-core coordinate induced by the exact full return; it agrees with the restricted-core return at $E=0$ and is analytic in $a,\zeta,E$. It is even in every signed external amplitude $z_h$, where $E_h=z_h^2$. Taylor's integral formula leaves only the displayed cubic intrinsic, $aE_h$, $E_hE_g$, and $\zeta E_h$ terms. The signed endpoint charts and the equality of symbolic and ordinary powers in Lemma~\ref{lem:4.3} keep their coefficients bounded through endpoints and factor collisions. In particular, the $O(\mathcal E)$ full/core transfer is not used to absorb the sign-bearing neutral linear terms: their limiting exact derivatives are those computed in Lemma~\ref{lem:8.3}.

A neutral mode on the pivot side of $\alpha$ satisfies $s\,\operatorname{sign}(h-\alpha)=-1$. Since $\rho_h^{\rm ext}\ge\rho_{\min}^{\rm ext}>0$ and $\zeta,a,\mathcal E\to0$, $$ d_\theta|a|^3\le\varepsilon d_\theta a^2,\qquad |a|\mathcal E+\mathcal E^2+\|\zeta\|\mathcal E \le\varepsilon\mathcal E $$ on a late enough tail. The stable linear terms are $O(E_{\rm st})=o(E_{\rm neu})$ by \eqref{eq:B.6.63}, and the same is true of the stable part of the displayed remainder. Absorbing them and taking $\varepsilon$ small gives the first assertion in \eqref{eq:B.8.24}. Completion convergence gives $q_n\to0$; an eventually nonincreasing sequence converging to zero is eventually nonnegative.

For the invariant half-line in the pivot normal, use $\mathfrak A$, not $a$. The unforced exact five-node return in Lemma~\ref{lem:7.3} gives $$ \mathfrak A^+=\mathfrak A\{1+O(\mathfrak A)\}. $$ Lemma~\ref{lem:8.3} and Taylor's integral formula give \begin{equation*} \mathfrak A^+ =R\mathfrak A+\sum_{h\in\mathcal N}\kappa_hE_h+\mathcal R_{\rm st}^{\rm piv} +O(\mathfrak A^2+|\mathfrak A|\mathcal E +\mathcal E^2+\|\zeta\|\mathcal E), \tag{B.8.25}\label{eq:B.8.25} \end{equation*} where $R\ge R_{\min}^{\rm piv}>0$ and $|\mathcal R_{\rm st}^{\rm piv}|\le C E_{\rm st}$. If $\theta$ is the lower second-kind root, then $D_J(\theta)<0$ and $\operatorname{sign}A_h(\theta)=\operatorname{sign}(h-\theta)$. If it is the upper root, $D_J(\theta)>0$ and $\operatorname{sign}A_h(\theta)=-\operatorname{sign}(h-\theta)$. Together with \eqref{eq:B.8.17b}, both cases give, for $h\in\mathcal N$, \begin{equation*} \operatorname{sign}\kappa_h=-\operatorname{sign}(h-\theta). \end{equation*} Neutral nodes on one side of $I$ have one common forcing sign. At the boundary of the corresponding closed half-line, that forcing is bounded below by a positive multiple of $E_{\rm neu}$; the stable terms and their nonlinear companions are $o(E_{\rm neu})$ by \eqref{eq:B.6.63}. Shrinking the neighbourhood therefore makes the half-line invariant under the exact recurrence \eqref{eq:B.8.25}, so $\mathfrak A$ crosses zero at most once.

The function $\overline B=d_\theta B$ extends as a positive analytic unit. For the pivot increment and its uniform remainder, define $$ F=\overline B(\theta^+-\theta) -d_\theta\frac{\mathfrak A+3\mathfrak A^+}{2}. $$ On the fixed component $F=0$. In the intrinsic $\mathfrak A$-normal, Lemma~\ref{lem:7.3} gives $$ d(\theta^+-\theta)=\frac{2}{B}\,d\mathfrak A,\qquad d\mathfrak A^+=d\mathfrak A, $$ so $dF=2d_\theta d\mathfrak A -d_\theta(4d\mathfrak A)/2=0$. In a neutral external direction at the limiting fixed component, \eqref{eq:B.8.21c} gives $$ B\,d(\theta^+-\theta)=\frac32\,d\mathfrak A^+, $$ and again $dF=0$. At a nearby tangential base point the corresponding neutral-source jet is $O(\|\zeta\|)$, while stable external first derivatives are uniformly bounded after division by $d_\theta$. Tangential derivatives along the fixed component vanish because the return is the identity there.

In one common signed endpoint chart, with coordinates $(y,\mathfrak A,z_h,\zeta)$, both terms defining $F$ are analytic functions of the same variables. Recovery gives $y^2=u d_\theta$ with $u$ a positive analytic unit. Each term vanishes at $y=0$, and their difference $F$ is even in $y$; hence $F(0,\cdot)=\partial_yF(0,\cdot)=0$ and $$ \frac{F(y,\cdot)}{y^2} =\int_0^1(1-t)\,\partial_y^2F(ty,\cdot)\,dt $$ is analytic. Multiplication by the unit $u$ proves that $F/d_\theta$ is analytic. The coefficient formulas for the return are denominator-free, and the symbolic-to-ordinary power argument applies in the collision charts. Hence this quotient has a uniform $C^2$ bound on a finite compact signed-chart cover, including simultaneous collisions. On the interval interior the first-derivative computations above show that the quotient has zero value and zero intrinsic first derivative at the limiting component; its neutral-external first derivatives are $O(\|\zeta\|)$ and its stable-external derivatives are uniformly bounded. Equality of analytic germs extends these statements to the endpoint, while the preceding integral formula supplies the uniform Hessian bound. Taylor's integral remainder gives \begin{equation*} \overline B(\theta^+-\theta) =d_\theta\left\{\frac{\mathfrak A+3\mathfrak A^+}{2}+R_2\right\}, \tag{B.8.27}\label{eq:B.8.27} \end{equation*} \begin{equation*} |R_2|\le C\bigl(\mathfrak A^2+\mathcal E^2+\|\zeta\|\mathcal E+E_{\rm st}\bigr). \tag{B.8.28}\label{eq:B.8.28} \end{equation*}

The exact full-state parity chord is uniformly coercive: \begin{equation*} \|Y^+-Y\|^2 \ge c\left(d_\theta\mathfrak A^2+\mathcal E\right). \tag{B.8.29}\label{eq:B.8.29} \end{equation*} To see the uniform Gram bound, let $\mathbf e_p$ be the disappearing retained signed-coordinate axis and let $\mathbf e_h$ be the mutually orthogonal external axes. In the interval interior the core-normal leading vector is an analytic vector $v_0$ with $\|v_0\|\ge c_0$. At an endpoint its scaled limit is $$ v_0^{\rm ep}=c_p\mathbf e_p+v_{\rm ret},\qquad |c_p|\ge c_0, $$ and it is multiplied by the generator $y\mathfrak A$, whose square is a positive unit times $d_\theta\mathfrak A^2$. For any surviving external signed coordinate $z_h$, the exact chord is $(\widetilde M_h-1)z_h\mathbf e_h$, with $\widetilde M_h\to m_h\ne1$. Dependence of the induced normalised core on $z_h$ starts with $E_h=z_h^2$, so its first derivative in $z_h$ vanishes at the fixed face. In the generator order $(\sqrt{d_\theta}\mathfrak A,(z_h)_h)$, the limiting Gram matrix is \begin{equation*} \Gamma_{\rm chord} =\operatorname{diag}\left(\|v_0^{\rm ep}\|^2,((m_h-1)^2)_h\right), \qquad \lambda_{\min}(\Gamma_{\rm chord}) \ge\min\left\{c_0^2,\min_h(m_h-1)^2\right\}>0. \end{equation*} If no external coordinate is present, the external block and its inner minimum are simply omitted. Neutral modes give the entry $4$ and stable modes give other positive entries. At zero external amplitude the reconstruction \eqref{eq:B.5.24} has the core derivative and the mutually orthogonal standard external axes, while \eqref{eq:B.5.25} makes all full/core cross terms higher order. The signed endpoint and collision generator changes have invertible positive-unit diagonal limits. Taylor's integral formula and a finite-cover singular-value minimum therefore prove \eqref{eq:B.8.29} uniformly through every endpoint and collision.

Lemma~\ref{lem:3.1} and \eqref{eq:B.8.29} imply \begin{equation*} \sum_n\left(d_{\theta,n}\mathfrak A_n^2+\mathcal E_n\right)<\infty. \tag{B.8.30}\label{eq:B.8.30} \end{equation*} After the possible single crossing, $\mathfrak A_n$ and $\mathfrak A_{n+1}$ have a common sign. The additional $E_{\rm st}$ term in \eqref{eq:B.8.28} is summable by \eqref{eq:B.8.30}. Equations \eqref{eq:B.8.27}--\eqref{eq:B.8.30} show that the total pivot variation opposite that sign is finite. Boundedness of $\theta_n\in[\ell,r]$ then makes the other directional variation finite as well. Hence $\theta_n$ converges.
\end{proof}

The single-crossing argument settles configurations in which the external forcing has one sign. Mixed configurations require a scalar for which the first-order intrinsic drift cancels; this is the purpose of the next two lemmas.

\begin{lemma}[Distinct-root Lagrange cocycle]\label{lem:8.5}
Assume $\alpha\notin J$. For each $z\in J\cup\{\alpha\}$, define \begin{equation*} \beta_z(h)=\frac{2C_*}{\mathscr D_*'(z)(z-h)}. \end{equation*} We abbreviate $\beta_i(h)=\beta_{\lambda_i}(h)$ for $\lambda_i\in J$. On a tail before any exact deletion, with all the displayed factors positive, define \begin{equation*} \mathcal S_h =E_h|g|^{\beta_\alpha(h)} \prod_{\lambda_i\in J}w_i^{\beta_{\lambda_i}(h)}. \tag{B.8.32}\label{eq:B.8.32} \end{equation*} Here $E_h$ is the exact full-state mass and $g,w_i$ are induced normalised-core quantities. Then \begin{equation*} \boxed{\displaystyle \log\frac{\mathcal S_h^+}{\mathcal S_h} =-\frac{4x}{C_*} +\frac{2D_J(h)}{C_*}e +O(x^2+e^2+a^2+\mathcal E).} \tag{B.8.33}\label{eq:B.8.33} \end{equation*} The remainder is uniform at both pivot endpoints and through every off-node factor collision.
\end{lemma}

\begin{proof}
Neutrality gives $\mathscr D_*(h)=-2C_*$. The quantity $$ \beta_z(h)=\frac{\mathscr D_*(h)}{\mathscr D_*'(z)(h-z)} $$ is the Lagrange cardinal at $h$ associated with the root $z$. Hence for every $R\in\mathcal P_5$, \begin{equation*} R(h)=\sum_{\lambda_i\in J}\beta_{\lambda_i}(h)R(\lambda_i) +\beta_\alpha(h)R(\alpha). \tag{B.8.34}\label{eq:B.8.34} \end{equation*}

Let $W_i$ be the full retained weights, $S=\sum_iW_i$, and $w_i=W_i/S$. If $\widetilde H_0,\widetilde H_1$ are the two exact full squared monic norms, the exact signed two-block multiplier is \begin{equation*} \widetilde M(t)=\frac{\widetilde P(t)\widetilde Q(t)}{\widetilde C} =1+\frac{\widetilde{\mathscr D}(t)}{\widetilde C},\qquad \widetilde C=\sqrt{\widetilde H_0\widetilde H_1}. \end{equation*} Therefore the exact full recurrence and the normalisation of the retained core give \begin{equation*} \begin{aligned} \log\frac{E_h^+\prod_i(w_i^+)^{\beta_i}} {E_h\prod_iw_i^{\beta_i}} ={}&2\log|\widetilde M(h)| +2\sum_i\beta_i\log|\widetilde M(\lambda_i)|\\ &+\left(\sum_i\beta_i\right)\log\frac S{S^+}. \end{aligned} \tag{B.8.36}\label{eq:B.8.36} \end{equation*} Here both sums are over $i\in J$. Since $S=1-\sum_kE_k$ and the finitely many full multipliers are bounded, the final term is $O(\mathcal E)$.

Set the external masses equal to zero and let $\dot{\mathscr D}$ denote a first variation of the untilded core completion, with $x=C-C_*$. At $E=0$ the full and core objects agree, and \eqref{eq:B.5.25} shows that their derivatives in core chart variables agree there. Since $\deg\dot{\mathscr D}\le5$, the core-direction derivative of the exact multiplier is \begin{equation*} \left.d\widetilde M(t)\right|_{E=0}=\frac{\dot{\mathscr D}(t)}{C_*} -\frac{\mathscr D_*(t)}{C_*^2}x. \end{equation*} At the five core nodes and at $\alpha$, $\widetilde M=1$, while at $h$, $\widetilde M=-1$. Differentiate \eqref{eq:B.8.36} only in these core directions and apply \eqref{eq:B.8.34} to $\dot{\mathscr D}$. The core evaluations cancel the evaluation at $h$, except for the missing cardinal at $\alpha$, and the scalar normalisation contributes the two equal radial terms. Explicitly, \begin{equation*} \left.d\mathcal B_h\right|_{E=0} =-\frac4{C_*}x -\frac{2\beta_\alpha(h)}{C_*}\dot{\mathscr D}(\alpha), \tag{B.8.38}\label{eq:B.8.38} \end{equation*} where $\mathcal B_h$ denotes the logarithm on the left of \eqref{eq:B.8.36}. Derivatives in an $E_k$-direction are the bounded Uvarov derivatives of Lemma~\ref{lem:8.3} and contribute to the displayed $O(\mathcal E)$ remainder in \eqref{eq:B.8.33}.

For the affine normal $U$, polynomial division at $x=0$ gives $$ \dot{\mathscr D}(t)=\operatorname{quo}_K(\mathscr D_*U)(t). $$ Since $\mathscr D_*=D_J(t)(t-\alpha)$, $$ \operatorname{quo}_K(\mathscr D_*U)(t) =\frac{D_J(t)U(t)-D_J(\theta)U(\theta)}{t-\theta}. $$ Using $U(\theta)=-a$ and $U(\alpha)=-ge$, \begin{equation*} \dot{\mathscr D}(\alpha) =-D_J(\alpha)e+\frac{D_J(\theta)}g\,a. \end{equation*} The $e$-term in \eqref{eq:B.8.38} is \begin{equation*} \frac{2\beta_\alpha(h)D_J(\alpha)}{C_*}e =\frac4{\alpha-h}e =\frac{2D_J(h)}{C_*}e, \end{equation*} where the last equality is $D_J(h)(h-\alpha)=-2C_*$. The remaining $a$-coefficient is \begin{equation*} -\frac{4D_J(\theta)} {D_J(\alpha)(\alpha-h)g}\,a. \tag{B.8.41}\label{eq:B.8.41} \end{equation*} The exact root shear is \begin{equation*} \theta^+-\theta =-\frac{2D_J(\theta)}{C_*}a+O(a^2). \end{equation*} Hence \begin{equation*} \begin{aligned} \beta_\alpha(h)\{\log|g^+|-\log|g|\} &=-\frac{\beta_\alpha(h)}{g}(\theta^+-\theta)+O(a^2)\\ &=\frac{4D_J(\theta)} {D_J(\alpha)(\alpha-h)g}\,a+O(a^2), \end{aligned} \tag{B.8.43}\label{eq:B.8.43} \end{equation*} which cancels \eqref{eq:B.8.41}.

For the uniform remainder, use the exact full multiplier expression \eqref{eq:B.8.36} together with the induced core recurrence \eqref{eq:B.5.26}. At a signed endpoint the normalised-core multiplier of the vanishing retained coordinate is near $+1$, so its logarithm is analytic; $g$ is a unit by Lemma~\ref{lem:8.1}. At a factor collision every full node evaluation $\widetilde P(\lambda_i),\widetilde Q(\lambda_i)$ stays nonzero because its product tends to $C_*>0$. Equations \eqref{eq:B.5.25}--\eqref{eq:B.5.26} make the hybrid full/core log-return $C^2$ on a finite compact signed-chart cover. Taylor's formula with integral remainder, after the explicitly cancelled core first derivatives \eqref{eq:B.8.38}--\eqref{eq:B.8.43}, and with the bounded external derivatives retained as $O(\mathcal E)$, gives $$ |R|\le C\{(|x|+|e|+|a|)^2+\mathcal E\}. $$ This is \eqref{eq:B.8.33}. After an exact deletion, the corresponding lower-support cocycle applies.
\end{proof}

Only now do we compare stable-mass deletion for this logarithmic quantity. Let $\mathcal B_h^{\rm Lag}$ denote the real-analytic extension of the logarithm of the positive return ratio in \eqref{eq:B.8.33}, and let $\mathcal B_{h,0}^{\rm Lag}$ be its value after the stable masses are set to zero and the remaining full weights are renormalised. This extended log-return is even in every stable signed amplitude and hence analytic in the corresponding squared masses. The mean-value theorem in those squared-mass variables on the finite compact cover gives
\begin{equation*} |\mathcal B_h^{\rm Lag}-\mathcal B_{h,0}^{\rm Lag}|\le C_hE_{\rm st}. \tag{B.8.43a}\label{eq:B.8.43a} \end{equation*}
This concerns the extended log-increment, not the generally unbounded logarithm of $\mathcal S_h$.

\begin{lemma}[Confluent Hermite cocycle]\label{lem:8.6}
Suppose $\alpha=c\in J$, and factor \begin{equation*} \mathscr D_*(t)=(t-c)^2R(t). \end{equation*} For a simple root $\lambda_i\ne c$, define $$ H_i(t)=\frac{\mathscr D_*(t)} {\mathscr D_*'(\lambda_i)(t-\lambda_i)}. $$ At the double root define \begin{equation*} H_{c,0}(t)= \left(1-(t-c)\frac{R'(c)}{R(c)}\right)\frac{R(t)}{R(c)}, \qquad H_{c,1}(t)=(t-c)\frac{R(t)}{R(c)}. \end{equation*} Define \begin{equation*} b_i(h)= \begin{cases} H_i(h),&\lambda_i\ne c,\\ H_{c,0}(h),&\lambda_i=c, \end{cases} \qquad \kappa_h=H_{c,1}(h) =\frac{2C_*}{R(c)(c-h)}. \end{equation*} Then the hybrid scalar \begin{equation*} \mathcal S_h^{\rm conf} =E_h\exp\!\left(\frac{\kappa_h}{c-\theta}\right) \prod_{i=1}^5w_i^{b_i(h)} \tag{B.8.47}\label{eq:B.8.47} \end{equation*} uses the exact full-state mass $E_h$ and induced normalised-core factors, and satisfies \begin{equation*} \boxed{\displaystyle \log\frac{(\mathcal S_h^{\rm conf})^+} {\mathcal S_h^{\rm conf}} =-\frac{4x}{C_*} +\frac{2D_J(h)}{C_*}e +O(x^2+e^2+a^2+\mathcal E).} \tag{B.8.48}\label{eq:B.8.48} \end{equation*}
\end{lemma}

\begin{proof}
For $\lambda_i\ne c$, $H_i$ is one at $\lambda_i$, zero at the other simple nodes, and has a double zero at $c$. Moreover $$ H_{c,0}(c)=1,\quad H_{c,0}'(c)=0,\qquad H_{c,1}(c)=0,\quad H_{c,1}'(c)=1, $$ and both functions vanish at the four simple nodes. For $F\in\mathcal P_5$, the difference between $F$ and $$ \sum_{\lambda_i\ne c}F(\lambda_i)H_i(t) +F(c)H_{c,0}(t)+F'(c)H_{c,1}(t) $$ has four simple zeros and a double zero at $c$, while its degree is at most five. It must be zero. Evaluation at $h$ gives \begin{equation*} F(h)=\sum_i b_i(h)F(\lambda_i)+\kappa_hF'(c). \tag{B.8.49}\label{eq:B.8.49} \end{equation*} Neutrality yields $(h-c)^2R(h)=-2C_*$, whence $$ H_{c,1}(h)=\frac{(h-c)R(h)}{R(c)} =\frac{2C_*}{R(c)(c-h)}. $$

Use the exact full multiplier identity \eqref{eq:B.8.36}, with the weights $b_i(h)$, set $E=0$, and apply \eqref{eq:B.8.49} to the core variation $\dot{\mathscr D}$. The displayed interpolation gives term by term \begin{equation*} \left.d\mathcal B_h\right|_{E=0}=-\frac4{C_*}x -\frac{2\kappa_h}{C_*}\dot{\mathscr D}'(c). \tag{B.8.50}\label{eq:B.8.50} \end{equation*} Let $g=c-\theta$. Since $K=(t-\theta)(t-c)$, $$ \operatorname{quo}_K(\mathscr D_*U) =\frac{F(t)-F(\theta)}{t-\theta}, \qquad F(t)=(t-c)R(t)U(t). $$ Now $$ F(c)=0,\qquad F'(c)=R(c)U(c)=-gR(c)e,\qquad F(\theta)=gR(\theta)a. $$ Differentiating the divided difference at $c$ gives \begin{equation*} \dot{\mathscr D}'(c) =\frac{F'(c)g-\{F(c)-F(\theta)\}}{g^2} =-R(c)e+\frac{R(\theta)}g\,a. \end{equation*} The term involving the remaining completion root in \eqref{eq:B.8.50} becomes \begin{equation*} \frac{2\kappa_hR(c)}{C_*}e =\frac4{c-h}e =\frac{2D_J(h)}{C_*}e. \end{equation*} The remaining normal coefficient is \begin{equation*} -\frac{4R(\theta)} {R(c)(c-h)(c-\theta)}a. \tag{B.8.53}\label{eq:B.8.53} \end{equation*} Because $D_J(\theta)=(\theta-c)R(\theta)$, the root shear is \begin{equation*} \theta^+-\theta =\frac{2(c-\theta)R(\theta)}{C_*}a+O(a^2). \end{equation*} Hence \begin{equation*} \frac{\kappa_h}{c-\theta^+}-\frac{\kappa_h}{c-\theta} =\frac{\kappa_h}{(c-\theta)^2}(\theta^+-\theta)+O(a^2) =\frac{4R(\theta)} {R(c)(c-h)(c-\theta)}a+O(a^2), \end{equation*} which cancels \eqref{eq:B.8.53}.

Lemma~\ref{lem:8.1} ensures $c\notin\{\ell,r\}$, so the exponential in \eqref{eq:B.8.47} is a bounded positive unit. At an endpoint $j$, $$ b_j(h)=\frac{2C_*}{\mathscr D_*'(j)(j-h)}, $$ so the endpoint powers coincide with those in Lemma~\ref{lem:8.5}. The exact full multiplier expression, the transfer \eqref{eq:B.5.25}, and the compact signed charts used after \eqref{eq:B.8.43} give the uniform remainder in \eqref{eq:B.8.48}, with all external derivatives included in $O(\mathcal E)$.
\end{proof}

Let $\mathcal B_h^{\rm conf}$ be the analytically extended log-return in \eqref{eq:B.8.48}, and append a subscript $0$ after deleting and renormalising the stable masses. This extension is even in every stable signed amplitude, hence analytic in the corresponding squared masses. The mean-value theorem in those squared masses gives
\begin{equation*} |\mathcal B_h^{\rm conf}-\mathcal B_{h,0}^{\rm conf}|\le C_hE_{\rm st}. \tag{B.8.55a}\label{eq:B.8.55a} \end{equation*}
Thus \eqref{eq:B.6.63} transfers every strict Lagrange or confluent log-increment estimate from the stable-free comparison state to the exact full state.

The distinct-root and confluent formulae have the same leading cocycle. Together with the single-crossing lemma, they now cover every possible position of the persistent neutral nodes.

\begin{lemma}[Convergence of five-node tails with a persistent neutral mode]\label{lem:8.7}
Every tail whose limiting completion has five nodes has a singleton squared-weight $\omega$-limit set. This includes a nonspectral or confluent remaining completion root, simple or simultaneous factor collisions, all positions of persistent neutral nodes, both pivot endpoints, and exact deletion into a four-node face.
\end{lemma}

\begin{proof}
Retain the stable modes in the exact recurrence; by \eqref{eq:B.6.63} their total forcing is $o(E_{{\rm neu},n})$. A neutral coordinate is either exactly deleted or remains strictly positive forever, because $$ E_{h,n+1}=\widetilde M_{h,n}^2E_{h,n},\qquad \widetilde M_{h,n}\to-1. $$ There are finitely many coordinates; pass beyond the last exact external deletion. If no neutral coordinate remains, Corollary~\ref{cor:7.4} applies.

At the pivot endpoints, \begin{equation*} \mathscr D_*'(\ell)>0,\qquad \mathscr D_*'(r)<0. \tag{B.8.56}\label{eq:B.8.56} \end{equation*} For the lower pivot interval, $D_J'(\ell)<0<D_J'(r)$ and $\ell-\alpha,r-\alpha<0$; for the upper interval, $D_J'(\ell)>0>D_J'(r)$ and $\ell-\alpha,r-\alpha>0$. This is \eqref{eq:B.8.56}. Since no neutral node lies in $I$, the endpoint exponents in Lemmas~\ref{lem:8.5}--\ref{lem:8.6} have signs \begin{equation*} \begin{array}{c|cc} &\text{exponent at }\ell&\text{exponent at }r\\ \hline h<\ell&+&-\\ h>r&-&+ \end{array}. \tag{B.8.57}\label{eq:B.8.57} \end{equation*}

On the nonstationary tail under consideration, with orbit subscripts suppressed in the following local bounds, Lemmas~\ref{lem:5.2}--\ref{lem:5.3} give \begin{equation*} x<0,\qquad -x\ge cX,\qquad |x|+|e|+|a|=O(X),\qquad \mathcal E=O(X^2). \end{equation*} Here $x,e,a$ are core coordinates and \eqref{eq:B.5.27} absorbs their full-state differences into $O(X^2)$. If $D_J(h)e\ge0$, \eqref{eq:B.8.33} or \eqref{eq:B.8.48}, together with the post-extension stable estimates \eqref{eq:B.8.43a} and \eqref{eq:B.8.55a}, implies \begin{equation*} \log\frac{\mathcal S_{h,n+1}}{\mathcal S_{h,n}} \ge c_hX_n>0 \end{equation*} on the exact full tail for all sufficiently large $n$.

At an interior cluster point every factor in $\mathcal S_h$, other than $E_h$, is bounded above and below, while $E_h\to0$. The scalar tends to zero along that subsequence. At a pivot endpoint with positive exponent, exactly one core weight tends to zero and every other factor is a unit, so the scalar again tends to zero. An eventually increasing positive sequence is bounded below by its positive value at the start of the monotone tail and cannot have a subsequence of either kind. Hence such a mode with $h<\ell$ forces $\theta_n\to r$, whereas one with $h>r$ forces $\theta_n\to\ell$.

For a pair $$ h_-<\alpha<h_+,\qquad q_-=\alpha-h_-,\qquad q_+=h_+-\alpha, $$ neutrality gives $$ D_J(h_-)=\frac{2C_*}{q_-},\qquad D_J(h_+)=-\frac{2C_*}{q_+}, $$ and hence \begin{equation*} q_-D_J(h_-)+q_+D_J(h_+)=0. \end{equation*} The paired scalar \begin{equation*} \mathcal T=\mathcal S_{h_-}^{q_-}\mathcal S_{h_+}^{q_+} \end{equation*} uses the confluent versions when $\alpha\in J$. Its transverse first-order term cancels; its radial term $$ -\frac{4(q_-+q_+)x}{C_*} $$ is positive and dominates the $O(X^2)$ remainder. The full/core transfer \eqref{eq:B.5.27} contributes only $O(X^2)$, while the stable-mass comparisons \eqref{eq:B.8.43a} and \eqref{eq:B.8.55a}, applied to the two already extended log-returns, contribute $o(E_{\rm neu})=o(X^2)$. Its exponent at a pivot endpoint $j$ is \begin{equation*} \begin{aligned} \gamma_j &=\frac{2C_*}{D_J'(j)(j-\alpha)} \left(\frac{q_-}{j-h_-}+\frac{q_+}{j-h_+}\right)\\ &=\frac{2C_*(h_+-h_-)} {D_J'(j)(j-h_-)(j-h_+)}. \end{aligned} \tag{B.8.62}\label{eq:B.8.62} \end{equation*} For the second equality, reduction of the numerator in parentheses to a common denominator gives $$ \frac{(h_+-h_-)(j-\alpha)} {(j-h_-)(j-h_+)}. $$ If $h_-<\ell<r<h_+$, the two endpoint exponents have opposite signs, so the paired scalar leaves exactly one possible pivot endpoint.

Suppose first that $\alpha>r$. Every neutral node belongs to exactly one of \begin{equation*} I_L=(-\infty,\ell),\qquad I_M=(r,\alpha),\qquad I_R=(\alpha,\infty). \end{equation*} If no $I_R$-mode persists, all modes are on the pivot side of $\alpha$. Here $s=1$, so Lemma~\ref{lem:8.4} gives $e\ge0$, and neutrality gives $$ D_J(h)=\frac{-2C_*}{h-\alpha}>0. $$ Every individual scalar is favourable. An $I_L$-mode selects $r$, while an $I_M$-mode selects $\ell$; the two families cannot both persist. The surviving family fixes the pivot limit.

If only $I_R$-modes persist, all nodes lie to the right of $I$, and the single-crossing conclusion of Lemma~\ref{lem:8.4} gives convergence. If modes persist on both sides of $\alpha$ and an $I_L$-mode is present, pair it with an $I_R$-mode; the pair straddles the pivot interval, so \eqref{eq:B.8.62} selects one endpoint. If no $I_L$-mode is present, all modes lie in $I_M\cup I_R$, to the right of $I$, and the single-crossing conclusion applies. These cases are disjoint and exhaustive.

If $\alpha<\ell$, use the reflected partition \begin{equation*} I_L=(-\infty,\alpha),\qquad I_M=(\alpha,\ell),\qquad I_R=(r,\infty). \end{equation*} When $I_L$ is absent, all modes are on the pivot side of $\alpha$. Now $s=-1$, so Lemma~\ref{lem:8.4} gives $e\le0$, while $D_J(h)<0$; hence $D_J(h)e\ge0$. Modes in $I_M$ and $I_R$ select opposite endpoints and cannot both persist. If only $I_L$-modes persist, all lie to the left of $I$, and the single-crossing estimate applies. In a configuration straddling $\alpha$, either an $I_R$-mode is present and pairs with an $I_L$-mode across the pivot interval, or every mode lies to the left of $I$ and the single-crossing estimate applies. This exhausts the reflected case.

Finally, consider exact deletion of a pivot-endpoint coordinate. Once that squared coordinate is zero, coordinate-factor invariance keeps it zero. Let $J_4$ be the remaining four completion nodes, and define $$ E_n=\|(\operatorname{Id}-\Pi_{J_4})Y_n\|^2,\qquad \widehat Y_n=\frac{\Pi_{J_4}Y_n}{\sqrt{1-E_n}}. $$ The exact two-block map on $J_4$ is the identity, while \eqref{eq:B.5.6} makes the total outside squared mass summable. Uniform Gram regularity and the mean-value theorem give $$ \|\widehat Y_{n+1}-\widehat Y_n\| \le C(E_n+E_{n+1}),\qquad \sum_nE_n<\infty. $$ Hence $\widehat Y_n$ converges and every outside coordinate tends to zero.

If no exact deletion occurs, the positional classification gives $\theta_n\to\theta_*$, while Lemmas~\ref{lem:5.1} and~\ref{lem:5.3} give $$ P_n\to P_*,\qquad Q_n\to Q_*,\qquad a_n^2\to C_*. $$ The five-node recovery formula $$ w_{n,i} =\frac{a_n^2D_J'(\lambda_i)^{-1} (\lambda_i-\theta_n)} {P_n(\lambda_i)} $$ then gives convergence of every retained squared weight. The denominator has nonzero limit because $P_*(\lambda_i)Q_*(\lambda_i)=C_*>0$. This also covers asymptotic approach to a pivot endpoint, where exactly one limiting weight is zero. Hence the squared-weight $\omega$-limit set is a singleton.
\end{proof}

\section{Neutral external modes at a six-node limit}\label{sec:9}

When $|S_*|=6$, the two completion roots range over a closed rectangle. A persistent neutral mode forces these roots with definite endpoint signs. Lemmas~\ref{lem:9.1}--\ref{lem:9.3} combine those signs with a logarithmic cocycle to prove convergence away from two support-boundary faces. Lemma~\ref{lem:9.4} treats those faces in coordinates compatible with the five-node limit and with exact support loss.

Let \begin{equation*} D(t)=\prod_{i=1}^6(t-\lambda_i),\qquad K=(t-\rho)(t-\sigma), \end{equation*} with \begin{equation*} (\rho,\sigma)\in [\lambda_2,\lambda_3]\times[\lambda_4,\lambda_5]. \end{equation*} The limiting completion is $\mathscr D_*=D$. Stable modes remain in the exact full state, but \eqref{eq:B.6.63} makes their forcing negligible relative to any persistent neutral mass. Every persistent nonstable external mode is neutral and satisfies \begin{equation*} D(h)=-2C_*. \end{equation*}

\begin{lemma}[External forcing at a six-node limit]\label{lem:9.1}
Let $k_P(\,\cdot\,,h)$, $k_Q(\,\cdot\,,h)$ be the degree-two reproducing kernels of the two positive fixed phase measures. The affine normal forced by inserting squared mass at $h$ is \begin{equation*} A_h(t)=4P(h)\operatorname{quo}_{D} \bigl(K(t)\{P(t)k_Q(t,h)-Q(t)k_P(t,h)\}\bigr). \tag{B.9.4}\label{eq:B.9.4} \end{equation*} It satisfies the exact endpoint signs \begin{equation*} \operatorname{sign}A_h(\rho)=\operatorname{sign}(h-\rho),\qquad \operatorname{sign}A_h(\sigma)=\operatorname{sign}(\sigma-h). \tag{B.9.5}\label{eq:B.9.5} \end{equation*}
\end{lemma}

\begin{proof}
The Uvarov differentiation is carried out coefficient by coefficient in Appendix B.A1.1. There $P,Q$ are the fixed normalised-core polynomials, while $\widetilde P_\varepsilon,\widetilde Q_\varepsilon$ are the polynomials of the state after an external squared mass has been inserted. Equations \eqref{eq:B.A1.1}--\eqref{eq:B.A1.5} derive $$ p=-P(h)k_P,\qquad q=3P(h)k_Q $$ from the two monic orthogonality systems. Equations \eqref{eq:B.A1.6}--\eqref{eq:B.A1.9} then differentiate the normalised retained core induced by the exact full return and give $$ p_{\rm out}=-4P(h)k_P,\qquad q_{\rm out}=4P(h)k_Q. $$ Finally, \eqref{eq:B.A1.10} is the differentiated factor identity; polynomial division by $D$ gives \eqref{eq:B.A1.11}, which is exactly \eqref{eq:B.9.4}. As a polynomial in coefficient coordinates, the quotient extends through factor collisions.

Appendix B.A1.2 proves the signs. The residue identities \eqref{eq:B.A1.12}--\eqref{eq:B.A1.13} identify $K$ as the common second-kind polynomial. The positive Jacobi recurrences are expanded in \eqref{eq:B.A1.14}--\eqref{eq:B.A1.17}, and Christoffel--Darboux is reduced explicitly in \eqref{eq:B.A1.18}--\eqref{eq:B.A1.25}. This yields the formulas with explicit signs $$ A_h(0)=\frac4h \left\{\frac{P(h)^2}{bc}+\frac{C_*}{Bf}\right\}, $$ $$ A_h(1)=-\frac4{h-1} \left\{\frac{P(h)^2}{b(1-c)} +\frac{C_*}{B(1-f)}\right\}. $$ These are \eqref{eq:B.A1.26}--\eqref{eq:B.A1.27}, and their braces are strictly positive. Undoing the affine map $t\mapsto(t-\rho)/(\sigma-\rho)$ gives \eqref{eq:B.9.5}. All Jacobi coefficients occurring in the denominators are strictly positive for the positive six-node measures, so these formulas remain finite at every off-node factor collision.
\end{proof}

The signs in Lemma~\ref{lem:9.1} describe the local forcing of the root coordinates. A global restriction on their accumulation set comes from the neutral cocycle below.

\begin{lemma}[Neutral cocycle at a six-node limit]\label{lem:9.2}
Define $c_i=1/D'(\lambda_i)$ and \begin{equation*} \beta_i(h)=\frac{2C_*c_i}{\lambda_i-h}. \tag{B.9.13}\label{eq:B.9.13} \end{equation*} As long as the six normalised-core coordinates and the exact external mass $E_h$ are positive, define the hybrid scalar \begin{equation*} \mathcal J_h =E_h\prod_{i=1}^6|K(\lambda_i)|^{\beta_i(h)}. \tag{B.9.14}\label{eq:B.9.14} \end{equation*} Here $K$ belongs to the normalised core induced by the exact full state. On every sufficiently late nonstationary tail, \begin{equation*} \log\frac{\mathcal J_{h,n+1}}{\mathcal J_{h,n}} \ge c_hX_n>0. \tag{B.9.15}\label{eq:B.9.15} \end{equation*} In particular, \begin{equation*} \prod_{i=1}^6|K_n(\lambda_i)|^{\beta_i(h)} \longrightarrow+\infty. \tag{B.9.16}\label{eq:B.9.16} \end{equation*} The logarithmic return in \eqref{eq:B.9.15} extends analytically through every support boundary and every simple or simultaneous factor collision.
\end{lemma}

\begin{proof}
Neutrality and the partial-fraction identity for $1/D$ give $$ \beta_i(h)=\frac{D(h)c_i}{h-\lambda_i}. $$ These are the six Lagrange cardinal values at $h$, so $\sum_i\beta_i(h)=1$.

The exact root and cocycle calculation is in Appendix B.A1.3. For an affine intrinsic normal $u$, equations \eqref{eq:B.A1.28}--\eqref{eq:B.A1.31} perform the Euclidean division $$ Du=Kq_u+r_u $$ and derive, by interpolation at $\rho,\sigma,h$, the two root velocities and the exact formula for $q_u(h)$. Equation \eqref{eq:B.A1.32} differentiates the exact full external multiplier at zero external mass in this core direction, and \eqref{eq:B.A1.33}--\eqref{eq:B.A1.36} integrate that one-form to $$ \Phi_h=-2C_*\sum_i\frac{c_i}{\lambda_i-h} \log|K(\lambda_i)|. $$ Equation \eqref{eq:B.A1.37} shows that the root product in \eqref{eq:B.9.14} is $e^{-\Phi_h}$; its first intrinsic normal variation cancels the full multiplier variation exactly at $E=0$.

To continue through a support boundary, equation \eqref{eq:B.A1.38} expresses each normalised-core root-factor quotient under the exact full return as a product of analytic squared-weight, polynomial-evaluation, and norm quotients. In a signed boundary chart, $w_i=\chi_i^2$ and $|K(\lambda_i)|=\chi_i^2u_i$ with $u_i>0$ analytic. The sign of $K(\lambda_i)$ is fixed on that chart, so the quotient has an analytic positive extension. The coefficient identity holds on the dense coprime locus and is polynomial after its nonvanishing node denominators are cleared; reducedness from Lemma~\ref{lem:4.3} extends it through simple and simultaneous collisions. Together with the exact full multiplier and the transfer \eqref{eq:B.5.25}, this makes the hybrid log-return $C^2$ on the finite signed-chart cover. Its external first derivatives are the bounded Uvarov derivatives and are retained in the $O(\mathcal E)$ term.

On a fixed component at level $C$, where full and core states agree, the exact logarithmic return is \eqref{eq:B.A1.39}: $$ 2\log\left|1-\frac{2C_*}{C}\right|. $$ Its radial derivative at $C_*$ is $-4/C_*$; the preceding calculation shows that all tangent and intrinsic normal derivatives vanish. Taylor's integral remainder on the finite compact signed-chart cover gives \eqref{eq:B.A1.40}, $$ \log\frac{\mathcal J_h^+}{\mathcal J_h} =-\frac{4(C-C_*)}{C_*} +O\{(C-C_*)^2+\|U\|^2+\mathcal E\}. $$ Here $C,U$ are normalised-core coordinates; \eqref{eq:B.5.25}--\eqref{eq:B.5.27} place the exact full return in these coordinates with a $C^2$ error $O(\mathcal E)$ in core directions, while its generally nonzero external derivatives remain in the displayed $O(\mathcal E)$ term.

Only after this real-analytic extension of the log-return has been established do we delete stable masses for comparison. Let $\mathcal B_h^{(6)}=\log(\mathcal J_h^+/\mathcal J_h)$ and let $\mathcal B_{h,0}^{(6)}$ denote its value after the stable masses are set to zero and the remaining exact weights are renormalised. The extended log-return is even in the stable signed amplitudes, so the mean-value theorem in their squared masses on the finite compact cover gives
\begin{equation*} |\mathcal B_h^{(6)}-\mathcal B_{h,0}^{(6)}|\le C_hE_{\rm st}. \tag{B.9.22a}\label{eq:B.9.22a} \end{equation*}
This concerns the log-return, not the generally unbounded logarithm of $\mathcal J_h$. By \eqref{eq:B.6.63}, the right-hand side is $o(E_{\rm neu})=o(X^2)$ whenever a neutral mode persists.

Suppressing orbit subscripts in the preceding estimate, Lemmas~\ref{lem:5.2}--\ref{lem:5.3} give $$ C-C_*<0,\quad -(C-C_*)\ge cX,\quad \|U\|=O(X),\quad \mathcal E=O(X^2). $$ The radial term therefore dominates on the exact full tail and proves \eqref{eq:B.9.15}. Since $E_{h,n}\to0$ while the positive scalar $\mathcal J_{h,n}$ is eventually bounded below by its value at the beginning of the monotone tail, $$ \prod_i|K_n(\lambda_i)|^{\beta_i(h)} =\frac{\mathcal J_{h,n}}{E_{h,n}}\longrightarrow+\infty. $$ Equation \eqref{eq:B.9.16} follows.
\end{proof}

If a root fails to converge, \eqref{eq:B.9.16} can diverge only through a vanishing boundary factor with negative exponent. The following sign analysis reduces that possibility to two faces.

\begin{lemma}[Classification away from the remaining boundary faces]\label{lem:9.3}
The three components of $\{D<0\}$ are \begin{equation*} I_L=(\lambda_1,\lambda_2),\quad I_M=(\lambda_3,\lambda_4),\quad I_R=(\lambda_5,\lambda_6). \end{equation*} For the four moving faces, the signs of \eqref{eq:B.9.13} are \begin{equation*} \begin{array}{c|rrrr} &\beta_2&\beta_3&\beta_4&\beta_5\\ \hline I_L&+&-&+&-\\ I_M&-&+&+&-\\ I_R&-&+&-&+. \end{array} \tag{B.9.24}\label{eq:B.9.24} \end{equation*} Every six-node tail with a persistent neutral mode has convergent roots except possibly \begin{equation*} \text{an $I_L$-mode with }\sigma_n\to\lambda_5, \quad\text{or}\quad \text{an $I_R$-mode with }\rho_n\to\lambda_2. \tag{B.9.25}\label{eq:B.9.25} \end{equation*}
\end{lemma}

\begin{proof}
In the formula $\beta_i(h)=2C_*/\{D'(\lambda_i)(\lambda_i-h)\}$, the signs of $D'(\lambda_i)$ alternate. This gives the sign table \eqref{eq:B.9.24}. The resulting cluster alternatives are listed in \eqref{eq:B.A1.47}: \begin{equation*} \begin{array}{c|c} h\in I_L&\rho=\lambda_3\ \hbox{or}\ \sigma=\lambda_5,\\ h\in I_M&\rho=\lambda_2\ \hbox{or}\ \sigma=\lambda_5,\\ h\in I_R&\rho=\lambda_2\ \hbox{or}\ \sigma=\lambda_4. \end{array} \end{equation*} Indeed, \eqref{eq:B.9.16} can diverge only when a vanishing root factor has negative exponent.

Appendix B.A1.4 derives the required one-sided convergence from the exact face relation $y^2=u_\sigma d_\sigma$ and analytic division along all fixed components. Equations \eqref{eq:B.A1.42}--\eqref{eq:B.A1.43} give $$ z_\rho^+ =z_\rho R_\rho+d_\sigma z_\sigma^2q_\rho +\sum_h E_h\,a_{h,\rho}, $$ where, on a sufficiently late tube, $$ R_\rho>0,\qquad q_\rho>0,\qquad a_{h,\rho}>0\quad\hbox{for persistent neutral $h>\rho$}. $$ The first two signs are \eqref{eq:B.4.14} and \eqref{eq:B.6.24}; the last is the neutral kernel sign \eqref{eq:B.9.5}. The coefficients belonging to stable modes are uniformly bounded and contribute $O(E_{\rm st})$. The root increment is the exact divided formula \eqref{eq:B.A1.45}, \begin{equation*} \rho^+-\rho=d_\rho \left\{A_\rho z_\rho+B_\rho d_\sigma z_\sigma^2 +\sum_hB_{h,\rho}E_h\right\}, \qquad A_\rho<0, \tag{B.9.25b}\label{eq:B.9.25b} \end{equation*} with bounded analytic $B$'s. Equation \eqref{eq:B.A1.46}, which combines \eqref{eq:B.6.38} with the exact full-state chord of every surviving external coordinate, gives \begin{equation*} \sum_n\{d_{\rho,n}z_{\rho,n}^2+ d_{\sigma,n}z_{\sigma,n}^2+\mathcal E_n\}<\infty. \tag{B.9.25c}\label{eq:B.9.25c} \end{equation*}

If no persistent neutral node lies in $I_L$, every persistent neutral node has $h>\rho$. Their contribution in \eqref{eq:B.A1.42} has a uniform positive sign. Stable modes remain in the exact formula, but their $O(E_{\rm st})=o(E_{\rm neu})$ contribution is eventually absorbed, so the half-line $z_\rho\ge0$ is forward invariant and $z_\rho$ crosses zero at most once. In \eqref{eq:B.9.25b}, the total variation opposite the sign of its leading term is bounded by the summable last two core terms in \eqref{eq:B.9.25c}, together with the geometrically summable stable forcing. Since $\rho$ is bounded, the other directional variation is also finite, exactly as in \eqref{eq:B.6.40}. Hence $\rho_n$ converges. Reversing the node order proves: if no persistent neutral node lies in $I_R$, then $\sigma_n$ converges.

Middle-only forcing makes both roots converge. If an $I_L$-mode persists but no $I_R$-mode does, $\sigma_n\to\sigma_*$. When $\sigma_*<\lambda_5$, the $I_L$-row of \eqref{eq:B.9.24}, together with \eqref{eq:B.9.16}, forces $\rho_n\to\lambda_3$. The only remaining case is $\sigma_n\to\lambda_5$, the first remaining boundary face in \eqref{eq:B.9.25}. Reversal gives the second boundary face.

Suppose both $h_L\in I_L$ and $h_R\in I_R$ persist. Their individual cocycles restrict every cluster point to $$ (\{\rho=\lambda_3\}\cup\{\sigma=\lambda_5\}) \cap (\{\rho=\lambda_2\}\cup\{\sigma=\lambda_4\}), $$ which consists of the inner corner $(\lambda_3,\lambda_4)$ and the outer corner $(\lambda_2,\lambda_5)$. Define $$ t_i=\frac{\lambda_i-h_L}{h_R-\lambda_i}. $$ The derivative of $(x-h_L)/(h_R-x)$ is $(h_R-h_L)/(h_R-x)^2>0$, so \begin{equation*} 0<t_2<t_3<t_4<t_5. \tag{B.9.28}\label{eq:B.9.28} \end{equation*} For $p,q>0$ and $r=p/q$, factor $2C_*/|D'(\lambda_i)|$ from $p\beta_i(h_L)+q\beta_i(h_R)$. The remaining numerator is $q(r-t_i)$ at $i=2,4$ and $q(t_i-r)$ at $i=3,5$. Hence \begin{equation*} p\beta_i(h_L)+q\beta_i(h_R)>0 \Longleftrightarrow \begin{cases} r>t_i,&i=2,4,\\ r<t_i,&i=3,5. \end{cases} \tag{B.9.29}\label{eq:B.9.29} \end{equation*} Choose $t_2<r<t_5$. At the outer corner the two vanishing root factors have positive combined exponents, while $E_{h_L}^pE_{h_R}^q\to0$; all other factors are units. The eventually increasing positive scalar $\mathcal J_{h_L}^p\mathcal J_{h_R}^q$ would tend to zero, contradicting its positive lower bound. Only the inner corner remains, so both roots converge to $(\lambda_3,\lambda_4)$.

If a middle mode also persisted, its exponents at the two vanishing inner-corner factors are both positive by the $I_M$-row of \eqref{eq:B.9.24}. Its own increasing scalar would tend to zero, another contradiction. The alternatives above are exhaustive and prove the lemma.
\end{proof}

It remains to analyse the two faces in \eqref{eq:B.9.25}, where the six-node coordinates meet the five-node strata.

\begin{lemma}[Reduction at the remaining boundary faces]\label{lem:9.4}
The two residual cases in \eqref{eq:B.9.25} have convergent roots and squared core weights. The assertion is uniform through the pivot endpoints, simple or simultaneous factor collisions, and exact deletion into a five-node face.
\end{lemma}

\begin{proof}
Consider first the top face $\sigma_n\to\lambda_5$. Define $$ \alpha=\lambda_5,\qquad J=\{\lambda_1,\lambda_2,\lambda_3,\lambda_4,\lambda_6\}, \qquad \theta=\rho. $$ Hence \begin{equation*} \mathscr D_*(t)=D_J(t)(t-\alpha),\qquad \lambda_2\le\theta\le\lambda_3. \end{equation*} The six-node recovery formula at the omitted node is \begin{equation*} w_\alpha= \frac{H_0c_\alpha(\alpha-\rho)(\alpha-\sigma)} {P(\alpha)}. \end{equation*} All factors except $\alpha-\sigma$ are nonzero analytic units. In a signed chart, with $w_\alpha=y^2$, \begin{equation*} \alpha-\sigma=u_y y^2,\qquad u_y>0\quad\hbox{analytic}. \tag{B.9.30b}\label{eq:B.9.30b} \end{equation*} For the affine completion normal, let \begin{equation*} a=-U(\theta),\qquad e=-\frac{U(\alpha)}{\alpha-\theta}. \end{equation*} This is a uniformly invertible analytic coordinate change because $\alpha-\theta\ge\lambda_5-\lambda_3>0$. On $y=0$ the five-node fixed component is $a=0$, with $e$ its transverse shear; on the positive six-node fixed component $a=e=0$, with $y$ free. Their union ideal is \begin{equation*} \mathcal I_{\rm hf}=(a,ye). \end{equation*} The quantitative remainder follows from the exact cocycle calculation below.

Define \begin{equation*} \delta=\alpha-\sigma=u_y y^2,\qquad g=\alpha-\theta,\qquad \Delta=\sigma-\theta=g-\delta. \end{equation*} Both $g$ and $\Delta$ are uniformly positive units. Let $x=C-C_*$, with $C$ the normalised-core product. For the hybrid distinct-root scalar \eqref{eq:B.8.32}, let $\mathcal B_h=\log(\mathcal S_h^+/\mathcal S_h)$; its external factor follows the exact full recurrence, while its retained factors use the induced normalised core. At zero external mass the two states agree. The exact calculation in \eqref{eq:B.A1.59}--\eqref{eq:B.A1.63} gives, on the positive six-node fixed component, \begin{equation*} d\mathcal B_h =-\frac4{C_*}\,dx +\frac{2\beta_\alpha(h)D_J(\sigma)}{C_*}\,de +\frac{2\beta_\alpha(h)\delta D_J(\sigma)} {C_*g\Delta}\,da +\sum_k\chi_{hk}(y)\,dE_k. \end{equation*} Here the $\chi_{hk}$ are the uniformly bounded exact Uvarov derivatives, not full/core error terms. Define \begin{equation*} R_h^{\rm Lag} =\mathcal B_h+\frac{4x}{C_*} -\frac{2D_J(h)}{C_*}e. \end{equation*} Neutrality and $\mathscr D_*(h)=D_J(h)(h-\alpha)=-2C_*$ give $$ \frac{2D_J(h)}{C_*} =\frac{2\beta_\alpha(h)D_J(\alpha)}{C_*}. $$ The two nonzero jets are \begin{equation*} A_h^{\rm coc}(y,\vartheta) :=\left.\partial_aR_h^{\rm Lag}\right|_{\rm six} =\frac{2\beta_\alpha(h)\delta D_J(\sigma)} {C_*g\Delta}, \end{equation*} \begin{equation*} B_h^{\rm coc}(y,\vartheta) :=\left.\partial_eR_h^{\rm Lag}\right|_{\rm six} =\frac{2\beta_\alpha(h)}{C_*} \{D_J(\sigma)-D_J(\alpha)\}. \end{equation*} Here $\vartheta$ denotes the compact tangential factor, root, and norm parameters. For $y\ne0$, these jets need not vanish, but their anisotropic size is quadratic in $y$. The identity \begin{equation*} D_J(\sigma)-D_J(\alpha) =-\delta\int_0^1D_J'(\alpha-t\delta)\,dt \end{equation*} together with \eqref{eq:B.9.30b} gives \begin{equation*} |A_h^{\rm coc}(y,\vartheta)|+|B_h^{\rm coc}(y,\vartheta)|\le C_hy^2. \tag{B.9.30k}\label{eq:B.9.30k} \end{equation*}

Subtracting these jets, define \begin{equation*} \widetilde R_h :=R_h^{\rm Lag}-A_h^{\rm coc}(y,\vartheta)a-B_h^{\rm coc}(y,\vartheta)e. \end{equation*} The exact full multiplier formula \eqref{eq:B.8.36}, the induced core recurrence \eqref{eq:B.5.26}, the $C^2$ transfer \eqref{eq:B.5.25}, and Taylor's integral formula on the finite signed-chart cover give \begin{equation*} |\widetilde R_h| \le C_h\{(|x|+|a|+|e|)^2+\mathcal E\}. \tag{B.9.30m}\label{eq:B.9.30m} \end{equation*} Appendix B.A1.5.3 proves \eqref{eq:B.9.30m}, including all mixed normal terms and every external squared weight, at both pivot endpoints and every simple or simultaneous factor collision. In particular, its constant is independent of $y\downarrow0$.

For the orbit, let $\mathfrak a_n=-U_n(\theta_n)$ denote the indexed core transverse coordinate and write $\mathcal E_n=E_{{\rm neu},n}+E_{{\rm st},n}$. On the nonstationary tail under consideration, Lemmas~\ref{lem:5.2}--\ref{lem:5.3} give \begin{equation*} y_n^2\to0,\qquad |x_n|+|\mathfrak a_n|+|e_n|\le C X_n,\qquad \mathcal E_n=O(X_n^2),\qquad -x_n\ge cX_n. \tag{B.9.30n}\label{eq:B.9.30n} \end{equation*} Equations \eqref{eq:B.9.30k}--\eqref{eq:B.9.30m} give \begin{equation*} R_{h,n}^{\rm Lag}=o(X_n). \tag{B.9.30o}\label{eq:B.9.30o} \end{equation*}

The transverse sign admits an estimate of the same order. Let $v=y^2$. Define the normal coordinate for the six-node normalised core at $\sigma$ by $$ d_\theta=(\theta-\lambda_2)(\lambda_3-\theta), \qquad q=-\frac{U(\sigma)}{\Delta} =e+\frac{\delta}{g\Delta}a. $$ The moment normal for the whole six-node core is \begin{equation*} \mathfrak A_y=\sum_{i\in J\cup\{\alpha\}} \frac{w_i}{\lambda_i-\theta}. \end{equation*} The $\alpha$-summand is essential when $y>0$. Appendix B.A1.5.4 constructs a bounded analytic function $\kappa$ and the adapted transverse coordinate \begin{equation*} \widehat q=q+\kappa(v,\vartheta)\,v q \end{equation*} such that \begin{equation*} \widehat q_{n+1}-\widehat q_n \le-c\{d_{\theta,n}\mathfrak a_n^2+E_{{\rm neu},n}\}+\mathfrak r_n, \qquad \sum_{k=n}^{\infty}|\mathfrak r_k|=o(X_n). \tag{B.9.30r}\label{eq:B.9.30r} \end{equation*} The exact law for $v$ uses the normalised-core signed multiplier defined in \eqref{eq:B.A1.78b}. At zero external mass its leading coefficients are \begin{equation*} \ell_v=-\frac{2D_J(\sigma)}{C_*}, \end{equation*} and $\ell_{v,a}=2D_J(\theta)/(C_*\Delta)$. The first is a positive analytic unit. If $f_2$ is the coefficient of $vq^2$ in the exact transverse return, $\kappa=-f_2/\ell_v$ on the fixed set, and \begin{equation*} v^+-v=v\{\ell_v q+\ell_{v,a} a+O(a^2+q^2+\mathcal E)\}. \end{equation*} Thus the increment of $\kappa vq$ cancels the entire $vq^2f_2$ term. Signed endpoint parity supplies the factor $d_\theta$ in every mixed $a$-term. The neutral external terms have the favourable sign; stable linear terms have tail $O(E_{{\rm st},n})$ and are included in $\mathfrak r_n$ using \eqref{eq:B.6.63}. The remaining terms are absorbed into $d_\theta a^2+E_{\rm neu}$ or have tail bounded by \begin{equation*} C\sup_{k\ge n}v_k \sum_{k=n}^{\infty} \{d_{\theta,k}\mathfrak a_k^2+v_kq_k^2+\mathcal E_k\} =o(X_n). \end{equation*} Appendix B.A1.5.5 proves the energy estimate \begin{equation*} \sum_{k=n}^{\infty} \{d_{\theta,k}\mathfrak a_k^2+v_kq_k^2+\mathcal E_k\} \le C X_n \end{equation*} directly from the exact full two-parity chord identity; here $\mathcal E$ includes both neutral and stable external masses. This also proves the uniformity of \eqref{eq:B.9.30r} at the endpoints and collisions.

Since $q_n\to0$, also $\widehat q_n\to0$. Sum \eqref{eq:B.9.30r} from $n$ to infinity: $$ \widehat q_n\ge-\sum_{k=n}^{\infty}|\mathfrak r_k|=-o(X_n). $$ Equation \eqref{eq:B.9.30n} and $q_n=e_n+O(v_n\mathfrak a_n)$ give $|q_n|+|\mathfrak a_n|\le CX_n$. Since $v_n\to0$, it follows explicitly that $$ |\widehat q_n-q_n|\le Cv_n|q_n|=o(X_n),\qquad |q_n-e_n|\le Cv_n|\mathfrak a_n|=o(X_n). $$ Hence \begin{equation*} e_n\ge-o(X_n). \tag{B.9.30w}\label{eq:B.9.30w} \end{equation*} Every persistent neutral mode at this upper boundary face lies below $\alpha$. Neutrality gives \begin{equation*} D_J(h)=-\frac{2C_*}{h-\alpha}>0. \end{equation*} Combining \eqref{eq:B.9.30o}, \eqref{eq:B.9.30w}, and \eqref{eq:B.9.30n}, and using the post-extension stable comparison \eqref{eq:B.A1.72a}, in the exact cocycle increment yields, for every persistent neutral $h$, \begin{equation*} \log\frac{\mathcal S_{h,n+1}}{\mathcal S_{h,n}} =-\frac{4x_n}{C_*} +\frac{2D_J(h)}{C_*}e_n+o(X_n) \ge c_hX_n>0 \tag{B.9.30y}\label{eq:B.9.30y} \end{equation*} on every sufficiently late nonstationary tail.

The endpoint exponents are those in \eqref{eq:B.8.57}. A mode in $I_L$ makes its increasing scalar tend to zero at every interior cluster point and at $\theta=\lambda_2$; hence it forces $\theta_n\to\lambda_3$. A mode in $I_M$ similarly forces $\theta_n\to\lambda_2$. The two types cannot both persist, since their two increasing positive scalars would force incompatible limits. At least one $I_L$-mode is present in the residual case \eqref{eq:B.9.25}, so $\theta_n$ converges. If exactly one core coordinate is deleted before this conclusion, the deletion is permanent and the restricted five-node analysis of Lemma~\ref{lem:8.7} applies. If two core coordinates are deleted in the same return, let $J_4$ be the four surviving core nodes. The exact parity map on $J_4$ is the identity, while every coordinate outside $J_4$ is either identically zero or has summable squared mass by \eqref{eq:B.5.6}; the projection argument in the proof of Lemma~\ref{lem:7.1} therefore gives convergence. Three or more simultaneous core deletions would leave an $\omega$-limit support of at most three nodes, contrary to Lemma~\ref{lem:3.2}. After a deletion, the corresponding lower-support argument applies.

For the reflected boundary face $\rho_n\to\lambda_2$, reverse the node order. Then $\alpha=\lambda_2$, $\theta=\sigma$, and the oriented normal coordinate satisfies $q_n=-e_n+O(v_n\mathfrak a_n)$. All displayed identities are unchanged after this sign convention. Every persistent mode lies above $\alpha$, so $D_J(h)<0$ and $D_J(h)e_n\ge-o(X_n)$. The reflected version of \eqref{eq:B.9.30y} again selects a unique pivot endpoint. Both roots therefore converge.

Finally Lemma~\ref{lem:5.1} gives convergence of $P_n,Q_n,C_n$, Lemma~\ref{lem:5.3} gives $U_n\to0$, and the recovery formula \eqref{eq:B.6.42} gives coordinatewise convergence of all six squared core weights, including a coordinate that vanishes only asymptotically. This completes both boundary reductions.
\end{proof}

The interior and boundary dynamics are now settled. The six-node case can therefore be assembled, including every possible exact support loss.

\begin{lemma}[Six-node completion tails]\label{lem:9.5}
If $|S_*|=6$, then $Y_n^2$ converges.
\end{lemma}

\begin{proof}
Unstable external modes are exactly deleted by Lemma~\ref{lem:5.2}. If no neutral mode persists, Lemma~\ref{lem:6.2} shadows the tail by an exact orbit of positive weights on the invariant coordinate face. Lemma~\ref{lem:6.1} makes that exact six-node orbit converge unless a core coordinate is deleted exactly. If exactly one coordinate is deleted, factor the unchanged completion as \begin{equation*} \mathscr D_*(t)=D_J(t)(t-\alpha) \tag{B.9.32}\label{eq:B.9.32} \end{equation*} on the remaining five nodes; the restricted five-node argument of Lemma~\ref{lem:7.3} applies. If two coordinates are deleted in the same return, let $J_4$ be the four surviving core nodes. The parity map on $J_4$ is the identity; deleted core coordinates vanish identically and the squared mass of every external coordinate is summable by \eqref{eq:B.5.6}, so the projection argument in Lemma~\ref{lem:7.1} applies. Three or more simultaneous deletions would leave an $\omega$-limit support of size at most three, contrary to Lemma~\ref{lem:3.2}. A later deletion from the five-node face enters the same four-node argument.

If a neutral mode persists and no core coordinate is exactly deleted, Lemmas~\ref{lem:9.3}--\ref{lem:9.4} show that $\rho_n,\sigma_n$ converge. Factor convergence is Lemma~\ref{lem:5.1}, $U_n\to0$ is Lemma~\ref{lem:5.3}, and \eqref{eq:B.6.42} gives coordinatewise convergence of all six squared core weights, even when one or two tend to zero. Every external weight tends to zero. The global coefficient charts of Lemmas~\ref{lem:4.3},~\ref{lem:9.1}, and~\ref{lem:9.2} extend the argument through factor collisions.

If an exact deletion occurs during this neutral analysis, it is permanent. After one deletion, equation \eqref{eq:B.9.32} is the exact five-node completion, and the restricted analysis of Lemma~\ref{lem:8.7} covers every external-node position, including a spectral or confluent $\alpha$. Two simultaneous deletions, or a later second deletion, are covered by the four-node projection argument above; three or more simultaneous deletions would contradict the four-node lower bound for a nonterminating $\omega$-limit support. These alternatives exhaust the finite support-loss chain. Hence every six-node completion tail has a singleton squared-weight $\omega$-limit set.
\end{proof}

Sections~\ref{sec:7}--\ref{sec:9} prove convergence of $Y_n^2$ for every possible limiting support. The remaining step is the recovery of signs, the opposite parity, and the original spectral representation.

\section{Signed recovery, the opposite parity, and the original problem}\label{sec:10}

\begin{lemma}[Squared weights determine the signed parity limit]\label{lem:10.1}
If $Y_n^2\to w_*$, then $Y_n$ converges as a signed Euclidean unit vector. The opposite monic parity $Z_n$ also converges.
\end{lemma}

\begin{proof}
Lemma~\ref{lem:3.2} gives at least four positive entries in $w_*$, so all solutions of the moment equations depend continuously there. The next squared weight converges to $w_*^{\rm o}=T(w_*)$, and $$ w_*=T(w_*^{\rm o}). $$ The two limiting weights have the same support. By Lemma~\ref{lem:3.4}, the limiting signed two-block multiplier is $+1$ on every positive limiting coordinate. Equivalently, for $i\in\operatorname{supp}w_*$, \begin{equation*} \frac{\widetilde P_n(\lambda_i)\widetilde Q_n(\lambda_i)}{\widetilde a_n\widetilde b_n}\longrightarrow1. \end{equation*} The exact full-state recurrence is \begin{equation*} Y_{n+1,i}= \frac{\widetilde P_n(\lambda_i)\widetilde Q_n(\lambda_i)}{\widetilde a_n\widetilde b_n}Y_{n,i}. \end{equation*} Every coordinate with $w_{*,i}>0$ has an eventually fixed sign and absolute value tending to $\sqrt{w_{*,i}}$. Coordinates with $w_{*,i}=0$ tend to zero regardless of their signs. There are finitely many coordinates, so $Y_n\to Y_*$ in Euclidean norm.

The map $L$ is continuous at $Y_*$, again because its support has at least four nodes. Hence \begin{equation*} Z_n=L(Y_n)\longrightarrow L(Y_*). \end{equation*} Recall from \eqref{eq:B.3.1} that $Z_n=-y_{2n+1}$, so the original odd signed directions converge as well.
\end{proof}

\begin{lemma}[Return to the original SPD instance]\label{lem:10.2}
Convergence in the active spectral coordinates is exactly convergence of the signed normalised residuals of the original restarted CG problem. The same conclusion holds for repeated eigenvalues, initially zero spectral components, and permanent later support loss.
\end{lemma}

\begin{proof}
The vectors $v_i$ in Lemma~\ref{lem:2.1} are orthonormal, so $$ \mathcal U:\mathbb R^N\longrightarrow \operatorname{span}\{v_1,\ldots,v_N\},\qquad \mathcal U(\eta)=\sum_i\eta_iv_i $$ is a fixed isometry. It follows that coefficient convergence is equivalent to Euclidean convergence of the original signed vectors.

A repeated eigenspace contributes the single fixed line generated by the initial spectral projection; polynomial residual updates never rotate inside it. An eigenspace whose initial projection is zero is never activated. A later zero coordinate stays zero. Translation by $x_*=A^{-1}b$ leaves every residual unchanged, and reversal of an orthogonal diagonalisation is another fixed isometry. These observations establish the result.
\end{proof}

\begin{proof}[Proof of the theorem for restart length three]
If $r_0=0$, the iteration terminates at restart index zero. If at any later boundary the grade is at most three, Lemma~\ref{lem:2.1} proves termination in the next block of restart length three. We may assume that the orbit is infinite.

Lemmas~\ref{lem:3.1}--\ref{lem:3.3} give finite parity energy and regular $\omega$-limit supports on four, five, or six nodes. Lemmas~\ref{lem:4.4} and~\ref{lem:5.1} determine a single completion polynomial and a set $S_*$ of four, five, or six limiting spectral nodes. All off-completion coordinates fall under Lemma~\ref{lem:5.2}.

If $|S_*|=4$, Lemma~\ref{lem:7.1} gives convergence of the even squared weights. If $|S_*|=5$, Corollary~\ref{cor:7.4} treats the absence of persistent neutral modes and Lemma~\ref{lem:8.7} treats their presence. If $|S_*|=6$, Lemma~\ref{lem:9.5}, using the boundary analysis in Lemma~\ref{lem:9.4}, gives the same conclusion. Hence $Y_n^2$ converges in every nonterminating case.

Lemma~\ref{lem:10.1} recovers convergence of the signed even parity and of the opposite signed parity. Lemma~\ref{lem:10.2} carries those limits through the fixed spectral isometry to the original SPD residuals, including repeated eigenvalues, initially zero spectral components, and every permanent later support loss. This completes the proof of Theorem~\ref{thm:main}.
\end{proof}

\section*{Appendix B. Exact algebra for Lemmas~\ref{lem:9.1}--\ref{lem:9.4}}

This appendix contains the calculations used in Section~\ref{sec:9}. Retain the notation $D(t)=\prod_{i=1}^6(t-\lambda_i)$, $K=(t-\rho)(t-\sigma)$, and $C_*$ from that section. Throughout the appendix, $h$ is a persistent neutral external node, so $D(h)=-2C_*$. In Subsection B.A1.2 only, make the affine change $x=(t-\rho)/(\sigma-\rho)$ and suppress bars on the transformed quantities there. All other subsections use the original variable $t$ and roots $\rho,\sigma$. Returning to the original variable multiplies each evaluated affine normal by the positive number $(\sigma-\rho)^2$ and hence does not alter any endpoint sign.

The calculations follow the order in which they are used: Uvarov derivatives determine the returned core, residue identities give the endpoint signs and the logarithmic cocycle, and the final subsections treat the boundary strata. Two points concerning continuation across collisions will be used repeatedly. Identities first obtained on the coprime locus extend as identities of analytic germs in a reduced quotient. Uniform estimates follow from either denominator-free coefficient formulae or removable singularities in signed charts; a finite compact cover then supplies the required $C^2$ bounds.

\subsection*{B.A1.1. Uvarov derivatives and the returned core}

At a fixed six-node state let \begin{equation*} \mu_P(i)=c_iK(\lambda_i)Q(\lambda_i),\qquad \mu_Q(i)=c_iK(\lambda_i)P(\lambda_i). \tag{B.A1.1}\label{eq:B.A1.1} \end{equation*} Lemma~\ref{lem:4.1} shows that these are probability measures, that their monic cubics are $P,Q$, and that \begin{equation*} P(\lambda_i)Q(\lambda_i)=C_*,\qquad P(h)Q(h)=-C_*. \tag{B.A1.2}\label{eq:B.A1.2} \end{equation*} For $R=P,Q$, the reproducing kernel is the unique quadratic satisfying \begin{equation*} \langle k_R(\,\cdot\,,h),f\rangle_{\mu_R}=f(h) \quad(f\in\mathcal P_2); \end{equation*} existence and uniqueness follow because the $3\times3$ moment Gram matrix is positive definite.

Insert mass into the exact state by $\widetilde\mu_{P,\varepsilon}=(1-\varepsilon)\mu_P+ \varepsilon\delta_h$, with the expansion $$ \widetilde P_\varepsilon=P+\varepsilon p+O(\varepsilon^2). $$ Differentiation of $\langle \widetilde P_\varepsilon,f\rangle_{\widetilde\mu_{P,\varepsilon}}=0$ gives, for every $f\in\mathcal P_2$, $$ \langle p,f\rangle_{\mu_P}+P(h)f(h)=0. $$ Hence \begin{equation*} p=-P(h)k_P(\,\cdot\,,h). \tag{B.A1.4}\label{eq:B.A1.4} \end{equation*} The scalar factor $1-\varepsilon$ contributes a multiple of the unperturbed zero moment and hence disappears.

Let $\widetilde Q_\varepsilon=Q+\varepsilon q+O(\varepsilon^2)$ be the exact partner monic cubic. Its unnormalised orthogonality equations are $$ \int \widetilde P_\varepsilon(t)^2\widetilde Q_\varepsilon(t)f(t)\, d\widetilde\mu_{P,\varepsilon}(t)=0. $$ Their derivative is $$ C_*\langle q,f\rangle_{\mu_Q} +2C_*\langle p,f\rangle_{\mu_P} +P(h)^2Q(h)f(h)=0. $$ Using \eqref{eq:B.A1.2}, \eqref{eq:B.A1.4}, and $P(h)^2Q(h)=-C_*P(h)$, one obtains \begin{equation*} \langle q,f\rangle_{\mu_Q}=3P(h)f(h),\qquad q=3P(h)k_Q(\,\cdot\,,h). \tag{B.A1.5}\label{eq:B.A1.5} \end{equation*}

Let $R=Qp+Pq$. Since $\mu_P(i)P_i=C_*c_iK_i$ and $\mu_Q(i)Q_i=C_*c_iK_i$, equations \eqref{eq:B.A1.4}--\eqref{eq:B.A1.5} give \begin{equation*} \sum_i c_iK_iR_i f_i =\langle p,f\rangle_{\mu_P} +\langle q,f\rangle_{\mu_Q} =2P(h)f(h). \tag{B.A1.6}\label{eq:B.A1.6} \end{equation*} Modulo the common retained-mass normalisation, the logarithmic variation of the returned normalised six-node core is \begin{equation*} \delta\log\mu_{P,i}^{\rm out} =\frac{2p_i}{P_i}+\frac{2q_i}{Q_i} =\frac{2R_i}{C_*}. \tag{B.A1.7}\label{eq:B.A1.7} \end{equation*} If $P_{\rm out}=P+\varepsilon p_{\rm out}+O(\varepsilon^2)$, its differentiated orthogonality and \eqref{eq:B.A1.6} give $$ \begin{aligned} \langle p_{\rm out},f\rangle_{\mu_P} &=-\frac{2}{C_*}\sum_i\mu_{P,i}P_iR_if_i\\ &=-2\sum_ic_iK_iR_if_i=-4P(h)f(h). \end{aligned} $$ Hence \begin{equation*} p_{\rm out}=-4P(h)k_P(\,\cdot\,,h). \tag{B.A1.8}\label{eq:B.A1.8} \end{equation*} For its restricted-core partner, use $Q_{\rm out}=Q+\varepsilon q_{\rm out}+O(\varepsilon^2)$. Differentiating $\int P_{\rm out}^2Q_{\rm out}f\,d\mu_P^{\rm out}=0$ and using \eqref{eq:B.A1.7} gives $$ \begin{aligned} 0={}&C_*\langle q_{\rm out},f\rangle_{\mu_Q} +2C_*\langle p_{\rm out},f\rangle_{\mu_P} +\frac{2}{C_*}\int P^2QRf\,d\mu_P\\ ={}&C_*\langle q_{\rm out},f\rangle_{\mu_Q} -8C_*P(h)f(h)+4C_*P(h)f(h). \end{aligned} $$ Hence \begin{equation*} q_{\rm out}=4P(h)k_Q(\,\cdot\,,h). \tag{B.A1.9}\label{eq:B.A1.9} \end{equation*}

At a fixed point, differentiate the general factor identity $$ KPQ=D(K+U)+H_1K_1. $$ For the returned state its first variation is $$ \delta K\,PQ+K(p_{\rm out}Q+Pq_{\rm out}) =D(\delta K+\delta U)+\delta H_1K+C_*\delta K_1. $$ Because $PQ=D+C_*$, cancellation of $D\delta K$ leaves \begin{equation*} K(p_{\rm out}Q+Pq_{\rm out}) =D\delta U+ \{\delta H_1K+C_*(\delta K_1-\delta K)\}. \tag{B.A1.10}\label{eq:B.A1.10} \end{equation*} The braces have degree at most two. Taking the polynomial quotient by the degree-six polynomial $D$ in \eqref{eq:B.A1.10}, and substituting \eqref{eq:B.A1.8}--\eqref{eq:B.A1.9}, proves \begin{equation*} \delta U(t)=4P(h)\operatorname{quo}_{D} \bigl(K(t)\{P(t)k_Q(t,h)-Q(t)k_P(t,h)\}\bigr). \tag{B.A1.11}\label{eq:B.A1.11} \end{equation*} Equation \eqref{eq:B.A1.11} is \eqref{eq:B.9.4}; equations \eqref{eq:B.A1.8}--\eqref{eq:B.A1.9} include the contribution from the factor coordinates.

\subsection*{B.A1.2. Endpoint signs}

After the affine normalisation, $K=x(x-1)$. Matching residues at all six nodes and the coefficient of $z^{-1}$ at infinity gives \begin{equation*} \mathcal C_P(z)=\sum_i\frac{\mu_P(i)}{z-\lambda_i} =\frac{K(z)Q(z)}{D(z)},\qquad \mathcal C_Q(z)=\frac{K(z)P(z)}{D(z)}. \tag{B.A1.12}\label{eq:B.A1.12} \end{equation*} Furthermore $$ \sum_i\frac{\mu_P(i)P(\lambda_i)}{z-\lambda_i} =\frac{C_*K(z)}{D(z)}. $$ Hence the third second-kind polynomial is \begin{equation*} P(z)\mathcal C_P(z)- \sum_i\frac{\mu_P(i)P(\lambda_i)}{z-\lambda_i}=K(z), \tag{B.A1.13}\label{eq:B.A1.13} \end{equation*} and the identical equation holds for $Q$.

For $\mu_P$, the monic recurrence is \begin{equation*} \pi_0=1,\quad \pi_1=x-a_0,\quad \pi_2=(x-a_1)\pi_1-\beta_1,\quad P=(x-a_2)\pi_2-\beta_2\pi_1, \tag{B.A1.14}\label{eq:B.A1.14} \end{equation*} where $\beta_1,\beta_2>0$. Equation \eqref{eq:B.A1.13} says $$ (x-a_2)(x-a_1)-\beta_2=x(x-1). $$ Hence $$ a_1+a_2=1,\qquad \beta_2=a_2(1-a_2)>0. $$ With $a=a_0$, $b=\beta_1$, and $c=a_2$, this gives \begin{equation*} P=(x-a)K-b(x-c),\qquad b>0,\quad0<c<1, \tag{B.A1.15}\label{eq:B.A1.15} \end{equation*} \begin{equation*} R_P=\pi_2=(x-1+c)(x-a)-b,\qquad H_P=bc(1-c). \end{equation*} Applying the displayed derivation to $\mu_Q$ gives \begin{equation*} Q=(x-A)K-B(x-f),\qquad R_Q=(x-1+f)(x-A)-B,\qquad H_Q=Bf(1-f), \tag{B.A1.17}\label{eq:B.A1.17} \end{equation*} where $B>0$ and $0<f<1$.

The Christoffel--Darboux identity, obtained by telescoping the three recurrence relations \eqref{eq:B.A1.14}, is \begin{equation*} k_P(x,h)= \frac{P(x)R_P(h)-R_P(x)P(h)}{H_P(x-h)}, \tag{B.A1.18}\label{eq:B.A1.18} \end{equation*} and similarly for $Q$. Define $$ S(x)=P(x)k_Q(x,h)-Q(x)k_P(x,h),\qquad s_5=[x^5]S,\quad s_4=[x^4]S. $$ Directly taking the two leading coefficients in \eqref{eq:B.A1.18} yields \begin{equation*} s_5=\frac{R_Q(h)}{H_Q}-\frac{R_P(h)}{H_P}, \end{equation*} \begin{equation*} s_4=(h+p_2+q_2)s_5+ \frac{P(h)}{H_P}-\frac{Q(h)}{H_Q}, \end{equation*} where $p_2=[x^2]P$ and $q_2=[x^2]Q$.

Let $\Theta=\operatorname{quo}_D(KS)$. Since $$ D=x^6+(p_2+q_2)x^5+\cdots,\qquad K=x^2-x, $$ coefficient comparison in $KS=D\Theta+\operatorname{rem}_D(KS)$ gives \begin{equation*} \Theta(0)=(h-1)s_5+\frac{P(h)}{H_P}-\frac{Q(h)}{H_Q}, \tag{B.A1.21}\label{eq:B.A1.21} \end{equation*} \begin{equation*} \Theta(1)=hs_5+\frac{P(h)}{H_P}-\frac{Q(h)}{H_Q}. \tag{B.A1.22}\label{eq:B.A1.22} \end{equation*} The two elementary consequences of \eqref{eq:B.A1.15}, \begin{equation*} (h-a)(h-1)-b=\frac{P(h)-bc}{h}, \quad h(h-a)-b=\frac{P(h)+b(1-c)}{h-1}, \end{equation*} and their $Q$-counterparts reduce \eqref{eq:B.A1.21}--\eqref{eq:B.A1.22} to \begin{equation*} \Theta(0)=\frac{1}{h} \left(\frac{P(h)}{bc}-\frac{Q(h)}{Bf}\right), \end{equation*} \begin{equation*} \Theta(1)=\frac{1}{h-1} \left(-\frac{P(h)}{b(1-c)} +\frac{Q(h)}{B(1-f)}\right). \tag{B.A1.25}\label{eq:B.A1.25} \end{equation*} Multiplying by $4P(h)$ and using $P(h)Q(h)=-C_*$ gives \begin{equation*} A_h(0)=\frac{4}{h} \left(\frac{P(h)^2}{bc}+\frac{C_*}{Bf}\right), \tag{B.A1.26}\label{eq:B.A1.26} \end{equation*} \begin{equation*} A_h(1)=-\frac{4}{h-1} \left(\frac{P(h)^2}{b(1-c)} +\frac{C_*}{B(1-f)}\right). \tag{B.A1.27}\label{eq:B.A1.27} \end{equation*} The braces are strictly positive. Equations \eqref{eq:B.A1.26}--\eqref{eq:B.A1.27}, followed by reversal of the affine normalisation, prove \eqref{eq:B.9.5}.

The endpoint signs describe the direct forcing of the two roots. The next calculation identifies a primitive for the intrinsic root variation and thereby produces the cocycle of Lemma~\ref{lem:9.2}.

\subsection*{B.A1.3. A primitive for the root variation and the exact logarithmic cocycle}

Let $u$ be an affine intrinsic normal and divide \begin{equation*} D(t)u(t)=K(t)q_u(t)+r_u(t),\qquad \deg r_u\le1. \tag{B.A1.28}\label{eq:B.A1.28} \end{equation*} Equation \eqref{eq:B.4.13} evaluated at the two simple roots of $K$ gives \begin{equation*} \dot\rho=-\frac{2r_u(\rho)}{C_*(\sigma-\rho)},\qquad \dot\sigma=\frac{2r_u(\sigma)}{C_*(\sigma-\rho)}. \tag{B.A1.29}\label{eq:B.A1.29} \end{equation*} Here the dot denotes the coefficient of the first-order normal perturbation.

Both $u$ and $r_u$ are affine, and $r_u(\rho)=D(\rho)u(\rho)$, $r_u(\sigma)=D(\sigma)u(\sigma)$. Interpolating each at $h$ and substituting in $q_u(h)=\{D(h)u(h)-r_u(h)\}/K(h)$ gives \begin{equation*} q_u(h)=L_\rho r_u(\rho)+L_\sigma r_u(\sigma), \end{equation*} where, with $g=\sigma-\rho$, \begin{equation*} L_\rho=\frac{2C_*+D(\rho)}{(h-\rho)D(\rho)g},\qquad L_\sigma=-\frac{2C_*+D(\sigma)}{(h-\sigma)D(\sigma)g}. \tag{B.A1.31}\label{eq:B.A1.31} \end{equation*} For example, the coefficient of $r_u(\rho)$ is $$ \frac{(\sigma-h)\{D(h)/D(\rho)-1\}/g}{(h-\rho)(h-\sigma)} =\frac{2C_*+D(\rho)}{(h-\rho)D(\rho)g}; $$ the other coefficient is obtained by the displayed interpolation with $\rho,\sigma$ exchanged.

The exact full multiplier is $\widetilde M_h=\widetilde P(h)\widetilde Q(h)/\widetilde C$, where $\widetilde C=\sqrt{\widetilde H_0\widetilde H_1}$. At zero external mass the tilded objects agree with the fixed normalised core. The polynomial identity \eqref{eq:B.5.22a} then shows that, in the pure core $u$-direction, the numerator changes by $q_u(h)$ plus the same scalar first variation as the denominator. Since $\widetilde M_h=-1$ at neutrality, \begin{equation*} \ell_h(u):=\left.d\log \widetilde M_h^2\right|_{E=0}=-\frac{2}{C_*}q_u(h). \tag{B.A1.32}\label{eq:B.A1.32} \end{equation*} Substitution of \eqref{eq:B.A1.29}--\eqref{eq:B.A1.31} gives \begin{equation*} \ell_h(u)=f_h(\rho)\dot\rho+f_h(\sigma)\dot\sigma, \qquad f_h(t)=\frac{2C_*+D(t)}{(h-t)D(t)}. \tag{B.A1.33}\label{eq:B.A1.33} \end{equation*}

The partial fraction identity \begin{equation*} \frac{1}{D(t)}=\sum_i\frac{c_i}{t-\lambda_i}, \qquad \sum_i\frac{c_i}{\lambda_i-h}=\frac{1}{2C_*}, \end{equation*} uses $D(h)=-2C_*$. It gives \begin{equation*} f_h(t)=-2C_*\sum_i\frac{c_i}{\lambda_i-h} \frac{1}{t-\lambda_i}. \end{equation*} Hence \begin{equation*} \Phi_h(\rho,\sigma)= -2C_*\sum_i\frac{c_i}{\lambda_i-h} \log|(\lambda_i-\rho)(\lambda_i-\sigma)| \tag{B.A1.36}\label{eq:B.A1.36} \end{equation*} satisfies $d\Phi_h[u]=\ell_h(u)$. Since \begin{equation*} \prod_i|K(\lambda_i)|^{\beta_i(h)}=e^{-\Phi_h}, \tag{B.A1.37}\label{eq:B.A1.37} \end{equation*} the intrinsic first normal variation of the two factors in $\mathcal J_h$ cancels exactly.

The recovery identity \eqref{eq:B.6.42} yields the uniform discrete remainder. Here $w,K,P,H_0$ are the current normalised-core objects, and their plus versions belong to the normalised retained core induced by the exact full return. Applied at both core states, the recovery identity gives the exact quotient \begin{equation*} \frac{K^+(\lambda_i)}{K(\lambda_i)} =\frac{w_i^+}{w_i} \frac{P^+(\lambda_i)}{P(\lambda_i)} \frac{H_0}{H_0^+}. \tag{B.A1.38}\label{eq:B.A1.38} \end{equation*} At a support boundary, $w_i=\chi_i^2$ and $|K(\lambda_i)|=\chi_i^2u_i$ for a positive analytic unit $u_i$. The sign of $K(\lambda_i)$ is fixed on each boundary chart; hence \eqref{eq:B.A1.38} is still an analytic, strictly positive quotient. All $P(\lambda_i)$ stay nonzero because $P(\lambda_i)Q(\lambda_i)\to C_*>0$. Combining this core quotient with the exact full external multiplier and \eqref{eq:B.5.25}--\eqref{eq:B.5.26}, the hybrid logarithmic return $\mathcal B_h^{(6)}:=\log(\mathcal J_h^+/\mathcal J_h)$ is a $C^2$ function in the signed face charts.

On the fixed component at level $C$ and at zero external mass, \begin{equation*} \mathcal B_h^{(6)}=2\log\left|1-\frac{2C_*}{C}\right|. \tag{B.A1.39}\label{eq:B.A1.39} \end{equation*} Its $C$-derivative at $C_*$ is $-4/C_*$; its derivatives in all tangent directions to the fixed component vanish. Equations \eqref{eq:B.A1.32}--\eqref{eq:B.A1.37} show that its two core-normal derivatives vanish as well. Other external squared masses enter with their bounded, generally nonzero Uvarov first derivatives. Taylor's formula with integral remainder on a finite compact face cover, using the $C^2$ full/core transfer in the core variables, gives \begin{equation*} \mathcal B_h^{(6)}=-\frac{4(C-C_*)}{C_*} +R_h,\qquad |R_h|\le C_{\rm rem}\{(C-C_*)^2+\|U\|^2+\mathcal E\}. \tag{B.A1.40}\label{eq:B.A1.40} \end{equation*} Thus no external first-order term is silently absorbed into a core transfer error; all such terms are bounded by the last displayed summand. The expressions \eqref{eq:B.A1.28}--\eqref{eq:B.A1.38} are rational in coefficient coordinates with denominators nonzero at spectral nodes. Clearing those denominators shows that the cancellation is a polynomial identity on the dense coprime locus; continuity extends it to every common-factor collision. This gives the boundary- and collision-uniform return estimate used in Lemma~\ref{lem:9.2}; only after this extension is the stable comparison \eqref{eq:B.9.22a} invoked. Lemma~\ref{lem:5.3} then proves \eqref{eq:B.9.15}--\eqref{eq:B.9.16}.

The cocycle can fail to control the root product only when a boundary factor vanishes. The divided recurrences below govern the resulting approach to a support face.

\subsection*{B.A1.4. Divisibility, the single-crossing estimate, and root variation}

Use signed amplitudes for all vanishing weights. From \eqref{eq:B.4.9a}, if the $\sigma$-root approaches an endpoint, its vanishing squared weight is \begin{equation*} y^2=u_\sigma d_\sigma,\qquad d_\sigma=(\sigma-\lambda_4)(\lambda_5-\sigma), \end{equation*} with $u_\sigma>0$ analytic; for the $\rho$-root we similarly use $x^2=u_\rho d_\rho$, where $d_\rho=(\rho-\lambda_2)(\lambda_3-\rho)$ and $u_\rho>0$ is analytic. The return is even in each signed amplitude, so external coordinates enter as $E_h=\xi_h^2$.

Let $z_\rho=U(\rho)$ and $z_\sigma=U(\sigma)$. At zero external mass, the ideal divisions \eqref{eq:B.6.18}--\eqref{eq:B.6.23} already prove $$ z_\rho^+=z_\rho R_\rho+d_\sigma z_\sigma^2q_\rho. $$ For nonzero external coordinates, the difference from this expression is an analytic function of the independent squared amplitudes $(E_h)_h$ and vanishes when they are all zero. The multivariable Hadamard identity $$ F(E)-F(0)=\sum_hE_h\int_0^1 \partial_{E_h}F(E_1,\ldots,E_{h-1},sE_h,0,\ldots,0)\,ds $$ gives the exact analytic factorisation \begin{equation*} z_\rho^+ =z_\rho R_\rho+d_\sigma z_\sigma^2q_\rho +\sum_h E_h\,a_{h,\rho}. \tag{B.A1.42}\label{eq:B.A1.42} \end{equation*} The analytic coefficients $a_{h,\rho}$ absorb the mixed external terms. Equations \eqref{eq:B.4.14} and \eqref{eq:B.A1.11} identify their values on the fixed component for a persistent neutral source: \begin{equation*} q_\rho\big|_{\rm fixed} =\frac{\Gamma_\rho}{(\sigma-\rho)^2d_\sigma}>0,\qquad a_{h,\rho}\big|_{\rm fixed}=A_h(\rho)>0 \tag{B.A1.43}\label{eq:B.A1.43} \end{equation*} whenever $h>\rho$. The endpoint zeros of $D(\sigma)$ and $d_\sigma$ are the same simple zero, so the first quotient in \eqref{eq:B.A1.43} has a finite positive endpoint limit. Compactness makes $R_\rho,q_\rho$, and the indicated neutral coefficients positive on a sufficiently late tube. All stable coefficients remain uniformly bounded and contribute $O(E_{\rm st})$.

At zero external mass, \eqref{eq:B.6.32} gives the first two terms of \eqref{eq:B.A1.45}. For arbitrary external mass the root $\rho$ is fixed when its endpoint weight is zero, so the increment is divisible by $d_\rho=x^2/u_\rho$. After division, its difference from the unforced expression is analytic in $(E_h)$, vanishes at $E=0$, and the displayed multivariable Hadamard identity makes it $\sum_h E_h\,B_{h,\rho}$. Differentiating $K^+(\rho^+)=0$, using \eqref{eq:B.4.13}, gives \begin{equation*} \left.\frac{\partial(\rho^+-\rho)}{\partial z_\rho} \right|_{\rm fixed} =\frac{2D(\rho)}{C_*(\rho-\sigma)} =d_\rho A_\rho,\qquad A_\rho<0. \end{equation*} Hence \begin{equation*} \rho^+-\rho=d_\rho \{A_\rho z_\rho+ B_\rho d_\sigma z_\sigma^2+ \sum_h B_{h,\rho}E_h\}, \tag{B.A1.45}\label{eq:B.A1.45} \end{equation*} with bounded analytic $B$'s. Equations \eqref{eq:B.A1.42}--\eqref{eq:B.A1.45} are the exact divided recurrences used in Lemma~\ref{lem:9.3}.

If no persistent neutral node lies in $I_L$, then every persistent neutral node satisfies $h>\rho$. On the boundary $z_\rho=0$, equations \eqref{eq:B.A1.42}--\eqref{eq:B.A1.43} therefore give a positive contribution bounded below by a constant times $E_{\rm neu}$. Stable modes remain in the exact recurrence; their contribution, including its mixed analytic terms, is $O(E_{\rm st})=o(E_{\rm neu})$ by \eqref{eq:B.6.63}. Thus the exact return eventually leaves the half-space $z_\rho\ge0$ invariant. A negative $z_\rho$ can cross its boundary only into that half-space and can do so at most once.

It remains to justify the external part of the energy estimate without deleting those stable modes. For every surviving external signed amplitude $\xi_h$, the exact full-state chord has leading external component $(\widetilde M_h-1)\xi_h\mathbf e_h$, where $\widetilde M_h\to m_h\ne1$. The dependence of the normalised core on $\xi_h$ starts with $E_h=\xi_h^2$. Hence the limiting Gram matrix is the core Gram matrix from \eqref{eq:B.6.38} direct-summed with $\operatorname{diag}_h\{(m_h-1)^2\}$. The full/core transfer \eqref{eq:B.5.25}, the positive-unit endpoint changes, and the collision charts preserve a uniform lower singular-value bound. The exact parity-chord identity \eqref{eq:B.3.4} consequently gives \begin{equation*} \sum_n(d_{\rho,n}z_{\rho,n}^2+ d_{\sigma,n}z_{\sigma,n}^2+\mathcal E_n)<\infty. \tag{B.A1.46}\label{eq:B.A1.46} \end{equation*} Here $\mathcal E$ is the total external mass, including the stable modes. After the possible crossing, the variation of $\rho$ opposite to the sign of the leading term in \eqref{eq:B.A1.45} is bounded by a constant times the remaining summable terms in \eqref{eq:B.A1.46}. Since $\rho$ is bounded, its other directional variation is finite as well. Hence $\rho_n$ has finite total variation. Reflection proves the statement for $\sigma_n$.

\subsection*{B.A1.5. Analysis of the remaining boundary faces}

At a cluster point, a factor in \eqref{eq:B.9.16} can be unbounded only if a vanishing support-boundary factor has negative exponent. The sign table \eqref{eq:B.9.24} yields exactly \begin{equation*} \begin{array}{c|c} h\in I_L&\rho=\lambda_3\ \text{or}\ \sigma=\lambda_5,\\ h\in I_M&\rho=\lambda_2\ \text{or}\ \sigma=\lambda_5,\\ h\in I_R&\rho=\lambda_2\ \text{or}\ \sigma=\lambda_4. \end{array} \tag{B.A1.47}\label{eq:B.A1.47} \end{equation*} Equations \eqref{eq:B.A1.42}--\eqref{eq:B.A1.46} prove the one-sided alternatives in Lemma~\ref{lem:9.3}. For $h_L\in I_L$, $h_R\in I_R$, let $$ t_i=\frac{\lambda_i-h_L}{h_R-\lambda_i}. $$ The derivative of $(x-h_L)/(h_R-x)$ is $(h_R-h_L)/(h_R-x)^2>0$, which proves \eqref{eq:B.9.28}. Factoring the common positive number $$ \frac{2C_*}{|D'(\lambda_i)|} $$ from $p\beta_i(h_L)+q\beta_i(h_R)$ gives \eqref{eq:B.9.29} directly. Choosing $t_2<p/q<t_5$ makes the two outer-corner exponents positive. At that corner all remaining root factors are nonzero, whereas $E_{h_L}^pE_{h_R}^q\to0$; the increasing product of cocycles would tend to zero, contradicting its positive lower bound at the start of the monotone tail. The inner corner is the unique cluster point. A middle mode there has positive exponents at both vanishing factors and leads to the identical contradiction.

For the upper boundary face, take $\alpha=\lambda_5$ and $$ J=\{\lambda_1,\lambda_2,\lambda_3,\lambda_4,\lambda_6\},\qquad \theta=\rho. $$ The general recovery formula gives \begin{equation*} w_\alpha=\frac{H_0c_\alpha(\alpha-\rho)(\alpha-\sigma)}{P(\alpha)}. \end{equation*} All factors except $\alpha-\sigma$ are nonzero units. With a signed amplitude $y$, \begin{equation*} y^2=w_\alpha,\qquad \alpha-\sigma=u_y y^2,\quad u_y>0\ \text{analytic}. \tag{B.A1.49}\label{eq:B.A1.49} \end{equation*} The evaluation change for the affine normal is \begin{equation*} U\longmapsto \left(a=-U(\theta), e=-\frac{U(\alpha)}{\alpha-\theta}\right); \tag{B.A1.50}\label{eq:B.A1.50} \end{equation*} its determinant is $1/(\alpha-\theta)$, a bounded nonzero unit because $\alpha-\theta\ge\lambda_5-\lambda_3>0$. Equations \eqref{eq:B.A1.49}--\eqref{eq:B.A1.50} give the analytic coordinate change from the positive six-node stratum to the five-node stratum. The fixed union has the two components \begin{equation*} (y,a)=0\quad\text{and}\quad(a,e)=0, \end{equation*} whose union ideal is $\mathcal I_{\rm hf}=(a,ye)$. Lemma~\ref{lem:4.3} proves that this is a radical complete intersection and that $\mathcal I_{\rm hf}^{(2)}=\mathcal I_{\rm hf}^2$. We compute the four first jets directly. Define $$ \delta=\alpha-\sigma=u_y y^2,\qquad g=\alpha-\theta,\qquad \Delta=\sigma-\theta=g-\delta. $$ The spectral gaps give uniform positive lower bounds for $g$ and $\Delta$. The affine interpolation formula \eqref{eq:B.8.11a} gives \begin{equation*} U(t)=\frac{a(t-\alpha)}g-e(t-\theta),\qquad U(\sigma)=-\frac{\delta}{g}a-\Delta e. \tag{B.A1.52}\label{eq:B.A1.52} \end{equation*}

\subsubsection*{B.A1.5.1. First-order expansions in the transverse and pivot variables}

Let $\Phi_q=q^+-q$, let $\mathcal N$ be the persistent neutral external nodes, and put $\varepsilon_k=\operatorname{sign}(k-\alpha)$. At a current base point on either fixed component, define the exact coefficient directly from the full return by \begin{equation*} \gamma_k^{\rm ext}(y,e,\vartheta) :=\left.\partial_{E_k}\Phi_q\right|_{\rm fixed}. \tag{B.A1.53}\label{eq:B.A1.53} \end{equation*} Here $\vartheta$ denotes the remaining fixed-core tangential parameters, including $C$ and the pivot data. The exact moment Gram systems, normalisation, and simple-root evaluation show that every $\gamma_k^{\rm ext}$ is analytic and uniformly bounded through the two pivot endpoints and every simple or simultaneous factor collision. At the incidence $x=e=0$, put $m_k=P(k)Q(k)/C_*$. General Uvarov differentiation there gives $$ p^{[k]}=-P(k)k_P(\,\cdot\,,k),\qquad q^{[k]}=P(k)(2-m_k)k_Q(\,\cdot\,,k), $$ $$ p_{\rm out}^{[k]}=-2(1-m_k)P(k)k_P(\,\cdot\,,k),\qquad q_{\rm out}^{[k]}=2(1-m_k)P(k)k_Q(\,\cdot\,,k), $$ and hence the affine completion response $$ A_k^{(y)}(t)=2(1-m_k)P(k)\operatorname{quo}_{\mathscr D_*} \bigl(K(t)\{P(t)k_Q(t,k)-Q(t)k_P(t,k)\}\bigr),\qquad \gamma_k^{\rm ext}=-\frac{A_k^{(y)}(\sigma)}\Delta. $$ If $k\in\mathcal N$, then $m_k=-1$ at this incidence, so the response reduces to \eqref{eq:B.9.4} and the signs \eqref{eq:B.9.5}, with their five-node boundary version \eqref{eq:B.8.15}, give $$ \gamma_k^{\rm ext}=\varepsilon_k\rho_k^{\rm ext},\qquad \rho_k^{\rm ext}>0. $$ Analyticity preserves this sign on a smaller neighbourhood of the two fixed components. No sign is asserted or needed for a stable coefficient.

Define the precise transverse remainder using the coefficient at that same base point by \begin{equation*} R_q=\Phi_q+c_\theta d_\theta a^2 -\sum_k\gamma_k^{\rm ext}(y,e,\vartheta)E_k. \tag{B.A1.54}\label{eq:B.A1.54} \end{equation*} On the five-node component, the exact unforced return fixes every state with $a=0$; its linear normal return has no $a$-term by \eqref{eq:B.8.11b}. On the six-node component, \eqref{eq:B.4.13} gives $U^+-U\in(U_0,U_1)^2$, so the linear return is the identity in both $a$ and $e$. Varying $x=C-C_*$ within either fixed core also leaves $q$ fixed. Finally, every quantity determined by the squared weights is even in the signed omitted amplitude $y$. Together with \eqref{eq:B.A1.53}, these facts give the first-jet table \begin{equation*} \begin{array}{c|rrrrrr} &R_q&\partial_yR_q&\partial_aR_q&\partial_eR_q& \partial_xR_q&\partial_{E_k}R_q\\ \hline (y,a)=0&0&0&0&0&0&0\\ (a,e)=0&0&0&0&0&0&0. \end{array} \tag{B.A1.55}\label{eq:B.A1.55} \end{equation*} Here the entries are evaluated at $x=\mathcal E=0$, with $e$ free on the first component and $y$ free on the second. The exact coefficient \eqref{eq:B.A1.53} exhausts the linear transverse expansion. Terms of the forms $y^2a$, $y^2e$, and $y^2x$ belong to the quadratic anisotropic estimate below. Multivariable Taylor expansion in the external squared weights leaves only terms containing two such weights, or one external weight multiplied by a normal variable.

For the pivot formula, use the moment normal of the whole current core: \begin{equation*} \mathfrak A_y=\sum_{i\in J\cup\{\alpha\}} \frac{w_i}{\lambda_i-\theta},\qquad B_y=\sum_{i\in J\cup\{\alpha\}} \frac{w_i}{(\lambda_i-\theta)^2},\qquad \overline B_y=d_\theta B_y. \tag{B.A1.56}\label{eq:B.A1.56} \end{equation*} At $y=0$ these are the quantities in \eqref{eq:B.8.17a}. On the six-node fixed component, the $\alpha$-summand is part of the moment identity and gives $\mathfrak A_y=0$.

With \eqref{eq:B.A1.56}, define \begin{equation*} F_\theta= \overline B_y(\theta^+-\theta) -d_\theta\frac{\mathfrak A_y+3\mathfrak A_y^+}{2}, \qquad R_\theta=F_\theta/d_\theta. \tag{B.A1.57}\label{eq:B.A1.57} \end{equation*} The division is analytic in the signed endpoint charts, exactly as in the derivation of \eqref{eq:B.8.27}. In an intrinsic direction, $$ d(\theta^+-\theta)=\frac2{B_y}\,d\mathfrak A_y, \qquad d\mathfrak A_y^+=d\mathfrak A_y, $$ and hence $$ dF_\theta =2d_\theta d\mathfrak A_y -d_\theta( d\mathfrak A_y+3d\mathfrak A_y)/2=0. $$ At $x=e=0$, the Uvarov calculation \eqref{eq:B.8.21b}--\eqref{eq:B.8.21c}, with the additional core term included in both sums, gives for $k\in\mathcal N$ $$ B_y\,d(\theta^+-\theta)=\frac32\,d\mathfrak A_y^+, \qquad d\mathfrak A_y=0, $$ and again $dF_\theta=0$. On the five-node fixed component the current multiplier satisfies $m_k(e)=-1+O(e)$, so the general relation following \eqref{eq:B.8.21c} makes the corresponding derivative of $R_\theta$ only $O(|e|)$. The variables $x$ and, respectively, $e$ or $y$, are tangent to the two fixed components, on which $F_\theta$ vanishes identically. Hence, for $k\in\mathcal N$, \begin{equation*} \begin{array}{c|rrrrrr} &R_\theta&\partial_yR_\theta&\partial_aR_\theta& \partial_eR_\theta&\partial_xR_\theta& \partial_{E_k}R_\theta\\ \hline (y,a)=0&0&0&0&0&0&O(|e|)\\ (a,e)=0&0&0&0&0&0&0. \end{array} \end{equation*} Restoring the radial variable, the neutral contribution omitted from the zero jet is $O((|x|+|e|)E_{\rm neu})$. Stable $E_k$-derivatives of $R_\theta$ need not vanish, but analyticity of the exact Gram systems and the endpoint division make them uniformly bounded; their contribution is $O(E_{\rm st})$. The endpoint factor $d_\theta$ is common to both sides of \eqref{eq:B.A1.57}, while $d_\theta B_y$ is a positive unit. This gives the stated uniformity of the neutral mixed term and of the stable bound. The cocycle first jet is computed in the next subsection.

\subsubsection*{B.A1.5.2. First-order expansion of the distinct-root cocycle}

Let $$ Q_U=\operatorname{quo}_K(\mathscr D_*U),\qquad R_U=\operatorname{rem}_K(\mathscr D_*U). $$ Because $R_U$ is affine and agrees with $\mathscr D_*U$ at $\theta$ and $\sigma$, \eqref{eq:B.A1.52} gives \begin{equation*} \begin{aligned} R_U(\alpha) &= \frac{\alpha-\sigma}{\theta-\sigma} \mathscr D_*(\theta)U(\theta) +\frac{\alpha-\theta}{\sigma-\theta} \mathscr D_*(\sigma)U(\sigma)\\ &=\frac{\delta}{\Delta} \{-gD_J(\theta)a +\delta D_J(\sigma)a +g\Delta D_J(\sigma)e\}. \end{aligned} \tag{B.A1.59}\label{eq:B.A1.59} \end{equation*} Since $\mathscr D_*(\alpha)U(\alpha)=0$ and $K(\alpha)=g\delta$, Euclidean division yields \begin{equation*} Q_U(\alpha) =\frac{gD_J(\theta)-\delta D_J(\sigma)} {g\Delta}\,a-D_J(\sigma)e. \tag{B.A1.60}\label{eq:B.A1.60} \end{equation*} The radial contribution is recorded by $x=C-C_*$. The formulas retain both this scalar variation and the affine-normal contribution; their extension through factor collisions is justified in the actual state ring below.

The first root return follows directly from $$ K^+-K+\frac2{C_*}\operatorname{rem}_K(\mathscr D_*U) \in(U_0,U_1)^2. $$ Evaluation at $\theta$, followed by differentiation of $K^+(\theta^+)=0$, gives \begin{equation*} d(\theta^+-\theta) =-\frac{2gD_J(\theta)}{C_*\Delta}\,da. \tag{B.A1.61}\label{eq:B.A1.61} \end{equation*}

For $\alpha\notin J$, use the hybrid scalar $$ \mathcal S_h =E_hg^{\beta_\alpha(h)} \prod_{i\in J}w_i^{\beta_i(h)}, $$ where $E_h$ is an exact full-state mass and $w_i,g$ belong to the induced normalised core. At a fixed six-node core and zero external mass, the exact full multiplier is $1$ at every root of $\mathscr D_*$ and $-1$ at $h$. Lagrange interpolation of the degree-at-most-five core variation of $\mathscr D$, exactly as in \eqref{eq:B.8.34}--\eqref{eq:B.8.38}, gives for the weight part $$ d\mathcal B_h^{\rm wt} =-\frac4{C_*}\,dx -\frac{2\beta_\alpha(h)}{C_*}Q_U(\alpha). $$ The factor $g^{\beta_\alpha(h)}$ contributes $-\beta_\alpha(h)d(\theta^+-\theta)/g$. Substitution of \eqref{eq:B.A1.60}--\eqref{eq:B.A1.61} gives \begin{equation*} \boxed{ d\mathcal B_h =-\frac4{C_*}\,dx +\frac{2\beta_\alpha(h)D_J(\sigma)}{C_*}\,de +\frac{2\beta_\alpha(h)\delta D_J(\sigma)} {C_*g\Delta}\,da +\sum_k\chi_{hk}(y)\,dE_k.} \tag{B.A1.62}\label{eq:B.A1.62} \end{equation*}

The exact external part in \eqref{eq:B.A1.62} is best written using both multiplier types. Put $$ M_*(t)=\frac{P(t)Q(t)}{C_*},\qquad m_k=M_*(k),\qquad M^{\rm c}(t)=\widetilde M(t)\sqrt{S/S^+}, $$ where $M^{\rm c}$ is the signed multiplier of a retained normalised-core coordinate. Then \begin{equation*} \chi_{hk}(y) =-2\dot{\widetilde M}_k(h) +2\sum_{i\in J}\beta_i(h)\dot M_k^{\rm c}(\lambda_i) -\frac{\beta_\alpha(h)}g\dot\theta_k. \tag{B.A1.63}\label{eq:B.A1.63} \end{equation*} Here the dots denote $E_k$-derivatives at zero external mass. For an arbitrary source $k$, $$ p^{[k]}=-P(k)k_P(\cdot,k),\qquad q^{[k]}=P(k)(2-m_k)k_Q(\cdot,k), $$ while the returned core derivatives are $$ p_{\rm out}^{[k]}=-2(1-m_k)P(k)k_P(\cdot,k),\qquad q_{\rm out}^{[k]}=2(1-m_k)P(k)k_Q(\cdot,k). $$ Consequently, $$ \dot{\widetilde M}_k(t)=M_*(t) \left\{\frac{p^{[k]}(t)}{P(t)} +\frac{q^{[k]}(t)}{Q(t)} -\frac{\dot{\widetilde C}_k}{C_*}\right\}. $$ The derivative $\dot M_k^{\rm c}$ also contains the derivative of the displayed positive normalisation factor $\sqrt{S/S^+}$, and $\dot\theta_k$ is obtained from the general returned derivatives above by differentiating the simple root of the returned core quadratic. When $k\in\mathcal N$, one has $m_k=-1$ on this six-node fixed component and these formulas reduce to \eqref{eq:B.A1.4}--\eqref{eq:B.A1.9}. All quantities are obtained from positive definite Gram systems. Their denominators remain nonzero on the compact signed charts, so $\sup_{h,k,y}|\chi_{hk}(y)|<\infty$. Hence the exact external first-jet terms, including the $y^2E_k$-terms, are $O(\mathcal E)$, and the stable part satisfies $$ \left|\sum_{k\notin\mathcal N}\chi_{hk}(y)E_k\right|\le C_hE_{\rm st}. $$

The restrictions to the two fixed components are explicit. On the five-node component, $a=y=\mathcal E=0$ and the remaining completion root is $\alpha+e$, so \begin{equation*} \left.\mathcal B_h\right|_{(y,a)=0} =2\log\left| 1-\frac{2C_*+D_J(h)e}{C_*+x}\right|. \tag{B.A1.64}\label{eq:B.A1.64} \end{equation*} On the six-node component, $a=e=\mathcal E=0$, and \begin{equation*} \left.\mathcal B_h\right|_{(a,e)=0} =2\log\left|1-\frac{2C_*}{C_*+x}\right|, \tag{B.A1.65}\label{eq:B.A1.65} \end{equation*} independently of the free tangent variable $y$. After subtraction of the displayed $-4x/C_*+2D_J(h)e/C_*$, equations \eqref{eq:B.A1.64}--\eqref{eq:B.A1.65} have zero value and zero first derivatives in $x,e$ at the incidence point, with a uniform $O(x^2+e^2)$ remainder. Equation \eqref{eq:B.A1.63} gives all derivatives in the $E_k$-directions.

After the same subtraction, \eqref{eq:B.A1.62} gives at every positive six-node fixed point \begin{equation*} \partial_aR_h^{\rm Lag} =\frac{2\beta_\alpha(h)\delta D_J(\sigma)} {C_*g(g-\delta)}, \qquad \partial_eR_h^{\rm Lag} =\frac{2\beta_\alpha(h)}{C_*} \{D_J(\sigma)-D_J(\alpha)\}. \tag{B.A1.66}\label{eq:B.A1.66} \end{equation*} The first coefficient is nonzero for $y\ne0$. We isolate both coefficients as explicit anisotropic terms.

\subsubsection*{B.A1.5.3. Anisotropic remainder for the distinct-root cocycle}

Define \begin{equation*} A_h^{\rm coc}(y,\vartheta) =\frac{2\beta_\alpha(h)\delta D_J(\sigma)} {C_*g\Delta}, \qquad B_h^{\rm coc}(y,\vartheta) =\frac{2\beta_\alpha(h)}{C_*} \{D_J(\sigma)-D_J(\alpha)\}. \tag{B.A1.67}\label{eq:B.A1.67} \end{equation*} The parameter $\vartheta$ contains the root, factor, and norm coordinates tangent to the fixed component. Since $\delta=u_yy^2$, while $g,\Delta,C_*$ are positive units, $$ |A_h^{\rm coc}|\le Cy^2. $$ Also, the fundamental theorem of calculus gives the exact identity \begin{equation*} D_J(\sigma)-D_J(\alpha) =-\delta\int_0^1D_J'(\alpha-t\delta)\,dt, \end{equation*} and hence \begin{equation*} |A_h^{\rm coc}|+|B_h^{\rm coc}|\le C_hy^2. \tag{B.A1.69}\label{eq:B.A1.69} \end{equation*}

For the remainder bound, work in a regular signed state chart $r$ and regard $x,a,e$ as analytic functions on it. At a factor collision these functions are not assumed to be independent coordinates. The exact full multiplier in \eqref{eq:B.8.36}, the normalised-core recurrence \eqref{eq:B.5.26}, and the transfer \eqref{eq:B.5.25} define an analytic log-return $\mathcal B_h(r,E)$.

Here is the intrinsic expansion in the state ring. Put $D_6(t)=D_J(t)(t-\alpha)$ and $\mathcal J=(a,e)=(U_0,U_1)$. On the six-node fixed component, with the norm $C=C_*+x$ allowed to vary, the exact cocycle value is
\[
G_h(C)=2\log\left|1+\frac{D_6(h)}C\right|,
\qquad D_6(h)=-2C_*.
\]
Its normal differential is obtained without inverting a factor resultant. The first variations of the consecutive squared monic norms agree on the fixed set. Linearising \eqref{eq:B.4.10} therefore gives $d(PQ-C)=Q_U$, where $Q_U=\operatorname{quo}_K(D_6U)$ and $U$ here denotes the affine normal variation. The full multiplier at $h$ is $(D_6(h)+C)/C$, the retained multipliers are $1$, and \eqref{eq:B.A1.61} gives the pivot contribution. Thus
\begin{equation*}
\begin{aligned}
d\mathcal B_h&=G_h'(C)\,dC+\mathcal L_{h,C}(U),\\
\mathcal L_{h,C}(U)
 &=\frac{2Q_U(h)}{D_6(h)+C}
   +\frac2C\sum_{i\in J}\beta_i(h)Q_U(\lambda_i)
   +\frac{2\beta_\alpha(h)D_J(\theta)}{C\Delta}\,a.
\end{aligned}
\tag{B.A1.69a}\label{eq:B.A1.69a}
\end{equation*}
At $C=C_*$, interpolation gives $\mathcal L_{h,C_*}(U)=2\beta_\alpha(h)D_J(\sigma)q/C_*$, while $G_h'(C_*)=-4/C_*$. Substitution of $q=e+\delta a/(g\Delta)$ recovers \eqref{eq:B.A1.62}. Polynomial division by the monic $K$ and the nonzero denominators in \eqref{eq:B.A1.69a} show that its coefficients are analytic functions of $C$ and the current root data, including at a factor collision.

Subtract $G_h(C)$ and this linear normal expression from the intrinsic log-return. On every generic coprime factor branch the difference has zero value and zero first normal jet, so it belongs to $\mathcal J^{(2)}$. Lemma~\ref{lem:4.3} gives $\mathcal J^{(2)}=\mathcal J^2$ in the actual signed state ring, with every collision branch included. Moreover,
\[
G_h(C)+4x/C_*=O(x^2),\qquad
\mathcal L_{h,C}(U)-\mathcal L_{h,C_*}(U)\in x\mathcal J.
\]
After the linear $e$ term in $R_h^{\rm Lag}$ and the two coefficients in \eqref{eq:B.A1.67} are subtracted, it follows that the intrinsic remainder belongs to $(x,a,e)^2$. Analytic ideal membership thus supplies the expansion on every signed state chart.

Consequently there are analytic coefficients $A_{ij}(r)$, with $N=(x,a,e)$, such that the full remainder has the exact representation
\begin{equation*}
\begin{aligned}
\widetilde R_h
 &:=R_h^{\rm Lag}-A_h^{\rm coc}a-B_h^{\rm coc}e\\
 &=\sum_{i,j=1}^3 N_iN_j A_{ij}(r)
   +\sum_k E_k B_{hk}(r,E).
\end{aligned}
\tag{B.A1.70}\label{eq:B.A1.70}
\end{equation*}
The external coefficients can be chosen by the exact integral formula
\[
B_{hk}(r,E)=\int_0^1\partial_{E_k}\widetilde R_h(r,tE)\,dt.
\]
At a fixed incidence and $E=0$ they equal the first jets $\chi_{hk}$ in \eqref{eq:B.A1.63}. Thus all external--core products and products of external squared masses are retained.

We justify uniform boundedness of these coefficients. At a pivot endpoint, $M_i^{\rm c}=\widetilde M(\lambda_i)\sqrt{S/S^+}$ is a positive analytic unit extending the logarithmic quotient even when $w_i=0$. At a factor collision every retained node evaluation of $P,Q$ remains nonzero, the moment Gram matrices stay invertible, and the root gaps $g,\Delta$ and the normalisations have positive lower bounds. The log-return, its external derivatives, and the intrinsic analytic ideal coefficients therefore exist on neighbourhoods of every point of the compact signed cover. Shrinking these neighbourhoods and taking a finite subcover bounds the chosen coefficients. Equation \eqref{eq:B.A1.70} now gives
\begin{equation*}
|\widetilde R_h|
 \le C_h\{(|x|+|a|+|e|)^2+\mathcal E\}.
\tag{B.A1.71}\label{eq:B.A1.71}
\end{equation*}
The only first-order $a,e$ terms are the isolated $A_h^{\rm coc}a,B_h^{\rm coc}e$. The fixed value $G_h(C)$ is independent of $y$, and its subtracted $x$-jet is zero, so no $y^2x$ term remains. Every term containing an external squared mass, including a $y^2E_k$ term, is $O(\mathcal E)$.

With the tail bounds \eqref{eq:B.9.30n}, \eqref{eq:B.A1.69}--\eqref{eq:B.A1.71} give \begin{equation*} |R_{h,n}^{\rm Lag}| \le Cy_n^2(|\mathfrak a_n|+|e_n|) +C\{(|x_n|+|\mathfrak a_n|+|e_n|)^2+\mathcal E_n\} =o(X_n). \tag{B.A1.72}\label{eq:B.A1.72} \end{equation*} The cocycle estimate is uniform in the neighbourhood. Only now, after the boundary extension has been proved, let $\mathcal B_{h,0}^{\rm Lag}$ denote the same log-return with stable masses set to zero and the remaining exact weights renormalised. The extended log-return is even in every stable signed amplitude and hence analytic in its squared mass. The mean-value theorem in those squared masses and \eqref{eq:B.6.63} give \begin{equation*} |\mathcal B_h^{\rm Lag}-\mathcal B_{h,0}^{\rm Lag}|\le C_hE_{\rm st}=o(E_{\rm neu}). \tag{B.A1.72a}\label{eq:B.A1.72a} \end{equation*} This is an estimate for the extended log-return, not for the logarithm of the scalar itself.

\subsubsection*{B.A1.5.4. An adapted transverse coordinate near a support boundary}

For \eqref{eq:B.9.30y}, use the coordinates $$ v=y^2,\qquad q=-\frac{U(\sigma)}{\Delta} =e+\frac{\delta}{g\Delta}a,\qquad d_\theta=(\theta-\lambda_2)(\lambda_3-\theta). $$ Hence $q-e=O(va)$. This change is uniformly analytic, agrees with the five-node transverse coordinate at $v=0$, and is the $\sigma$-root normal while the six normalised-core weights are positive. At a pivot endpoint let $p$ be the signed amplitude of the disappearing pivot coordinate. Recovery gives \begin{equation*} p^2=u_p d_\theta,\qquad u_p>0 \tag{B.A1.73}\label{eq:B.A1.73} \end{equation*} with $u_p$ analytic. Use the moment normal for the whole six-node core \begin{equation*} \mathfrak A_y=\sum_{i\in J\cup\{\alpha\}} \frac{w_i}{\lambda_i-\theta},\qquad B_y=\sum_{i\in J\cup\{\alpha\}} \frac{w_i}{(\lambda_i-\theta)^2}. \end{equation*} At the fixed component, combine the intrinsic identities $$ d(\theta^+-\theta)=\frac2{B_y}\,d\mathfrak A_y $$ and \eqref{eq:B.A1.61}. This gives \begin{equation*} d\mathfrak A_y =-\frac{B_ygD_J(\theta)}{C_*\Delta}\,da. \tag{B.A1.75}\label{eq:B.A1.75} \end{equation*} Now $d_\theta B_y$ and $D_J(\theta)/d_\theta$ are nonzero analytic units at either endpoint. Hence the coefficient in \eqref{eq:B.A1.75} is a nonzero analytic unit. Thus, on the fixed component, the conormal differentials $d\mathfrak A_y$ and $da$ differ by a uniformly nonvanishing analytic factor. This is a first-order statement on that component. It does not identify the two functions up to a unit on the whole neighbourhood; mixed and quadratic terms away from the fixed component must be retained.

\begin{appendixBlemma}[Generators for the fixed strata after two support losses]
\label{lem:BA1-endpoint-generators}
Work in a regular signed state chart $r=(\zeta,p,y)$, where $\zeta$ denotes three independent retained-state coordinates. The physical quantities $a,q$ are analytic functions on this chart; they need not be independent coordinates at a factor collision. Let $\lambda_p$ be the disappearing pivot node, and let $\eta_p,\eta_y$ be the bracket generators of Lemma~\ref{lem:4.3} for the nodes $\lambda_p,\alpha$. There is an analytic matrix $\mathbf M$ with determinant $\alpha-\lambda_p\ne0$ such that
\begin{equation*}
\binom{p\eta_p}{y\eta_y}=\mathbf M\binom{pa}{yq}.
\tag{B.A1.75a}\label{eq:B.A1.75a}
\end{equation*}
Consequently the union of the endpoint fixed strata has radical complete-intersection ideal $\mathcal I_{\rm fix}=(pa,yq)$, and $\mathcal I_{\rm fix}^{(m)}=\mathcal I_{\rm fix}^m$ for every $m\ge1$.
\end{appendixBlemma}

\begin{proof}
Write $K=(t-\theta)(t-\sigma)$, $\Delta=\sigma-\theta$, $g=\alpha-\theta$, $\delta=\alpha-\sigma$, and $b=\theta-\lambda_p$. The exact recovery formula \eqref{eq:B.4.9a} gives
\[
\delta=y^2 u_y,\qquad b=p^2 u_b,\qquad
u_y=\frac{P(\alpha)}{H_0c_\alpha g},\qquad
u_b=\frac{P(\lambda_p)}{H_0c_{\lambda_p}(\sigma-\lambda_p)}.
\]
Here $c_{\lambda_p}$ denotes the barycentric coefficient at the pivot node. These are analytic units on the chosen chart; $u_y>0$, while the sign of $u_b$ depends on which pivot endpoint is used.
For $U=U_1t+U_0$, the definitions $a=-U(\theta)$ and $q=-U(\sigma)/\Delta$ give $U_1=a/\Delta-q$. Since $\bar K_\gamma=t+\gamma-\theta-\sigma$, the bracket definitions give the exact identities
\begin{equation*}
\begin{aligned}
\eta_p&=-U(\theta-\delta)
       =\frac g\Delta a-y^2u_yq,\\
\eta_y&=-U(\sigma+b)
       =-p^2\frac{u_b}{\Delta}a+(\sigma-\lambda_p)q.
\end{aligned}
\tag{B.A1.75c}\label{eq:B.A1.75c}
\end{equation*}
Thus
\begin{equation*}
\mathbf M=\begin{pmatrix}
g/\Delta&-pyu_y\\
-pyu_b/\Delta&\sigma-\lambda_p
\end{pmatrix},\qquad
\det\mathbf M
 =\frac{g(\sigma-\lambda_p)-\delta b}{\Delta}
 =\alpha-\lambda_p.
\tag{B.A1.75d}\label{eq:B.A1.75d}
\end{equation*}
This calculation is an identity in the affine coefficients of $U$, valid whether or not their differentials are independent. It uses only spectral gaps, positive norms, and nonzero node evaluations as denominators, and hence extends through every allowed factor collision. The invertible generator change transfers the ideal and power conclusions of Lemma~\ref{lem:4.3} to $(pa,yq)$.
\end{proof}

We now derive the intrinsic transverse remainder in the same actual state ring. After the external linear coefficients \eqref{eq:B.A1.53} are subtracted, put $R_q^{\rm int}=R_q|_{\mathcal E=0}$, with $R_q$ defined by \eqref{eq:B.A1.54}.

The support-type ideals are
\begin{equation*}
(p,y),\qquad(p,q),\qquad(a,y),\qquad(a,q).
\tag{B.A1.75e}\label{eq:B.A1.75e}
\end{equation*}
At a collision, their individual minimal primes must be distinguished.

\begin{appendixBlemma}[Transverse square through the collision branches]
\label{lem:BA1-collision-square}
The intrinsic remainder satisfies
\begin{equation*}
R_q^{\rm int}\in\mathcal I_{\rm fix}^{(2)}
 =\mathcal I_{\rm fix}^2=(pa,yq)^2.
\tag{B.A1.75f}\label{eq:B.A1.75f}
\end{equation*}
\end{appendixBlemma}

\begin{proof}
The fixed union is the union of the four support-type zero sets in \eqref{eq:B.A1.75e}. At a collision, each support-type ideal is decomposed into its individual factor-allocation minimal primes, as in Lemma~\ref{lem:4.3}.

At a generic point of the four-node component, the return is the identity on that face and separate evenness kills both absent-amplitude derivatives. Hence $R_q^{\rm int}$ has zero value and zero first differential there. On a generic five-node component $(a,y)=0$, the exact face return and its first normal jet are those in \eqref{eq:B.A1.55}; the absent-amplitude derivative vanishes by evenness. The other five-node component $(p,q)=0$ uses the same five-node calculation with the omitted nodes interchanged. The added term $c_\theta d_\theta a^2$ has zero first jet on both components, since $a=0$ on the first and $d_\theta$ is a unit times $p^2$ on the second. Thus the remainder has zero value and first differential on both generic five-node components. On every coprime six-node branch, $U^+-U\in(U_0,U_1)^2$ and the root-motion formulas imply the same zero first jet for the physical $q$ return. These statements hold for every fixed value of the norm and factor parameters, so their tangential derivatives vanish as well.

For clarity, consider separately a proper common quadratic factor. By \eqref{eq:B.4.17b2} and the signed recovery substitution \eqref{eq:B.4.19a}, the six-node fixed quotient is the node ring with the two free signed variables $p,y$. Lift its two node coordinates to $u,v$ in the regular state ring and straighten a nonzero linear generator to $z$. The quotient cotangent space has basis $p,y,u,v$, so together with $z$ these form an analytic state chart by the inverse-function theorem. In these coordinates
\[
(a,q)=(z,uv)=(z,u)\cap(z,v).
\]
The two ideals $(z,u)$ and $(z,v)$ are the two distinct factor-allocation branches. They each have dense coprime points; the zero first jets just computed therefore extend to each branch. In these coordinates ordinary analytic division gives
\[
(z,u)^2\cap(z,v)^2
 =(z^2,zuv,u^2v^2)=(z,uv)^2.
\]
Indeed, the constant coefficient in $z$ is divisible by both $u^2$ and $v^2$, while the coefficient of $z$ is divisible by both $u$ and $v$; the remaining terms are divisible by $z^2$. This explains the local square statement without treating $a,q$ as independent coordinates or treating their common zero set as one prime component. A common linear factor has the regular normal form in \eqref{eq:B.4.17b2}. In the present upper-face application the neutral external mode has $P(h)Q(h)=-C_*<0$, so the self-dual case $P=Q$ is excluded.

Finally, let $\mathfrak p$ range over \emph{all} minimal primes of $\mathcal I_{\rm fix}$. The preceding generic calculations give
\[
R_q^{\rm int}\in(\mathcal I_{\rm fix}R_{\mathfrak p})^2
\quad\hbox{for every such }\mathfrak p.
\]
Lemma~\ref{lem:BA1-endpoint-generators} makes the fixed ideal a radical complete intersection; its square has no embedded primary component. Intersecting these localised squares, using \eqref{eq:B.4.17e}, proves \eqref{eq:B.A1.75f} in the original signed state ring. The argument refines every support component into its collision branches before localising.
\end{proof}

The cocycle's nonzero first jet remains isolated in \eqref{eq:B.A1.66}--\eqref{eq:B.A1.67}.

Write the membership in \eqref{eq:B.A1.75f} as $$ R_q^{\rm int}=A(pa)^2+B(pa)(yq)+C(yq)^2. $$ This representation is not unique, so we choose it equivariantly. Average the identity over $(p,y)\mapsto(\varepsilon_pp,\varepsilon_yy)$, where $\varepsilon_p,\varepsilon_y\in\{\pm1\}$. Since $R_q^{\rm int}$ is even in $p$ and $y$, the identity is preserved if $A,C$ are replaced by their even--even averages and $B$ by $$ B_{--}(p,y,\ldots)=\frac14\sum_{\varepsilon_p,\varepsilon_y=\pm1}\varepsilon_p\varepsilon_y B(\varepsilon_pp,\varepsilon_yy,\ldots). $$ The analytic coefficient $B_{--}$ is odd in each amplitude, and hence $B_{--}=py\,\widetilde B$ with $\widetilde B$ analytic and even in both. The mixed generator is therefore $p^2y^2aq\,\widetilde B=u_pd_\theta v aq\,\widetilde B$, and the positive unit $u_p$ is absorbed into the coefficient below. Restore the external variables in \eqref{eq:B.A1.54} by multivariable Taylor division. In the upper-face case every $k\in\mathcal N$ lies below $\alpha$, so the exact neutral first jets from \eqref{eq:B.A1.53} have negative sign. Using \eqref{eq:B.A1.73} gives \begin{equation*} \begin{aligned} q^+-q={}&-c(v,\vartheta)d_\theta a^2 -\sum_{k\in\mathcal N} \rho_k^{\rm ext}(v,q,\vartheta)E_k\\ &+d_\theta v\,a q\,f_1(v,\vartheta) +v q^2f_2(v,\vartheta)+\mathcal R_E+\mathcal R_{\rm st}. \end{aligned} \tag{B.A1.76}\label{eq:B.A1.76} \end{equation*} The coefficients are analytic and bounded functions of the signed state; their displayed arguments suppress the remaining state dependence and do not assert independence of $a,q$ at a collision. The function $\rho_k^{\rm ext}(v,q,\vartheta)$ is the analytic re-expression of $\rho_k^{\rm ext}(y,e,\vartheta)$ under the coordinate change above. Equations \eqref{eq:B.8.11b} and \eqref{eq:B.A1.53} give, after the neighbourhood is shrunk, \begin{equation*} c(v,\vartheta)\ge c_0>0,\qquad \rho_k^{\rm ext}(v,q,\vartheta)\ge \rho_{\min}^{\rm ext}>0 \tag{B.A1.77}\label{eq:B.A1.77} \end{equation*} for every $k\in\mathcal N$. Terms in $\mathcal R_E$ contain either two external squared weights or one external weight multiplied by a normal variable tending to zero; the exact linear and mixed contributions of stable modes are placed in $\mathcal R_{\rm st}$. Along the tail, \begin{equation*} |\mathcal R_E|\le\epsilon_n\mathcal E,\qquad |\mathcal R_{\rm st}|\le C E_{\rm st}, \qquad \epsilon_n\longrightarrow0. \end{equation*} Tangential norm and factor parameters occur only in the bounded coefficients in \eqref{eq:B.A1.76}, since both fixed strata are fixed for every value of those parameters. Consequently, $x$ enters \eqref{eq:B.A1.76} only through its bounded coefficients.

We remove the sign-indefinite $vq^2$ term. Substitution of $e=q-\delta a/(g\Delta)$ in \eqref{eq:B.A1.60} gives the exact simplification \begin{equation*} Q_U(\alpha)=\frac{D_J(\theta)}{\Delta}a-D_J(\sigma)q. \end{equation*} Since $v=w_\alpha$ is a normalised-core weight, its exact law under the full return is \begin{equation*} v^+=v\widetilde M(\alpha)^2\frac S{S^+}=v\{M_\alpha^{\rm c}\}^2,\qquad M_\alpha^{\rm c}:=\widetilde M(\alpha)\sqrt{S/S^+}. \tag{B.A1.78b}\label{eq:B.A1.78b} \end{equation*} The positive square root is used. At zero external mass on the six-node fixed component, $M_\alpha^{\rm c}=1$. The exact Uvarov derivatives and the $C^2$ transfer \eqref{eq:B.5.25} therefore give \begin{equation*} v^+-v =v\{\ell_v q+\ell_{v,a} a +O(a^2+q^2+\mathcal E)\}, \tag{B.A1.79}\label{eq:B.A1.79} \end{equation*} where $\ell_v,\ell_{v,a}$ are analytic along the fixed component, including the norm parameter, and at $C=C_*$ satisfy \begin{equation*} \ell_v=-\frac{2D_J(\sigma)}{C_*}>0,\qquad \ell_{v,a}=\frac{2D_J(\theta)}{C_*\Delta}. \tag{B.A1.80}\label{eq:B.A1.80} \end{equation*} The first coefficient remains a uniformly positive analytic unit because $\sigma$ stays in the fixed positive gap below $\alpha$. Incorporating the norm parameter into these coefficients accounts for every mixed radial term. At a collision, subtract the displayed intrinsic first jet and apply the ordinary-square statement for $(U_0,U_1)=(a,q)$ in Lemma~\ref{lem:4.3}; this gives the same $O(a^2+q^2)$ remainder in the actual state ring. The external linear terms remain in $O(\mathcal E)$.

Retain the separately even analytic coefficient $f_2$ supplied by the actual state-ring decomposition in \eqref{eq:B.A1.76}. Use the analytic extension $\ell_v=-2D_J(\sigma)/C$, where $C$ is the current core norm coordinate. This agrees with the fixed-set coefficient and remains a positive unit. Define \begin{equation*} \kappa=-\frac{f_2}{\ell_v},\qquad \widehat q=q+\kappa vq. \end{equation*} Both $\kappa$ and its first derivatives are bounded, and $1+\kappa v$ is a positive unit on a smaller neighbourhood.

The increment satisfies \begin{equation*} v^+-v=v\{\ell_v q+\ell_{v,a} a+R_v\},\qquad |R_v|\le C(a^2+q^2+\mathcal E). \tag{B.A1.81a}\label{eq:B.A1.81a} \end{equation*} Because $D_J(\theta)=d_\theta u_\theta$ and all denominators in \eqref{eq:B.A1.80} are units, \begin{equation*} |\ell_{v,a}|\le Cd_\theta. \end{equation*} The increment of $\kappa$ can be controlled in the genuine signed state chart. Every component of the intrinsic signed return minus the identity vanishes on the fixed union, and hence belongs to $(pa,yq)$ by Lemma~\ref{lem:BA1-endpoint-generators}. For a retained component, separate evenness in $p,y$ turns this representation into $p^2a A+y^2q B$ with analytic coefficients. For the two omitted components, the corresponding parity projections give
\[
p^+-p=p\{a A_p+y^2q B_p\},\qquad
 y^+-y=y\{q A_y+p^2a B_y\}.
\]
It follows that the increments of the retained state coordinates, $p^2$, and $y^2$ are bounded by $C(p^2|a|+y^2|q|)$ on a smaller chart. The full/core transfer adds $O(\mathcal E)$. Both $f_2$ and $\ell_v$, hence $\kappa$, are analytic in these retained coordinates and in $p^2,y^2$. The mean-value theorem and \eqref{eq:B.A1.73} therefore give
\begin{equation*}
|\kappa^+-\kappa|
 \le C\{d_\theta|a|+v|q|+\mathcal E\}
 \le C\{d_\theta|a|+v|q|
       +d_\theta a^2+vq^2+\mathcal E\}.
\end{equation*}
This bound uses separate parity and the fixed ideal, and remains valid when $da,dq$ are dependent.

The two exact product identities are \begin{equation*} (vq)^+-vq=(v^+-v)q+v(q^+-q) +(v^+-v)(q^+-q), \end{equation*} \begin{equation*} \widehat q^+-\widehat q =q^+-q+\kappa\{(vq)^+-vq\} +(\kappa^+-\kappa)v^+q^+. \end{equation*} Substitute \eqref{eq:B.A1.76} and \eqref{eq:B.A1.81a} into these identities. The term $\kappa v\ell_v q^2=-vf_2q^2$ cancels the displayed $vf_2q^2$ in \eqref{eq:B.A1.76}. Every uncancelled term is bounded by \begin{equation*} \begin{aligned} C\{& v(d_\theta a^2+\mathcal E)+v^2q^2 +vd_\theta|aq|\\ &+(|a|+|q|+\mathcal E) (d_\theta a^2+vq^2+\mathcal E)\} +|\mathcal R_E|+|\mathcal R_{\rm st}|. \end{aligned} \end{equation*} For example, the first line comes respectively from $v(q^+-q)$, from the $vq^2$-part of that product and the $v|q|\,|\kappa^+-\kappa|$ term, and from $v\ell_{v,a}aq$. The product $(v^+-v)(q^+-q)$ and the remaining change of $\kappa$ give the second line. Thus all mixed terms, including the external--core products generated by the exact normalisation factor in \eqref{eq:B.A1.78b}, are accounted for.

Use $$ vd_\theta|aq| \le\varepsilon d_\theta a^2 +C_\varepsilon d_\theta v^2q^2 \le\varepsilon d_\theta a^2+C_\varepsilon v^2q^2. $$ Along the orbit, after reducing $c_0$ if necessary to the minimum in \eqref{eq:B.A1.77}, \begin{equation*} \widehat q^+-\widehat q \le-(c_0-\varepsilon)(d_\theta a^2+E_{\rm neu}) +C_\varepsilon v^2q^2 +\epsilon_n\{d_\theta a^2+vq^2+\mathcal E\}+CE_{\rm st}, \tag{B.A1.83}\label{eq:B.A1.83} \end{equation*} where $\epsilon_n\to0$. This estimate is uniform in the sign of $f_2$ and keeps the stable forcing in the exact recurrence.

\subsubsection*{B.A1.5.5. Tail energy and endpoint selection}

The adapted transverse coordinate has the desired sign. Its remainder is summable by the parity chord identity, whose leading signed directions admit a quantitative lower bound. In one endpoint/collision signed state chart put $r=(\zeta,p,y)$, with $a,q$ regarded as analytic functions of $r$, let $\xi=(\xi_h)_h$ be the surviving external signed amplitudes, and write the pair of exact full-state parity chords as $$ \Xi_{\rm full}(r,\xi)=\bigl(Y^+-Y,Z^+-Z\bigr). $$ At $\xi=0$ every component of $\Xi_{\rm full}$ vanishes on the union of fixed strata. Lemma~\ref{lem:BA1-endpoint-generators} identifies its analytic vanishing ideal as $\mathcal I_{\rm fix}=(pa,yq)$. Componentwise membership in this analytic ideal, followed by successive Hadamard division in the external amplitudes, therefore gives the exact vector identity \begin{equation*} \Xi_{\rm full}(r,\xi)=pa\,V_a(r,\xi)+yq\,V_q(r,\xi)+\sum_h\xi_hW_h(r,\xi) \tag{B.A1.83a}\label{eq:B.A1.83a} \end{equation*} with analytic full-state vector coefficients. This is an identity on a neighbourhood, not a Taylor expansion in the degenerate products $pa,yq$.

For each such chart let $\Pi$ select, from the first chord $Y^+-Y$, the actual signed spectral coordinates at the disappearing pivot node $\lambda_p$, at $\alpha$, and at every surviving external node. At $\xi=0$, coordinate-factor invariance, the exact divisions in \eqref{eq:B.A1.75c}--\eqref{eq:B.A1.75d}, and the separate signed parities give \begin{equation*} \binom{(Y^+-Y)_p}{(Y^+-Y)_\alpha} =C_{\rm core}\binom{pa}{yq},\qquad C_{\rm core}=\begin{pmatrix}c_a&py\,b_{12}\\py\,b_{21}&c_q\end{pmatrix},\qquad \sigma_{\min}(C_{\rm core})\ge d/2. \tag{B.A1.83b}\label{eq:B.A1.83b} \end{equation*} Here $b_{12},b_{21}$ are bounded and analytic. The five-node signed return \eqref{eq:B.4.24f} and \eqref{eq:B.A1.75} give, on the limiting fixed slice $x=0$ of the fixed union, $$ c_a=\frac{gD_J(\theta)}{C_*\Delta(\lambda_p-\theta)}, $$ whose apparent endpoint singularity cancels to a nonzero product of spectral gaps. Equations \eqref{eq:B.A1.79}--\eqref{eq:B.A1.80} similarly give $$ c_q=\ell_v/2=-D_J(\sigma)/C_*. $$ Thus both diagonal coefficients are analytic units. Compactness and a shrink of the chart give a common $d>0$ for their moduli and make the off-diagonal operator norm at most $d/2$, proving the last inequality in \eqref{eq:B.A1.83b}.

The core entries of $W_h(r,0)$ vanish because every reconstructed core coordinate is even in $\xi_h$. The exact $h$-coordinate of the first chord is $\xi_h(\widetilde M_h-1)$, and coordinate-factor invariance kills its derivatives in the other external rows at $\xi=0$. Hence, for the exact coefficient matrix $$ \mathcal C=\bigl(\Pi V_a,\Pi V_q,(\Pi W_h)_h\bigr), $$ its restriction to $\xi=0$ on the limiting fixed slice $x=0$ of the fixed union is $C_{\rm core}\oplus\operatorname{diag}_h(m_h-1)$. Every surviving limiting external multiplier satisfies $m_h\ne1$; neutral modes have $m_h=-1$, and the finitely many stable modes are uniformly separated from $1$. The generator change in \eqref{eq:B.A1.75d} has the nonzero determinant $\alpha-\lambda_p$, so it and its inverse have uniformly bounded norms on the compact chart. The projected output coordinates are actual signed spectral coordinates, and the coefficient identities extend analytically through the allowed factor collisions. Continuity, followed by a finite-cover minimum and a further shrink, therefore gives \begin{equation*} \sigma_{\min}(\mathcal C)\ge\sigma_0>0. \tag{B.A1.83c}\label{eq:B.A1.83c} \end{equation*} Applying $\Pi$ to the exact identity \eqref{eq:B.A1.83a} now yields $$ \|\Xi_{\rm full}\|\ge\|\Pi\Xi_{\rm full}\| \ge\sigma_0\left(|pa|^2+|yq|^2+\sum_h|\xi_h|^2\right)^{1/2}. $$ Finally, \eqref{eq:B.A1.73} gives $p^2=u_pd_\theta$ with $u_p$ bounded below by a positive constant, while $y^2=v$ and $\xi_h^2=E_h$. Squaring proves the exact estimate \begin{equation*} \|Y^+-Y\|^2+\|Z^+-Z\|^2 \ge c\{d_\theta a^2+vq^2+\mathcal E\} \end{equation*} uniformly through positive-unit endpoints and simple or simultaneous collisions. Here $Y,Z$ and both chords are exact full states, while $a,q,v$ are normalised-core coordinates.

The exact two-parity chord identity \eqref{eq:B.3.4}, with the exact full monic norms bounded above and below, and the definition $X_n=\tau-\widetilde a_n$, now give \begin{equation*} \sum_{k=n}^{\infty} \{d_{\theta,k}\mathfrak a_k^2+v_kq_k^2+\mathcal E_k\} \le CX_n. \tag{B.A1.85}\label{eq:B.A1.85} \end{equation*} Stable coordinates remain in the exact state and are already included in $\mathcal E_k$ and in the diagonal block of \eqref{eq:B.A1.83c}; no shadow deletion is used in this chord estimate.

Choose $\varepsilon<c_0/2$ in \eqref{eq:B.A1.83} and put $$ H_k=d_{\theta,k}\mathfrak a_k^2+v_kq_k^2+\mathcal E_k. $$ Since $v_n\to0$, \eqref{eq:B.A1.85} yields \begin{equation*} \sum_{k=n}^{\infty}v_k^2q_k^2 \le\left(\sup_{k\ge n}v_k\right) \sum_{k=n}^{\infty}v_kq_k^2=o(X_n), \end{equation*} and, more generally, every remainder carrying an extra factor $v_k$ has tail bounded by $C\sup_{k\ge n}v_k\sum_{k\ge n}H_k=o(X_n)$. Likewise, \begin{equation*} \sum_{k=n}^{\infty}|\epsilon_k|H_k \le\left(\sup_{k\ge n}|\epsilon_k|\right)\sum_{k=n}^{\infty}H_k=o(X_n). \end{equation*} By Lemma~\ref{lem:5.2}, specifically \eqref{eq:B.5.15} and its proof, $$ \sum_{k=n}^{\infty}E_{{\rm st},k}=O(E_{{\rm st},n}). $$ Equation \eqref{eq:B.6.63} gives $E_{{\rm st},n}=o(E_{{\rm neu},n})$, while \eqref{eq:B.5.13}, summed over the finitely many persistent neutral modes, gives $E_{{\rm neu},n}=O(X_n^2)=o(X_n)$. Hence $$ \sum_{k=n}^{\infty}E_{{\rm st},k}=o(X_n). $$ Thus \eqref{eq:B.A1.83} has the form \begin{equation*} \widehat q_{k+1}-\widehat q_k \le-c(d_{\theta,k}\mathfrak a_k^2+E_{{\rm neu},k})+\mathfrak r_k, \qquad \sum_{k=n}^{\infty}|\mathfrak r_k|=o(X_n). \tag{B.A1.87}\label{eq:B.A1.87} \end{equation*} Because $U_n\to0$, both $q_n$ and $\widehat q_n$ tend to zero. Summing \eqref{eq:B.A1.87} from $n$ to infinity gives $$ \widehat q_n\ge-o(X_n). $$ Equation \eqref{eq:B.9.30n} and $q_n=e_n+O(v_n\mathfrak a_n)$ give $|q_n|+|\mathfrak a_n|\le CX_n$. Since $v_n\to0$, $$ |\widehat q_n-q_n|\le C v_n|q_n|\le Cv_nX_n=o(X_n),\qquad |q_n-e_n|\le Cv_n|\mathfrak a_n|\le Cv_nX_n=o(X_n). $$ Therefore \begin{equation*} e_n\ge-o(X_n). \tag{B.A1.88}\label{eq:B.A1.88} \end{equation*}

Together, equations \eqref{eq:B.A1.72} and \eqref{eq:B.A1.88} yield \eqref{eq:B.9.30y}. They remain valid as $y\downarrow0$, through every pivot endpoint and factor collision. After exact deletion, coordinate-factor invariance reduces the argument to the corresponding lower-support lemma.

At the upper boundary face, $\alpha=\lambda_5\notin J$, and the distinct-root cocycle applies throughout the positive six-node neighbourhood. After exact deletion, if the remaining completion root $c$ coincides with a retained node, then $\mathscr D_*=(t-c)^2R$ has one double root at that retained node and four other simple retained roots. The Hermite scalar \eqref{eq:B.8.47} then applies, with value and derivative cardinals $H_{c,0},H_{c,1}$; its return is \eqref{eq:B.8.48}. If no such coincidence occurs, the distinct-root scalar \eqref{eq:B.8.32} remains the applicable formula. Reversal of the node order proves the lower-boundary version, with the oriented normal coordinate $q=-e+O(va)$.

Together with Subsections B.A1.1--B.A1.4, this completes the algebraic and boundary estimates invoked in Lemmas~\ref{lem:9.1}--\ref{lem:9.4}.

\section*{Part C. Counterexamples from restart length four onward} \gdef\proofpart{C} \renewcommand{\thesection}{C.\arabic{section}} \setcounter{section}{0} \setcounter{theorem}{0}

Part C proves the negative assertion in Theorem~\ref{thm:sharp}.

\section{Counterexamples at restart lengths at least four}\label{sec:counterexamples}

\begin{theorem}[Counterexamples at every restart length at least four]\label{thm:negative}
For every integer $s\ge4$, there exist $$ n=s+4,\qquad A\in\mathbb R^{n\times n},\qquad b,x_0\in\mathbb R^n, $$ such that $A$ is real symmetric positive definite, the restarted conjugate-gradient iteration with restart length $s$, carried out in exact arithmetic, never terminates, and its even normalised residuals do not converge in Euclidean norm. \end{theorem}

\begin{proposition}[Certified configuration at restart length four]\label{prop:s4-seed}
There is an eight-dimensional real symmetric positive definite restarted CG problem at restart length four whose orbit is nonterminating and whose even normalised residuals do not converge. The matrix and initial data are specified uniquely by the algebraic isolating conditions and the two contractions below. The dyadic constants in these contractions are \begin{equation*} \varepsilon_0=2^{-500000000},\qquad M=2^{600000000},\qquad N=M^2=2^{1200000000}. \tag{C.1.1}\label{eq:C.1.1} \end{equation*}
\end{proposition}

\subsection{Outline of the construction}

The proof has two stages. At restart length four, Sections~\ref{restarted-cg-as-an-exact-spectral-map}--\ref{exact-transfer-of-the-hopf-point-to-the-raw-field} construct a transverse Hopf point in scalar factor coordinates and identify it with the leading field of the rescaled squared-weight map. Sections~\ref{literal-selection-of-one-periodic-orbit}--\ref{the-one-fixed-exact-spd-datum} select a nonconstant periodic profile, shadow it by an exact orbit, and reconstruct a fixed eight-dimensional SPD instance. Sections~\ref{sec:C10}--\ref{sec:C18} then adjoin one distant core node at a time, preserve the Hopf crossing, and repeat the analytic periodic-orbit and shadowing arguments at each higher restart length. Appendix C identifies the finite assertions needed for the restart-four seed; degree elevation is analytic.

\subsection{Norm and coordinate conventions}

In Part C, finite coefficient and state vectors carry the $\ell^\infty$ norm unless a subscript $2$ or an explicitly Euclidean chord is shown. For a matrix, $\|M\|_\infty$ is the induced row-sum norm, whereas $|M|_{\max}=\max_{ij}|M_{ij}|$ is the entrywise maximum. Multilinear maps carry the operator norm induced by these vector norms, and finite product spaces use the maximum of their factor norms. For a map on a box $\mathcal B$, \begin{equation*} \|F\|_{C^k(\mathcal B)} =\max_{0\le|\gamma|\le k} \sup_{z\in\mathcal B}\|\partial^\gamma F(z)\|. \end{equation*} All boxes, including \eqref{eq:C.4.6}, are coordinatewise product boxes, equivalently $\ell^\infty$ balls with possibly different coordinate radii.

For $2\pi$-periodic functions use \begin{equation*} \begin{aligned} \widehat h_k &=\frac1{2\pi}\int_0^{2\pi}h(\theta)e^{-ik\theta}\,d\theta,\\ \|h\|_{L^2}^2 &=2\pi\sum_k|\widehat h_k|_2^2,\\ \|h\|_{H^1}^2 &=2\pi\sum_k(1+k^2)|\widehat h_k|_2^2. \end{aligned} \end{equation*} The Lyapunov--Schmidt domain $H^1_{\rm per}\times\mathbb R^2$ uses the maximum product norm and the range uses $L^2_{\rm per}$. These conventions also fix every operator norm in \eqref{eq:C.6.4}--\eqref{eq:C.6.14} and every sequence norm in Section~\ref{literal-selection-of-one-exact-map-shadow}.

\section{The spectral map for restarted CG}\label{restarted-cg-as-an-exact-spectral-map}

This section recasts the spectral reduction from the common core in the notation needed for the eight-node seed and isolates the nontermination criterion.

Let $A\in\mathbb R^{n\times n}$ be symmetric positive definite, let $r_k=b-Ax_k$, and fix a restart length $s$. At one restart step, $x_{k+1}$ is the unique vector in $x_k+\mathcal K_s(A,r_k)$ such that
$$ x_*-x_{k+1}\perp_A\mathcal K_s(A,r_k), \qquad x_*=A^{-1}b, $$
where $\mathcal K_s(A,r)=\operatorname{span}\{r,Ar,\ldots,A^{s-1}r\}$. Equivalently, it minimises the $A$-norm of the error over that affine Krylov space. We denote the residual at the next restart point by $r_{k+1}$.

Let $0<\lambda_1<\cdots<\lambda_m$, let $A=\operatorname{diag}(\lambda_1,\ldots,\lambda_m)$, and let a normalised residual have coordinates
\begin{equation*} y=\sum_i\eta_i\mathbf e_i,\qquad w_i=\eta_i^2,\qquad \sum_iw_i=1. \end{equation*} Here $\mathbf e_i$ are the standard coordinate vectors.

For restart length four, one CG block has residual
\begin{equation*} r^+=p_w(A)r,\qquad p_w(t)=1+\sum_{\ell=1}^4c_\ell t^\ell, \end{equation*}
where Galerkin orthogonality is exactly
\begin{equation*} \sum_iw_i p_w(\lambda_i)\lambda_i^j=0, \qquad 0\le j\le3. \tag{C.2.3}\label{eq:C.2.3} \end{equation*}
The coefficient matrix in \eqref{eq:C.2.3}, after shifting the column index, is the Gram matrix of $1,t,t^2,t^3$ for the positive measure $\sum_iw_i\lambda_i\delta_{\lambda_i}$. It is positive definite whenever at least four active nodes are distinct. Hence $p_w$ is unique. If the resulting residual is nonzero, then after normalisation the squared spectral weights obey
\begin{equation*} (T_4w)_i= \frac{w_i p_w(\lambda_i)^2} {\sum_jw_jp_w(\lambda_j)^2}. \tag{C.2.4}\label{eq:C.2.4} \end{equation*}

Translation of all nodes by five replaces a centred residual factor $R$ by $t\mapsto R(t-5)/R(-5)$; the common normalising constants cancel from \eqref{eq:C.2.4}.

In Section~\ref{restarted-cg-as-an-exact-spectral-map} only, the nodes are listed in increasing spectral order. From Section~\ref{the-exact-algebraic-hopf-seed} onward we use the fixed chart order ``six core nodes, then the two external nodes''. That order need not be increasing. Formula \eqref{eq:C.2.4}, grade, and every conclusion below are invariant under this one fixed permutation.

If every $w_i>0$ and every $p_w(\lambda_i)\ne0$, the next residual still has all $m$ active nodes. For $m=8>4$, its grade is eight and a block with restart length four cannot terminate. It suffices to construct a positive forward orbit of $T_4$ on eight nodes, with no zero block factor, whose even squared weights do not converge.

\section{The certified Hopf point at restart length four}\label{the-exact-algebraic-hopf-seed}

The seed consists of a six-node two-cycle of the core squared-weight map and two external nodes whose squared two-block multipliers equal one. The parameter $q$ moves the associated leading field through a transverse Hopf crossing. The following rational polynomials realise this arrangement.

Every terminating decimal in this section denotes the corresponding exact rational number. Define
\begin{equation*} \begin{aligned} P(t)&=\prod_{r\in R_P}(t-r),\\ R_P&=\{-3.812196375334383,\,-0.3408184613601878,\\ &\hspace{2.5em}1.5054785540322202,\,2.190808973956946\},\\ Q_q(t)&=(t-q)\prod_{r\in R_Q}(t-r),\\ R_Q&=\{-2.3304189906225,\,0.6733921795528344,\\ &\hspace{2.5em}4.323444838404325\},\\ C&=0.31355365301369736, \end{aligned} \tag{C.3.1}\label{eq:C.3.1} \end{equation*}
and
\begin{equation*} \mathfrak f_-(t,q)=P(t)Q_q(t)-C,\qquad \mathfrak f_+(t,q)=P(t)Q_q(t)+C. \end{equation*}

Define
\begin{equation*} Q_H=[-1.851888497035,-1.851888497033]. \tag{C.3.3}\label{eq:C.3.3} \end{equation*}

The certificate described in Lemma~\ref{lemma-3.1-certified-critical-branch-and-transverse-hopf-point} gives the following rational isolating intervals; in fact they are uniform on the larger interval $[-1.858,-1.851]\supset Q_H$: \begin{equation*} \begin{array}{ccl|ccl} &\mathfrak f_-&\text{role}&&\mathfrak f_+&\text{role}\\ \hline 1&[-3.812224,-3.812222]&\xi_1&1&[-3.812171,-3.812168]&--\\ 2&[-2.329779,-2.329767]&\xi_2&2&[-2.331069,-2.331056]&h_-\\ 3&[-1.859053,-1.852044]&\eta_-&3&[-1.856952,-1.849955]&--\\ 4&[-0.339465,-0.339456]&\xi_3&4&[-0.342178,-0.342169]&--\\ 5&[0.671419,0.671428]&\xi_4&5&[0.675358,0.675367]&--\\ 6&[1.507016,1.507022]&\eta_+&6&[1.503932,1.503938]&--\\ 7&[2.190298,2.190302]&\xi_5&7&[2.191315,2.191319]&h_+\\ 8&[4.323453,4.323455]&\xi_6&8&[4.323434,4.323437]&-- \end{array} \tag{C.3.3a}\label{eq:C.3.3a} \end{equation*} Endpoint Sturm counts give one root in every isolating interval; strict opposite boundary signs and a derivative interval excluding zero continue each label for every $q$ in the stated interval. The six selected core intervals have pairwise separation greater than one, and the combined six-core/two- external label set has minimum separation greater than $2^{-11}$. Hence every factor in $|\Delta'(\xi_i)|$ exceeds one. These assertions are certified by the programs described in Appendix C. Let $\eta_-,\eta_+$ be the third and sixth ordered roots of $\mathfrak f_-$, and define
\begin{equation*} H(t)=(t-\eta_-)(t-\eta_+),\qquad \Delta(t)=\mathfrak f_-(t,q)/H(t), \tag{C.3.4}\label{eq:C.3.4} \end{equation*}
and denote the six roots of $\Delta$ by $\xi_1<\cdots<\xi_6$. Let $h_-,h_+$ be the second and seventh ordered roots of $\mathfrak f_+$. The two phases of the fixed-point family have weights
\begin{equation*} w_i^P(a)=\frac{( \xi_i-a)Q_q(\xi_i)}{\Delta'(\xi_i)}, \qquad w_i^Q(a)=\frac{( \xi_i-a)P(\xi_i)}{\Delta'(\xi_i)}. \tag{C.3.5}\label{eq:C.3.5} \end{equation*}

The degree-five Lagrange identity $$ \sum_{\Delta(\xi)=0}\frac{R(\xi)}{\Delta'(\xi)}=[t^5]R(t), \qquad \deg R\le5, $$ shows that both vectors sum to one. It also gives, for $0\le j\le3$, $$ \sum_i w_i^P(a)P(\xi_i)\xi_i^j=0, \qquad \sum_i w_i^Q(a)Q_q(\xi_i)\xi_i^j=0. $$ The same identity with $(t-a)P(t)$ and $(t-a)Q_q(t)$ gives squared monic norm $C$ in the two phases. At each core node $P(\xi_i)Q_q(\xi_i)=C$, and therefore $$ T_4w^P(a)=w^Q(a),\qquad T_4w^Q(a)=w^P(a). $$ Thus \eqref{eq:C.3.5} is an exact two-cycle. At the external nodes $P(h_\pm)Q_q(h_\pm)=-C$, so their squared two-block multipliers are equal to one.

For an external node $g$, define
\begin{equation*} A_g(a,q)=\frac{4P(g)}{g-a} \left\{\frac{P(g)}{P(a)}-\frac{Q_q(g)}{Q_q(a)}\right\}, \qquad \ell_g(a,q)=-\frac{2}{C} \frac{\Delta(g)-\Delta(a)}{g-a}. \tag{C.3.6}\label{eq:C.3.6} \end{equation*}

The tangent forcing vectors $V_g(a,q)\in\mathbb R^3$ are obtained by the following finite polynomial calculation. For $R=P,Q_q$, divide $\Delta=A_RR+B_R$, set
\begin{equation*} \beta_R=\frac{B_R(a)}{R(a)},\qquad U_R(t)=\frac{\beta_RR(t)-B_R(t)}{C(t-a)}, \tag{C.3.7}\label{eq:C.3.7} \end{equation*}
and use the divided kernel
\begin{equation*} K_R(t,g)=\frac{R(t)U_R(g)-U_R(t)R(g)}{t-g}, \tag{C.3.8}\label{eq:C.3.8} \end{equation*}
with its derivative value when $t=g$. These are the reproducing kernels of $\mathcal P_3$ for the two measures in \eqref{eq:C.3.5} determined by the phases of the fixed-point family. This follows from the degree-five Lagrange coefficient identity after substituting $\Delta=A_RR+B_R$ in \eqref{eq:C.3.7}.

The polynomial division used below is exact. Indeed, if $x$ is a root of $\Delta$, substitution of \eqref{eq:C.3.7}--\eqref{eq:C.3.8}, using $P(x)Q_q(x)=C$, gives \begin{equation*} \begin{aligned} P(x)K_{Q_q}(x,g)-Q_q(x)K_P(x,g) ={}&P(g)-Q_q(g)\\ &-\frac{\Delta(a)}{(g-a)(x-a)} \left\{\frac{P(g)}{P(a)}-\frac{Q_q(g)}{Q_q(a)}\right\}. \end{aligned} \tag{C.3.8a}\label{eq:C.3.8a} \end{equation*} All denominators are nonzero on the certified parameter box. Differentiating the orthogonality equations after inserting an external squared mass at $g$ gives the derivatives of the two phase polynomials $$ p_g=-4P(g)K_P(\,\cdot\,,g),\qquad q_g=4P(g)K_{Q_q}(\,\cdot\,,g). $$ This is the Uvarov formula, obtained by pairing the differentiated moment equations with an arbitrary polynomial in $\mathcal P_3$.

Define $(V_g)_C=4(P(g)^2+C)$ and $$ \mathcal N_g(t)=p_g(t)Q_q(t)+P(t)q_g(t)-(V_g)_C -A_g\frac{\Delta(t)-\Delta(a)}{t-a}. $$ Since $P(g)Q_q(g)=-C$, equation \eqref{eq:C.3.8a} gives $\mathcal N_g(x)=0$ at each of the six distinct roots of $\Delta$. Moreover, $\deg\mathcal N_g\le7$. Hence there is a unique affine polynomial $(V_g)^H$ such that \begin{equation*} \boxed{\quad p_gQ_q+Pq_g-(V_g)_C -A_g\frac{\Delta(t)-\Delta(a)}{t-a} =\Delta(t)(V_g)^H(t).\quad} \tag{C.3.9}\label{eq:C.3.9} \end{equation*} Equation \eqref{eq:C.3.9} therefore defines the remaining two components of $V_g$ by polynomial division, without interpolation at selected external labels.

The same calculation identifies these quantities with the derivative of the two-block scalar map. Differentiating $$ SAA_+=\Delta(SH+U)+C_+S_+,\qquad S=t-a, $$ in the direction obtained by inserting an external mass gives $$ S\{p_gQ_q+Pq_g-(V_g)_C-\Delta(V_g)^H\} =\Delta\,\delta U-C(\delta a^+-\delta a). $$ After \eqref{eq:C.3.9} is substituted, comparison of the coefficient of $\Delta$ and the constant remainder yields \begin{equation*} \delta U=A_g,\qquad \delta a^+-\delta a=\frac{A_g\Delta(a)}C. \tag{C.3.9a}\label{eq:C.3.9a} \end{equation*} The second term has weighted degree two under the scaling in Section~\ref{a-complete-physical-raw-chart-and-its-exact-expansion} and so does not enter the leading $a$-equation. Finally,
\begin{equation*} d_g(q)[\delta C,\delta H] =-\frac{2}{C}\{\Delta(g)\delta H(g)+2\delta C\}. \tag{C.3.10}\label{eq:C.3.10} \end{equation*}
For $h,g\in\{h_-,h_+\}$, define $\mathsf T=(\mathsf T_{hg})_{h,g}$ immediately from this functional by $\mathsf T_{hg}=d_h[V_g]$. Substitution of \eqref{eq:C.3.9} gives the equivalent expanded formula
\begin{equation*} \boxed{\quad \mathsf T_{hg}=-\frac{2}{C}\left\{p_g(h)Q_q(h)+P(h)q_g(h)+(V_g)_C-A_g\frac{\Delta(h)-\Delta(a)}{h-a}\right\}.\quad} \tag{C.3.10a}\label{eq:C.3.10a} \end{equation*}
Thus the single occurrence of $(V_g)_C$ in \eqref{eq:C.3.10a} is the sum of the term $-(V_g)_C$ supplied by \eqref{eq:C.3.9} and the term $2(V_g)_C$ in \eqref{eq:C.3.10}. This is the same identity used in \eqref{eq:C.13.2}. These quantities are rational functions of the isolated algebraic roots, with every denominator separated from zero by the certificate. We identify $B=(B_C,B_0,B_1)\in\mathbb R^3$ with $(\delta C,\delta H)$ by $B_C=\delta C$ and $\delta H(t)=B_1t+B_0$. We also use $e=(e_-,e_+)$ and $z=(a,u,B,e)$.

These formulae define the algebraic ingredients of the vector field below. Sections~\ref{a-complete-physical-raw-chart-and-its-exact-expansion} and~\ref{exact-transfer-of-the-hopf-point-to-the-raw-field} show that, after weighted rescaling and passage to the simplex chart, this field is the leading term of the exact two-block map. The next lemma contains the finite input needed for that analytic passage.

\begin{lemma}[Certified critical branch and transverse Hopf point]\label{lemma-3.1-certified-critical-branch-and-transverse-hopf-point}
There is a unique analytic function
\begin{equation*} a=\alpha(q)\in[0.2952193,0.29521955],\qquad q\in Q_H. \end{equation*}

There are also positive analytic amplitudes $c_-(q),c_+(q)$. With these functions, define the seven-dimensional weighted field by
\begin{equation*} \begin{aligned} \dot a&=\kappa(a)u,\qquad \kappa(a)=2\Delta(a)/C,\\ \dot u&=u+\beta(a,q)u^2+\sum_{g=h_-,h_+}A_g(a,q)e_g,\\ \dot B&=B+\mathbf q(a,q)u^2+\sum_gV_g(a,q)e_g,\\ \dot e_g&=e_g\{2+\ell_g(a,q)u+d_g(q)[B]\}. \end{aligned} \tag{C.3.12}\label{eq:C.3.12} \end{equation*}
Denote the right-hand side of \eqref{eq:C.3.12} by $G_q(z)$. This field has the equilibrium
\begin{equation*} z_*(q)= \left(\alpha(q),0,-\sum_gc_g(q)V_g(\alpha(q),q), c_-(q),c_+(q)\right). \tag{C.3.13}\label{eq:C.3.13} \end{equation*}

Here $\beta,\mathbf q$ are defined exactly by the following three-phase coefficient recursion. At one phase, the factor relation is \begin{equation*} SAA_+=\Delta\Sigma+C_+S_+, \qquad S=t-a, \tag{C.3.14}\label{eq:C.3.14} \end{equation*} Coefficient vectors are ordered by ascending powers in the basis $(1,t,\ldots,t^8)$, with the parametrisations \begin{equation*} \begin{aligned} A_+(t)&=t^4+\sum_{j=0}^3a_j^+t^j,\\ \Sigma(t)&=t^3+\sum_{j=0}^2\sigma_jt^j,\\ R_+(t)&:=C_+S_+(t)=\varrho_0+\varrho_1t. \end{aligned} \tag{C.3.14a}\label{eq:C.3.14a} \end{equation*} The nine unknowns are ordered as \begin{equation*} x=(a_0^+,a_1^+,a_2^+,a_3^+, \sigma_0,\sigma_1,\sigma_2,\varrho_0,\varrho_1)^T. \tag{C.3.14b}\label{eq:C.3.14b} \end{equation*} For a polynomial $f$, let $\operatorname{sh}_j(f)$ be the coefficient vector of $t^jf$ padded to length nine. The common phase matrix is \begin{equation*} \begin{aligned} \mathcal M(A,S,\Delta) &=[\,\mathcal A(A,S)\mid\mathcal D(\Delta)\mid\mathcal R\,],\\ \mathcal A(A,S) &=[\,\operatorname{sh}_0(SA)\ \operatorname{sh}_1(SA)\ \operatorname{sh}_2(SA)\ \operatorname{sh}_3(SA)\,],\\ \mathcal D(\Delta) &=[\,-\operatorname{sh}_0(\Delta)\ -\operatorname{sh}_1(\Delta)\ -\operatorname{sh}_2(\Delta)\,],\\ \mathcal R&=[\,-\mathbf e_0\ -\mathbf e_1\,]. \end{aligned} \tag{C.3.14c}\label{eq:C.3.14c} \end{equation*} Here $\mathbf e_0,\mathbf e_1$ are the first two standard coordinate vectors in $\mathbb R^9$; the displayed objects, from left to right, are the nine columns. The order-zero right-hand side is the coefficient vector of $-SAt^4+\Delta t^3$. Expand every input and output through order two in an intrinsic parameter $\epsilon$. At order $m=1,2$, use the order-zero matrix \eqref{eq:C.3.14c}; the right-hand side is $$ [\epsilon^m]\{-SAt^4+\Delta t^3\} $$ minus every convolution containing a previously computed lower-order coefficient of $A_+,\Sigma$, or $R_+$. Thus the same order-zero matrix is used at orders $m=1,2$, with the lower-order convolution terms moved to the right-hand side.

After a solve, monicity gives $C_+=\varrho_1$ and $S_+=R_+/\varrho_1=t+\varrho_0/\varrho_1$. Recover the coordinates of the \emph{incoming} phase by dividing by its incoming factor $S$:
\begin{equation*}
\Sigma=SH_{\rm in}+U_{\rm in},\qquad \deg U_{\rm in}<1.
\tag{C.3.14d}\label{eq:C.3.14d}
\end{equation*}
The division is performed coefficientwise at each order. The scalar norm coordinate attached to this incoming phase is the partner norm $C_+$ obtained in the same solve. More explicitly, index consecutive inputs by $(A^{(j)},S^{(j)})$, $j=0,1,2$. The $j$th solve gives $(A^{(j+1)},S^{(j+1)},\Sigma^{(j)},C^{(j)})$; define $H^{(j)},U^{(j)}$ by $\Sigma^{(j)}=S^{(j)}H^{(j)}+U^{(j)}$. The state coordinates are $(S^{(j)},C^{(j)},H^{(j)},U^{(j)})$. Starting with $[\epsilon]U^{(0)}=1$ and zero first-order $(C^{(0)},H^{(0)})$ increments, perform the consecutive $P,Q,P$ solves. The input and returned states are those with indices $0$ and $2$. The third solve recovers the coordinates of the state reached after two blocks; it does not add a third block to the return. In the scalar coordinate order $B=(B_C,B_0,B_1)$, define \begin{equation*} \beta=[\epsilon^2](U^{(2)}-U^{(0)}), \qquad \mathbf q=[\epsilon^2](B^{(2)}-B^{(0)}). \tag{C.3.14e}\label{eq:C.3.14e} \end{equation*} Here $B^{(j)}$ consists of the increments of $C^{(j)}$ and the constant and linear coefficients of $H^{(j)}$ from their fixed-state values. Equations \eqref{eq:C.3.14a}--\eqref{eq:C.3.14e} determine the nonlinear coefficients in this incoming-phase convention. The intrinsic-lift routine listed in Appendix~C constructs the same shifted columns, convolution right-hand sides, monic recovery, and quotient/remainder recursion, and certifies both phase matrices uniformly.

There is a unique $q_*$ in the interior of $Q_H$ for which $J_*^{\rm sc}=D_zG_{q_*}(z_*(q_*))$ has spectrum
\begin{equation*} -1,\ 1,\ 2,\ \pm i\omega_*,\ 1\pm i\omega_*, \qquad 3.72<\omega_*<3.724. \tag{C.3.15}\label{eq:C.3.15} \end{equation*}
All seven eigenvalues are simple, no other one is imaginary, and
\begin{equation*} \left|\frac{d}{dq}\operatorname{Re}\lambda_H(q_*)\right|>1.94. \tag{C.3.16}\label{eq:C.3.16} \end{equation*}

The Hopf eigenvector has a nonzero $a$-component. Uniformly on the certified branch,
\begin{equation*} \min_i(w_i^P,w_i^Q)>2^{-19},\qquad c_->2^{-15},\qquad c_+>1/4, \end{equation*}
and $q_*$ is more than $2^{-41}$ from either endpoint of $Q_H$. \end{lemma}

\begin{proof}
Equating the two external tangent rates in \eqref{eq:C.3.6}--\eqref{eq:C.3.10} gives the critical equation. Its derivative with respect to $a$ lies between $-161517.159$ and $-161516.321$, and its values on the two $a$-faces have opposite strict signs. The interval implicit-function theorem gives the unique analytic branch. The normal-balance and two tangent-rate equations then determine its two positive amplitudes explicitly.

The characteristic polynomial follows from an explicit linearisation. At an equilibrium, let the logarithmic external variations be $\varphi_g=\delta e_g/c_g$, and define \begin{equation*} \begin{gathered} k=\kappa(\alpha),\quad p_{\rm lin}=\sum_gc_gA_g'(\alpha),\quad \mathbf A_c=(c_gA_g(\alpha))_g,\\ q_{\rm lin}=\sum_gc_gV_g'(\alpha),\quad \mathsf V_c=(c_gV_g(\alpha))_g,\quad \ell=(\ell_g(\alpha))_g,\quad D B=(d_g[B])_g. \end{gathered} \tag{C.3.17a}\label{eq:C.3.17a} \end{equation*} The linearisation of \eqref{eq:C.3.12} is \begin{equation*} \lambda a=ku,\qquad \lambda u=p_{\rm lin}a+u+\mathbf A_c\varphi,\qquad \lambda B=q_{\rm lin}a+B+\mathsf V_c\varphi,\qquad \lambda\varphi=\ell u+DB. \tag{C.3.17b}\label{eq:C.3.17b} \end{equation*} With $\mathbf1=(1,1)^T$, the equilibrium equations are \begin{equation*} \mathbf A_c\mathbf1=0,\qquad D\mathsf V_c\mathbf1=2\mathbf1. \tag{C.3.17c}\label{eq:C.3.17c} \end{equation*} The first is the normal balance; $B_*=-\mathsf V_c\mathbf1$ and $0=2\mathbf1+DB_*$ give the second.

For $\lambda\ne0,1$, let $\theta=\lambda(\lambda-1)$ and $\zeta=\lambda\varphi$. Eliminating $a,B$ from \eqref{eq:C.3.17b} gives \begin{equation*} (\theta-kp_{\rm lin})u-\mathbf A_c\zeta=0, \qquad (\theta I_2-D\mathsf V_c)\zeta-(\theta\ell+kDq_{\rm lin})u=0. \tag{C.3.17d}\label{eq:C.3.17d} \end{equation*} With $\chi=\zeta-\ell u$, this is the ordinary eigenvalue problem \begin{equation*} \theta\binom u\chi= \mathcal M\binom u\chi,\qquad \mathcal M= \begin{pmatrix} kp_{\rm lin}+\mathbf A_c\ell&\mathbf A_c\\ D\mathsf V_c\ell+kDq_{\rm lin}&D\mathsf V_c \end{pmatrix}. \end{equation*} Equation \eqref{eq:C.3.17c} gives $\mathcal M(0,\mathbf1)^T=2(0,\mathbf1)^T$. Moreover the general determinant reduction of \eqref{eq:C.3.17b} is \begin{equation*} \det(\lambda I-J_*^{\rm sc})=(\lambda-1) \det\{\lambda(\lambda-1)I_3-\mathcal M\}. \tag{C.3.17f}\label{eq:C.3.17f} \end{equation*} To verify \eqref{eq:C.3.17f}, perform the eliminations leading to \eqref{eq:C.3.17d} for $\lambda\ne0,1$, multiply back the two divisions by $\lambda$ and $\lambda-1$, and compare the resulting monic polynomials of degree seven. Equality on that infinite set of $\lambda$'s proves the polynomial identity, including at $0,1$. The explicit factor $\lambda-1$ is also seen directly: $D:\mathbb R^3\to\mathbb R^2$ has a nonzero kernel, and $(a,u,\varphi,B)=(0,0,0,B_{\ker})$, $DB_{\ker}=0$, is a $\lambda=1$ vector.

If $q_\theta$ denotes the monic quadratic factor obtained after removing the eigenvalue $2$ of $\mathcal M$, the characteristic polynomial factors exactly as
\begin{equation*} (\lambda+1)(\lambda-1)(\lambda-2) q_\theta(\lambda(\lambda-1)). \end{equation*}

The two scalar coefficients evaluated by the certificate are precisely those of this quadratic: $n_g=A_g(\alpha)$, $c=(c_-,c_+)^T$, $\mathsf T=D(V_-,V_+)$, $\varrho_{\rm bal}=c_-n_-=-c_+n_+$, $g_{\rm lin}=\mathsf T'(\alpha)c=Dq_{\rm lin}$, $\delta_\ell=\ell_- -\ell_+$, $\delta_g=(g_{\rm lin})_--(g_{\rm lin})_+$, $\sigma=c_+\mathsf T_{-,+}+c_-\mathsf T_{+,-}$, $d=2-\sigma$. Because $\mathsf T\operatorname{diag}(c)\mathbf1=2\mathbf1$, its two eigenvalues are $2,d$, and $(1,-1)$ is a left eigenvector for $d$. Hence \begin{equation*} \operatorname{tr}\mathcal M-2 =kp_{\rm lin}+\varrho_{\rm bal}\delta_\ell+d=:\tau. \end{equation*} For $d\ne0$, block elimination and that left eigenvector give \begin{equation*} \frac{\det\mathcal M}{2} =k\{dp_{\rm lin}-\varrho_{\rm bal}\delta_g\}=:K_\theta; \tag{C.3.20}\label{eq:C.3.20} \end{equation*} explicitly, $\det(D\mathsf V_c)=2d$, $\mathbf A_c(D\mathsf V_c)^{-1}=\mathbf A_c/d$, and the term $\mathbf A_c\ell=\varrho_{\rm bal}\delta_\ell$ cancels from the Schur complement. Both sides of \eqref{eq:C.3.20} are polynomial in the displayed entries, so the identity also holds when $d=0$. Hence \begin{equation*} q_\theta(\theta)=\theta^2-\tau\theta+K_\theta. \tag{C.3.21}\label{eq:C.3.21} \end{equation*} This is the $(\tau,q_0)$ calculation in the certificate, where the certificate denotes $K_\theta$ by $q_0$.

Define \begin{equation*} \mathscr D_\theta=\tau^2-4K_\theta. \end{equation*} On the whole branch $\mathscr D_\theta<-55.40$. At the two rational endpoints of $Q_H$, the Hopf functional $\mathcal H_H=\tau/2-\mathscr D_\theta/4$ has opposite strict signs; its derivative along the critical branch lies in $(109.89,110.82)$, and hence has one simple zero. At that zero the two $\theta$-roots are $-\omega_*^2\pm i\omega_*$, whose preimages are precisely the four nonreal eigenvalues in \eqref{eq:C.3.15}. If the imaginary-axis branch is chosen compatibly, differentiating $\theta=\lambda(\lambda-1)$ at $\lambda=i\omega_*$ gives \begin{equation*} \frac{d}{dq}\operatorname{Re}\lambda_H =-\frac{\mathcal H_H'}{1+4\omega_*^2}. \end{equation*} The certified intervals $3.72<\omega_*<3.724$ and $109.89<\mathcal H_H'<110.82$ give \eqref{eq:C.3.16}.

If the Hopf mode had $u=0$, \eqref{eq:C.3.17d} would give $\mathbf A_c\zeta=0$. Here $\mathbf A_c=\varrho_{\rm bal}(1,-1)$ with $\varrho_{\rm bal}\ne0$, so $\zeta$ is proportional to $\mathbf1$. A nonzero eigenvector has $\zeta\ne0$, and the second equation in \eqref{eq:C.3.17d}, together with $D\mathsf V_c\mathbf1=2\mathbf1$, then forces $\theta=2$. This contradicts the certified nonreal $\theta$-root. Since $\dot a=\kappa u$ and $\kappa\ne0$, the $a$-component is nonzero.

The rational interval and Sturm computations establishing these inequalities are described in Appendix C. Equations \eqref{eq:C.3.17a}--\eqref{eq:C.3.21} give the characteristic reduction by linear algebra. \end{proof}

\section{Weight coordinates and the return-map expansion}\label{a-complete-physical-raw-chart-and-its-exact-expansion}

The certificate of Section~\ref{the-exact-algebraic-hopf-seed} is formulated in scalar polynomial coordinates. To construct a restarted CG orbit, we need a chart on the weight simplex and a return-map expansion that is uniform under the singular scaling used below.

Translate every centred node by five, so all eight nodes lie in $(1,10)$. On the $P$-phase of the fixed-point family, let $w_i^*(a)=w_i^P(a)$ and
\begin{equation*} s_1=\partial_aw_1^*(a)>1/16. \end{equation*} Formula \eqref{eq:C.3.5} is affine in $a$, so $s_1$ is independent of $a$; the displayed lower bound holds uniformly for every $q\in Q_H$.

For a nearby eight-component weight vector define
\begin{equation*} \begin{aligned} a_{\rm raw}&=\frac{w_1-w_1^*(0)}{s_1},\\ \zeta_j&=w_{j+1}-w_{j+1}^*(a_{\rm raw}),\quad1\le j\le4,\\ E_-&=w_7,\qquad E_+=w_8. \end{aligned} \tag{C.4.2}\label{eq:C.4.2} \end{equation*}
The sixth core weight is fixed by $\sum_iw_i=1$. We denote the affine weight-to-coordinate map in \eqref{eq:C.4.2} by $\chi_q$ and abbreviate $E=(E_-,E_+)$. For each fixed $q$, $\chi_q$ and its inverse are affine in the state variables, with coefficients analytic in $q$. In the rescaled weight chart we relabel $a_{\rm raw}$ as $a$; Section~\ref{exact-transfer-of-the-hopf-point-to-the-raw-field} temporarily denotes it by $a_{\rm wt}$ when comparing it with the scalar recovery root.

For each basis perturbation $p_j=t^j$, $0\le j\le3$, the first-phase order-one system \eqref{eq:C.3.14a}--\eqref{eq:C.3.14d} has an output column ordered as $(\delta C,\delta H_0,\delta H_1,\delta U)^T$. These four columns form $J_{\rm lift}$. On the full certified parameter box its interval matrix is contained in {\small \setlength{\arraycolsep}{3pt} \begin{equation*} J_{\rm lift}\in \begin{pmatrix} [12,13]&[-222,-221]&[452,453]&[-2582,-2581]\\ [-729,-728]&[1693,1694]&[-9071,-9070]&[28934,28935]\\ [485,486]&[-1124,-1123]&[6038,6039]&[-19226,-19225]\\ [-4,-3]&[0,1]&[-26,-25]&[45,46] \end{pmatrix}. \end{equation*} } The midpoint-residual Neumann defect is below $1/4$, and the certified induced row-sum bound is $\|J_{\rm lift}^{-1}\|_\infty<2^{20}$.

This lift identifies the normal coordinate used in the weighted scaling. For $\mathbf n=(u,B_C,B_0,B_1)$, define \begin{equation*} p=J_{\rm lift}^{-1}(B_C,B_0,B_1,u)^T,\qquad \delta w_i=w_i^*(a) \left\{\frac{B_C}{C}-\frac{p(\xi_i)}{P(\xi_i)}\right\}, \tag{C.4.2a}\label{eq:C.4.2a} \end{equation*} and \begin{equation*} (K_{\rm wt}(a)\mathbf n)_j=\delta w_{j+1} -\frac{\partial_aw_{j+1}^*(a)}{s_1}\delta w_1, \qquad1\le j\le4. \end{equation*} This $K_{\rm wt}(a)$ is invertible. If $K_{\rm wt}(a)\mathbf n=0$, then $\sum_i\delta w_i=\sum_i\partial_aw_i^*=0$ shows that $\delta w$ is a scalar multiple of the tangent to the fixed-point family on all six core coordinates. Differentiating the four moment equations first gives $$ \sum_i\delta w_iP(\xi_i)\xi_i^j+ \sum_iw_i^*p(\xi_i)\xi_i^j=0,\qquad0\le j\le3. $$ The first sum is zero because it is the derivative of the same moment equation in the $a$-direction tangent to the fixed-point family, where $P$ is fixed. Substitution of $p(t)=\sum_{\ell=0}^3p_\ell t^\ell$ gives \begin{equation*} \sum_{\ell=0}^3\left(\sum_iw_i^*\xi_i^{j+\ell}\right)p_\ell=0, \qquad0\le j\le3. \end{equation*} The matrix is a positive definite Gram matrix for cubics on six distinct nodes with positive weights, so $p=0$, and invertibility of $J_{\rm lift}$ gives $\mathbf n=0$. At the scalar equilibrium let \begin{equation*} \mathbf n_*=(0,B_*),\qquad B_*=-\sum_gc_gV_g(\alpha,q),\qquad v_*(q)=K_{\rm wt}(\alpha(q))\mathbf n_*. \end{equation*} Estimating the lift and barycentric formulas gives \begin{equation*} |v_*|+|c_-|+|c_+|<2^{500001}. \tag{C.4.2e}\label{eq:C.4.2e} \end{equation*} Hence $v_*$ is obtained algebraically from the certified data, with its dependence on $q$ included in the derivative estimates below.

The scalar normal is ordered as $\mathbf n=(u,B_C,B_0,B_1)$, while the right-hand side in \eqref{eq:C.4.2a} follows the lift row order $(B_C,B_0,B_1,u)$.

\begin{lemma}[Uniform bounds for the exact weight map]\label{lemma-4.1-uniform-exact-raw-map-bounds}
For $q\in Q_H$, on
\begin{equation*} |a_{\rm raw}-\alpha(q)|\le2^{-24},\qquad |\zeta|+|E_-|+|E_+|\le2^{-1000000}, \tag{C.4.3}\label{eq:C.4.3} \end{equation*}
the signed analytic extension of the exact two-block map $F_q=\chi_qT_4^2\chi_q^{-1}$ has the joint parameter/state bound
\begin{equation*} \max_{|\gamma|\le5} \left\|\partial_{(q,a_{\rm raw},\zeta,E)}^\gamma F_q\right\| <2^{1000000}. \tag{C.4.4}\label{eq:C.4.4} \end{equation*}

Let $v_*(q)$ be the representation of \eqref{eq:C.3.13} in the coordinates of weighted degree one. Use the uncentred scaled variables
\begin{equation*} \zeta=tV,\qquad E_g=t^2\mathcal E_g,\qquad t^+=\frac{t}{1+t}. \tag{C.4.5}\label{eq:C.4.5} \end{equation*}
Here $\mathcal E=(\mathcal E_-,\mathcal E_+)$ and $z=(a,V,\mathcal E)$.

For
\begin{equation*} |a-\alpha(q)|\le2^{-24},\quad |V-v_*(q)|\le1,\quad |\mathcal E_g-c_g(q)|\le2^{-20},\quad |t|\le2^{-2000000}, \tag{C.4.6}\label{eq:C.4.6} \end{equation*}
the exact rescaled map has the two-sided analytic expansion
\begin{equation*} \widehat F_q(t,z) =z+tG_q^{\rm wt}(z)+t^2R_q(t,z), \tag{C.4.7}\label{eq:C.4.7} \end{equation*}
Here $G^{\rm wt}(q,z)=G_q^{\rm wt}(z)$ and $R(t,q,z)=R_q(t,z)$. These functions satisfy the joint bounds
\begin{equation*} \|G^{\rm wt}\|_{C^3(Q_H\times\mathcal Z)}<2^{3000000}, \qquad \sup_{|t|\le2^{-2000000}} \|R(t,\cdot)\|_{C^1(Q_H\times\mathcal Z)}<2^{3000000}, \tag{C.4.8}\label{eq:C.4.8} \end{equation*}
where $\mathcal Z$ is the state box in \eqref{eq:C.4.6}. All points reconstructed from $t>0$ in \eqref{eq:C.4.6} have positive core weights and $\mathcal E_g>0$. At every such point in a slightly smaller tube, both successive one-block factors are nonzero at all eight nodes and have modulus greater than $2^{-31}$. \end{lemma}

\begin{proof}
For the $q$-derivatives, define $\overline Q(t)=\prod_{r\in R_Q}(t-r)$, $\mathscr U(t)=P(t)\overline Q(t)$, and $$ \Phi_\pm(t,q)=\mathscr U(t)(t-q)\pm C. $$ For any selected root branch $x(q)$, differentiation of $\Phi_\pm(x(q),q)=0$ gives, for $1\le k\le5$, \begin{equation*} \begin{aligned} x^{(k)}=-\frac{1}{\partial_t\Phi_\pm}\Bigg\{& \sum_{r=2}^k\partial_t^r\Phi_\pm B_{k,r}(x',\ldots,x^{(k-r+1)})\\ &{}-k\sum_{r=0}^{k-1}\mathscr U^{(r)} B_{k-1,r}(x',\ldots,x^{(k-r)})\Bigg\}. \end{aligned} \tag{C.4.9}\label{eq:C.4.9} \end{equation*} Here $B_{j,r}$ are the partial Bell polynomials, $B_{0,0}=1$, and an empty first sum is zero. Formula \eqref{eq:C.4.9} is the rearrangement of $$ \frac{d^k}{dq^k}\{\Phi_\pm(x(q),q)\} =\sum_{r=1}^k\partial_t^r\Phi_\pm B_{k,r} -k\sum_{r=0}^{k-1}\mathscr U^{(r)}B_{k-1,r}=0. $$ The root certificate gives $|\partial_t\Phi_\pm|>2^{-14}$ and bounds every displayed derivative of $\Phi,\mathscr U$ by $2^{44}$. Substitution in \eqref{eq:C.4.9} gives \begin{equation*} |x'|<2^{58},\quad |x''|<2^{175},\quad |x'''|<2^{294},\quad |x''''|<2^{420},\quad |x^{(5)}|<2^{550}. \tag{C.4.10}\label{eq:C.4.10} \end{equation*} For the last estimate the largest Bell monomial is $5x'x''''$, of binary exponent below $58+420+3=481$; after the $2^{44}$ coefficient bound, the fewer than $2^9$ total terms, and $1/|\partial_t\Phi_\pm|<2^{14}$, the exponent is below $548$.

Differentiating the elementary symmetric functions through order five now bounds every coefficient of $P,Q_q,\Delta,H$ and its $q$-derivatives by $2^{1024}$. Since $\Delta'(\xi_i)$ is a product of five gaps greater than one, five quotient differentiations in \eqref{eq:C.3.5} give \begin{equation*} |\partial_q^kw_i^*(a)|+ |\partial_q^k\partial_aw_i^*(a)|<2^{10000}, \quad0\le k\le5. \end{equation*} Together with $s_1>1/16$, the explicit affine formulas \eqref{eq:C.4.2} give \begin{equation*} \|D^k\chi_q\|+\|D^k\chi_q^{-1}\|<2^{10000(k+1)}, \quad0\le k\le5, \tag{C.4.12}\label{eq:C.4.12} \end{equation*} for joint $q$- and state derivatives.

The rational CG map requires a uniform bound that includes its normalisation denominator. Denote by $\Gamma$ the moment matrix in \eqref{eq:C.2.3}, namely the Gram matrix for $\sum_iw_i\lambda_i\delta_{\lambda_i}$. Four core nodes have weights greater than $2^{-19}$ and mutual gaps greater than one. Cauchy--Binet and the trace bound give
\begin{equation*} \det\Gamma>2^{-76},\qquad \operatorname{tr}\Gamma<2^{26}, \qquad \|\Gamma^{-1}\|_2<2^{159},\quad \|\Gamma^{-1}\|_\infty\le2\|\Gamma^{-1}\|_2<2^{160}. \end{equation*}
The signed external neighbourhood changes $\Gamma$ by less than $2^{-999972}$, so the Neumann lemma gives the weaker uniform inverse bound $2^{164}$. The coefficient equation $\Gamma\boldsymbol c_{\rm res}=-\boldsymbol m$ involves only polynomials in the independent weight variables and nodes. Holding those weights independent at this stage, \eqref{eq:C.4.10} bounds its joint weight/node/$q$ derivatives through order five by $2^{2000}$. Repeated differentiation of $\Gamma^{-1}\Gamma=I$ gives sums of words \begin{equation*} \Gamma^{-1}(D^{\gamma_1}\Gamma)\Gamma^{-1}\cdots (D^{\gamma_r}\Gamma)\Gamma^{-1},\qquad \gamma_1+\cdots+\gamma_r=\gamma, \end{equation*} with $r\le|\gamma|\le5$ and fewer than $6!\,2^5$ terms. All coefficient and evaluated-factor derivatives through order five are below $2^{20000}$.

Before differentiating \eqref{eq:C.2.4}, we bound its denominator. On either phase of the fixed-point family, the product of the two monic factors at every selected core or external node is $\pm C$. Each monic quartic has modulus below $2^{14}$ on the centred node interval and $C>1/4$; hence each monic factor has modulus above $2^{-16}$. Normalising at the shifted origin costs at most another factor $2^{14}$, so every normalised one-block factor on either phase of the fixed-point family has modulus above $2^{-30}$. As a core weight is greater than $2^{-19}$, the squared normalisation \begin{equation*} Z(w)=\sum_iw_ip_w(\lambda_i)^2 \end{equation*} satisfies $Z>2^{-79}$ on the fixed-point family. Differentiating the moment equations gives \begin{equation*} \|D\,p_w(\lambda_i)\|<2^{20000}. \end{equation*} Comparing a point in \eqref{eq:C.4.3} with the point of the fixed-point family having the same $(a,q)$ changes every first-phase factor by less than $2^{20000-1000000}<2^{-980000}$, and changes $Z$ by less than $2^{-900000}$. It follows that $Z>2^{-80}$. Its intermediate state is within $2^{-900000}$ of the $Q$-phase of the fixed-point family, where the identical estimate applies. Subtracting either factor variation from the corresponding fixed-point margin proves \begin{equation*} |p_w(\lambda_i)|>2^{-31},\qquad |p_{T_4w}(\lambda_i)|>2^{-31}, \quad1\le i\le8. \end{equation*}

The preceding factor bound gives $\|D^jZ\|<2^{40020}$ for $0\le j\le5$. Every term in $D^j(1/Z)$ is $Z^{-r-1}$ times a product of $r\le j$ positive-order derivatives of $Z$; hence it is below $2^{40020j+500}$, including fewer than $720$ terms. Combining this with the numerator $w_ip_w(\lambda_i)^2$, whose derivatives are below $2^{40020}$, gives the order-specific one-block bounds \begin{equation*} \|D^jT_4\|<2^{50000(j+1)},\qquad0\le j\le5. \tag{C.4.18}\label{eq:C.4.18} \end{equation*} In a fifth derivative of $T_4^2$, every chain-rule term has an outer derivative of order $r\le5$ and inner derivative orders summing to five. Equation \eqref{eq:C.4.18} bounds its binary exponent by $50000(5+2r+1)\le800000$. Both weight charts are affine in the state: higher derivatives in \eqref{eq:C.4.12} arise only when $q$-derivatives hit their coefficients. At total order five the inner-chart factors add at most $10000(5+r)\le100000$, the outer chart adds at most $60000$, and the Bell count adds less than $20000$. The sum is strictly below $1000000$, and lower orders are smaller. This gives \eqref{eq:C.4.4}, with the normalisation division already justified.

The first derivative of the return map used in the weighted rescaling is obtained as follows. At a point $F_q(a,0,0)=(a,0,0)$ of the fixed-point family, the tangent directions to the fixed set return exactly. For its one remaining intrinsic direction, start from \begin{equation*} SPQ=\Delta(SH)+CS,\qquad S=t-a. \end{equation*} Let $u=dU$ be a scalar intrinsic perturbation and divide \begin{equation*} \Delta u=SR_{\rm div}+\varrho,\qquad \deg R_{\rm div}\le5,\quad \varrho\in\mathbb R. \tag{C.4.20}\label{eq:C.4.20} \end{equation*} The certified $P,Q$ are coprime. Hence the map $$ (p_{\rm var},q_{\rm var})\longmapsto Qp_{\rm var}+Pq_{\rm var},\qquad \mathcal P_3\times\mathcal P_3\longrightarrow\mathcal P_7, $$ is injective and, by equal dimensions, bijective. Choose its unique pair with $Qp_{\rm var}+Pq_{\rm var}=R_{\rm div}$. Define $ (A_0,A_1,A_2)=(P,Q,P), $ $ (p_0,p_1,p_2)=(p_{\rm var},q_{\rm var},p_{\rm var}),\qquad \dot S_j=-j\varrho/C\quad(0\le j\le2). $ For $j=0,1$, direct substitution gives \begin{equation*} C\dot S_j+S(p_jA_{j+1}+A_jp_{j+1}) =-j\varrho+SR_{\rm div}=\Delta u-(j+1)\varrho =\Delta u+C\dot S_{j+1}. \tag{C.4.21}\label{eq:C.4.21} \end{equation*} This is exactly the coefficientwise differential of \begin{equation*} S_jA_jA_{j+1}=\Delta(S_jH+U_j)+CS_{j+1}. \tag{C.4.22}\label{eq:C.4.22} \end{equation*} Here $dS_j=\dot S_j$, $dU_j=u$, and $dC=dH=0$. In the notation of \eqref{eq:C.3.14c}, the phase-zero coefficient matrix is $\mathcal M(P,S,\Delta)$ and the phase-one matrix is $\mathcal M(Q_q,S_+,\Delta)$; on the fixed-point family $S_+=S$. These are exactly the two $9\times9$ matrices certified uniformly over the stated $(q,a)$ parameter box, with the same column and unknown ordering as \eqref{eq:C.3.14b}. Their invertibility shows that \eqref{eq:C.4.21} is the derivative of the return map. After two blocks $p_2=p_{\rm var}$, $dU_2=dU_0$, and $dC,dH$ also return, while the recovery root has only the displayed shear. Since $(a,C,H,U)$ are local coordinates by the invertible maps $J_{\rm lift}$ and $K_{\rm wt}$ proved above, and the directions tangent to the fixed-point family return exactly, the full intrinsic normal block is the identity.

Passing to the affine weight chart does not change that block. If to first order $\zeta=K_{\rm wt}(a)\mathbf n$, then $\mathbf n^+=\mathbf n$ and $K_{\rm wt}(a^+)\mathbf n^+=K_{\rm wt}(a)\mathbf n+O(|\mathbf n|^2)$, even though $a^+-a$ has a linear intrinsic shear. Hence \begin{equation*} D_\zeta(F_q)_\zeta(a,0,0)=I_4. \end{equation*} Coordinate divisibility of \eqref{eq:C.2.4} also gives \begin{equation*} (F_q)_{E_g}=E_gm_g,\qquad m_g(a,0,0)=1. \tag{C.4.24}\label{eq:C.4.24} \end{equation*} The two-block external multiplier is exactly $$ m_g(a,0,0)=\{P(g)Q_q(g)/C\}^2=(-C/C)^2=1. $$

For $t\ne0$, define the outgoing uncentred coordinates by \begin{equation*} V^+=(1+t)(F_q)_\zeta/t,\qquad \mathcal E_g^+=(1+t)^2(F_q)_{E_g}/t^2. \tag{C.4.25}\label{eq:C.4.25} \end{equation*} Hadamard division in $\zeta$ and the exact factor in \eqref{eq:C.4.24} remove both apparent singularities. Taylor's integral formula gives \eqref{eq:C.4.7}. Denote the $a$-component of $F_q$ by $(F_q)_a$ and its normal component by $(F_q)_\zeta$. Then, with every derivative on the right evaluated at $(a,0,0)$, \begin{equation*} \begin{aligned} G_a&=D_\zeta(F_q)_a[V],\\ G_V&=V+\tfrac12D_{\zeta\zeta}^2(F_q)_\zeta[V,V] +\sum_gD_{E_g}(F_q)_\zeta\,\mathcal E_g,\\ G_{\mathcal E_g}&=\mathcal E_g \{2+D_\zeta m_g[V]\}. \end{aligned} \end{equation*} Three joint $q,z$ derivatives of $G$ use at most five joint derivatives of $F_q$, which is why \eqref{eq:C.4.4} was proved through order five. The value of the integral remainder and each of its first $q,z$ derivatives contain at most three factors among $V,\mathcal E$. Only derivatives with respect to $q,z$ are needed. By \eqref{eq:C.4.2e}, the resulting exponent is below $1000000+3(500001)+100<3000000$. This proves the joint estimates \eqref{eq:C.4.8}, with the uncentred variables of \eqref{eq:C.4.5}. The $t$-cutoff places \eqref{eq:C.4.5} far inside \eqref{eq:C.4.3}, and $c_->2^{-15}>2^{-20}$ and $c_+>1/4$ prove $\mathcal E_g>0$, hence positivity. \end{proof}

Appendix C describes the certification of the finite margins used here.

\section{The Hopf point in weight coordinates}\label{exact-transfer-of-the-hopf-point-to-the-raw-field}

The expansion is in weight coordinates, whereas the Hopf spectrum was certified in the scalar chart. The following lemma identifies the two linearisations and retains the nonzero $a$-component.

\begin{lemma}[Tangent conjugacy in weight coordinates]\label{lemma-5.1-weighted-tangent-conjugacy}
At the critical branch, the scalar leading vector field and the leading vector field in weight coordinates are analytically conjugate. Their tangent map $L$ is invertible, \begin{equation*} J_*^{\rm wt}=LJ_*^{\rm sc}L^{-1},\qquad \|L\|_2<2^{5000}. \tag{C.5.1}\label{eq:C.5.1} \end{equation*} It maps $\delta a$ to the same $\delta a$. Hence the Hopf vector in weight coordinates can be normalised by $v_a=1$, and then \begin{equation*} \|v\|_2<2^{30000}. \tag{C.5.2}\label{eq:C.5.2} \end{equation*} \end{lemma}

\begin{proof}
For six core nodes, an incoming monic quartic $A$, monic affine recovery factor $S$, and scalar squared monic norm $h$, orthogonality is equivalent to \begin{equation*} w_i=\frac{hS(\xi_i)}{A(\xi_i)\Delta'(\xi_i)}. \tag{C.5.3}\label{eq:C.5.3} \end{equation*} Indeed, multiplying by $A(\xi_i)\xi_i^j$ and summing reduces the assertion to the degree-five Lagrange coefficient identity for $S(t)t^j$, $0\le j\le3$. The normalisation is given explicitly by \begin{equation*} h^{-1}=\sum_{i=1}^6\frac{S(\xi_i)}{A(\xi_i)\Delta'(\xi_i)}. \end{equation*} The weights in \eqref{eq:C.5.3} sum to one, and uniqueness of the positive moment solution proves the asserted equivalence.

For an incoming perturbation $P+\epsilon p$, coefficient comparison in \eqref{eq:C.3.14} gives a $9\times9$ system. Extraction of the incoming-phase coordinates $(\delta C,\delta H_0,\delta H_1,\delta U)$ by \eqref{eq:C.3.14d} defines a $4\times4$ lift matrix $J_{\rm lift}$, and the exact certificate proves \begin{equation*} \|J_{\rm lift}^{-1}\|_\infty<2^{20}. \end{equation*} Hence, in the scalar order $\mathbf n=(u,B_C,B_0,B_1)$, \begin{equation*} p=J_{\rm lift}^{-1}(B_C,B_0,B_1,u)^T. \tag{C.5.5}\label{eq:C.5.5} \end{equation*} Differentiating \eqref{eq:C.5.3} at fixed $a$ gives \begin{equation*} \delta w_i=w_i^*(a) \left\{\frac{B_C}{C}-\frac{p(\xi_i)}{P(\xi_i)}\right\}. \end{equation*} If $s_i=\partial_aw_i^*(a)$, the normal in weight coordinates is \begin{equation*} (K_{\rm wt}(a)\mathbf n)_j=\delta w_{j+1}-\frac{s_{j+1}}{s_1}\delta w_1, \qquad1\le j\le4. \tag{C.5.7}\label{eq:C.5.7} \end{equation*}

If $K_{\rm wt}(a)\mathbf n=0$, normalisation implies $\delta w$ is a multiple of the tangent vector $(s_i)_{i=1}^6$ to the fixed-point family. Explicitly, $\sum_i\delta w_i=\sum_is_i=0$; hence the relations for components $2,\ldots,5$ in \eqref{eq:C.5.7} also force the same relation for the sixth core component. Differentiating the four moment equations then gives \begin{equation*} \sum_{\ell=0}^3 \left(\sum_iw_i^*\xi_i^{j+\ell}\right)p_\ell=0, \qquad0\le j\le3. \tag{C.5.8}\label{eq:C.5.8} \end{equation*} The matrix in \eqref{eq:C.5.8} is positive definite because a nonzero cubic cannot vanish at six distinct positive-weight nodes. It follows that $p=0$, and \eqref{eq:C.5.5} gives $\mathbf n=0$. Hence $K_{\rm wt}(a)$ is invertible.

To relate the field \eqref{eq:C.3.12} to the exact two-block map, consider a nearby eight-node weight, retain the two external masses $E=(E_-,E_+)$ and normalise the six-node core to have mass one. Use the incoming-phase convention of \eqref{eq:C.3.14d}: solve for the partner, attach its norm coordinate to the input state, and divide the quotient by the input recovery factor. Equations \eqref{eq:C.5.3}--\eqref{eq:C.5.5} and the inverse-function theorem then give analytic scalar coordinates $$ (a,U,B,E),\qquad B=(\delta C,\delta H_0,\delta H_1). $$ Let $\mathscr S_q$ be the exact two-block map in these coordinates. Assign weight one to $U,B$ and weight two to each $E_g$; the notation $O_{\rm w}(j)$ denotes an analytic remainder whose monomials have weighted degree at least $j$. Its weighted Taylor expansion at the fixed-point set is \begin{equation*} \begin{aligned} a^+-a&=\kappa(a)U+O_{\rm w}(2),\\ U^+-U&=\beta(a,q)U^2+\sum_gA_g(a,q)E_g+O_{\rm w}(3),\\ B^+-B&=\mathbf q(a,q)U^2+\sum_gV_g(a,q)E_g+O_{\rm w}(3),\\ \log m_g&=\ell_g(a,q)U+d_g(q)[B]+O_{\rm w}(2), \end{aligned} \tag{C.5.8a}\label{eq:C.5.8a} \end{equation*} where $E_g^+=E_gm_g$ and $m_g$ is the analytic external multiplier. The remainders are uniform in the $C^1$ norm on the compact certified parameter box.

We justify each coefficient in \eqref{eq:C.5.8a}. The intrinsic calculation \eqref{eq:C.4.20}--\eqref{eq:C.4.22} gives $\kappa=2\Delta(a)/C$. The normal derivative of the return map is the identity at every point of the fixed-point family, so differentiation in a direction tangent to that family rules out terms $UB$ in the two middle equations. Their remaining quadratic intrinsic coefficients are exactly those obtained from the three-phase recursion \eqref{eq:C.3.14a}--\eqref{eq:C.3.14e}. Equations \eqref{eq:C.3.9}--\eqref{eq:C.3.9a} give the external coefficients. Finally, evaluation of the differentiated factor identity at $g$ gives $$ \delta(PQ)(g)=\frac{\Delta(g)-\Delta(a)}{g-a}U \quad\hbox{and}\quad \delta(PQ)(g)=\Delta(g)\delta H(g)+\delta C $$ in the intrinsic direction and in directions tangent to the fixed-point family, respectively. Differentiating $\log\{(P(g)Q_q(g)/C)^2\}$ at its fixed value yields precisely $\ell_gU+d_g[B]$ from \eqref{eq:C.3.6} and \eqref{eq:C.3.10}, which establishes \eqref{eq:C.5.8a}.

For $t$, let $$ \Sigma_t(a,u,B,e)=(a,tu,tB,t^2e),\qquad t^+=\frac{t}{1+t}, $$ and let the exact scalar rescaled map be $$ \widehat{\mathscr S}_q(t,\cdot) =\Sigma_{t^+}^{-1}\circ\mathscr S_q\circ\Sigma_t. $$ Equation \eqref{eq:C.5.8a} gives \begin{equation*} \widehat{\mathscr S}_q(t,z) =z+tG_q(z)+t^2R_q^{\rm sc}(t,z), \qquad \|R_q^{\rm sc}\|_{C^1}\le C_{\rm rem}. \tag{C.5.8b}\label{eq:C.5.8b} \end{equation*} Hence \eqref{eq:C.3.12} is the leading vector field of the scalar map.

Let $\Lambda_t(a,V,e)=(a,tV,t^2e)$ be the weight-coordinate scaling and let $\mathfrak C_q$ be the unscaled scalar-to-weight coordinate change. The map $$ \mathfrak C_{q,t}=\Lambda_t^{-1}\circ\mathfrak C_q\circ\Sigma_t $$ extends analytically to $t=0$, since $a_{\rm wt}=a+t\delta w_1/s_1+O(t^2)$, and \begin{equation*} \mathfrak C_{q,0}(a,u,B,e) =\bigl(a,K_{\rm wt}(a)(u,B),e\bigr). \end{equation*} For $t\ne0$, coordinate covariance gives the exact identity \begin{equation*} \widehat F_q(t,\mathfrak C_{q,t}z) =\mathfrak C_{q,t^+}\bigl(\widehat{\mathscr S}_q(t,z)\bigr). \tag{C.5.10}\label{eq:C.5.10} \end{equation*} Both sides extend analytically to $t=0$. Comparing their coefficients of $t$ gives \begin{equation*} G_q^{\rm wt}(\mathfrak C_{q,0}z) =D\mathfrak C_{q,0}(z)G_q(z). \tag{C.5.10a}\label{eq:C.5.10a} \end{equation*} At the critical equilibrium its differential is \begin{equation*} L(\delta a,\delta\mathbf n,\delta e) =\bigl(\delta a,K_{\rm wt}(\alpha(q_*))\delta\mathbf n +\partial_aK_{\rm wt}(\alpha(q_*))\mathbf n_*\delta a,\delta e\bigr). \tag{C.5.10b}\label{eq:C.5.10b} \end{equation*} Its diagonal blocks are $1,K_{\rm wt}(\alpha(q_*)),I_2$, so $L$ is invertible. Equation \eqref{eq:C.5.10a} gives the similarity \eqref{eq:C.5.1}; \eqref{eq:C.5.10b} also exhibits the shear that would be lost by treating the uncentred equilibrium $\mathbf n_*$ as zero.

Direct estimates in \eqref{eq:C.5.3}--\eqref{eq:C.5.7} and \eqref{eq:C.5.10b}, including the shear $(\partial_aK_{\rm wt})\mathbf n_*$, are certified in Appendix C and give $\|L\|_2<2^{5000}$.

For the quantitative eigenvector estimate, normalise the Hopf vector in \eqref{eq:C.3.17b} by $a=1$, take $\lambda=i\omega_*$, $\theta=\lambda(\lambda-1)$, and use its first equation to get $u=\lambda/k$. Eliminating $B$ from the last two equations gives \begin{equation*} (\theta I_2-D\mathsf V_c)\varphi=(\lambda-1)\ell u+Dq_{\rm lin},\qquad B=(\lambda-1)^{-1}(q_{\rm lin}+\mathsf V_c\varphi). \tag{C.5.11}\label{eq:C.5.11} \end{equation*} The real matrix $D\mathsf V_c$ has eigenvalues $2$ and $\operatorname{tr}(D\mathsf V_c)-2$, by \eqref{eq:C.3.17c}. The certified Hopf value is $\theta=-\omega_*^2-i\omega_*$, so for either real eigenvalue $\rho$, \begin{equation*} |\theta-\rho|\ge|\operatorname{Im}\theta|=\omega_*, \qquad |\det(\theta I_2-D\mathsf V_c)|\ge\omega_*^2>1. \tag{C.5.12}\label{eq:C.5.12} \end{equation*} The interval certificate reconstructs every block in \eqref{eq:C.5.11} from \eqref{eq:C.3.6}--\eqref{eq:C.3.10} and proves \begin{equation*} \begin{gathered} |k|>2^6,\quad 3.72<\omega_*<3.724,\quad |A_g|,|A_g'|<2^{10},\quad |V_g|,|V_g'|<2^{17},\\ |\ell_g|<2^7,\quad |D|_{\max}<2^4,\quad \|\mathsf V_c\|_\infty<2^3,\quad \|q_{\rm lin}\|_\infty<2^2. \end{gathered} \tag{C.5.13}\label{eq:C.5.13} \end{equation*} These are closed rational interval comparisons; see Appendix C. Since $D$ has three columns, its action satisfies \begin{equation*} \|DB\|_\infty\le3|D|_{\max}\|B\|_\infty <2^6\|B\|_\infty. \end{equation*} We use this dimension factor throughout the remaining estimates.

The bound for $V_g$ is obtained by using $(V_g)_C=4(P(g)^2+C)$ and the two equations $d_{h_\pm}[V_g]=\mathsf T_{h_\pm,g}$ to solve for its affine $H$-part; differentiating the same identities gives the bound for $V_g'$.

Now $|u|<1$. Every entry of $\theta I_2-D\mathsf V_c$ is below $2^{10}$; its two-by-two adjugate has row sum below $2^{11}$, so \eqref{eq:C.5.12} gives $$ \|(\theta I_2-D\mathsf V_c)^{-1}\|_\infty<2^{11}. $$ The right side of the first equation in \eqref{eq:C.5.11} is below $2^{11}$. Hence $\|\varphi\|_\infty<2^{22}$, and the second equation gives $\|B\|_\infty<2^{27}$. With seven complex components, \begin{equation*} \|v_{\rm scalar}\|_2<2^{29},\qquad \|Lv_{\rm scalar}\|_2<2^{5029}<2^{30000}. \tag{C.5.14}\label{eq:C.5.14} \end{equation*} This is \eqref{eq:C.5.2}. \end{proof}

Lemmas~\ref{lemma-3.1-certified-critical-branch-and-transverse-hopf-point} and~\ref{lemma-5.1-weighted-tangent-conjugacy} now give a simple transverse Hopf point for the leading vector field in weight coordinates, and its Hopf eigenvector has nonzero $a$-component.

\section{Construction of a periodic orbit}\label{literal-selection-of-one-periodic-orbit}

Continue in the uncentred weight variables $z=(a,V,\mathcal E)$ of \eqref{eq:C.4.5}. Their equilibrium curve is \begin{equation*} z_*^{\rm wt}(q)=(\alpha(q),v_*(q),c_-(q),c_+(q)). \end{equation*} We translate this equilibrium curve to zero only in the Lyapunov--Schmidt equation. Let \begin{equation*} f_\mu(x)=G^{\rm wt}_{q_*+\mu} \bigl(z_*^{\rm wt}(q_*+\mu)+x\bigr). \end{equation*} Set $J^{\rm wt}(q)=D_zG_q^{\rm wt}(z_*^{\rm wt}(q))$; at the critical parameter, $J_*^{\rm wt}=J^{\rm wt}(q_*)=LJ_*^{\rm sc}L^{-1}$ by \eqref{eq:C.5.1}. Choose the $a$-normalised eigenvector $v$ from Lemma~\ref{lemma-5.1-weighted-tangent-conjugacy} and define \begin{equation*} \psi(\theta)=ve^{i\theta}+\overline v e^{-i\theta}. \end{equation*} Let $X=H^1_{\rm per}([0,2\pi],\mathbb R^7)$ and $Y=L^2_{\rm per}$. The complex equation $(\widehat h_1)_a=0$ means the two real conditions $$ \operatorname{Re}(\widehat h_1)_a=0, \qquad \operatorname{Im}(\widehat h_1)_a=0; $$ let $X_\perp$ denote this closed real codimension-two subspace of $X$. For $\mathbf u=(h,\mu,\nu)\in X_\perp\times\mathbb R^2$, define \begin{equation*} \mathcal H(\mathbf u;\varepsilon)= (\omega_*+\nu)(\psi+h)' -\int_0^1Df_\mu(s\varepsilon(\psi+h))(\psi+h)\,ds. \tag{C.6.3}\label{eq:C.6.3} \end{equation*} For $\varepsilon\ne0$, its zeros are the periodic solutions satisfying the phase and amplitude conditions, and \eqref{eq:C.6.3} supplies the removable value at $\varepsilon=0$.

\begin{lemma}[Quantitative Hopf contraction]\label{lemma-6.1-quantitative-hopf-contraction}
The linearisation $D_0=D_{\mathbf u}\mathcal H(0;0):X_\perp\times\mathbb R^2\to Y$ is invertible and \begin{equation*} \|D_0^{-1}\|<M_H=2^{70000000}. \tag{C.6.4}\label{eq:C.6.4} \end{equation*} On the explicitly fixed product box $$ \mathcal B_H=\{(\mathbf u,\varepsilon):\|\mathbf u\|\le2^{-350000000},\ |\varepsilon|\le\varepsilon_0\}, $$ one has \begin{equation*} \|\partial_\varepsilon\mathcal H(0;\varepsilon)\|<B_H,\qquad \|D_{\mathbf u}^2\mathcal H\|+\|\partial_\varepsilon D_{\mathbf u}\mathcal H\|<L_H, \quad B_H=L_H=2^{70000000}. \tag{C.6.5}\label{eq:C.6.5} \end{equation*} For $\varepsilon_0$ in \eqref{eq:C.1.1}, define \begin{equation*} R_0=2M_HB_H\varepsilon_0=2^{-359999999}. \tag{C.6.6}\label{eq:C.6.6} \end{equation*} Then $\mathcal H(\mathbf u;\varepsilon_0)=0$ has exactly one solution in the closed $R_0$-ball, \begin{equation*} \mathbf u_0=(h_0,\mu_0,\nu_0),\qquad\|\mathbf u_0\|\le R_0. \tag{C.6.7}\label{eq:C.6.7} \end{equation*} \end{lemma}

\begin{proof}
On Fourier mode $k$, the $h$-operator is $ik\omega_*I-J_*^{\rm wt}$. By \eqref{eq:C.3.15}, it is invertible except at $k=\pm1$. Those two complex-conjugate modes form one real two-dimensional kernel. The two real gauge equations defining $X_\perp$ remove the corresponding kernel in the domain. Choose an adjoint Hopf vector $v^\dagger$ with $\langle v^\dagger,v\rangle=1$. Pairing the $k=1$ equation, and its conjugate, with $v^\dagger$ identifies the remaining real cokernel with the parameter and frequency columns. In the ordered real basis their matrix is \begin{equation*} \begin{pmatrix} \partial_q\operatorname{Re}\lambda_H&0\\ \partial_q\operatorname{Im}\lambda_H&-1 \end{pmatrix}, \tag{C.6.8}\label{eq:C.6.8} \end{equation*} which is invertible by \eqref{eq:C.3.16}. This also fixes the adjoint normalisation used in the quantitative inverse.

From \eqref{eq:C.4.8}, every spectral projector $\Pi_\lambda$ of $J_*^{\rm wt}$ is bounded by \begin{equation*} P_0=2^{18100000}, \end{equation*} because its product formula has six numerator factors below $2^{3000001}$ and eigenvalue gaps at least one. For $|k|\ge2$, \eqref{eq:C.3.15} and $\omega_*>3.72$ give $$ \min_{\lambda\in\operatorname{spec}J_*^{\rm wt}}|ik\omega_*-\lambda| \ge1+|k|. $$ The spectral resolution yields the high-mode estimate \begin{equation*} \|(ik\omega_*I-J_*^{\rm wt})^{-1}\|_2 \le\frac{7P_0}{1+|k|}, \qquad |k|\ge2. \end{equation*} Since $(1+k^2)^{1/2}\le1+|k|$, Parseval's identity shows that the inverse on all high modes maps $L^2$ to $H^1$ with norm at most $7P_0$. The $k=0$ inverse and the reduced $k=\pm1$ inverses obey the same coarser bound, because all remaining spectral gaps are at least one.

The equilibrium derivatives satisfy \begin{equation*} \|(z_*^{\rm wt})'\|<2^{21200000},\qquad \|(z_*^{\rm wt})''\|<2^{64000000}. \tag{C.6.10}\label{eq:C.6.10} \end{equation*} These follow by differentiating $G_q^{\rm wt}(z_*^{\rm wt}(q))=0$ and expanding $(J_*^{\rm wt})^{-1}$ in the seven projectors. It follows that \begin{equation*} \|(J^{\rm wt})'(q_*)\|<2^{24200001}. \end{equation*} Pairing with the fixed adjoint normalisation, solving \eqref{eq:C.6.8}, applying the low- and high-mode reduced resolvents, and imposing the $a$-component normalisation is bounded by \begin{equation*} P_0^2\,2^{3000000}(1+2^{21200000})\,2^{400000} <2^{60800001}. \tag{C.6.12}\label{eq:C.6.12} \end{equation*} Estimate \eqref{eq:C.6.12} proves \eqref{eq:C.6.4}.

Two derivatives of \eqref{eq:C.6.3} use the $C^3$ bound in \eqref{eq:C.4.8}. Sobolev multiplication and embedding, \eqref{eq:C.5.2}, and \eqref{eq:C.6.10} give the upper exponent \begin{equation*} 20+3000000+64000000+3(200001) =67600023<70000000. \end{equation*} The remaining estimates contribute at most three further binary powers, which gives \eqref{eq:C.6.5}.

For $\mathcal T(\mathbf u)=\mathbf u-D_0^{-1}\mathcal H(\mathbf u;\varepsilon_0)$, the centre displacement is at most $M_HB_H\varepsilon_0=R_0/2$, while \begin{equation*} M_HL_H(2M_HB_H+1)\varepsilon_0 <2^{-219999998}<1/4. \tag{C.6.14}\label{eq:C.6.14} \end{equation*} Since $R_0<2^{-350000000}$, every segment used in the mean-value estimates lies in $\mathcal B_H$. Accordingly, $\mathcal T$ is a strict contraction of the closed $R_0$-ball into itself. Banach's theorem proves \eqref{eq:C.6.7}, with uniqueness in that ball. \end{proof}

Define \begin{equation*} \begin{aligned} q_{\rm per}&=q_*+\mu_0,\qquad \Omega_0=\omega_*+\nu_0,\\ \Gamma_0(\theta)&=z_*^{\rm wt}(q_{\rm per}) +\varepsilon_0\{\psi(\theta)+h_0(\theta)\}. \end{aligned} \tag{C.6.15}\label{eq:C.6.15} \end{equation*} Then \begin{equation*} \Omega_0\Gamma_0'=G^{\rm wt}_{q_{\rm per}}(\Gamma_0), \qquad \operatorname{Re}(\widehat{\Gamma_0}_1)_a=\varepsilon_0,\quad \operatorname{Im}(\widehat{\Gamma_0}_1)_a=0. \end{equation*} Hence $\Gamma_0$ is nonconstant. Its state displacement from the certified critical point is less than $2^{-330000000}$, so it lies in the positive chart and remains in \eqref{eq:C.4.6}. Moreover \begin{equation*} |\mu_0|\le R_0<2^{-41}, \end{equation*} while Lemma~\ref{lemma-3.1-certified-critical-branch-and-transverse-hopf-point} places $q_*$ more than $2^{-41}$ from both parameter faces. Hence $q_{\rm per}\in Q_H$. The rational ball and conditions in \eqref{eq:C.6.6}--\eqref{eq:C.6.7} uniquely determine the periodic profile and its parameter, frequency, and phase.

\section{Shadowing by an orbit of the exact map}\label{literal-selection-of-one-exact-map-shadow}

An orbit of the squared-weight return map shadowing the periodic curve $\Gamma_0$ will yield a restarted CG counterexample.

Let \begin{equation*} \gamma_0(\tau)=\Gamma_0(\Omega_0\tau),\qquad T_{\rm per}=\frac{2\pi}{\Omega_0}. \end{equation*} Then $\dot\gamma_0=G^{\rm wt}_{q_{\rm per}}(\gamma_0)$ and $3/2<T_{\rm per}<2$.

Let $t_n=1/n$. Since $t_{n+1}=t_n/(1+t_n)$, the scaling law \eqref{eq:C.4.5} turns iteration of the fixed map $T_4^2$ into the nonautonomous coordinate recurrence \begin{equation*} z_{n+1}=\widehat F_{q_{\rm per}}(1/n,z_n). \tag{C.7.1a}\label{eq:C.7.1a} \end{equation*}

\begin{lemma}[Quantitative nonhyperbolic shadowing]\label{lemma-7.1-quantitative-nonhyperbolic-shadowing}
There is an exact forward solution of \eqref{eq:C.7.1a} such that, for \begin{equation*} \tau_n=\sum_{k=N}^{n-1}\frac{1}{k}, \tag{C.7.2}\label{eq:C.7.2} \end{equation*} one has \begin{equation*} z_n-\gamma_0(\tau_n)\longrightarrow0. \tag{C.7.3}\label{eq:C.7.3} \end{equation*} With the Floquet error $x$ and projections defined in the proof, this orbit is the unique one in the closed Banach ball $$ \sup_{n\ge N}n^{1/2}|x_n|\le1 $$ that satisfies the zero initial stable condition $\Pi_{\rm st}x_N=0$ and the boundary condition at infinity in \eqref{eq:C.7.13}. Every reconstructed weight is positive. \end{lemma}

\begin{proof}
The variational coefficient along $\gamma_0$ differs from $J_*^{\rm wt}$ by less than \begin{equation*} \epsilon_{\rm var}=2^{-300000000}. \end{equation*} Let $\Pi_{\rm st}^*$ be the spectral projector of $J_*^{\rm wt}$ for $-1$, let $\Pi_{\rm cu}^*=I-\Pi_{\rm st}^*$, and take $C_0=2^{18100003}$. The projector formula gives \begin{equation*} \|e^{\tau J_*^{\rm wt}}\Pi_{\rm st}^*\|\le C_0e^{-\tau}\quad(\tau\ge0), \qquad \|e^{\tau J_*^{\rm wt}}\Pi_{\rm cu}^*\|\le C_0\quad(\tau\le0). \end{equation*} For $\mathsf A(\tau)=DG_{q_{\rm per}}^{\rm wt}(\gamma_0(\tau))$ and $\mathsf E=\mathsf A-J_*^{\rm wt}$, the stable and complementary whole-line integral operators have norm at most \begin{equation*} C_0\epsilon_{\rm var}(8+8/7)<10C_0\epsilon_{\rm var} <2^{-281899990}<1/4. \end{equation*} Their contractions define periodic invariant stable and complementary bundles, with forward rate $7/8$ and backward rate $1/8$. This separates only the continuation of $-1$ from its entire complement.

Let $\mathsf X(0)=I$ for the fundamental variational matrix. Variation of constants on $0\le\tau\le T_{\rm per}<2$ gives \begin{equation*} \|\mathsf X(T_{\rm per})-e^{T_{\rm per}J_*^{\rm wt}}\| <4(2^{18100009})^2\epsilon_{\rm var} <2^{-263799978}. \tag{C.7.4c}\label{eq:C.7.4c} \end{equation*} Let $\mathscr C_-,\mathscr C_0,\mathscr C_1,\mathscr C_2$ be the positively oriented circles of radii $1/64,1/16,1/16,1/16$ about $e^{-T_{\rm per}},1,e^{T_{\rm per}},e^{2T_{\rm per}}$, respectively. They are disjoint and contained in the open right half-plane. Including the already bounded frequency displacement, the comparison monodromy has algebraic multiplicities $$ 1,\quad2,\quad3,\quad1 $$ inside these four contours. They correspond, respectively, to the exponents $-1$, $\pm i\omega_*$, $1$ together with $1\pm i\omega_*$, and $2$. On the smaller stable circle, spectral resolution bounds the comparison resolvent by $7P_0/(1/64)<2^{18100009}$. The other three circles obey this same bound, including the frequency displacement; \eqref{eq:C.7.4c} makes the resolvent Neumann correction below $2^{-245000000}$. Riesz projection ranks are unchanged, so the monodromy has the same four multiplicities and no spectrum outside their union.

The union of the four disks avoids the closed negative real axis. Hence the principal matrix logarithm exists. Since $\mathsf X(T_{\rm per})$ is real and the principal scalar logarithm is invariant under conjugation on this domain, its holomorphic functional calculus is real. Define \begin{equation*} B_{\rm Fl}=T_{\rm per}^{-1}\operatorname{Log}\mathsf X(T_{\rm per}),\qquad \mathsf P(\tau)=\mathsf X(\tau)e^{-\tau B_{\rm Fl}}. \end{equation*} The definitions give real $B_{\rm Fl}$ and $\mathsf P$, with $\mathsf P(\tau+T_{\rm per})=\mathsf P(\tau)$; Dunford's formula on the same contours gives \begin{equation*} \|B_{\rm Fl}\|<2^{18200000}. \end{equation*} On one period, $\|\mathsf X\|+\|\mathsf X^{-1}\|<2^{18100010}$ and $\|e^{\pm\tau B_{\rm Fl}}\|<2^{18200010}$. Hence $\mathsf P,\mathsf P^{-1}$ have exponent below $36300021$. The identities \begin{equation*} \mathsf P'=\mathsf A\mathsf P-\mathsf PB_{\rm Fl},\qquad \mathsf P''=\mathsf A'\mathsf P+\mathsf A\mathsf P'-\mathsf P'B_{\rm Fl} \end{equation*} bound their first two derivative exponents by $54500030$ and $72700050$, respectively. We have the real Floquet representation \begin{equation*} \mathsf X(\tau)=\mathsf P(\tau)e^{\tau B_{\rm Fl}} \end{equation*} and \begin{equation*} \|\mathsf P\|+\|\mathsf P^{-1}\|+\|\mathsf P'\|+\|\mathsf P''\| <K_{\rm Fl}=2^{80000000}. \end{equation*}

Let $\Pi_{\rm st}$ be the Riesz projection of $B_{\rm Fl}$ corresponding to $\mathscr C_-$ under the logarithm, and define $\Pi_{\rm cu}=I-\Pi_{\rm st}$. For $3/2<T_{\rm per}<2$, the stable contour satisfies
\begin{equation*}
|z|\le e^{-T_{\rm per}}+1/64<e^{-7T_{\rm per}/8},
\end{equation*}
whereas the complementary contours satisfy $|z|\ge15/16>e^{-T_{\rm per}/8}$. The resolvent Neumann estimate bounds the exact resolvent by $2^{18100010}$. Dunford's formula applied to $z^{\tau/T_{\rm per}}$ therefore gives forward rate $7/8$ on $\Pi_{\rm st}$ and backward rate $1/8$ on $\Pi_{\rm cu}$. Multiplying the resolvent bound by the contour radii gives prefactors at most $2^{18100004}$ and $3\cdot2^{18100006}$, respectively, both below $2^{20000000}$. For $j\ge N$, let $\mathsf A_j=I+B_{\rm Fl}/j$; the evolution operators used below are \begin{equation*} \begin{aligned} \Phi_{\rm st}(n,k+1)&=\mathsf A_{n-1}\cdots \mathsf A_{k+1}\Pi_{\rm st}, &&N\le k<n,\\ \Phi_{\rm cu}(n,k+1)&=\mathsf A_n^{-1}\cdots \mathsf A_k^{-1}\Pi_{\rm cu}, &&k\ge n. \end{aligned} \end{equation*} An empty stable product is the identity on $\Pi_{\rm st}$. For $j\ge N$, $\|B_{\rm Fl}\|/j<1/2$, every $\mathsf A_j$ is invertible, and all these factors commute. The matrix-logarithm series gives, for any forward product, \begin{equation*} \begin{aligned} \prod_{j=m}^{n-1}(I+B_{\rm Fl}/j) &=\exp\left\{B_{\rm Fl}\sum_{j=m}^{n-1}j^{-1}+\mathsf R_{n,m}\right\},\\ \|\mathsf R_{n,m}\| &\le2\|B_{\rm Fl}\|^2\sum_{j=m}^\infty j^{-2}<1. \end{aligned} \end{equation*} The same estimate holds for inverse complementary products. On a generalised eigenspace, a Jordan block of size $d\le7$ contributes at most $(1+s)^{d-1}$ at logarithmic time $s$. The elementary inequality $(1+s)^6\le C e^{s/8}$ absorbs every such factor. Decreasing the stable exponent from $7/8$ to $3/4$ and increasing the complementary backward exponent from $1/8$ to $1/4$, then comparing the harmonic sum with $\log(n/(k+1))$, gives \begin{equation*} \begin{aligned} \|\Phi_{\rm st}(n,k+1)\| &\le K_{\rm prod}\left(\frac{k+1}{n}\right)^{3/4}, &&N\le k<n,\\ \|\Phi_{\rm cu}(n,k+1)\| &\le K_{\rm prod}\left(\frac{k+1}{n}\right)^{1/4}, &&k\ge n,\\ K_{\rm prod}&=2^{100000000}. \end{aligned} \tag{C.7.7}\label{eq:C.7.7} \end{equation*}

To derive the nonlinear recurrence, define \begin{equation*} \bar z_n=\gamma_0(\tau_n),\qquad z_n=\bar z_n+\mathsf P(\tau_n)x_n. \end{equation*} Taylor's integral formula gives \begin{equation*} G_{q_{\rm per}}^{\rm wt}(\bar z_n+e)-G_{q_{\rm per}}^{\rm wt}(\bar z_n)=\mathsf A(\tau_n)e+\mathcal Q_n(e), \end{equation*} where, with $D_{\rm quad}=2^{3000000}$, \begin{equation*} |\mathcal Q_n(e)|\le D_{\rm quad}|e|^2,\qquad |\mathcal Q_n(e)-\mathcal Q_n(f)| \le D_{\rm quad}(|e|+|f|)|e-f|. \end{equation*} Apply \eqref{eq:C.4.7} at $t=t_n$, expand $\gamma_0(\tau_n+t_n)$ and $\mathsf P(\tau_n+t_n)$ through second order, and use $\mathsf P'=\mathsf A\mathsf P-\mathsf PB_{\rm Fl}$. The order-$t_n$ linear expression is $I+t_nB_{\rm Fl}$. The transformed $\mathcal Q_n$ is $t_n$ times a quadratic function of $x$; every other remaining value and first state derivative has a factor $t_n^2$. The map in Floquet error coordinates takes the form \begin{equation*} x_{n+1}=(I+B_{\rm Fl}/n)x_n+n^{-1}\mathcal N_n(x_n) +n^{-2}\mathcal R_n(x_n), \tag{C.7.8}\label{eq:C.7.8} \end{equation*} on the radius \begin{equation*} \delta_x=2^{-81000000}, \end{equation*} where \begin{equation*} \begin{aligned} |\mathcal N_n(x)|&\le L_{\rm nl}|x|^2,& |\mathcal N_n(x)-\mathcal N_n(y)| &\le L_{\rm nl}(|x|+|y|)|x-y|,\\ |\mathcal R_n(x)|&\le R_{\rm rem},& |\mathcal R_n(x)-\mathcal R_n(y)|&\le R_{\rm rem}|x-y|, \end{aligned} \tag{C.7.10}\label{eq:C.7.10} \end{equation*} with \begin{equation*} L_{\rm nl}=2^{250000000},\qquad R_{\rm rem}=2^{410000000}. \end{equation*} One inverse output factor and two input factors of $\mathsf P$ give \begin{equation*} 2^{20}K_{\rm Fl}^3D_{\rm quad}<2^{250000000}=L_{\rm nl}. \end{equation*} The exact-map defect and the second-order time/Floquet terms have at most one inverse output factor and four factors among $\mathsf P,\mathsf P',\mathsf P''$, so \begin{equation*} 2^{20}K_{\rm Fl}^5D_{\rm quad}<2^{410000000}=R_{\rm rem}. \end{equation*} These estimates hold for both the values and first state derivatives and prove \eqref{eq:C.7.10}.

On the Banach ball \begin{equation*} \sup_{n\ge N}n^{1/2}|x_n|\le1, \tag{C.7.12}\label{eq:C.7.12} \end{equation*} define \begin{equation*} \begin{aligned} (\mathcal Lx)_n^{\rm st} &=\sum_{k=N}^{n-1}\Phi_{\rm st}(n,k+1)f_k(x_k),\\ (\mathcal Lx)_n^{\rm cu} &=-\sum_{k=n}^\infty\Phi_{\rm cu}(n,k+1)f_k(x_k), \end{aligned} \tag{C.7.13}\label{eq:C.7.13} \end{equation*} where $f_k=k^{-1}\mathcal N_k+k^{-2}\mathcal R_k$. Integral comparison in \eqref{eq:C.7.7}, with exponents $3/4$ and $1/4$, gives the common convolution constant $44K_{\rm prod}/3$. The four sufficient inequalities are \begin{equation*} \begin{aligned} M&>\delta_x^{-1},\\ M&>\frac{44K_{\rm prod}}{3}(L_{\rm nl}+R_{\rm rem}),\\ M&>\frac{352K_{\rm prod}}{3}L_{\rm nl},\\ M^2&>\frac{176K_{\rm prod}}{3}R_{\rm rem}. \end{aligned} \tag{C.7.14}\label{eq:C.7.14} \end{equation*} For the values of $M,N$ in \eqref{eq:C.1.1}, their right-hand binary exponents are respectively below \begin{equation*} 81000000,\quad510000005,\quad350000009,\quad510000008, \tag{C.7.15}\label{eq:C.7.15} \end{equation*} whereas the left exponents are $600000000,600000000,600000000,1200000000$. Hence \eqref{eq:C.7.13} is a contraction of \eqref{eq:C.7.12} into itself. For $$ \|x\|_*=\sup_{n\ge N}n^{1/2}|x_n|, \qquad C_{\rm conv}=44K_{\rm prod}/3, $$ one has \begin{equation*} \|\mathcal Lx\|_* \le\frac{C_{\rm conv}(L_{\rm nl}+R_{\rm rem})}{M}<1, \end{equation*} and \begin{equation*} \|\mathcal Lx-\mathcal Ly\|_* \le\left( \frac{2C_{\rm conv}L_{\rm nl}}{M} +\frac{C_{\rm conv}R_{\rm rem}}{M^2} \right)\|x-y\|_* <\frac12\|x-y\|_*. \tag{C.7.15b}\label{eq:C.7.15b} \end{equation*} The last two inequalities in \eqref{eq:C.7.14} make the two terms in \eqref{eq:C.7.15b} separately less than $1/4$. The unique fixed point solves \eqref{eq:C.7.8} and satisfies $|x_n|\le n^{-1/2}$, which gives \eqref{eq:C.7.3}.

The stable sum in \eqref{eq:C.7.13} is empty at $n=N$, so $\Pi_{\rm st}x_N=0$. Conversely, any solution of \eqref{eq:C.7.8} in the ball \eqref{eq:C.7.12} with this stable initial condition satisfies the stable variation-of-constants sum. Iterating its complementary component backwards from an index $m$ leaves a boundary term bounded by $C(m/n)^{1/4}m^{-1/2}$, which tends to zero. It therefore satisfies the complementary sum in \eqref{eq:C.7.13}, hence is the same fixed point. This is precisely the local uniqueness stated in the lemma.

Since $1/N<2^{-2000000}$, every evaluation lies in \eqref{eq:C.4.6}, and $K_{\rm Fl}/\sqrt N<2^{-100}$ keeps the reconstructed tail in the positive tube. The corresponding state in weight coordinates is, by definition, \begin{equation*} z_n=\gamma_0(\tau_n)+\mathsf P(\tau_n)x_n, \tag{C.7.16}\label{eq:C.7.16} \end{equation*} where the periodic Floquet factor is evaluated modulo $T_{\rm per}$. \end{proof}

The complementary block in \eqref{eq:C.7.13} contains every neutral or unstable Floquet direction, together with every weakly negative direction outside the chosen stable subspace. The construction accommodates nonhyperbolic phase dynamics and every value of the first Lyapunov coefficient.

\section{Accumulation along the periodic orbit}\label{harmonic-times-see-the-whole-periodic-orbit}

Shadowing gives proximity to a periodic curve. Nonconvergence follows once the harmonic time mesh is shown to accumulate at every phase of that curve.

The sequence $\tau_n$ is increasing, unbounded, and $\tau_{n+1}-\tau_n=1/n\to0$. Given any phase $\theta\pmod {T_{\rm per}}$, let $n_j$ be the first index with $\tau_{n_j}\ge jT_{\rm per}+\theta$. Its overshoot is at most $1/(n_j-1)$, which tends to zero. Hence $\tau_n\bmod T_{\rm per}$ has every phase as an accumulation point.

The $a$-coordinate of $\Gamma_0$ has first Fourier coefficient $\varepsilon_0$, so it is nonconstant. Choose two phases with distinct $a$-values and subsequences approaching them. By \eqref{eq:C.7.3}, the weight-coordinate sequence $a_n$ has two distinct accumulation values.

At the limiting $P$-phase of the fixed-point family, \begin{equation*} w_i^P(a)=\frac{( \xi_i-a)Q_{q_{\rm per}}(\xi_i)}{\Delta'(\xi_i)}. \end{equation*} Every slope is nonzero, so $a\mapsto(w_1^P(a),\ldots,w_6^P(a))$ is injective. Define the reconstructed squared-weight vector $\widetilde W_n=\chi_{q_{\rm per}}^{-1}(a_n,n^{-1}V_n,n^{-2}\mathcal E_n)$. Equation \eqref{eq:C.4.5} gives uniformly \begin{equation*} \widetilde W_{n,i}=w_i^P(a_n)+O(n^{-1}),\qquad1\le i\le6. \end{equation*} The sequence $(\widetilde W_n)$ has at least two distinct accumulation points.

\begin{figure}[t]
\centering
\includegraphics[width=0.72\textwidth]{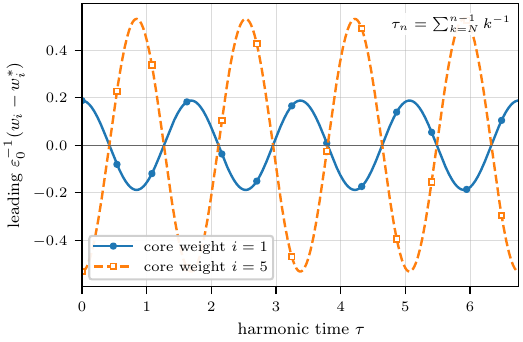}
\caption{Magnified leading periodic profile for the certified restart-four counterexample. The Hopf eigenvector is normalised by $v_a=1$, so the first harmonic of the recovery-root coordinate is $2\varepsilon_0\cos(\omega_*\tau)$. Since the fixed-family weights \eqref{eq:C.3.5} are affine in $a$, the displayed curves are $2\partial_a w_1^P\cos(\omega_*\tau)$ and $2\partial_a w_5^P\cos(\omega_*\tau)$, evaluated at the certified Hopf point, with $w_i^*=w_i^P(\alpha(q_*))$. The correction $h_0$ in \eqref{eq:C.6.15} has norm at most $2^{-359999999}$ by \eqref{eq:C.6.6}--\eqref{eq:C.6.7}. The exact orbit shadows the profile at $\tau_n=\sum_{k=N}^{n-1}k^{-1}$ as in \eqref{eq:C.7.2}--\eqref{eq:C.7.3}; these times accumulate at every phase modulo the period. Because both displayed slopes are nonzero, the even squared weights, and hence the even normalised residuals, have more than one accumulation point.}
\label{fig:s4-nonconvergence}
\end{figure}

Figure~\ref{fig:s4-nonconvergence} illustrates the leading periodic modulation behind this nonconvergence.

\section{Reconstruction of the SPD counterexample}\label{the-one-fixed-exact-spd-datum}

The positive nonconvergent orbit in weight coordinates obtained in Sections~\ref{literal-selection-of-one-periodic-orbit}--\ref{harmonic-times-see-the-whole-periodic-orbit} determines a fixed eight-dimensional SPD problem.

Fix the parameter $q_{\rm per}$ from \eqref{eq:C.6.15}. In the fixed rational isolating intervals, let $\xi_1<\cdots<\xi_6$ be the six selected roots of $PQ_{q_{\rm per}}-C$, and let $\xi_7,\xi_8$ be the selected second and seventh roots of $PQ_{q_{\rm per}}+C$.

These are uniquely specified real numbers. Define \begin{equation*} \lambda_i=\xi_i+5,\qquad A=\operatorname{diag}(\lambda_1,\ldots,\lambda_8). \tag{C.9.1}\label{eq:C.9.1} \end{equation*} The isolating intervals give $1<\lambda_i<10$, so $A$ is real symmetric positive definite.

Let $z_N$ be its uncentred state in weight coordinates from \eqref{eq:C.7.16}. Define $$ \Psi_{q,n}(a,V,\mathcal E) :=\chi_q^{-1}\bigl(a,n^{-1}V,n^{-2}\mathcal E\bigr). $$ This is the weight reconstructed from the scaled weight coordinate at $t=1/n$. The definitions \eqref{eq:C.4.5} and \eqref{eq:C.4.25} give the covariance identity $$ \Psi_{q,n+1}\!\left(\widehat F_q(1/n,z)\right) =T_4^2\Psi_{q,n}(z). $$ Hence define \begin{equation*} W_0=\Psi_{q_{\rm per},N}(z_N). \end{equation*} This is the squared-weight vector reconstructed from \eqref{eq:C.4.2} and \eqref{eq:C.4.5}. With $z_N=(a_N,V_N,\mathcal E_N)$, its weight coordinates satisfy \begin{equation*} \zeta=N^{-1}V_N,\qquad E_g=N^{-2}\mathcal E_{g,N}, \end{equation*} with the six core weights recovered affinely by \eqref{eq:C.4.2}. Lemma~\ref{lemma-7.1-quantitative-nonhyperbolic-shadowing} proves $W_{0,i}>0$ and $\sum_iW_{0,i}=1$.

Finally define \begin{equation*} x_0=0,\qquad b=(\sqrt{W_{0,1}},\ldots,\sqrt{W_{0,8}})^T, \tag{C.9.4}\label{eq:C.9.4} \end{equation*} using the positive square roots. Equations \eqref{eq:C.3.1}--\eqref{eq:C.3.3a}, \eqref{eq:C.6.3}, \eqref{eq:C.6.6}--\eqref{eq:C.6.7}, \eqref{eq:C.7.12}--\eqref{eq:C.7.13}, and \eqref{eq:C.9.1}--\eqref{eq:C.9.4} uniquely specify $A,b,x_0$. The critical parameter $q_*$ and the critical scalar configuration are algebraic, being selected by rational polynomial isolation. The displayed contractions and isolation conditions select $\mu_0$, $q_{\rm per}$, the periodic and shadowing profiles, the roots at $q_{\rm per}$, $W_0$, and its square roots as computable real numbers.

For every $j\ge0$, let $W_j=\Psi_{q_{\rm per},N+j}(z_{N+j})$. Coordinate covariance gives \begin{equation*} T_4^2W_j=W_{j+1}. \tag{C.9.5}\label{eq:C.9.5} \end{equation*} Here $W_j$ is the squared-weight vector at restart index $2j$, and $T_4W_j$ is the vector at restart index $2j+1$. The spectral parameter $q_{\rm per}$, all eight eigenvalues, and $A$ are fixed once and for all; the rescaling index $N+j$ is only a coordinate used to describe the orbit.

The positive tube preserves every core and external weight and every individual first- and second-phase block factor. Hence both $W_j$ and $T_4W_j$ have all eight positive coordinates for every $j$. Every residual at a restart point has grade eight, so no block with restart length four terminates.

By Section~\ref{harmonic-times-see-the-whole-periodic-orbit}, $W_j$ has two distinct accumulation points. If the normalised residuals $y_{2j}$ converged in Euclidean norm, each continuous coordinate square $\langle y_{2j},\mathbf e_i\rangle^2=W_{j,i}$ would converge, contradicting the two distinct accumulation points of $W_j$. Hence the even residual subsequence does not converge. Proposition~\ref{prop:s4-seed} follows.

\section*{Degree elevation and arbitrary restart length}

The construction at restart length four is complete. To pass from length $r$ to length $r+1$, we isolate the algebraic data determining the Hopf point at length $r$, adjoin one distant core node, and prove that the critical branch and transverse crossing persist. The return-map expansion and the two functional-analytic arguments are then repeated at the elevated length.

Throughout Sections~\ref{sec:C10}--\ref{sec:C18}, every $O(R^{-k})$, $o(1)$, and differentiated asymptotic estimate is uniform on the fixed compact parameter set in the displayed $(q,a)$ variables and on each fixed bounded set of evaluation variables. Its constant may depend on the unelevated configuration, the fixed separation $d$, the chosen compact parameter sets, and the stated derivative order, but not on $R$ or on a variable ranging in those sets. Before a limiting argument is used to preserve a strict inequality, the parameter set is shrunk once and a positive uniform margin is fixed.

\section{Scalar two-cycle configurations with two external nodes}\label{sec:C10}

Fix an integer $r\ge4$. A scalar two-cycle configuration with two external nodes at restart length $r$ consists of monic degree-$r$ polynomials $P_q,Q_q$, analytic in one real parameter $q$, a constant $C>0$, and a factorisation into monic polynomials \begin{equation*} P_qQ_q-C=\Delta_qH_q, \qquad \deg\Delta_q=r+2, \qquad \deg H_q=r-2. \end{equation*} All roots of $P_q,Q_q,\Delta_q,H_q$ are assumed real and simple, $P_q,Q_q$ are coprime, and the roots of $\Delta_q$ and $H_q$ are disjoint. The $r+2$ roots $\xi_i$ of $\Delta_q$ are the core spectral nodes. With \begin{equation*} S_a(t)=t-a, \end{equation*} the two phases of the fixed-point family have weights \begin{equation*} w_i^P(a)=\frac{S_a(\xi_i)Q_q(\xi_i)}{\Delta_q'(\xi_i)}, \qquad w_i^Q(a)=\frac{S_a(\xi_i)P_q(\xi_i)}{\Delta_q'(\xi_i)}. \tag{C.10.3}\label{eq:C.10.3} \end{equation*} The roots of $P_q$ and $Q_q$ occupy the same $r$ of the $r+1$ open gaps between successive roots of $\Delta_q$, with one root of each polynomial in each occupied gap; $a$ ranges over the remaining gap. Consequently all weights in \eqref{eq:C.10.3} are positive.

Here and below, positivity refers to the phase weights. Section~\ref{sec:C17} translates the spectral nodes into $(0,\infty)$ without changing the squared-weight dynamics.

Let $T_r$ denote the squared-weight map for one block of restart length $r$. The Lagrange identity \begin{equation*} \sum_{\Delta_q(\xi)=0}\frac{R(\xi)}{\Delta_q'(\xi)} =[t^{r+1}]R(t),\qquad \deg R\le r+1, \end{equation*} shows that both weight vectors sum to one. At every core node, $P_q(\xi_i)Q_q(\xi_i)=C$. Hence $P_q$ is the monic degree-$r$ orthogonal polynomial for $w^P$, its squared monic norm is $C$, and similarly for $Q_q,w^Q$. It follows exactly that \begin{equation*} T_rw^P=w^Q,\qquad T_rw^Q=w^P. \end{equation*}

Choose two bounded simple roots $h_-,h_+$ of \begin{equation*} P_q(t)Q_q(t)+C=0. \tag{C.10.6}\label{eq:C.10.6} \end{equation*} At either node the signed two-block multiplier is $-1$ and its square is $+1$. These two nodes are the external modes.

Let $B=(\delta C,\delta H)\in\mathbb R^{r-1}$, let $u$ be the single intrinsic core-normal coordinate, and let $\varphi\in\mathbb R^2$ be the logarithmic external-amplitude coordinate. Let $A_g,V_g,\ell_g,d_g,\kappa$ denote the higher-degree analogues defined explicitly in \eqref{eq:C.13.1}--\eqref{eq:C.13.3}. At a critical point at which the recovery root $a$ is constant, let $c_g$ denote the equilibrium external amplitudes and set $ p_{\rm lin}=\sum_gc_g\partial_aA_g,\quad \mathbf A_c=(c_gA_g)_g,\quad q_{\rm lin}=\sum_gc_g\partial_aV_g, $ $ \mathsf V_c=(c_gV_g)_g,\qquad \ell=(\ell_g)_g,\qquad DB=(d_g[B])_g. $ The higher-degree analogues of \eqref{eq:C.3.6}--\eqref{eq:C.3.10} give the following linearisation of the leading vector field. Section~\ref{sec:C14} derives the same coefficients from the exact two-block map: \begin{equation*} \begin{aligned} \lambda a&=\kappa u,\\ \lambda u&=p_{\rm lin}a+u+\mathbf A_c\varphi,\\ \lambda B&=q_{\rm lin}a+B+\mathsf V_c\varphi,\\ \lambda\varphi&=\ell u+DB. \end{aligned} \end{equation*} Here $\mathbf A_c$ is a $1\times2$ row, $\mathsf V_c$ is an $(r-1)\times2$ matrix, and $D$ is a $2\times(r-1)$ matrix. For $\lambda\ne0,1$, set $\theta=\lambda(\lambda-1)$, $\zeta=\lambda\varphi$, and $\chi=\zeta-\ell u$. Eliminating $a,B$ gives the three-dimensional eigenvalue problem \begin{equation*} \theta\binom{u}{\chi} =\mathcal M\binom{u}{\chi}, \qquad \mathcal M= \begin{pmatrix} \kappa p_{\rm lin}+\mathbf A_c\ell&\mathbf A_c\\ D\mathsf V_c\ell+\kappa Dq_{\rm lin}&D\mathsf V_c \end{pmatrix}. \tag{C.10.8}\label{eq:C.10.8} \end{equation*} Let $J_{\rm sc}$ denote the Jacobian of this leading field at the critical point. Multiplying back the eliminated factors and comparing the resulting monic polynomials gives the exact characteristic identity \begin{equation*} \det(\lambda I-J_{\rm sc}) =(\lambda-1)^{r-3} \det\!\left(\lambda(\lambda-1)I_3-\mathcal M\right). \tag{C.10.9}\label{eq:C.10.9} \end{equation*} This identity is valid also at $\lambda=0,1$, by polynomial continuation.

At a critical point the normal balance and the two external rate equations are \begin{equation*} \mathbf A_c\mathbf1=0,\qquad D\mathsf V_c\mathbf1=2\mathbf1. \end{equation*} Hence $\mathcal M(0,\mathbf1)^T=2(0,\mathbf1)^T$: the $\theta$-eigenvalue $2$ is universal. We call the configuration a positive transverse Hopf configuration if the remaining quadratic in $\theta$ has negative discriminant and, at a parameter $q_*$, its four $\lambda$-preimages include a simple pair $\pm i\omega_*$, with nonzero crossing speed, and the corresponding eigenvector has a nonzero $a$-component.

The construction at restart length four, including the rational and Sturm certificates given in Appendix C, supplies such a configuration at $r=4$. The analytic argument below raises its degree one step at a time.

\section{One-step degree elevation}\label{sec:C11}

\begin{lemma}[Root continuation and positivity]\label{lem:C11-root-continuation}
Suppose a scalar two-cycle configuration with two external nodes exists at restart length $r$. Fix $d>0$. For a large positive number $R$, define \begin{equation*} \begin{aligned} P_{R,q}(t)&=(t-R)P_q(t),\\ Q_{R,q}(t)&=(t-R-d)Q_q(t),\\ C_R&=R(R+d)C. \end{aligned} \end{equation*} For every sufficiently large $R$, uniformly for $q$ in a small neighbourhood of $q_*$, this is the phase-polynomial part of a positive scalar two-cycle configuration at restart length $r+1$. The roots of the polynomial $P_qQ_q-C$ from the unelevated configuration continue, and the two additional roots are \begin{equation*} \begin{aligned} x_R^-&=R-\frac{C}{d}R^{2-2r}\{1+O(R^{-1})\},\\ x_R^+&=R+d+\frac{C}{d}R^{2-2r}\{1+O(R^{-1})\}. \end{aligned} \tag{C.11.2}\label{eq:C.11.2} \end{equation*} The unelevated bounded roots of $P_qQ_q+C$, including $h_-,h_+$, continue to bounded neutral external nodes. \end{lemma}

\begin{proof}
Dividing $P_{R,q}Q_{R,q}-C_R$ by $R(R+d)$ gives \begin{equation*} \frac{P_{R,q}(t)Q_{R,q}(t)-C_R}{R(R+d)} =P_q(t)Q_q(t) \left(1-\frac{t}{R}\right) \left(1-\frac{t}{R+d}\right)-C. \tag{C.11.3}\label{eq:C.11.3} \end{equation*} On bounded $t$-sets this converges, with every prescribed finite number of $q,t$ derivatives, to $P_qQ_q-C$. All $2r$ unelevated simple roots have unique bounded analytic continuations.

Let $t=R+R^{2-2r}z$. Uniformly for bounded $z$, $$ 1-\frac{t}{R}=-R^{1-2r}z, \qquad 1-\frac{t}{R+d}=\frac{d+O(R^{2-2r})}{R+d}, $$ and $P_q(t)Q_q(t)=R^{2r}\{1+O(R^{-1})\}$. Equation \eqref{eq:C.11.3} becomes \begin{equation*} -dz-C+O(R^{-1})(1+|z|)=0. \end{equation*} The implicit-function theorem gives the first expansion in \eqref{eq:C.11.2}. The substitution $t=R+d+R^{2-2r}z$ similarly gives \begin{equation*} dz-C+O(R^{-1})(1+|z|)=0, \end{equation*} and hence the second expansion. The expansions and their first two $q$-derivatives are uniform. Together with the $2r$ bounded roots, these two roots exhaust the degree $2r+2$ polynomial. Hence all its roots are real and simple for large $R$.

Let the monic polynomial $\Delta_R$ have the continued core roots of the unelevated configuration and $x_R^+$, and let the monic polynomial $H_R$ have the continued unelevated $H$-roots and $x_R^-$. Then \begin{equation*} P_{R,q}Q_{R,q}-C_R=\Delta_RH_R, \quad \deg\Delta_R=r+3, \quad \deg H_R=r-1. \end{equation*} The new roots $R,R+d$ of $P_{R,q},Q_{R,q}$ lie in the same new core gap. Positivity follows quantitatively: every bounded phase weight converges in $C^2(q,a)$ to the corresponding strictly positive unelevated weight, while direct substitution at $X=x_R^+$ gives \begin{equation*} w_{X,R}^P=\frac{C}{d}R^{1-2r}\{1+O_{C^2}(R^{-1})\}, \qquad w_{X,R}^Q=\frac{d}{R}\{1+O_{C^2}(R^{-1})\}. \tag{C.11.7}\label{eq:C.11.7} \end{equation*} Hence every phase weight is positive for all sufficiently large $R$. Applying the bounded-root argument to $P_{R,q}Q_{R,q}+C_R$ continues $h_-,h_+$. At either continued root, $P_{R,q}Q_{R,q}/C_R=-1$, so its squared two-block multiplier remains exactly $+1$. \end{proof}

Coprimality also persists: $R$ and $R+d$ are distinct, both lie beyond all unelevated roots, and neither is an unelevated root of the opposite phase polynomial.

\section{The distant spectral atom}\label{sec:C12}

The disparity in \eqref{eq:C.11.7} is the remaining continuity issue: in one phase a weight tends to zero much faster than its node tends to infinity. The next lemma controls the high-degree moment direction created by this far atom.

\begin{lemma}[Limit of the reproducing kernel]\label{lem:C12-kernel-limit}
Let $\mu_R$ be either elevated phase measure, let $X=x_R^+$, and let $K_R^{\rm rep}$ be the reproducing kernel of $\mathcal P_r$ in $L^2(\mu_R)$. Let $K^{\rm rep}$ be the unelevated reproducing kernel of $\mathcal P_{r-1}$. For bounded arguments, \begin{equation*} K_R^{\rm rep}(t,g)\longrightarrow K^{\rm rep}(t,g) \end{equation*} uniformly with two derivatives in $(q,a)$, and also locally uniformly with the required derivatives in the bounded evaluation variables $t,g$. Apart from the common $O(R^{-1})$ perturbation of the bounded atoms, the complementary kernel term is $O(R^{-1})$ in the $P$-phase and $O(R^{1-2r})$ in the $Q$-phase. \end{lemma}

\begin{proof}
Use the basis \begin{equation*} v_{j,R}(t)=\left(1-\frac{t}{X}\right)t^j, \quad0\le j<r, \qquad \phi_R(t)=\left(\frac{t}{X}\right)^r. \end{equation*} The $v_{j,R}$ vanish at $X$, span an $r$-dimensional subspace, and converge on bounded sets to $1,t,\ldots,t^{r-1}$. Evaluation at $X$ shows that their span is complementary to $\phi_R$.

Define \begin{equation*} \Gamma_R=(\langle v_{i,R},v_{j,R}\rangle_{\mu_R})_{i,j}, \qquad b_R=(\langle v_{i,R},\phi_R\rangle_{\mu_R})_i. \tag{C.12.3}\label{eq:C.12.3} \end{equation*} The far atom contributes exactly zero to both quantities. The bounded atoms and weights converge in $C^2(q,a)$, so $\Gamma_R$ converges to the unelevated positive moment Gram matrix $\Gamma$ and $\Gamma_R^{-1}\to\Gamma^{-1}$. Moreover $b_R=O_{C^2}(R^{-r})$. Orthogonalise the last vector: \begin{equation*} \rho_R(t)=\phi_R(t)-v_R(t)^T\Gamma_R^{-1}b_R, \qquad \sigma_R=\|\rho_R\|_{\mu_R}^2. \end{equation*} For bounded $t$, $\rho_R(t)=O_{C^2}(R^{-r})$, while $\rho_R(X)=1$. It follows that \begin{equation*} \sigma_R=w_{X,R}+O_{C^2}(R^{-2r}). \end{equation*} The kernel has the orthogonal decomposition \begin{equation*} K_R^{\rm rep}(t,g)=v_R(t)^T\Gamma_R^{-1}v_R(g) +\frac{\rho_R(t)\rho_R(g)}{\sigma_R}. \tag{C.12.6}\label{eq:C.12.6} \end{equation*} The first term tends to $K^{\rm rep}(t,g)$. In the $P$-phase the second term is \begin{equation*} \frac{O(R^{-2r})}{(C/d)R^{1-2r}(1+o(1))}=O(R^{-1}), \end{equation*} and in the $Q$-phase it is $O(R^{1-2r})$. Two derivatives of \eqref{eq:C.12.3}--\eqref{eq:C.12.6} satisfy the same estimates because the leading coefficients in \eqref{eq:C.11.7} are uniformly separated from zero. \end{proof}

\section{Persistence of the critical point and Hopf crossing}\label{sec:C13}

Lemma~\ref{lem:C11-root-continuation} supplies a positive elevated configuration, and Lemma~\ref{lem:C12-kernel-limit} controls the only kernel contribution from the distant node. These estimates continue the critical point and its transverse Hopf crossing. In this section and the next, the parameter $q$ is suppressed in polynomial subscripts: $P=P_q$, $P_R=P_{R,q}$, and similarly for $Q,\Delta,H$. Let $K_P$ and $K_Q$ denote the two phase reproducing kernels obtained from the analogues of \eqref{eq:C.3.7}--\eqref{eq:C.3.8}.

For a bounded external source node $g$, define \begin{equation*} \begin{aligned} A_g(a)&=\frac{4P(g)}{g-a} \left(\frac{P(g)}{P(a)}-\frac{Q(g)}{Q(a)}\right),\\ p_g(t)&=-4P(g)K_P(t,g),\\ q_g(t)&= 4P(g)K_Q(t,g),\\ (V_g)_C&=4\{P(g)^2+C\}. \end{aligned} \tag{C.13.1}\label{eq:C.13.1} \end{equation*} The polynomial component $(V_g)^H\in\mathcal P_{r-3}$ is determined by the differentiated factor identity \begin{equation*} p_g(t)Q(t)+P(t)q_g(t)-(V_g)_C-A_g(a)\frac{\Delta(t)-\Delta(a)}{t-a}=\Delta(t)(V_g)^H(t). \end{equation*} Thus $V_g=((V_g)_C,(V_g)^H)$. For any bounded external node $h$, define the linear functional $ d_h[\delta C,\delta H]=-\frac{2}{C}\{\Delta(h)\delta H(h)+2\delta C\}. $ For a source node $g$, the linear coupling coefficient is \begin{equation*} \mathsf T_{hg}=d_h[V_g] =-\frac{2}{C}\left\{ p_g(h)Q(h)+P(h)q_g(h)+(V_g)_C -A_g\frac{\Delta(h)-\Delta(a)}{h-a}\right\}. \tag{C.13.2}\label{eq:C.13.2} \end{equation*} Let $\mathsf T=(\mathsf T_{hg})_{h,g}$. The remaining coefficients used in \eqref{eq:C.10.8} are \begin{equation*} \kappa=\frac{2\Delta(a)}{C}, \qquad \ell_h=-\frac{2}{C}\frac{\Delta(h)-\Delta(a)}{h-a}. \tag{C.13.3}\label{eq:C.13.3} \end{equation*}

For the elevated family use the bounded normalisations \begin{equation*} \begin{aligned} \bar P_R&=-P_R/R,& \bar Q_R&=-Q_R/(R+d),\\ \bar\Delta_R&=-\Delta_R/(R+d),& \bar H_R&=-H_R/R,\\ \bar C_R&=C_R/[R(R+d)]=C. \end{aligned} \end{equation*} All bounded evaluations of these polynomials converge in $C^2(q,a)$ to the unelevated data. These normalisations arise from an analytic change of scalar-field coordinates. If $u_R$ is the original intrinsic normal and $B_R=(\delta C_R,\delta H_R)$, define \begin{equation*} \bar u=-\frac{u_R}{R},\qquad \delta\bar C=\frac{\delta C_R}{R(R+d)},\qquad \delta\bar H=-\frac{\delta H_R}{R}. \tag{C.13.4a}\label{eq:C.13.4a} \end{equation*} The external scaled amplitudes and the recovery root are unchanged. Since $\kappa_R=-\bar\kappa_R/R$, this change turns $\dot a=\kappa_Ru_R$ into $\dot a=\bar\kappa_R\bar u$; the remaining block scalings below arise from this same analytic coordinate change. Define \begin{equation*} \begin{aligned} \bar\kappa_R&=\frac{2\bar\Delta_R(a)}{C},\\ \bar A_{g,R}&=-A_{g,R}/R,\\ \bar\ell_{g,R}&=-\frac{2}{C} \frac{\bar\Delta_R(g)-\bar\Delta_R(a)}{g-a},\\ \bar p_{g,R}&=-p_{g,R}/R, \qquad \bar q_{g,R}=-q_{g,R}/(R+d),\\ (\bar V_{g,R})_C&=(V_{g,R})_C/[R(R+d)],\\ (\bar V_{g,R})^H&=-(V_{g,R})^H/R. \end{aligned} \end{equation*} Lemma~\ref{lem:C12-kernel-limit} gives the following convergences in $C^2$: \begin{equation*} \begin{aligned} \bar\kappa_R&\to\kappa, & \bar A_{g,R}&\to A_g, & \bar\ell_{g,R}&\to\ell_g,\\ \bar p_{g,R}&\to p_g, & \bar q_{g,R}&\to q_g, & (\bar V_{g,R})_C&\to4\{P(g)^2+C\}. \end{aligned} \tag{C.13.6}\label{eq:C.13.6} \end{equation*}

The convergence of the remaining $H$-component follows from the differentiated factor identity \begin{equation*} \bar p_{g,R}\bar Q_R+\bar P_R\bar q_{g,R} -(\bar V_{g,R})_C -\bar A_{g,R}\frac{\bar\Delta_R(t)-\bar\Delta_R(a)}{t-a} =\bar\Delta_R(\bar V_{g,R})^H. \tag{C.13.7}\label{eq:C.13.7} \end{equation*} Evaluate \eqref{eq:C.13.7} at $r-1$ fixed bounded points that remain away from the roots of $\Delta$. Polynomial interpolation in $\mathcal P_{r-2}$ proves coefficientwise $C^2$ convergence of $(\bar V_{g,R})^H$ to the unelevated $(V_g)^H$; its new highest coefficient tends to zero. The normalised external functional is \begin{equation*} \bar d_{h,R}[\delta\bar C,\delta\bar H] =-\frac{2}{C}\{\bar\Delta_R(h)\delta\bar H(h) +2\delta\bar C\}. \tag{C.13.8}\label{eq:C.13.8} \end{equation*} Hence $\mathsf T_{hg,R}=\bar d_{h,R}[\bar V_{g,R}]\to \mathsf T_{hg}$ in $C^2(q,a)$, including all mixed derivatives needed below. Let $\mathsf T_R=(\mathsf T_{hg,R})_{h,g}$.

\begin{lemma}[Critical branch and transverse Hopf persistence]\label{lem:C13-hopf-persistence}
For every sufficiently large finite $R$, the elevated family of restart length $r+1$ has a positive transverse Hopf point, whose Hopf eigenvector has a nonzero $a$-component, near the unelevated point. \end{lemma}

\begin{proof}
Let $c=(c_-,c_+)^T$ be the positive scaled external amplitudes and let $A=(A_-,A_+)$. Let $\bar A_R=(\bar A_{-,R},\bar A_{+,R})$. The critical equations are \begin{equation*} \mathcal C_R(a,c,q)= \begin{pmatrix} \bar A_R(a,q)c\\ (\mathsf T_R(a,q)c)_--2\\ (\mathsf T_R(a,q)c)_+-2 \end{pmatrix}=0. \end{equation*} Denote the unelevated map by $\mathcal C$. At the unelevated critical point, $D_{(a,c)}\mathcal C$ is invertible. For its determinant, define $$ \mathbf A_c=(c_jA_j)_j, \qquad \mathsf T_c=(\mathsf T_{ij}c_j)_{ij}, \qquad g_{\rm lin}=(\sum_jc_j\partial_a\mathsf T_{ij})_i. $$ Multiplying the two amplitude columns of $D_{(a,c)}\mathcal C$ by $c_-,c_+$ gives \begin{equation*} J_c= \begin{pmatrix} \sum_jc_jA_j'&\mathbf A_c\\ g_{\rm lin}&\mathsf T_c \end{pmatrix}, \qquad \det J_c=c_-c_+\det D_{(a,c)}\mathcal C. \end{equation*} The corresponding $\theta$-matrix is \begin{equation*} \mathcal M= \begin{pmatrix} \kappa\sum_jc_jA_j'+\mathbf A_c\ell&\mathbf A_c\\ \mathsf T_c\ell+\kappa g_{\rm lin}&\mathsf T_c \end{pmatrix}. \end{equation*} Its first column is $\kappa$ times the first column of $J_c$, plus $\ell_-$ and $\ell_+$ times the last two columns of $J_c$; its last two columns equal those of $J_c$. Multilinearity of the determinant gives the orientation identity \begin{equation*} \det\mathcal M=\kappa c_-c_+\det D_{(a,c)}\mathcal C. \end{equation*} At the Hopf point $\mathcal M$ has eigenvalues $2,\theta,\bar\theta$, none zero, and $\kappa c_-c_+\ne0$. Hence the critical Jacobian is invertible. Equations \eqref{eq:C.13.6}--\eqref{eq:C.13.8} and the implicit-function theorem produce a unique nearby positive critical branch for all sufficiently large $R$.

Let $y_R(q)=(a_R(q),c_{-,R}(q),c_{+,R}(q))$ denote the critical branch just obtained. We distinguish the elevated $\theta$-matrix $\mathcal M_R$ in the original scalar coordinates from its matrix $\bar{\mathcal M}_R$ in the barred coordinates \eqref{eq:C.13.4a}. Along this branch, define $$ \begin{gathered} \bar{\mathbf A}_{c,R}=(c_{j,R}\bar A_{j,R})_{j=-,+},\qquad (\mathsf T_{c,R})_{ij}=\mathsf T_{ij,R}c_{j,R},\\ (\bar g_{{\rm lin},R})_i=\sum_{j=-,+}c_{j,R}\partial_a\mathsf T_{ij,R}, \qquad \bar\ell_R=(\bar\ell_{-,R},\bar\ell_{+,R})^T. \end{gathered} $$ The reduced matrix computed directly in these coordinates is \begin{equation*} \bar{\mathcal M}_R= \begin{pmatrix} \bar\kappa_R\sum_jc_{j,R}\partial_a\bar A_{j,R} +\bar{\mathbf A}_{c,R}\bar\ell_R&\bar{\mathbf A}_{c,R}\\ \mathsf T_{c,R}\bar\ell_R+\bar\kappa_R\bar g_{{\rm lin},R} &\mathsf T_{c,R} \end{pmatrix}. \end{equation*} This rescaling is a coordinate change. The full linear change is $$ (a,u_R,\delta C_R,\delta H_R,\varphi) \longmapsto \left(a,-\frac{u_R}{R},\frac{\delta C_R}{R(R+d)}, -\frac{\delta H_R}{R},\varphi\right). $$ Using the same amplitude-column change $\varphi\mapsto\operatorname{diag}(c_{-,R},c_{+,R})\varphi$ in both systems and substituting in the reduced eigenvalue equations gives \begin{equation*} \boxed{\quad \bar{\mathcal M}_R=\mathcal D_R^{-1}\mathcal M_R\mathcal D_R, \qquad \mathcal D_R=\operatorname{diag}(-R,1,1). \quad} \tag{C.13.13a}\label{eq:C.13.13a} \end{equation*} In particular, the upper-right and lower-left blocks of $\mathcal M_R$ are $-R$ and $-R^{-1}$ times the corresponding barred blocks.

Let $y(q)=(a(q),c_-(q),c_+(q))$ denote the unelevated critical branch. Fix a small closed $q$-interval $I_q$ about the unelevated critical parameter. Equations \eqref{eq:C.13.6}--\eqref{eq:C.13.8}, including their mixed $(q,a)$ derivatives, and the parameter-dependent implicit-function theorem give \begin{equation*} \|y_R-y\|_{C^1(I_q)}\longrightarrow0, \qquad \|\bar{\mathcal M}_R-\mathcal M\|_{C^1(I_q)}\longrightarrow0. \tag{C.13.13b}\label{eq:C.13.13b} \end{equation*} At every critical point, $(0,1,1)^T$ is an eigenvector of $\bar{\mathcal M}_R$ with eigenvalue $2$. Define $$ \tau_R=\operatorname{tr}\bar{\mathcal M}_R-2,\qquad K_{\theta,R}=\tfrac12\det\bar{\mathcal M}_R,\qquad \mathscr D_{\theta,R}=\tau_R^2-4K_{\theta,R}. $$ Then $$ \det(\theta I_3-\bar{\mathcal M}_R) =(\theta-2)(\theta^2-\tau_R\theta+K_{\theta,R}), $$ and \eqref{eq:C.13.13a} gives the same factorisation for $\mathcal M_R$. Define \begin{equation*} \mathcal H_R=\frac{\tau_R}{2}-\frac{\mathscr D_{\theta,R}}4. \end{equation*} These quantities converge in $C^1(I_q)$ to their unelevated counterparts.

To verify transversality, let $\theta_R=\sigma_{\theta,R}+i\upsilon_R$ be one root of the quadratic and choose $\lambda_R$ so that $\lambda_R(\lambda_R-1)=\theta_R$. At a zero of $\mathcal H_R$, label the conjugate roots so that $\lambda_R=i\omega_R$, $\omega_R>0$. Then $\theta_R=-\omega_R^2-i\omega_R$ and $\mathcal H_R=\sigma_{\theta,R}+\upsilon_R^2$. Differentiating $(2\lambda_R-1)\lambda_R'=\theta_R'$ gives \begin{equation*} \left.\frac{d}{dq}\operatorname{Re}\lambda_R(q)\right|_{q=q_R} =-\frac{\mathcal H_R'(q_R)}{1+4\omega_R^2}. \tag{C.13.14a}\label{eq:C.13.14a} \end{equation*} The unelevated functional has a simple zero. Equation \eqref{eq:C.13.13b} gives, for every sufficiently large $R$, a unique nearby $q_R$ for which \begin{equation*} \mathcal H_R(q_R)=0,\qquad \mathcal H_R'(q_R)\ne0, \qquad \mathscr D_{\theta,R}(q_R)<0. \end{equation*} Equation \eqref{eq:C.13.14a} proves the transverse Hopf crossing, and every strict positivity margin persists.

Let $J_{{\rm sc},R}$ denote the elevated scalar-field Jacobian at the critical point. For restart length $r+1$, \eqref{eq:C.10.9} becomes \begin{equation*} \det(\lambda I-J_{{\rm sc},R}) =(\lambda-1)^{r-2} \det\!\left(\lambda(\lambda-1)I_3-\mathcal M_R\right). \end{equation*} Elevation adds one further real eigenvalue $+1$. The universal $\theta=2$ gives $\lambda=-1,2$, while the nonreal $\theta$-pair gives $\lambda=\pm i\omega_R,1\pm i\omega_R$. No other eigenvalue lies on the imaginary axis.

Finally, suppose the Hopf eigenvector had zero intrinsic normal component $u$. In the $\theta$-eigenproblem, its two-dimensional external component would be a nonzero eigenvector of the real matrix $\mathsf T_{c,R}$ with nonreal eigenvalue $\theta$. But $\mathsf T_{c,R}$ already has the real eigenvalue $2$, and a real $2\times2$ matrix with one real eigenvalue has both eigenvalues real. This is impossible. Hence $u\ne0$; since $\lambda a=\bar\kappa_Ru$, the $a$-component is nonzero. \end{proof}

\section{The two-block return map at arbitrary restart length}\label{sec:C14}

Lemma~\ref{lem:C13-hopf-persistence} supplies a Hopf point in scalar coordinates. To obtain a restarted CG orbit, that scalar field must be identified with the leading term of the exact two-block map. Sections~\ref{a-complete-physical-raw-chart-and-its-exact-expansion}--\ref{exact-transfer-of-the-hopf-point-to-the-raw-field} established this identification quantitatively at restart length four; the following analytic form is valid at every restart length.

\begin{lemma}[First-order two-block return map]\label{lem:C14-first-order}
At every positive coprime scalar configuration of restart length $r\ge4$, the core simplex has analytic local coefficient variables $(a,C,H,U)$, where $U\in\mathbb R$ is the single intrinsic normal. Centre $(C,H)$ at their values on the fixed-point family. In any affine normal coordinate $\zeta\in\mathbb R^r$ in weight space formed from $(\delta C,\delta H,U)$, the exact two-block map $\mathscr R_q$ satisfies \begin{equation*} \mathscr R_q(a,0,0)=(a,0,0), \qquad D_\zeta(\mathscr R_q)_\zeta(a,0,0)=I_r. \tag{C.14.1}\label{eq:C.14.1} \end{equation*} For either external node $h$, \begin{equation*} (\mathscr R_q)_{E_h}=E_hm_h, \qquad m_h(a,0,0)=\{P(h)Q(h)/C\}^2=1. \tag{C.14.2}\label{eq:C.14.2} \end{equation*} \end{lemma}

\begin{proof}
Coprimality makes the map \begin{equation*} (p_{\rm var},q_{\rm var})\longmapsto Qp_{\rm var}+Pq_{\rm var}: \mathcal P_{r-1}\times\mathcal P_{r-1}\longrightarrow\mathcal P_{2r-1} \end{equation*} bijective, since $Qp_{\rm var}+Pq_{\rm var}=0$ implies $P\mid p_{\rm var}$; since $\deg p_{\rm var}<\deg P$, $p_{\rm var}=q_{\rm var}=0$, and the dimensions are equal.

At three consecutive phases the scalar coefficient chart is \begin{equation*} S_jA_jA_{j+1} =\Delta(S_jH+U_j)+CS_{j+1}, \qquad (A_0,A_1,A_2)=(P,Q,P), \quad j=0,1. \tag{C.14.4}\label{eq:C.14.4} \end{equation*} On the fixed-point family, $U_j=0$ and $S_j=S=t-a$. For a pure normal variation $u=dU_0$, divide \begin{equation*} \Delta u=SR_{\rm div}+\varrho, \qquad R_{\rm div}\in\mathcal P_{r+1},\quad\varrho\in\mathbb R, \tag{C.14.5}\label{eq:C.14.5} \end{equation*} and solve uniquely \begin{equation*} Qp_{\rm var}+Pq_{\rm var}=R_{\rm div}, \qquad p_{\rm var},q_{\rm var}\in\mathcal P_{r-1}. \end{equation*} Define \begin{equation*} (p_0,p_1,p_2)=(p_{\rm var},q_{\rm var},p_{\rm var}), \quad \dot S_j=-j\varrho/C. \end{equation*} For $j=0,1$, \begin{equation*} \begin{aligned} C\dot S_j+S(p_jA_{j+1}+A_jp_{j+1}) &=-j\varrho+SR_{\rm div}\\ &=\Delta u-(j+1)\varrho =\Delta u+C\dot S_{j+1}. \end{aligned} \tag{C.14.8}\label{eq:C.14.8} \end{equation*} This is exactly the differential of \eqref{eq:C.14.4}, after cancelling the common $\Delta H\,dS_j$ term. After two blocks the phase-polynomial variation returns from $p_{\rm var}$ to $p_{\rm var}$, the normal returns from $u$ to $u$, and $(C,H)$ return; only the $a$-coordinate has the displayed linear shear.

These coefficient variables are local spectral coordinates. For an incoming monic phase polynomial $A$, the weights are \begin{equation*} w_i=\frac{\nu S(\xi_i)}{A(\xi_i)\Delta'(\xi_i)}, \qquad \nu^{-1}=\sum_i\frac{S(\xi_i)}{A(\xi_i)\Delta'(\xi_i)}. \tag{C.14.9}\label{eq:C.14.9} \end{equation*} If a degree-below-$r$ phase variation $v$ induced no first-order weight variation, differentiation of the $r$ orthogonality equations would give \begin{equation*} \sum_iw_iv(\xi_i)\xi_i^j=0, \qquad0\le j<r. \end{equation*} Taking the linear combination with the coefficients of $v$ gives \begin{equation*} \sum_iw_iv(\xi_i)^2=0. \end{equation*} Every weight is positive, hence $v=0$. If also the weights are fixed, differentiating \eqref{eq:C.14.9} gives $$ 0=\frac{\delta\nu}{\nu}-\frac{\delta a}{\xi_i-a} $$ at every core node; subtraction at two nodes gives $\delta a=0$ and then $\delta\nu=0$. The coefficient-to-weight differential is injective; because both $(a,A)$ and the mass-one core simplex have dimension $r+1$, it is an analytic local diffeomorphism.

To relate the phase coefficients to $(C,H,U)$, linearise \eqref{eq:C.14.4} at one phase while fixing $a$. If the variations $(\delta C,\delta H,\delta U)$ vanish, evaluation at the root of $S$ first forces the outgoing variation of $a$ to vanish, and the remaining homogeneous equation is $Qp_{\rm var}+Pq_{\rm var}=0$. Coprimality and the degree bounds then give $p_{\rm var}=q_{\rm var}=0$. The phase coefficient system is square, so its Jacobian is invertible. Equivalently, the $r-1$ directions $(C,H)$ tangent to the fixed-point family, together with the independent $U$-direction explicitly constructed in \eqref{eq:C.14.5}--\eqref{eq:C.14.8}, span all $r$ monic phase-polynomial directions. Hence $(a,C,H,U)\mapsto(a,A)$ has an invertible differential. Composing with the preceding coefficient-to-weight map and applying the analytic inverse-function theorem proves that $(a,C,H,U)$ is an analytic chart on the core simplex.

Passing to affine weight coordinates does not alter the normal identity. A linear $a$-coordinate shear changes the coefficient-to-weight differential by an $a$-derivative term proportional to $\delta a$, but this is multiplied by the first-order normal coordinate and is quadratic; this gives \eqref{eq:C.14.1}. Coordinate divisibility of the squared-weight map and \eqref{eq:C.10.6} give \eqref{eq:C.14.2}.
\end{proof}

We fix the weighted-order notation before using it. Let $\mathfrak B_{qa}$ be a compact subchart in the unscaled variables $(q,a)$, let $y=(U,\mathbf b-\mathbf b_*(q))$, and choose $\delta_0>0$ so that the signed polydiscs $$ \mathcal B_\delta =\{(q,a,y,E):(q,a)\in\mathfrak B_{qa},\ \|y\|_\infty\le\delta,\ \|E\|_\infty\le\delta^2\}, \qquad0<\delta\le\delta_0, $$ remain in the analytic chart. Assign weight one to every component of $y$ and weight two to every component of $E$. Hence the monomial $$ U^i(\mathbf b-\mathbf b_*)^\alpha E^\gamma $$ has weighted degree $i+|\alpha|+2|\gamma|$. For an analytic scalar- or vector-valued function $R$, the notation $R=O_{\rm w}(m)$ means that every Taylor monomial of weighted degree below $m$ vanishes, uniformly in $(q,a)\in\mathfrak B_{qa}$. Equivalently, after decreasing $\delta_0$ if necessary, there is one constant $K_{\rm w}$ such that \begin{equation*} \sup_{\mathcal B_\delta} \bigl\|D_y^\alpha D_E^\gamma R\bigr\| \le K_{\rm w}\delta^{m-|\alpha|-2|\gamma|} \end{equation*} whenever $|\alpha|+2|\gamma|\le m$. This derivative formulation is the one used below; derivatives of larger weighted order are merely required to be uniformly bounded on the fixed polydisc.

\begin{lemma}[Weighted Taylor form]\label{lem:C14-weighted-taylor}
Let $\mathbf b\in\mathbb R^{r-1}$ denote the coefficient vector $(C,H)$ along the fixed-point family, and let $\mathbf b_*(q)$ be its value at the critical equilibrium. The coefficients $\kappa,\beta,\mathbf q_2,A_h,V_h,\ell_h,d_h$ below are evaluated at the fixed-family value $\mathbf b=\mathbf b_*(q)$; their dependence on $q$ is suppressed, and their displayed arguments record any dependence on $a$. The function $D_h$ retains its displayed dependence on $\mathbf b$. Dependence of the unreduced analytic coefficients on $\mathbf b-\mathbf b_*(q)$ is included in the weighted remainders. On a compact positive chart the analytic two-block map has the form \begin{equation*} \begin{aligned} a^+-a&=\kappa(a)U+R_a,\\ U^+-U&=\beta(a)U^2+\sum_hA_h(a)E_h+R_U,\\ \mathbf b^+-\mathbf b&=\mathbf q_2(a)U^2+\sum_hV_h(a)E_h+R_{\mathbf b},\\ \log(E_h^+/E_h)&=D_h(\mathbf b)+\ell_h(a)U+R_{L,h}, \end{aligned} \tag{C.14.12}\label{eq:C.14.12} \end{equation*} where $D_h(\mathbf b_*)=0$, $\nabla_{\mathbf b}D_h(\mathbf b_*)=d_h$, and, assigning weight one to $U,\mathbf b-\mathbf b_*$ and weight two to each $E_h$, \begin{equation*} R_a=O_{\rm w}(2), \qquad R_U=O_{\rm w}(3), \qquad R_{\mathbf b}=O_{\rm w}(3), \qquad R_{L,h}=O_{\rm w}(2). \tag{C.14.13}\label{eq:C.14.13} \end{equation*} The remainders are uniform in the $C^1$ norm on compact subcharts. \end{lemma}

\begin{proof}
The two-block map is the identity on the entire fixed-point family $U=E=0$. Lemma~\ref{lem:C14-first-order} says that its first derivative in the normal and tangent variables of weighted degree one is the identity at every point of the fixed-point family, except for the explicitly retained $a$-coordinate shear. Analytic Hadamard division rules out pure tangent terms and mixed terms linear in $U$: after the displayed linear shear, the core normal and tangent increments lie in the ideal $(U^2,E_h)$. Their quadratic $U^2$ and linear $E_h$ coefficients are, by definition, $\beta,\mathbf q_2,A_h,V_h$; the remaining terms have the weighted degrees stated in \eqref{eq:C.14.13}. Equation \eqref{eq:C.14.2} makes each external coordinate its own analytic factor. The multiplier tangent to the fixed-point family is one at $\mathbf b=\mathbf b_*$; Taylor expansion of its logarithm gives the last line. Taylor's integral formula on a fixed polydisc gives the uniform $C^1$ statement. \end{proof}

\begin{theorem}[Uniform rescaling]\label{thm:C14-rescaling}
Define \begin{equation*} U_n=\frac{u_n}{n}, \qquad \mathbf b_n=\mathbf b_*(q)+\frac{B_n}{n}, \qquad E_{h,n}=\frac{e_{h,n}}{n^2}. \tag{C.14.14}\label{eq:C.14.14} \end{equation*} Let $e_n=(e_{h,n})_h$ and $z_n=(a_n,u_n,B_n,e_n)$. Let $\chi_q$ denote the analytic weight-to-coordinate chart and define the weight reconstruction \begin{equation*} \Psi_{q,n}(a,u,B,e) =\chi_q^{-1}\!\left(a,\frac{u}{n},\mathbf b_*(q)+\frac{B}{n}, \frac{e}{n^2}\right). \end{equation*} Follow one exact two-block map by the inverse scaling with $n+1$, or equivalently define \begin{equation*} F_{n,q}=\Psi_{q,n+1}^{-1}\circ T_r^2\circ\Psi_{q,n}. \tag{C.14.14b}\label{eq:C.14.14b} \end{equation*} The apparent nonautonomy below is solely a change of coordinates for one fixed squared-weight map $T_r^2$. On every compact set of scaled variables in a positive chart, there is an integer $N_0$ such that, for $n\ge N_0$, the resulting exact maps satisfy \begin{equation*} z_{n+1}=z_n+\frac{1}{n}G_q(z_n)+\frac{1}{n^2}R_{n,q}(z_n), \qquad \sup_{n\ge N_0,\,q}\|R_{n,q}\|_{C^1}<\infty, \tag{C.14.15}\label{eq:C.14.15} \end{equation*} where \begin{equation*} \begin{aligned} G_a&=\kappa(a)u,\\ G_u&=u+\beta(a)u^2+\sum_hA_h(a)e_h,\\ G_B&=B+\mathbf q_2(a)u^2+\sum_hV_h(a)e_h,\\ G_{e_h}&=e_h\{2+\ell_h(a)u+d_h[B]\}. \end{aligned} \tag{C.14.16}\label{eq:C.14.16} \end{equation*} \end{theorem}

\begin{proof}
An analytic monomial of weighted degree $m$ in \eqref{eq:C.14.12} becomes $O(n^{-m})$ under \eqref{eq:C.14.14}. Differentiation with respect to a scaled variable of weighted degree one removes one weighted unit and contributes the chain factor $1/n$; differentiation with respect to $e_h$ removes two units and contributes $1/n^2$. The order remains $O(n^{-m})$ uniformly in $C^1$. Also \begin{equation*} D_h(\mathbf b_*+B/n)=n^{-1}d_h[B]+O_{C^1}(n^{-2}). \end{equation*} Multiplying the outgoing $U,\mathbf b-\mathbf b_*$ by $n+1$, and the outgoing external weights by $(n+1)^2$, produces the $+u,+B,+2e_h$ terms in \eqref{eq:C.14.16}. Finally, $$ \begin{aligned} &(1+n^{-1})^2 \exp\{n^{-1}(d_h[B]+\ell_hu)+O_{C^1}(n^{-2})\}\\ &\hspace{5em}=1+n^{-1}\{2+d_h[B]+\ell_hu\}+O_{C^1}(n^{-2}). \end{aligned} $$ This proves the final equation. All remaining terms have a uniform $n^{-2}$ factor with a bounded first derivative. \end{proof}

The coefficients $\kappa,A_h,V_h,\ell_h,d_h$ in \eqref{eq:C.14.16} are precisely those used in \eqref{eq:C.13.1}--\eqref{eq:C.13.8}. Thus the Hopf point continued in scalar coordinates is a Hopf point of the leading field of the exact rescaled map.

\section{Selection of a nonconstant periodic orbit}\label{sec:C15}

The restart-four seed required explicit majorants. At an elevated length only local existence is needed, and the following qualitative form of the same reduction suffices. The preceding section reduces the exact squared-weight map to the leading vector field $G_q$; its simple transverse Hopf point yields the periodic profile that the exact map will shadow.

\begin{theorem}[Periodic orbits from a simple transverse Hopf point]\label{thm:C15-hopf}
Let $G_\mu$ be a $C^4$ family of vector fields with a $C^3$ equilibrium curve $z_*(\mu)$. Suppose at $\mu=0$ the linearisation has a simple pair $\pm i\omega_0$, no other imaginary eigenvalue, and \begin{equation*} \frac{d}{d\mu}\operatorname{Re}\lambda_H(0)\ne0. \end{equation*} Then, for every sufficiently small prescribed nonzero amplitude, there exist a nearby parameter and frequency, together with a nonconstant periodic solution satisfying the phase condition. The conclusion holds independently of the first Lyapunov coefficient. \end{theorem}

\begin{proof}
Translate the equilibrium curve to zero: $$ \widetilde G_\mu(x) =G_\mu(z_*(\mu)+x). $$ Introduce an unknown frequency $\Omega$ and seek $2\pi$-periodic solutions of \begin{equation*} \Omega x'=\widetilde G_\mu(x). \tag{C.15.2}\label{eq:C.15.2} \end{equation*} Use the real Banach spaces $$ X=H^1_{\rm per}([0,2\pi],\mathbb R^d), \qquad Y=L^2_{\rm per}([0,2\pi],\mathbb R^d). $$ Let $A_0=D\widetilde G_0(0)$, choose a complex Hopf eigenvector $v$ and an adjoint eigenvector $v^\dagger$ with $\langle v^\dagger,v\rangle=1$, and let $\Pi_X:X\to X$ and $\Pi_Y:Y\to Y$ be the real rank-two Fourier projections onto $$ \operatorname{span}_{\mathbb R} \{ve^{i\theta}+\overline v e^{-i\theta}, i ve^{i\theta}-i\overline v e^{-i\theta}\} $$ and the corresponding adjoint cokernel. The operator $ \mathscr L_0=\omega_0\partial_\theta-A_0:X\to Y $ is Fredholm of index zero. Simplicity of the Hopf pair and absence of other imaginary eigenvalues make $$ (I-\Pi_Y)\mathscr L_0:(I-\Pi_X)X\longrightarrow(I-\Pi_Y)Y $$ an isomorphism.

Decompose $$ x=\zeta_{\rm H}ve^{i\theta}+\overline{\zeta_{\rm H}}\,\overline v e^{-i\theta}+x_\perp,\qquad x_\perp\in(I-\Pi_X)X. $$ The implicit-function theorem applied to the $(I-\Pi_Y)$ equation gives a unique $$ x_\perp=x_\perp(\zeta_{\rm H},\overline{\zeta_{\rm H}},\mu,\Omega)\in(I-\Pi_X)X $$ near the origin. Time translation by $\phi$ acts as $\zeta_{\rm H}\mapsto e^{i\phi}\zeta_{\rm H}$ and preserves \eqref{eq:C.15.2}. Uniqueness of the complementary solution makes $x_\perp$ equivariant under this action. Pairing the remaining $\Pi_Y$ equation with $v^\dagger$ gives an $S^1$-equivariant complex scalar equation. Its smooth equivariant factorisation is \begin{equation*} \zeta_{\rm H}\Phi(|\zeta_{\rm H}|^2,\mu,\Omega)=0. \end{equation*} At $\zeta_{\rm H}=0$ the complementary solution is $x_\perp(0,0,\mu,\Omega)=0$. Differentiate the linearised complementary and reduced equations at $(\mu,\Omega)=(0,\omega_0)$. Every term containing a derivative of the complementary correction lies in the range of $ \mathscr L_0=\omega_0\partial_\theta-A_0 $ and vanishes after pairing with the adjoint cokernel vector $v^\dagger$. Using $\langle v^\dagger,v\rangle=1$ and the standard derivative formula for a simple eigenvalue, one fixed nonzero complex normalisation of the reduced equation gives \begin{equation*} \partial_\mu\Phi(0,0,\omega_0)=\lambda_H'(0), \qquad \partial_\Omega\Phi(0,0,\omega_0)=-i. \end{equation*} Hence \begin{equation*} D_{(\mu,\Omega)} \binom{\operatorname{Re}\Phi}{\operatorname{Im}\Phi} (0,0,\omega_0) = \begin{pmatrix} \operatorname{Re}\lambda_H'(0)&0\\ \operatorname{Im}\lambda_H'(0)&-1 \end{pmatrix}. \end{equation*} Its determinant is $-\operatorname{Re}\lambda_H'(0)\ne0$. For every sufficiently small prescribed amplitude $\varepsilon>0$, fix the phase by taking $\zeta_{\rm H}=\varepsilon$ and apply the implicit-function theorem to $\Phi(\varepsilon^2,\mu,\Omega)=0$. The theorem determines $\mu=\mu(\varepsilon^2)$, $\Omega=\Omega(\varepsilon^2)$, and the unique complementary correction $x_\perp$. The resulting periodic profile satisfies the phase condition and, since $\varepsilon\ne0$, is nonconstant. Thus transversality in the parameter direction suffices, independently of the first Lyapunov coefficient. \end{proof}

For an elevated configuration, Lemma~\ref{lem:C13-hopf-persistence} supplies the hypotheses. Because the Hopf eigenvector has a nonzero $a$-component, the first harmonic of $a$ on the selected periodic orbit is nonzero. Taking the amplitude small also keeps every core weight and both scaled external amplitudes strictly positive.

\section{Shadowing with harmonic time steps}\label{sec:C16}

Section~\ref{literal-selection-of-one-exact-map-shadow} supplied quantitative shadowing for the certified seed. After degree elevation, the uniform remainder in \eqref{eq:C.14.15} allows the initial index to be chosen after the configuration and periodic profile. The following qualitative statement accommodates saddle periodic orbits and neutral radial Floquet multipliers.

\begin{theorem}[Shadowing without hyperbolicity]\label{thm:C16-shadowing}
Let $\mathcal U\subset\mathbb R^d$ be open, let $G\in C^3(\mathcal U,\mathbb R^d)$, and let $\gamma$ be a periodic orbit of $\dot z=G(z)$. Fix an integer $N_0\ge1$. On one fixed neighbourhood $\mathcal V\Subset\mathcal U$ of $\gamma$, suppose that $F_n,\mathscr R_n\in C^1(\mathcal V,\mathbb R^d)$ are given for every $n\ge N_0$ and satisfy \begin{equation*} F_n(z)=z+n^{-1}G(z)+n^{-2}\mathscr R_n(z), \qquad \sup_{n\ge N_0}\|\mathscr R_n\|_{C^1(\mathcal V)}<\infty. \tag{C.16.1}\label{eq:C.16.1} \end{equation*} The periodic orbit may have arbitrarily many unstable directions and arbitrary neutral Jordan blocks. Let $\beta_{\rm st}$ be the smallest absolute value of the negative real parts of the stable Floquet exponents, with $\beta_{\rm st}=\infty$ if there is no stable exponent. For every \begin{equation*} 0<\rho<\min(1,\beta_{\rm st}) \end{equation*} there is an integer $N_1\ge N_0$ such that, for every $N\ge N_1$, there is a forward orbit of $F_n$ in $\mathbb R^d$ satisfying \begin{equation*} z_n-\gamma(\tau_n+\phi)=O(n^{-\rho}), \qquad \tau_N=0,\qquad \tau_n=\sum_{k=N}^{n-1}\frac{1}{k}\quad(n>N), \tag{C.16.3}\label{eq:C.16.3} \end{equation*} for a prescribed phase $\phi$. \end{theorem}

\begin{proof}
Let $\gamma_n=\gamma(\tau_n+\phi)$ and $z_n=\gamma_n+\eta_n$. Taylor expansion of \eqref{eq:C.16.1} and of $\gamma(t+1/n)$ gives \begin{equation*} \eta_{n+1}=\{I+n^{-1}\mathsf A(\tau_n+\phi)\}\eta_n +n^{-1}\mathcal N_n(\eta_n)+n^{-2}\mathcal R_n(\eta_n), \end{equation*} where $\mathsf A(t)=DG(\gamma(t))$, \begin{equation*} |\mathcal N_n(e)|\le C|e|^2, \quad |\mathcal N_n(e)-\mathcal N_n(f)|\le C(|e|+|f|)|e-f|, \quad \sup_{n\ge N_0}\|\mathcal R_n\|_{C^1}<\infty. \tag{C.16.5}\label{eq:C.16.5} \end{equation*} All phase and radial neutral couplings remain in the linear variational equation.

Let $T_\gamma$ be a period of $\gamma$ and let $\mathsf X(0)=I$ be the real fundamental matrix of $x'=\mathsf A(t)x$. Its monodromy $\mathscr M=\mathsf X(T_\gamma)$ is real and invertible. The square $\mathscr M^2=\mathsf X(2T_\gamma)$ has a real logarithm. To see this, recall that the real-logarithm criterion \cite[Theorem~1]{Culver1966} requires invertibility and, for each negative real eigenvalue, an even number of Jordan blocks of every size. Real eigenvalues of $\mathscr M$ square to positive numbers. Nonreal eigenvalues occur in conjugate pairs; a negative eigenvalue of $\mathscr M^2$ can arise only from a purely imaginary pair, whose two Jordan blocks have the same size after squaring because the derivative of $z\mapsto z^2$ is nonzero there. The criterion is satisfied, so we may choose a real matrix $B_{\rm Fl}$ with $$ e^{2T_\gamma B_{\rm Fl}}=\mathscr M^2 $$ and define $\mathsf P(t)=\mathsf X(t)e^{-tB_{\rm Fl}}$. Then $\mathsf P$ is real, invertible, and $2T_\gamma$-periodic. This is the real period-doubled Floquet representation $$ \mathsf X(t)=\mathsf P(t)e^{tB_{\rm Fl}}. $$ The real parts of the spectrum of $B_{\rm Fl}$ are the usual Floquet growth rates, so period doubling does not change $\beta_{\rm st}$.

Since $\mathsf P'=\mathsf A\mathsf P-\mathsf PB_{\rm Fl}$ and $\mathsf P(t+1/n)=\mathsf P(t)+n^{-1}\mathsf P'(t)+O(n^{-2})$, set $\eta_n=\mathsf P(\tau_n+\phi)y_n$. The transformed real error satisfies \begin{equation*} y_{n+1}=(I+B_{\rm Fl}/n)y_n +n^{-1}\widetilde{\mathcal N}_n(y_n) +n^{-2}\widetilde{\mathcal R}_n(y_n), \end{equation*} with the same bounds as \eqref{eq:C.16.5}. Split the real generalised spectral subspaces of $B_{\rm Fl}$ according to the real parts of their complexified spectra: \begin{equation*} E^{\rm st}=\bigoplus_{\operatorname{Re}\lambda<0}E_\lambda^{\mathbb R}, \qquad E^{\rm cu}=\bigoplus_{\operatorname{Re}\lambda\ge0}E_\lambda^{\mathbb R}. \end{equation*} After slightly decreasing the stable gap, the products of $I+B_{\rm Fl}/j$ satisfy \begin{equation*} \|\Phi_{\rm st}(n,k)\|\le C(n/k)^{-\beta}, \qquad n\ge k, \end{equation*} for some $\rho<\beta<\beta_{\rm st}$. For every $\varepsilon>0$, the inverse products on $E^{\rm cu}$ satisfy \begin{equation*} \|\Phi_{\rm cu}(n,k)\| \le C_\varepsilon(k/n)^\varepsilon, \qquad k\ge n. \end{equation*} Positive-real-part blocks decay backwards. A zero-real-part Jordan block contributes only a fixed power of $\log(k/n)$, which is bounded by $C_\varepsilon(k/n)^\varepsilon$. Choose \begin{equation*} 0<\varepsilon<\min(\rho,1-\rho). \tag{C.16.10}\label{eq:C.16.10} \end{equation*}

Let $\mathsf A_j=I+B_{\rm Fl}/j$, and let $\Pi_{\rm st},\Pi_{\rm cu}$ be the spectral projections of $B_{\rm Fl}$. The evolution operators used below are \begin{equation*} \begin{aligned} \Phi_{\rm st}(n,k+1)&=\mathsf A_{n-1}\cdots \mathsf A_{k+1}\Pi_{\rm st}, &&N\le k<n,\\ \Phi_{\rm cu}(n,k+1)&=\mathsf A_n^{-1}\cdots \mathsf A_k^{-1}\Pi_{\rm cu}, &&k\ge n. \end{aligned} \end{equation*} Choose $N$ large enough that every $\mathsf A_j$ is invertible. On the Banach space with norm $\|y\|_\rho=\sup_{n\ge N}n^\rho|y_n|$, define \begin{equation*} \begin{aligned} (\mathcal Ly)_n^{\rm st} &=\sum_{k=N}^{n-1}\Phi_{\rm st}(n,k+1)f_k(y_k),\\ (\mathcal Ly)_n^{\rm cu} &=-\sum_{k=n}^{\infty}\Phi_{\rm cu}(n,k+1)f_k(y_k), \end{aligned} \tag{C.16.11}\label{eq:C.16.11} \end{equation*} where \begin{equation*} f_k(y)=k^{-1}\widetilde{\mathcal N}_k(y)+k^{-2}\widetilde{\mathcal R}_k(y). \end{equation*} On a fixed ball, \begin{equation*} |f_k(y_k)|\le C(k^{-1-2\rho}+k^{-2}). \end{equation*} Integral comparison gives \begin{equation*} \begin{aligned} \sum_{k=N}^{n-1}(n/k)^{-\beta} (k^{-1-2\rho}+k^{-2})&=O(n^{-\rho}),\\ \sum_{k=n}^{\infty}(k/n)^\varepsilon k^{-1-2\rho} &=O(n^{-2\rho}),\\ \sum_{k=n}^{\infty}(k/n)^\varepsilon k^{-2} &=O(n^{-1}). \end{aligned} \end{equation*} The last two quantities are $O(n^{-\rho})$ by \eqref{eq:C.16.10}. The difference estimate gains the factor $O(k^{-1-\rho}+k^{-2})$, so the Lipschitz constant of $\mathcal L$ is $O(N^{-\rho}+N^{-1})$. For large $N$, $\mathcal L$ maps a fixed ball into itself and is a contraction. Its fixed point solves the recurrence and satisfies \eqref{eq:C.16.3}. Every coefficient, projection, evolution operator, and nonlinear term is real, and the contraction acts on the real Banach space of real sequences. Its fixed point, and hence the resulting sequence $z_n$, is real. The stable component has been given a fixed initial value (zero in \eqref{eq:C.16.11}); the boundary condition at infinity determines every unstable, phase, neutral radial, and resonant coordinate.
\end{proof}

The additional $+1$ eigenvalues introduced by degree elevation give unstable Floquet directions in the finite-dimensional complementary subspace covered by Theorem~\ref{thm:C16-shadowing}.

\section{Positivity, nontermination, and nonconvergence}\label{sec:C17}

Sections~\ref{sec:C14}--\ref{sec:C16} supply, at each elevated restart length, an orbit shadowing a nonconstant periodic profile. The final step preserves the strict positivity margins, excludes termination, and translates nonconvergence of the squared weights into nonconvergence of the residuals.

Apply Theorems~\ref{thm:C14-rescaling},~\ref{thm:C15-hopf}, and~\ref{thm:C16-shadowing} to a fixed elevated configuration. Denote the parameter of the resulting periodic orbit by $q_{\rm per}$. Throughout this section, $P,Q,\Delta$ and the node labels refer to their values at $q_{\rm per}$. Choose the periodic orbit in a compact tube on which all limiting core weights and scaled external amplitudes have positive lower bounds, while the first- and second-phase block factors are bounded away from zero. For sufficiently large $N$, the shadowing error and the factors $n^{-1},n^{-2}$ keep every reconstructed even and odd weight positive and every block factor nonzero. There are $r+2$ core nodes and two external nodes, hence $r+4>r$ active nodes at every restart point. The grade is $r+4$, so no block with restart length $r$ terminates.

The harmonic times $\tau_n$ are increasing and unbounded, while $\tau_{n+1}-\tau_n=1/n\to0$. For any phase modulo the period, the first $n$ for which $\tau_n$ passes that phase in the $j$-th period has overshoot at most $1/(n-1)$; hence every phase is an accumulation point modulo the period. The recovery root on the periodic orbit is nonconstant, and on the limiting $P$-phase of the fixed-point family \begin{equation*} w_i^P(a)=\frac{(\xi_i-a)Q(\xi_i)}{\Delta'(\xi_i)}. \end{equation*} Every slope is nonzero, so this affine map is injective. The even squared-weight sequence has at least two distinct accumulation points. If the even normalised residuals converged in Euclidean norm, every coordinate square would converge, contrary to the two accumulation points. Hence the even residual subsequence is nonconvergent.

The spectral nodes constructed so far need not be positive. Label the combined node set in the fixed chart order by $(\widehat\xi_i)_{i=1}^{r+4}=(\xi_1,\ldots,\xi_{r+2},h_-,h_+)$. Choose an integer $L$ with $-L<\min_{1\le i\le r+4}\widehat\xi_i$, and let $\lambda_i=\widehat\xi_i+L>0$. If $P_w$ is the unshifted monic degree-$r$ orthogonal polynomial, all its zeros lie strictly inside the convex hull of the real support of the positive measure $\sum_{i=1}^{r+4}w_i\delta_{\widehat\xi_i}$. Hence $P_w(-L)\ne0$. Define \begin{equation*} \widetilde\pi_w(t)=\frac{P_w(t-L)}{P_w(-L)}. \tag{C.17.2}\label{eq:C.17.2} \end{equation*} Then $\widetilde\pi_w(0)=1$, and for $0\le j<r$, \begin{equation*} \sum_{i=1}^{r+4}w_i\widetilde\pi_w(\lambda_i)\lambda_i^j =\frac{1}{P_w(-L)} \sum_{i=1}^{r+4}w_iP_w(\widehat\xi_i)(\widehat\xi_i+L)^j=0. \end{equation*} Uniqueness makes \eqref{eq:C.17.2} the restarted CG residual polynomial for the shifted nodes. Its node values differ from the unshifted values by one common nonzero scalar, which cancels from the normalised squared-weight map. Positivity, nontermination, and nonconvergence are unchanged by the shift.

Let $z_N$ be the Lyapunov--Perron fixed point at the chosen initial index and define \begin{equation*} W_0=\Psi_{q_{\rm per},N}(z_N). \end{equation*} The covariance \eqref{eq:C.14.14b} gives, for every $j\ge0$, \begin{equation*} T_r^{2j}W_0=\Psi_{q_{\rm per},N+j}(z_{N+j}). \end{equation*} Hence the reconstructed sequence is an orbit of the fixed squared-weight map $T_r^2$. The intermediate states are $T_r^{2j+1}W_0$, and the strict two-phase factor margins prove that all of their coordinates remain positive as well.

Define $A$ to be diagonal with the distinct positive nodes $\lambda_1,\ldots,\lambda_{r+4}$ as its entries, let $x_0=0$, and take $b_i=\sqrt{W_{0,i}}$ with the positive roots. Since $\sum_iW_{0,i}=1$, the initial residual has exactly those squared spectral weights. This reconstructs one restarted CG instance with a real SPD matrix and proves nontermination for all iterates and nonconvergence of the even subsequence.

\section{Induction and exact existence}\label{sec:C18}

Fix an integer $s\ge4$. For $s=4$, Proposition~\ref{prop:s4-seed} supplies the counterexample. If $s>4$, the construction in Sections~\ref{the-exact-algebraic-hopf-seed}--\ref{exact-transfer-of-the-hopf-point-to-the-raw-field} supplies the positive transverse Hopf configuration at restart length four. Apply Lemmas~\ref{lem:C11-root-continuation}--\ref{lem:C13-hopf-persistence} exactly $s-4$ times, with $d=1$ and at each step a sufficiently large integer scale. The final configuration has $s+2$ core nodes; the two bounded external nodes give exactly \begin{equation*} n=s+4. \end{equation*} Section~\ref{sec:C17} converts it into the required SPD orbit. The elevation scales, isolating neighbourhoods, periodic amplitude, and shadowing index may each depend on the prescribed value of $s$.

At each elevation, choose a rational neighbourhood on which the preceding critical point is the unique solution of the critical equations and all its strict root, weight, factor, and transversality inequalities hold. Lemmas~\ref{lem:C11-root-continuation}--\ref{lem:C13-hopf-persistence} show that, for every sufficiently large integer $R$, the elevated equations have a unique critical point in the corresponding continuation neighbourhood and the same strict inequalities hold. Choose any such integer. Repeating this finite construction gives a real configuration at every prescribed $s\ge4$.

All estimates required in the subsequent Lyapunov--Schmidt and Lyapunov--Perron arguments are finite on a sufficiently small compact neighbourhood of the selected critical point. Choose a sufficiently small nonzero dyadic amplitude $2^{-m}$ and then an initial index $N$ large enough for the contraction and positivity conditions of Theorem~\ref{thm:C16-shadowing}.

Once these finite choices are fixed, the local implicit-function construction selects a periodic profile and the Lyapunov--Perron contraction has a unique fixed point. Together they determine an exact real SPD instance.

\begin{proof}[Proof of Theorem~\ref{thm:negative}] Proposition~\ref{prop:s4-seed} gives the assertion for $s=4$. For $s>4$, the induction in this section produces a positive transverse Hopf configuration of restart length $s$, and Section~\ref{sec:C17} turns it into a fixed SPD instance in dimension $s+4$ whose iteration is nonterminating and whose even normalised residuals do not converge. \end{proof}

\begin{proof}[Proof of Theorem~\ref{thm:sharp}] Assertion \textup{(i)} follows from Theorems~\ref{s2:thm:main} and~\ref{thm:main}, and assertion \textup{(ii)} is Theorem~\ref{thm:negative}. \end{proof}

\section{An exact lift for restart lengths divisible by four}\label{sec:C19}

When the restart length is divisible by four, there is also an independent exact conjugacy. Let $d\ge1$, let $T_d^{\rm Ch}$ be the Chebyshev polynomial of the first kind, and define \begin{equation*} \Phi_d(x)=16\{1-(-1)^dT_d^{\rm Ch}(x-1)\}. \end{equation*} Then $\Phi_d(0)=0$. Every spectral node in the base configuration lies in $(1,10)$, and each equation $\Phi_d(x)=\lambda_i$ has exactly $d$ distinct roots $x_{i,a}\in(0,2)$. Define \begin{equation*} (\mathcal L_dw)_{i,a}=w_i/d, \qquad (J_d\eta)_{i,a}=\eta_i/\sqrt d. \end{equation*} For \begin{equation*} M_j(y)=\frac{1}{d}\sum_{\Phi_d(x)=y}x^j, \end{equation*} Newton's identities give \begin{equation*} M_j\in\mathbb R[y], \qquad \deg M_j\le\lfloor j/d\rfloor. \end{equation*} After making $\Phi_d(x)-y$ monic, only its constant coefficient depends on $y$, and it depends affinely.

If $p_{4,w}$ is the residual polynomial of the base map, then for $0\le j<4d$, \begin{equation*} \begin{aligned} \sum_{i,a}\frac{w_i}{d} p_{4,w}(\Phi_d(x_{i,a}))x_{i,a}^j &=\sum_iw_ip_{4,w}(\lambda_i)M_j(\lambda_i)=0, \end{aligned} \end{equation*} because $\deg M_j\le3$. The composition has degree $4d$ and value one at zero, so uniqueness gives \begin{equation*} p_{4d,\mathcal L_dw}=p_{4,w}\circ\Phi_d. \end{equation*} If $\mathcal F_s$ denotes the signed normalised-residual map for one block of restarted CG with restart length $s$, then \begin{equation*} T_{4d}\mathcal L_d=\mathcal L_dT_4, \qquad \mathcal F_{4d}J_d=J_d\mathcal F_4. \end{equation*} The map $J_d$ is an isometry, all $8d$ components remain active, and the grade is $8d>4d$. The base parity nonconvergence transfers under this conjugacy, giving an independent dimension-$2s$ family for $s=4d$.

\section*{Consequences and further questions}

Together with Akaike's theorem for steepest descent, Theorem~\ref{thm:sharp} shows that universal convergence along the even and odd subsequences holds precisely for $s\in\{1,2,3\}$. For each $s\ge4$, the counterexample constructed here shows that the corresponding universal statement is false. The error norms still decrease; the change at $s=4$ occurs in the normalised spectral dynamics.

The counterexamples have dimension $s+4$. Natural questions include the smallest dimension in which nonconvergence can occur at each restart length, and whether the Hopf mechanism persists when the spectra are subject to additional structure. The examples constructed here are selected by algebraic isolating boxes and contraction mappings; more convenient closed-form representatives would also be useful.

For $s=2$, directional convergence yields the vector Aitken-type postprocessing described in \cite{ColbrookStepaniantsTownsend2026}. The theorem for restart length three suggests an analogous acceleration, although its finite-precision behaviour and the effect of inexact inner solves remain to be studied. For $s\ge4$, any extrapolation theory must accommodate recurrent normalised directions rather than a universal two-cycle.

% Keep repository paths readable by breaking them at punctuation, not letters.
\begingroup
\def\UrlBreaks{\do\/\do\-\do\_\do\.}
\Urlmuskip=0mu
\section*{Appendix C. Exact certification of the configuration at restart length four}
\label{app:cap}

The computer-assisted numerical part of the proof is confined to the finite verification of the configuration at restart length four. The ancillary material also contains supplementary finite regression checks for the higher-degree identities proved analytically in the text. For general background on validated numerics, see \cite{Tucker2011}; for a related computer-assisted approach to the rigorous verification of Hopf bifurcations, see \cite{VanDenBergLessardQueirolo2021}. The current manuscript, computer-assisted proof, and Lean formalisation are available at \url{https://github.com/sgstepaniants/Forsythe}. Paths below are relative to the root of this repository; the computer-assisted proof is in \path{computer-assisted-proof/}, and the Lean formalisation is described in \hyperref[app:lean]{Appendix D}.

The eight programs listed below are in \path{computer-assisted-proof/evidence/launch/}; their filenames within that directory are:

\begin{enumerate} \def\labelenumi{\arabic{enumi}.} \item \path{all_s_s4_transverse_hopf_exact.py}

This verifies the algebraic root branches, the unique positive critical branch, the unique transverse Hopf parameter, the kernel evaluations, the tangent-rate and positivity bounds, and the $\theta$-quadratic invariants. Equations \eqref{eq:C.3.17a}--\eqref{eq:C.3.21} derive the full characteristic polynomial from these invariants. \item \path{all_s_s4_coupling_identity_exact.py}

This is a separate interval consistency and strict-sign check for the coupling formulas. It constructs $p_g$ and $q_g$, obtains $(V_g)^H$ by the literal polynomial division in \eqref{eq:C.3.9}, and then compares \eqref{eq:C.3.10}, \eqref{eq:C.3.10a}, and \eqref{eq:C.13.2}, including their five retained jet components. Its division residuals are enclosed in intervals containing zero, and its compared jet intervals overlap. These tests establish consistency, not exact equality; the exact identities follow from the polynomial derivations in \eqref{eq:C.3.9}--\eqref{eq:C.3.10a} and \eqref{eq:C.13.2}. It does not use the stored entries of $\mathsf T$ to construct $(V_g)^H$. The comparison uses target order $(h_-,h_+)$ for the rows and source order $(h_-,h_+)$ for the columns, and verifies the value-level sign pattern $[-,+;+,-]$. \item \path{all_s_s4_local_majorants_exact.py}

This verifies the root separations, the core-weight and factor margins, the inverse of the phase Gram matrix, and the inverse bounds for both $9\times9$ phase matrices. \item \path{all_s_s4_intrinsic_lift_exact.py},\\ \path{all_s_s4_scalar_hopf_vector_exact.py}, and\\ \path{all_s_s4_scalar_to_raw_tangent_exact.py}

These enclose the $4\times4$ intrinsic lift and its quadratic coefficients, verify the weighted tangent reconstruction and the hypotheses needed for its invertibility, and establish the scalar block and $a$-component eigenvector bounds in \eqref{eq:C.5.13}--\eqref{eq:C.5.14}. The algebraic recursions \eqref{eq:C.3.14a}--\eqref{eq:C.3.14e} supply the coefficient identities underlying these enclosures. \item \path{all_s_s4_explicit_counterexample_checker.py}

This combines the transverse-Hopf, local-majorant, intrinsic-lift, and scalar-to-weight tangent calculations, and verifies the retained binary-exponent ledger for Lemmas~\ref{lemma-6.1-quantitative-hopf-contraction} and~\ref{lemma-7.1-quantitative-nonhyperbolic-shadowing}. The \path{computer-assisted-proof/verify.sh} driver separately invokes the coupling and scalar checks before the aggregate calculation. \item \path{independent_all_s_symbolic_checks.py}

This supplementary SymPy calculation checks the characteristic reduction on seeded exact rational test instances with $r-1=3,\ldots,6$, corresponding to restart lengths $r=4,\ldots,7$, and checks the Newton-sum degree bound for elevation degrees $d=1,\ldots,8$. These are finite regression tests, not proofs for arbitrary degree; the algebraic proofs of the general statements appear in Sections~\ref{sec:C10} and~\ref{sec:C19}.

\end{enumerate}

The interface with the proof has three parts. First, the polynomial coefficients, terminating decimal seed data, and interval endpoints are exact rational inputs. Sturm counts and the ordered rational intervals in \eqref{eq:C.3.3}--\eqref{eq:C.3.3a} enclose and label the algebraic roots used in \eqref{eq:C.3.4}--\eqref{eq:C.3.5}. Rational interval propagation then encloses the derived algebraic quantities and establishes the denominator signs, positivity margins, transverse Hopf crossing, and quantitative majorants invoked in Sections~\ref{the-exact-algebraic-hopf-seed}--\ref{literal-selection-of-one-exact-map-shadow} at restart length four. Root ordering, branch choice, nonzero denominators, and strict containment are fixed before any division or contraction.

Second, once those strict finite bounds are available, the coordinate changes in Sections~\ref{a-complete-physical-raw-chart-and-its-exact-expansion}--\ref{exact-transfer-of-the-hopf-point-to-the-raw-field}, the two Banach contractions in Sections~\ref{literal-selection-of-one-periodic-orbit}--\ref{literal-selection-of-one-exact-map-shadow}, and the reconstruction in Sections~\ref{harmonic-times-see-the-whole-periodic-orbit}--\ref{the-one-fixed-exact-spd-datum} are analytic consequences of these finite premises. Equations \eqref{eq:C.4.2}, \eqref{eq:C.5.3}, \eqref{eq:C.5.8b}, \eqref{eq:C.5.10}, and \eqref{eq:C.9.1}--\eqref{eq:C.9.5} carry the certified scalar data to squared spectral weights and then to the fixed restart-four CG orbit for the explicitly defined family.

Third, Sections~\ref{sec:C10}--\ref{sec:C18} are analytic: they contain the degree elevation, qualitative periodic-orbit construction, shadowing theorem, and final SPD reconstruction. The final item supplies only the finite supplementary tests just described.

The programs for restart length four use integers, \texttt{Fraction} arithmetic, closed rational intervals, Euclidean polynomial division, and Sturm counts. Each terminating decimal in the input is parsed digit by digit as its base-ten rational value, so no input hypothesis is routed through binary floating-point arithmetic. Every interval division is performed only after its denominator interval has been shown to exclude zero, and all sign and root-count assertions are exact. The supplementary all-$s$ checks use exact SymPy expressions.

The file \path{computer-assisted-proof/CERTIFICATE_LEDGER.md} maps the retained groups of finite premises in \eqref{eq:C.3.3}--\eqref{eq:C.7.15} to the corresponding programs and direct assertions, and records the dependency and independence status, coordinate order, sign convention, and numerical bounds used by the analytic proof. The four stages of \path{computer-assisted-proof/verify.sh} separately run the coupling-identity, scalar-vector, aggregate restart-four, and supplementary all-$s$ checks; the driver first rejects disabled Python assertions and checks the pinned symbolic dependencies when they are required. The files \path{computer-assisted-proof/verification-lock.json} and \path{computer-assisted-proof/VERIFICATION_ENVIRONMENT.md} record the machine-readable requirements and tested environments, while \path{computer-assisted-proof/verification-transcript.txt} records a complete replay. The internal \path{computer-assisted-proof/SHA256SUMS} manifest covers exactly the seven restart-four programs whose calculations supply finite premises to the proof; documentation, replay and archive tooling, environment metadata, transcripts, and the supplementary all-$s$ regression program are intentionally outside it. The file \path{computer-assisted-proof/build_archive.py} creates the ancillary ZIP from deterministic directory and file lists with fixed timestamps and fixed POSIX modes.

To replay the four stages from the repository root, install the dependencies in \path{computer-assisted-proof/requirements-verification.txt} and run
\begin{verbatim}
bash computer-assisted-proof/verify.sh
\end{verbatim}
The setup instructions are in \path{computer-assisted-proof/README.md}; \path{computer-assisted-proof/checks/README.md} describes the additional mathematical sanity checks. The finite calculations establish the algebraic inequalities used above. Banach's theorem gives the unique fixed points in the balls and gauges specified in Sections~\ref{literal-selection-of-one-periodic-orbit}--\ref{literal-selection-of-one-exact-map-shadow}; Sections~\ref{harmonic-times-see-the-whole-periodic-orbit}--\ref{the-one-fixed-exact-spd-datum} then prove phase accumulation, SPD reconstruction, and nonconvergence.

\section*{Appendix D. Lean formalisation and verification}
\label{app:lean}

The directory \href{https://github.com/sgstepaniants/Forsythe/tree/main/lean-proof}{\nolinkurl{lean-proof/}} in the same repository contains a Lean proof of Theorem~\ref{thm:sharp}. The proof includes the analytic deductions from the algebraic Hopf input to exact CG trajectories. Its entry point, \path{lean-proof/Solution.lean}, exports the following five declarations, all in the namespace \path{ProofProject}:
\begin{enumerate}
\item \path{forsytheConjecture_restartLengthTwo} proves Theorem~\ref{s2:thm:main}.
\item \path{forsytheConjecture_restartLengthThree} proves Theorem~\ref{thm:main}.
\item \path{restartFourAlgebraicHopfCertificate} proves the exact algebraic Hopf certificate associated with Lemma~\ref{lemma-3.1-certified-critical-branch-and-transverse-hopf-point}, including the characteristic reduction in \eqref{eq:C.3.17a}--\eqref{eq:C.3.21}.
\item \path{counterexamples_restartLengthAtLeastFour} proves Theorem~\ref{thm:negative}, with dimension exactly $s+4$ for every $s\ge4$.
\item \path{forsytheSharpClassification} proves Theorem~\ref{thm:sharp}, including its stronger diagonal conclusion.
\end{enumerate}

The positive conclusions concern the signed residual directions, normalised in Euclidean norm, separately at even and odd restart boundaries, with finite termination allowed. The negative conclusions include an actual nonterminating exact restarted-CG trajectory whose even normalised residuals have no Euclidean limit. The formal proof uses the qualitative periodic-orbit, shadowing, and degree-elevation route of Sections~\ref{sec:C10}--\ref{sec:C18}. Its scope does not include a separate proof of Akaike's $s=1$ theorem or the particular $\varepsilon_0,M,N$ specification of Proposition~\ref{prop:s4-seed}; these are not assertions of Theorem~\ref{thm:sharp}.

The mathematical statement can be reviewed through two Lean files and one short configuration file:
\begin{enumerate}
\item \path{lean-proof/Challenge.lean} states the five theorem signatures. Its deliberate \path{sorry} bodies are placeholders for the statements; the proved declarations are supplied independently by \path{lean-proof/Solution.lean}.
\item \path{lean-proof/ProofProject/Definitions.lean} defines the exact affine-Krylov energy minimisation, termination, signed Euclidean normalisation, parity convergence, and counterexample predicates. These definitions occupy lines 18--127. The remaining definitions, through line 414, state the separate algebraic Hopf certificate, including its exact rational inputs, domains, selected roots, identities, and strict margins.
\item \path{lean-proof/comparator.json} selects the five fully qualified theorem names and permits only \path{propext}, \path{Classical.choice}, and \path{Quot.sound}.
\end{enumerate}

Comparator compares the declarations and their statement dependencies between the challenge and solution, checks transitive axiom dependencies, and replays the exported proof declarations in an initially empty Lean kernel environment, including a quotient-consistency check. All five exports depend only on the three axioms listed above: propositional extensionality, classical choice, and quotient soundness. The numerical certificate is proved within Lean; successful execution of an external Python program is not a premise of the formal theorem. Given the kernel check and the standard Lean/Mathlib definitions, the remaining human task is to confirm that the specified statements express the intended mathematics. It does not require reading all of the proof implementation files.

The project retains Lean 4.33.1 and its dependency selection in \path{lean-proof/lean-toolchain}, \path{lean-proof/lakefile.toml}, and \path{lean-proof/lake-manifest.json}. From the repository root, the ordinary build is
\begin{verbatim}
cd lean-proof
lake build Solution
\end{verbatim}
The Linux prerequisites and negative controls for the full isolated Comparator check are documented in \path{lean-proof/reproduction/README.md}. The driver is \path{lean-proof/scripts/verify.sh}.

The successful 5 September 2026 run is recorded in \path{lean-proof/VERIFICATION.md}, with its complete transcript in \path{lean-proof/verification-logs/linux-comparator-final.log}.

\endgroup

\section*{Acknowledgements}

M. J. C. and A. T. thank the organisers of the 2025 Simons workshop at Berkeley, where they learned of the Forsythe conjecture. A. T. acknowledges support from the Defense Advanced Research Projects Agency through The Right Space (TRS) Disruption Opportunity 25 (DARPA-PA-24-04-07) and from NSF CAREER grant DMS-2045646. G. S. acknowledges support from the NSF Mathematical Sciences Postdoctoral Research Fellowship under award 2402074.
MJC and GS also thank the Isaac Newton Institute for Mathematical Sciences, Cambridge, for its support and hospitality during the programme \emph{Operator Methods for Dynamical Systems}, and in particular for the freedom afforded by many late nights at the Institute and infinite supply of coffee, during which this solution first took shape! This work was supported by EPSRC grant no EP/Z000580/1.

\bibliographystyle{amsalpha} \bibliography{forsythe}

\end{document}